Pontificia Universidad Católica del Perú

Escuela de Posgrado

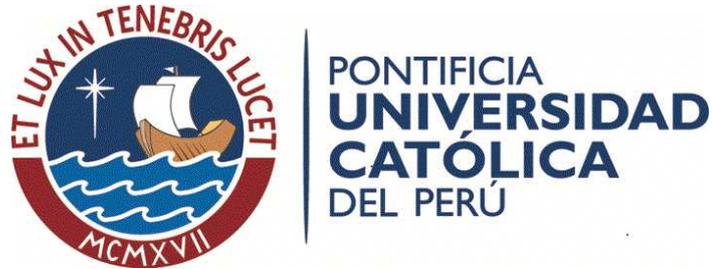

# Clasificación de planos torcidos graduados

Tesis para optar el grado académico de
Doctor en Matemáticas

Autor

Ricardo Manuel Bances Hernández

Asesor

Christian Holger Valqui Haase

Jurado

Dr. Jorge Alberto Guccione
Dr. Jack Denne Arce Flores
Dr. Gabriel Armando Muñoz Márquez
Dr. Guillermo Cortiñas

Lima - Perú

2021

# Resumen

En esta tesis se obtiene una clasificación casi completa de todos los productos tensoriales torcidos graduados de $K[x]$ con $K[y]$. Para ello se usa una representación de un producto tensorial torcido graduado de $K[x]$ con $K[y]$ en el álgebra $L(K^{\mathbb{N}_0})$, la cual está inmersa en el conjunto de matrices infinitas con entradas en $K$. De esta manera el problema de clasificar a los productos tensoriales torcidos graduados de $K[x]$ con $K[y]$ se traduce en el problema de clasificar a las matrices infinitas con entradas en $K$ que satisfacen ciertas condiciones. Con este método se logra clasificar a los productos tensoriales graduados de $K[x]$ con $K[y]$ en un ejemplo particular y tres casos principales: álgebras cuadráticas, clasificadas por Conner y Goetz por métodos diferentes, una familia llamada $A(n,d,a)$ con la propiedad de $n+1-$ extensión para cualquier $n \geq 2$ y un tercer caso no completamente clasificado, para el cual se describen los cálculos iniciales que ilustran cómo se puede alcanzar la clasificación de las posibles aplicaciones de torcimiento con una cantidad creciente de cálculo computacional. Además, en este tercer caso, se obtiene una familia de productos tensoriales torcidos graduados $B(a,L)$ parametrizada por una familia de sucesiones casi-balanceadas. Los miembros de la familia $B(a,L)$ no tienen la propiedad de $m$- extensión, para ningún $m$.

**Palabras clave:** producto tensorial torcido graduado, aplicación de torcimiento, álgebra cuadrática, familia A(n,d,a), sucesiones cuasi-balanceadas.

# Abstract

In this thesis, an almost complete classification of all graduated twisted tensorial products of K[x] with K[y] is obtained. For this purpose, a representation of a graduated twisted tensor product of $K[x]$ with $K[y]$ in the algebra $L(K^{N_0})$, which is immersed in the set of infinite matrices with entries in $K$, it is used. Thus the problem of classifying the graduated twisted tensor products of $K[x]$ with $K[y]$ results in the problem of classifying infinite matrices with inputs in K that satisfy certain conditions. With this method it is possible to classify the graduated tensor products of $K[x]$ with $K[y]$ in a particular example and three main cases: quadratic algebras, classified by Conner and Goetz by different methods, a family called $A(n,d,a)$ with the property of $n+1$ - extension for $n \geq 2$, and a third case not fully classified for which there are shown initial calculations illustrating how classification of possible twisting applications with an increasing amount of computational calculation can be achieved. Furthermore, in this third case, a family of products graduated twisted tensor $B(a,L)$ parametrized by a family of quasi-balanced sequences is obtained. Members of $B(a,L)$ family do not have the $m$-extension property, for no $m$.

*A mi esposa Ada*
*A mi hija Diana*

# Agradecimientos

# Índice general

# Introducción

En [3] los autores introdujeron la noción de producto tensorial torcido de $K$-álgebras unitarias, donde $K$ es un anillo unitario.

Asumiremos que $K$ es un cuerpo, y consideraremos el problema básico de clasificar a todos los productos tensoriales torcidos de A con B, para un par de álgebras A y B dado.

En general este problema está lejos de lograr resolverse, aunque se han obtenido algunos resultados, principalmente para álgebras de dimensión finita (ver [1], [2], [4], [7],[9] y [10]).

En particular, en [7] fueron encontradas algunas familias de productos tensoriales torcidos de $K[x]$ con $K[y]$. La clasificación completa de estos productos tensoriales parece estar lejos de poder alcanzarse, pero en [6](ver también [5]) los productos tensoriales torcidos graduados de $K[x]$ con $K[y]$ que son álgebras cuadráticas, fueron clasificados completamente.

Por otro lado el producto tensorial torcido de $K[x]/\langle x^n \rangle$ con un álgebra A puede ser representada en $M_n(A)$ (ver [8], Theorem 1.10). Esta representación puede ser generalizada a álgebras de dimensión finita (ver [1]).

En este trabajo empezamos con una representación de un producto tensorial torcido graduado de $K[x]$ con $K[y]$ en $L(K[x]^{\mathbb{N}_0})$, la cual es muy similar a la representación en [8, Theorem 1.10]. En el caso graduado esta representación puede ser simplificada aún más a una representación del producto tensorial torcido graduado en el álgebra $L(K^{\mathbb{N}_0})$, inmersa en el álgebra de matrices infinitas con entradas en $K$. De esta manera logramos traducir el problema de clasificar a todos los productos tensoriales torcidos graduados de $K[x]$ con $K[y]$ al problema de clasificar a las matrices infinitas con entradas en $K$ que satisfacen ciertas condiciones (ver Proposición 1.4.1). Con este método demostramos que todos los productos tensoriales torcidos graduados de $K[x]$ con $K[y]$ pueden ser clasificados en tres casos principales, excepto un ejemplo particular (ver Proposición 3.1.4). El primer caso es el caso de álgebras cuadráticas, ya clasificadas en [6], el segundo caso produce una familia llamada $A(n,d,a)$ y en el tercer caso solamente tenemos algunos resultados de clasificación parcial y obtenemos una familia llamada $B(a,L)$.

Uno puede describir un producto tensorial torcido graduado de $K[x]$ con $K[y]$ determinando cómo conmuta $y^k$ con $x$, lo cual significa determinar los coeficientes $a_i's$ en

$$y^k x = \sum_{i=0}^{k+1} a_i x^{k+1-i} y^i.$$

Por ejemplo la relación de conmutación $yx = qxy$ nos da el plano cuántico, y en la subsección 2.1.1 exploramos el caso $yx = bxy + cy^2$.

Decimos que el producto tensorial torcido graduado $A = K[x] \otimes K[y]$ es $m$−extensible, o que tiene la propiedad de $m$−extensión, si las relaciones $y^k x = \sum_{i=0}^{k+1} a_i x^{k+1-i} y^i$, $k = 1, 2, \cdots, m-1$ determinan completamente la multiplicación en $A$.
En particular $A$ es 2-extensible o cuadrática si la relación

$$yx = ax^2 + bxy + cy^2$$

determina completamente la multiplicación en $A$.
En general la relación $yx = ax^2 + bxy + cy^2$ produce un álgebra cuadrática, siempre que $ac \neq 1$ y que $(b, ac)$ no es una raíz de ningún miembro de una cierta familia de polinomios $Q_n(b, ac)$. En este caso nuestros resultados coinciden con los resultados de [6], los cuales son obtenidos con métodos muy diferentes.

Si $a \neq 0$, entonces un tal producto tensorial es equivalente a uno con $a = 1$, y de este modo nos enfocaremos en el caso $yx = x^2 + bxy + cy^2$. En el caso $c = 1$ uno puede demostrar que necesariamente $b = -1$ (ver Lema 2.1.5) y las álgebras resultantes no son cuadráticas, es decir, no tienen la propiedad de 2-extensión. Obtenemos un álgebra particular con $y^k x = x^{k+1} - x^k y + y^{k+1}$ para todo $k \in \mathbb{N}$ (ver Proposición 3.1.4). Esta álgebra no tiene la propiedad de $m$-extensión para ningún $m$.

Para todo producto tensorial torcido graduado con $yx = x^2 - xy + y^2$, el cual no es el caso particular mencionado antes, existe $n \geq 2$ tal que

$$y^k x = x^{k+1} - x^k y + y^{k+1}, \quad \text{para todo } k < n, \quad \text{y} \quad y^n x \neq x^{n+1} - x^n y + y^{n+1}.$$

Un resultado central es la Proposición 3.1.6, la cual demuestra que tenemos exactamente dos posibilidades para $y^n x$. En el primer caso

$$y^n x = dx^{n+1} - dx^n y - axy^n + (a+1)y^{n+1},$$

para $a, d$ en $K$ satisfaciendo ciertas condiciones, a saber, $(a, d)$ no es raíz de ningún elemento de una cierta familia de polinomios $R_j(a, d)$. Esto produce una familia $A(n, d, a)$ de productos tensoriales torcidos con la propiedad de $(n+1)$-extensión (ver secciones 3.1 y 3.2), lo cual significa que la multiplicación es determinada por las relaciones de conmutación hasta el grado $n + 1$.

La segunda posibilidad se analiza en el Capítulo 4 y la relación de conmutación en grado $n+1$ es
$$y^n x = d x^{n+1} - x^n y + (a+1) y^{n+1}$$
donde $d(a+1) = 1$. Describiremos todas las posibles relaciones de conmutación hasta $y^{3n+2} x$. Aunque en este caso no se alcanza la clasificación completa nuestros cálculos ilustran cómo se puede alcanzar la clasificación de todas las aplicaciones de torcimiento posibles hasta cualquier grado con una cantidad creciente de cálculo computacional.

Además, logramos encontrar una familia de productos tensoriales torcidos a la cual llamamos $B(a,L)$, parametrizada por $a \in K \setminus \{0,-1\}$ y $L \in \mathscr{L}$, donde $\mathscr{L}$ es el conjunto de sucesiones casi-balanceadas (ver Definición 4.3.1). Estas sucesiones son interesantes de por sí, por ejemplo ellas muestran una conexión sorprendente a la función de Euler $\varphi$. Toda sucesión casi-balanceada truncada puede ser continuada de varias maneras, lo cual implica que todos los miembros de la familia $B(a,L)$ no tienen la propiedad de $m$- extensión para ningún $m$.

La siguiente tabla contiene todas las posibles aplicaciones de torcimiento graduadas de $K[x]$ con $K[y]$. Las únicas aplicaciones de torcimiento que no han sido totalmente clasificadas están en la última fila. Para las familias con $yx = x^2 - xy + y^2$ formulamos las relaciones de conmutación con respecto a las matrices infinitas $Y = \psi(y)$, $M = \psi(x)$ y $\widetilde{M} = \psi(x-y)$, obtenidas de la representación fiel (ver Observación 1.4.2)
$$\psi : K[x] \otimes_\sigma K[y] \longrightarrow L(K^{\mathbb{N}_0}).$$
Por ejemplo la relación $yx = x^2 - xy + y^2$ corresponde a $Y\widetilde{M} = M\widetilde{M}$ y $Y^k \widetilde{M} = M^k \widetilde{M}$ representa $y^k x = x^{k+1} - x^k y + y^{k+1}$.

| Relaciones de conmutación | Clasificación y Parámetros | Referencia | Propiedad de m-extensión |
|---|---|---|---|
| $yx = bxy + cy^2$ | $b, c \in K$ | Subsección 2.1.1 | Cuadrática |
| $yx = x^2 + bxy + cy^2$ | $b, c \in K$ $Q_k(b,c) \neq 0, \forall k \in \mathbb{N}$ | Teorema 2.2.7 | Cuadrática |
| $Y^k \widetilde{M} = M^k \widetilde{M}, \forall k \geq 1$ | Ejemplo particular | Proposición 3.1.4 | No es $m$-extensible para ningún $m$ |
| $Y^k \widetilde{M} = M^k \widetilde{M}, \forall k < n$ $Y^n \widetilde{M} = dM^n \widetilde{M} - a\widetilde{M} Y^n$ | Familia $A(n,d,a)$ $a, d \in K, n \in \mathbb{N}$ $R_k(a,d) \neq 0 \, \forall k \in \mathbb{N}$ | Teorema 3.2.3 y Corolario 3.2.7 | Propiedad de n+1-exrensión |
| $Y^k \widetilde{M} = M^k \widetilde{M}, \forall k < n$ $Y^k \widetilde{M} = d^r M^k \widetilde{M},$ si $L_r < k < L_{r+1}$ $Y^k \widetilde{M} = d^r M^{k+1} - d^{r-1} M^k Y$ $+ a Y^{k+1}$, si $k = L_r$ | $a \in K \setminus \{-1, 0\}$ $L \in \mathscr{L}$ Subfamilia del caso que sigue. | Proposición 4.3.5 y Proposición 4.4.8 | No es $m$-extensible para ningún $m$. |
| $Y^k \widetilde{M} = M^k \widetilde{M}, \forall k < n$ $Y^n \widetilde{M} = dM^{n+1} - M^n Y$ $+ a Y^{n+1}$ | No totalmente clasificado | Proposición 4.2.2 y Proposición 4.2.4 | Conjetura: No es $m$-extensible. |

# Capítulo 1

# Preliminares

## 1.1. Productos tensoriales torcidos

**Definición 1.1.1.** *Sea K un anillo conmutativo y sean A, B, K-álgebras unitarias. Un producto tensorial torcido de A con B es una estructura de álgebra en el K-módulo $A \otimes B$ tal que*
*(1) las aplicaciones canónicas*

$$i_A : A \hookrightarrow A \otimes B, \quad i_B : B \hookrightarrow A \otimes B$$
$$a \mapsto a \otimes 1_B, \quad b \mapsto 1_A \otimes b$$

*son homomorfismos de álgebras y*
*(2)* $\quad\quad\quad\quad\quad\quad\quad \mu \circ (i_A \otimes i_B) = Id_{A \otimes B}$
*donde $\mu$ denota la aplicación multiplicación del producto tensorial torcido.*

**Observación 1.1.2.** La condición (2) es equivalente a la condición

$(2')\quad\quad (a \otimes 1_B) \cdot (1_A \otimes b) = a \otimes b, \quad \forall a \in A, b \in B.$

En efecto, el diagrama conmutativo para la condición (2) es el siguiente:

$$\begin{array}{ccc} A \otimes B & \xrightarrow{i_A \otimes i_B} & A \otimes B \otimes A \otimes B \\ & \searrow{\scriptstyle Id_{A \otimes B}} & \downarrow{\scriptstyle \mu} \\ & & A \otimes B \end{array}$$

Sean $a \in A, b \in B$ arbitrarios,

$$\mu \circ (i_A \otimes i_B) = Id_{A \otimes B} \Longleftrightarrow (\mu \circ (i_A \otimes i_B))(a \otimes b) = a \otimes b$$

$$\Longleftrightarrow \mu(a \otimes 1_B \otimes 1_A \otimes b) = a \otimes b \Longleftrightarrow (a \otimes 1_B) \cdot (1_A \otimes b) = a \otimes b.$$

Para caracterizar a los productos tensoriales torcidos partimos de la definición usual de la multiplicación en el producto tensorial $A \otimes B$:

$$(a_1 \otimes b_1) \cdot (a_2 \otimes b_2) = (a_1 a_2) \otimes (b_1 b_2).$$

Considerando el morfismo de multiplicación en un álgebra $C$ como una aplicación $\mu_C : C \otimes C \longrightarrow C$, la definición anterior nos da:

$$\mu_{A \otimes B} = (\mu_A \otimes \mu_B) \circ (Id_A \otimes \tau \otimes Id_B)$$

donde $\tau : B \otimes A \longrightarrow A \otimes B$ es el flip, que simplemente intercambia los factores. Es evidente que esta multiplicación da una estructura de producto tensorial torcido en $A \otimes B$.

Generalizando este caso analizaremos las aplicaciones $K$-lineales $\sigma : B \otimes A \longrightarrow A \otimes B$ que definen una estructura de álgebra asociativa en $A \otimes B$ a través de la aplicación multiplicación $\mu := (\mu_A \otimes \mu_B) \circ (Id_A \otimes \sigma \otimes Id_B)$.

**Proposición 1.1.3.** *Sea $\sigma : B \otimes A \longrightarrow A \otimes B$ una aplicación $K$-lineal. Si $\sigma$ define una estructura de producto tensorial torcido en $A \otimes B$ vía la multiplicación*

$$\mu := (\mu_A \otimes \mu_B) \circ (Id_A \otimes \sigma \otimes Id_B),$$

*entonces $\sigma$ satisface*

  *a)* $\sigma(1_B \otimes a) = a \otimes 1_B,$

  *b)* $\sigma(b \otimes 1_A) = 1_A \otimes b,$

  *c)* $\sigma \circ (Id_B \otimes \mu_A) = (\mu_A \otimes Id_B) \circ (Id_A \otimes \sigma) \circ (\sigma \otimes Id_A),$

  *d)* $\sigma \circ (\mu_B \otimes Id_A) = (Id_A \otimes \mu_B) \circ (\sigma \otimes Id_B) \circ (Id_B \otimes \sigma),$

  *para todo $a \in A$, todo $b \in B$ y donde $\mu_A$, $\mu_B$ son las multiplicaciones en $A$, $B$ respectivamente.*

**Demostración.** Sea $\sigma : B \otimes A \longrightarrow A \otimes B$ una aplicación $K$-lineal tal que la aplicación $\mu := (\mu_A \otimes \mu_B) \circ (Id_A \otimes \sigma \otimes Id_B)$ es una multiplicación que determina una estructura de producto tensorial torcido en $A \otimes B$.

**Demostración de a)**

$$\begin{aligned}
a \otimes 1_B &= i_A(1_A a) = i_A(1_A) \cdot i_A(a) = \mu(i_A(1_A) \otimes i_A(a)) \\
&= (\mu_A \otimes \mu_B) \circ (Id_A \otimes \sigma \otimes Id_B)(i_A(1_A) \otimes i_A(a)) \\
&= (\mu_A \otimes \mu_B) \circ (Id_A \otimes \sigma \otimes Id_B)(1_A \otimes 1_B \otimes a \otimes 1_B) \\
&= (\mu_A \otimes \mu_B)(1_A \otimes \sigma(1_B \otimes a) \otimes 1_B) \\
&= (\mu_A \otimes \mu_B)(1_A \otimes (\Sigma_i a_i \otimes b_i) \otimes 1_B) = (\mu_A \otimes \mu_B)(\Sigma_i 1_A \otimes a_i \otimes b_i \otimes 1_B) \\
&= \Sigma_i a_i \otimes b_i = \sigma(1_B \otimes a).
\end{aligned}$$

**Demostración de b)**
La demostración es similar a la de a).
**Demostración de c)**
Como $\mu$ es asociativo se cumple

$$i_B(b) \cdot (i_A(a) \cdot i_A(a')) = (i_B(b) \cdot i_A(a)) \cdot i_A(a').$$

Evaluando el primer miembro

$$\begin{aligned}
i_B(b) \cdot (i_A(a) \cdot i_A(a')) &= (1_A \otimes b) \cdot i_A(aa') \\
&= (1_A \otimes b) \cdot (aa' \otimes 1_B) \\
&= (\mu_A \otimes \mu_B) \circ (Id_A \otimes \sigma \otimes Id_B)(1_A \otimes b \otimes aa' \otimes 1_B) \\
&= (\mu_A \otimes \mu_B)(1_A \otimes \sigma(b \otimes aa') \otimes 1_B) = \sigma(b \otimes aa') \\
&= \sigma \circ (Id_B \otimes \mu_A)(b \otimes a \otimes a').
\end{aligned}$$

Desarrollando el segundo miembro

$$\begin{aligned}
(i_B(b) \cdot i_A(a)) \cdot i_A(a') &= [(1_A \otimes b) \cdot (a \otimes 1_B)] \cdot (a' \otimes 1_B) \\
&= \sigma(b \otimes a) \cdot (a' \otimes 1_B) = \mu[\sigma(b \otimes a) \otimes (a' \otimes 1_B)] \\
&= [(\mu_A \otimes \mu_B) \circ (Id_A \otimes \sigma \otimes Id_B)](\sigma(b \otimes a) \otimes a' \otimes 1_B) \\
&= [(\mu_A \otimes \mu_B) \circ (Id_A \otimes \sigma \otimes Id_B) \circ (\sigma \otimes Id_A \otimes Id_B)](x \otimes 1_B),
\end{aligned}$$

donde $x = b \otimes a \otimes a'$.
Pero
$$(Id_A \otimes \sigma \otimes Id_B) \circ (\sigma \otimes Id_A \otimes Id_B) = \{(Id_A \otimes \sigma) \circ (\sigma \otimes Id_A)\} \otimes Id_B.$$

Luego

$$\begin{aligned}
(i_B(b) \cdot i_A(a)) \cdot i_A(a') &= [(\mu_A \otimes \mu_B) \circ \{(Id_A \otimes \sigma) \circ (\sigma \otimes Id_A)\} \otimes Id_B](x \otimes 1_B) \\
&= (\mu_A \otimes \mu_B)(\{(Id_A \otimes \sigma) \circ (\sigma \otimes Id_A)\}(x) \otimes 1_B) \\
&= [(\mu_A \otimes Id_B) \circ (Id_A \otimes \sigma) \circ (\sigma \otimes Id_A)](x).
\end{aligned}$$

Igualando ambos miembros se obtiene:
$\sigma \circ (Id_B \otimes \mu_A) = (\mu_A \otimes Id_B) \circ (Id_A \otimes \sigma) \circ (\sigma \otimes Id_A)$.
**Demostración de d)**
Similarmente, usando el hecho de que

$$i_B(b) \cdot (i_B(b') \cdot i_A(a)) = (i_B(b) \cdot i_B(b')) \cdot i_A(a),$$

se demuestra d). □

**Observación 1.1.4.** Este tipo de demostraciones engorrosas se pueden evitar utilizando el conocido cálculo gráfico empleado en ciertas categorias. Consiste en usar diagramas donde los morfismos paralelos a la misma altura corresponden a los factores en algún morfismo producto tensorial y las composiciones se grafican bajando nivel a nivel.
La aplicación identidad de un álgebra será representada por una linea vertical.
Dada un álgebra $C$, el siguiente diagrama representa a la multiplicación $\mu_C$.
La aplicación $\sigma$ será representada por el diagrama

Por ejemplo la igualdad c) de la Proposición 1.1.3

$$\sigma \circ (Id_B \otimes \mu_A) = (\mu_A \otimes Id_B) \circ (Id_A \otimes \sigma) \circ (\sigma \otimes Id_A)$$

se vería así:

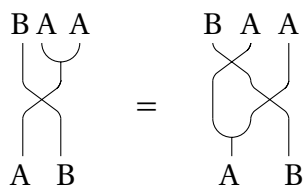

La igualdad d)

$$\sigma \circ (\mu_B \otimes Id_A) = (Id_A \otimes \mu_B) \circ (\sigma \otimes Id_B) \circ (Id_B \otimes \sigma)$$

se vería así:

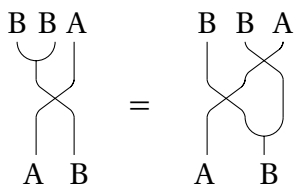

La asociatividad de $\mu$ implica las igualdades

$$i_B(b) \cdot (i_A(a) \cdot i_A(a')) = (i_B(b) \cdot i_A(a)) \cdot i_A(a'),$$

$$i_B(b) \cdot ((i_B(b') \cdot i_A(a)) = (i_B(b) \cdot i_B(b')) \cdot i_A(a),$$

que se han usado en la demostración de c) y d).
El recíproco de esta proposición también se cumple.

**Proposición 1.1.5.** *Sea $\sigma : B \otimes A \longrightarrow A \otimes B$ una aplicación K-lineal que satisface a), b), c) y d) de la Proposición 1.1.3. Entonces $\sigma$ define una estructura de producto tensorial torcido vía la multiplicación*

$$\mu_\sigma := (\mu_A \otimes \mu_B) \circ (Id_A \otimes \sigma \otimes Id_B).$$

*Esta álgebra será denotada por $A \otimes_\sigma B$ y a la aplicación $\sigma$ se le llama aplicación de torcimiento o "twist".*

**Demostración.** Para demostrar esta proposición usaremos el cálculo gráfico.
En primer lugar la relación $\mu_\sigma = (\mu_A \otimes \mu_B) \circ (Id_A \otimes \sigma \otimes Id_B)$ se representa así:

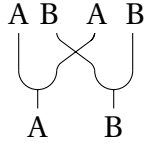

La asociatividad de $\mu_\sigma$ significa

$$\mu_\sigma \circ (\mu_\sigma \otimes Id_{A \otimes B}) = \mu_\sigma \circ (Id_{A \otimes B} \otimes \mu_\sigma).$$

Demostraremos gráficamente que ambos lados son iguales a

$\Lambda := (\mu_A \otimes \mu_B) \circ (Id_A \otimes \mu_A \otimes \mu_B \otimes Id_B) \circ (Id_A \otimes Id_A \otimes \sigma \otimes Id_B \otimes Id_B) \circ (Id_A \otimes \sigma \otimes \sigma \otimes Id_B).$

En efecto, el lado izquierdo $\mu_\sigma \circ (\mu_\sigma \otimes Id_{A \otimes B})$ está representado por el diagrama

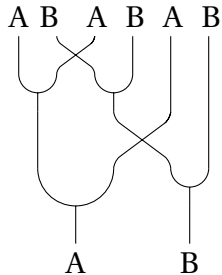

Ahora la expresión

$\Lambda := (\mu_A \otimes \mu_B) \circ (Id_A \otimes \mu_A \otimes \mu_B \otimes Id_B) \circ (Id_A \otimes Id_A \otimes \sigma \otimes Id_B \otimes Id_B) \circ (Id_A \otimes \sigma \otimes \sigma \otimes Id_B)$

se representa gráficamente así:

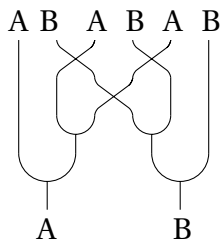

Vemos que
$$\mu_\sigma \circ (\mu_\sigma \otimes Id_{A\otimes B}) = \Lambda.$$

Por otro lado el lado derecho $\mu_\sigma \circ (Id_{A\otimes B} \otimes \mu_\sigma)$ está representado por el diagrama

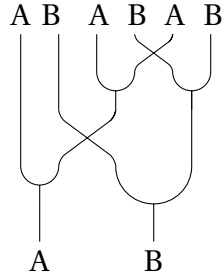

Vemos también que $\mu_\sigma \circ (Id_{A\otimes B} \otimes \mu_\sigma) = \Lambda$.

Así hemos demostrado que $\mu_\sigma \circ (\mu_\sigma \otimes Id_{A\otimes B}) = \mu_\sigma \circ (Id_{A\otimes B} \otimes \mu_\sigma)$.

Demostraremos ahora que $1_A \otimes 1_B$ es el elemento unidad de $A \otimes_\sigma B$.

En efecto:

$$\begin{aligned}
(1_A \otimes 1_B) \cdot (a \otimes b) &= \mu_\sigma[(1_A \otimes 1_B) \otimes (a \otimes b)] \\
&= [(\mu_A \otimes \mu_B) \circ (Id_A \otimes \sigma \otimes Id_B)](1_A \otimes 1_B \otimes a \otimes b) \\
&= (\mu_A \otimes \mu_B)(1_A \otimes \sigma(1_B \otimes a) \otimes b) \\
&= (\mu_A \otimes \mu_B)(1_A \otimes a \otimes 1_B \otimes b) \\
&= \mu_A(1_A \otimes a) \otimes \mu_B(1_B \otimes b) \\
&= (1_A a) \otimes (1_B b) = a \otimes b.
\end{aligned}$$

Análogamente se demuestra que

$$(a \otimes b) \cdot (1_A \otimes 1_B) = a \otimes b.$$

Ahora demostraremos que $i_A, i_B$ son homomorfismos de álgebras.

Para demostrar que $i_A$ es un homomorfismo de álgebras es suficiente demostrar que $i_A(aa') = i_A(a) \cdot i_A(a')$.

En efecto:

$$\begin{aligned}
i_A(a) \cdot i_A(a') &= (a \otimes 1_B) \cdot (a' \otimes 1_B) = \mu_\sigma(a \otimes 1_B \otimes a' \otimes 1_B) \\
&= (\mu_A \otimes \mu_B) \circ (Id_A \otimes \sigma \otimes Id_B)(a \otimes 1_B \otimes a' \otimes 1_B) \\
&= (\mu_A \otimes \mu_B)(a \otimes \sigma(1_B \otimes a') \otimes 1_B) \\
&= \mu_A \otimes \mu_B)(a \otimes a' \otimes 1_B \otimes 1_B) \\
&= \mu_A(a \otimes a') \otimes \mu_B(1_B \otimes 1_B) = aa' \otimes 1_B = i_A(aa').
\end{aligned}$$

Análogamente se demuestra que $i_B$ es un homomorfismo de álgebras .

Finalmente demostraremos que $\mu_\sigma$ cumple $\mu_\sigma \circ (i_A \otimes i_B) = Id_{A \otimes B}$.

En efecto:

$$\begin{aligned}
\mu_\sigma \circ (i_A \otimes i_B)(a \otimes b) &= (\mu_A \otimes \mu_B) \circ (Id_A \otimes \sigma \otimes Id_B) \circ (i_A \otimes i_B)(a \otimes b) \\
&= (\mu_A \otimes \mu_B) \circ (Id_A \otimes \sigma \otimes Id_B)(a \otimes 1_B \otimes 1_A \otimes b) \\
&= (\mu_A \otimes \mu_B)(a \otimes \sigma(1_B \otimes 1_A) \otimes b) \\
&= (\mu_A \otimes \mu_B)(a \otimes 1_A \otimes 1_B \otimes b) \\
&= a \otimes b = Id_{A \otimes B}(a \otimes b). \qquad \square
\end{aligned}$$

**Observación 1.1.6.** Si $\sigma : B \otimes A \longrightarrow A \otimes B$ es una aplicación $K$-lineal que define una estructura de producto tensorial torcido en $A \otimes B$ vía la multiplicación $\mu_\sigma := (\mu_A \otimes \mu_B) \circ (Id_A \otimes \sigma \otimes Id_B)$ entonces $\sigma = \mu_\sigma \circ (i_B \otimes i_A)$.

En efecto:

$$\begin{aligned}
\mu_\sigma \circ (i_B \otimes i_A)(b \otimes a) &= \mu_\sigma(1_A \otimes b \otimes a \otimes 1_B) \\
&= (\mu_A \otimes \mu_B) \circ (Id_A \otimes \sigma \otimes Id_B)(1_A \otimes b \otimes a \otimes 1_B) \\
&= (\mu_A \otimes \mu_B)(1_A \otimes \sigma(b \otimes a) \otimes 1_B) = \sigma(b \otimes a).
\end{aligned}$$

El producto tensorial torcido $A \otimes_\sigma B$ tiene la siguiente propiedad universal.

**Proposición 1.1.7** (Propiedad universal del producto tensorial torcido)**.** *Dados dos morfismos de álgebras $f : A \longrightarrow C$, $g : B \longrightarrow C$ tales que*

$$\mu_C \circ (g \otimes f) = \mu_C \circ (f \otimes g) \circ \sigma$$

*entonces existe un morfismo de álgebras $h : A \otimes_\sigma B \longrightarrow C$ que satisface $f = h \circ i_A$, $g = h \circ i_B$.*

**Demostración.** En efecto, $h = \mu_C \circ (f \otimes g)$ cumple: $f = h \circ i_A$, $g = h \circ i_B$.

$$\begin{aligned}
(h \circ i_A)(a) &= h(i_A(a)) = h(a \otimes 1_B) \\
&= [\mu_C \circ (f \otimes g)](a \otimes 1_B) \\
&= \mu_C[f(a) \otimes g(1_B)] = \mu_C(f(a) \otimes 1_C) \\
&= f(a) 1_C = f(a)
\end{aligned}$$

$$\begin{aligned}
(h \circ i_B)(b) &= h(i_B(b)) = h(1_A \otimes b) \\
&= [\mu_C \circ (f \otimes g)](1_A \otimes b) \\
&= \mu_C[f(1_A) \otimes g(b)] = \mu_C(1_C \otimes g(b)) \\
&= 1_C g(b) = g(b).
\end{aligned}$$

Ahora demostraremos que $h$ es un morfismo de álgebras.

Es evidente que $h$ es $K$-lineal. Resta demostrar que:

$$h((a_1 \otimes b_1).(a_2 \otimes b_2)) = h(a_1 \otimes b_1)h(a_2 \otimes b_2),$$

donde el producto de la izquierda es el producto torcido $\mu_\sigma$ dado por

$$\mu_\sigma := (\mu_A \otimes \mu_B) \circ (Id_A \otimes \sigma \otimes Id_B),$$

y el producto de la derecha corresponde a la multiplicación $\mu_C$.

En efecto:

$$
\begin{aligned}
h((a_1 \otimes b_1).(a_2 \otimes b_2)) &= h(\mu_\sigma((a_1 \otimes b_1) \otimes (a_2 \otimes b_2))) \\
&= h((\mu_A \otimes \mu_B) \circ (Id_A \otimes \sigma \otimes Id_B)(a_1 \otimes b_1 \otimes a_2 \otimes b_2)) \\
&= h((\mu_A \otimes \mu_B)(a_1 \otimes \sigma(b_1 \otimes a_2) \otimes b_2)) \\
&= h((\mu_A \otimes \mu_B)(a_1 \otimes (\Sigma_i a_{2i} \otimes b_{1i}) \otimes b_2)) \\
&= h(\Sigma_i \mu_A(a_1 \otimes a_{2i}) \otimes \mu_B(b_{1i} \otimes b_2)) = h(\Sigma_i a_1 a_{2i} \otimes b_{1i} b_2) \\
&= \mu_C \circ (f \otimes g)(\Sigma_i a_1 a_{2i} \otimes b_{1i} b_2) \\
&= \Sigma_i \mu_C(f(a_1)f(a_{2i}) \otimes g(b_{1i})g(b_2)) \\
&= \Sigma_i f(a_1)f(a_{2i})g(b_{1i})g(b_2) \\
&= f(a_1)(\Sigma_i f(a_{2i})g(b_{1i}))g(b_2) \\
&= f(a_1)(\Sigma_i \mu_C(f(a_{2i}) \otimes g(b_{1i})))g(b_2) \\
&= f(a_1)(\Sigma_i \mu_C \circ (f \otimes g)(a_{2i} \otimes b_{1i}))g(b_2) \\
&= f(a_1)(\mu_C \circ (f \otimes g)(\Sigma_i a_{2i} \otimes b_{1i}))g(b_2) \\
&= f(a_1)((\mu_C \circ (f \otimes g) \circ \sigma)(b_1 \otimes a_2))g(b_2) \\
&= f(a_1)((\mu_C \circ (g \otimes f))(b_1 \otimes a_2))g(b_2) \\
&= f(a_1)g(b_1)f(a_2)g(b_2) \\
&= (\mu_C \circ (f \otimes g))(a_1 \otimes b_1)(\mu_C \circ (f \otimes g))(a_2 \otimes b_2) \\
&= h(a_1 \otimes b_1)h(a_2 \otimes b_2).
\end{aligned}
$$

$\square$

Las aplicaciones de torcimiento son objetos de una categoría.

**Definición 1.1.8.** *Sean $\sigma : B \otimes A \longrightarrow A \otimes B$, $\tau : D \otimes C \longrightarrow C \otimes D$ dos aplicaciones de torcimiento. Un morfismo $(f, g) : \sigma \longrightarrow \tau$ es un par de morfismos de álgebras $f : A \longrightarrow C$, $g : B \longrightarrow D$ tal que $\tau \circ (g \otimes f) = (f \otimes g) \circ \sigma$. Es decir el siguiente diagrama es conmutativo.*

$$
\begin{array}{ccc}
B \otimes A & \xrightarrow{\sigma} & A \otimes B \\
{\scriptstyle g \otimes f} \downarrow & & \downarrow {\scriptstyle f \otimes g} \\
D \otimes C & \xrightarrow{\tau} & C \otimes D
\end{array}
$$

**Definición 1.1.9.** *Dos aplicaciones de torcimiento $\sigma$, $\tau : B \otimes A \longrightarrow A \otimes B$ se dice que son equivalentes si son isomorfas. Es decir, si existen automorfismos $f : A \longrightarrow A$ y $g : B \longrightarrow B$ tales que*

$$\tau = (f^{-1} \otimes g^{-1}) \circ \sigma \circ (g \otimes f).$$

*En otras palabras, si el siguiente diagrama es conmutativo.*

$$\begin{array}{ccc} B \otimes A & \xrightarrow{\sigma} & A \otimes B \\ {\scriptstyle g \otimes f} \uparrow & & \downarrow {\scriptstyle f^{-1} \otimes g^{-1}} \\ B \otimes A & \xrightarrow{\tau} & A \otimes B \end{array}$$

## 1.2. Productos tensoriales torcidos de extensiones polinomiales

Ahora veamos el caso en que el álgebra $B$ es $K[y]$, el álgebra de polinomios sobre $y$. Al producto tensorial torcido $A \otimes K[y]$ se le llama **extensión polinomial** de $A$. En este caso todo morfismo lineal $\sigma : K[y] \otimes A \longrightarrow A \otimes K[y]$ está determinado por los morfismos $\gamma_j^r : A \longrightarrow A$ definidos por la ecuación

$$\sigma(y^r \otimes a) = \sum_{j \geq 0} \gamma_j^r(a) \otimes y^j. \tag{1.2.1}$$

Las condiciones para que $\sigma$ sea una aplicación de torcimiento están dadas por la siguiente proposición.

**Proposición 1.2.1.** *Un conjunto de aplicaciones $\gamma_j^r : A \longrightarrow A$ define una aplicación de torcimiento $\sigma$ vía* (1.2.1) *si y solamente si se cumplen las siguientes propiedades:*

1) $\gamma_j^0 = \delta_{0,j} Id_A$

2) $\gamma_j^r(1) = \delta_{r,j}$

3) *Para todo $r, j$ tenemos*

$$\gamma_j^r(ab) = \sum_{k \geq 0} \gamma_k^r(a) \gamma_j^k(b)$$

   *para todo $a, b \in A$*

4) *Para todo $r, j$, y todo $i < r$*

$$\gamma_j^r = \sum_{l=0}^{j} \gamma_l^i \circ \gamma_{j-l}^{r-i}$$

**Demostración.** Por la Proposición 1.1.3 y la Proposición 1.1.5 es suficiente demostrar que las condiciones 1), 2), 3), 4) corresponden a las siguientes condiciones sobre $\sigma$:

$a)$ $\sigma(1_{K[y]} \otimes a) = a \otimes 1_{K[y]}$
$b)$ $\sigma(y^r \otimes 1_A) = 1_A \otimes y^r$
$c)$ $\sigma \circ (Id_{K[y]} \otimes \mu_A) = (\mu_A \otimes Id_{K[y]}) \circ (Id_A \otimes \sigma) \circ (\sigma \otimes Id_A)$
$d)$ $\sigma \circ (\mu_{K[y]} \otimes Id_A) = (Id_A \otimes \mu_{K[y]}) \circ (\sigma \otimes Id_{K[y]}) \circ (Id_{K[y]} \otimes \sigma)$
para todo $a \in A$ y todo $r \in \mathbb{N}$.

**Demostración de 1)**

$$\sigma(1 \otimes a) = a \otimes 1 \Leftrightarrow \sum_{j \geq 0} \gamma_j^0(a) \otimes 1 = a \otimes 1 \Leftrightarrow \gamma_j^0(a) = \delta_{0,j} Id_A(a).$$

**Demostración de 2)**

$$\sigma(y^r \otimes 1) = 1 \otimes y^r \Leftrightarrow \sum_{j \geq 0} \gamma_j^r(1) \otimes y^j = 1 \otimes y^r \Leftrightarrow \gamma_j^r(1) = \delta_{r,j}.$$

**Demostración de 3)**
Por definición de $\gamma_j^r$:
$$\sigma(y^r \otimes ab) = \sum_{j \geq 0} \gamma_j^r(ab) \otimes y^j. \qquad (1.2.2)$$

Por otro lado, usando la condición c) sobre $\sigma$ se tiene

$$\begin{aligned}
\sigma(y^r \otimes ab) &= [\sigma \circ (Id_{K[y]} \otimes \mu_A)](y^r \otimes a \otimes b) \\
&= [(\mu_A \otimes Id_{K[y]}) \circ (Id_A \otimes \sigma) \circ (\sigma \otimes Id_A)](y^r \otimes a \otimes b) \\
&= [(\mu_A \otimes Id_{K[y]}) \circ (Id_A \otimes \sigma)](\sigma(y^r \otimes a) \otimes b) \\
&= [(\mu_A \otimes Id_{K[y]}) \circ (Id_A \otimes \sigma)]\left(\sum_{k \geq 0} \gamma_k^r(a) \otimes y^k \otimes b\right) \\
&= (\mu_A \otimes Id_{K[y]})\left(\sum_{k \geq 0} \gamma_k^r(a) \otimes \sigma(y^k \otimes b)\right) \\
&= (\mu_A \otimes Id_{K[y]})\left(\sum_{k \geq 0} \gamma_k^r(a) \otimes \sum_{j \geq 0} \gamma_j^k(b) \otimes y^j\right) \\
&= \sum_{k \geq 0} \sum_{j \geq 0} \mu_A(\gamma_k^r(a) \otimes \gamma_j^k(b)) \otimes y^j \\
&= \sum_{k \geq 0} \sum_{j \geq 0} (\gamma_k^r(a) \gamma_j^k(b)) \otimes y^j.
\end{aligned}$$

Es decir
$$\sigma(y^r \otimes ab) = \sum_{k \geq 0} \sum_{j \geq 0} (\gamma_k^r(a) \gamma_j^k(b)) \otimes y^j. \qquad (1.2.3)$$

Comparando (1.2.2) y (1.2.3) se obtiene

$$\gamma_j^r(ab) = \sum_{k \geq 0} \gamma_k^r(a) \gamma_j^k(b).$$

**Demostración de 4)**

Usando la condición d) sobre $\sigma$, para $k < r$, se tiene

$$\begin{aligned}
\sigma(y^r \otimes a) &= \sigma(y^{r-k} y^k \otimes a) = \sigma[\mu_{K[y]}(y^{r-k} \otimes y^k) \otimes Id_A(a)] \\
&= [\sigma \circ (\mu_{K[y]} \otimes Id_A)](y^{r-k} \otimes y^k \otimes a) \\
&= [(Id_A \otimes \mu_{K[y]}) \circ (\sigma \otimes Id_{K[y]}) \circ (Id_{K[y]} \otimes \sigma)](y^{r-k} \otimes y^k \otimes a) \\
&= [(Id_A \otimes \mu_{K[y]}) \circ (\sigma \otimes Id_{K[y]})](y^{r-k} \otimes \sigma(y^k \otimes a)) \\
&= [(Id_A \otimes \mu_{K[y]}) \circ (\sigma \otimes Id_{K[y]})]\left(y^{r-k} \otimes \sum_{l \geq 0} \gamma_l^k(a) \otimes y^l\right) \\
&= (Id_A \otimes \mu_{K[y]})\left(\sum_{l \geq 0} \sigma(y^{r-k} \otimes \gamma_l^k(a)) \otimes y^l\right) \\
&= (Id_A \otimes \mu_{K[y]})\left(\sum_{l \geq 0} \sum_{t \geq 0} (\gamma_t^{r-k}(\gamma_l^k(a)) \otimes y^t) \otimes y^l\right) \\
&= \sum_{l \geq 0} \sum_{t \geq 0} (\gamma_t^{r-k} \circ \gamma_l^k)(a) \otimes y^{t+l}.
\end{aligned}$$

Haciendo $j = t + l$ ($\Leftrightarrow t = j - l$)

$$\sigma(y^r \otimes a) = \sum_{j \geq 0} \left(\sum_{l=0}^{j} (\gamma_{j-l}^{r-k} \circ \gamma_l^k)(a)\right) \otimes y^j.$$

Por definición de $\gamma_j^r$:

$$\sigma(y^r \otimes a) = \sum_{j \geq 0} \gamma_j^r(a) \otimes y^j.$$

Comparando se concluye

$$\gamma_j^r = \sum_{l=0}^{j} \gamma_{j-l}^{r-k} \circ \gamma_l^k, \quad k < r. \qquad \square$$

## 1.3. Planos torcidos

**Definición 1.3.1.** *Sea $K$ un cuerpo y sea $K[x]$ el $K$-álgebra de polinomios en $x$. Un plano torcido es una extensión polinomial de $K[x]$. Es decir un plano torcido es un producto tensorial torcido de $K[x]$ con $K[y]$, al cual denotaremos con $K[x] \otimes_\sigma K[y]$ donde $\sigma$ es la aplicación de torcimiento correspondiente.*

Clasificar a los productos tensoriales $K[x] \otimes_\sigma K[y]$ es equivalente a clasificar a las aplicaciones de torcimiento $\sigma : K[y] \otimes K[x] \longrightarrow K[x] \otimes K[y]$.

Sabemos que una aplicación lineal $\sigma : K[y] \otimes K[x] \longrightarrow K[x] \otimes K[y]$ determina y es determinada por aplicaciones lineales

$$\gamma_j^r : K[x] \longrightarrow K[x], \ r, j \in \mathbb{N}_0,$$

tal que $\gamma_j^r(a) = 0$ para $r, a$ fijos, $j$ suficientemente grande, a través de la fórmula

$$\sigma(y^r \otimes a) = \sum_{j \geq 0} \gamma_j^r(a) \otimes y^j.$$

Además $\sigma$ es una aplicación de torcimiento si se cumplen las condiciones 1), 2), 3), 4) de la Proposición 1.2.1.

Si $\sigma$ es una aplicación de torcimiento, entonces las aplicaciones $\gamma_j^r$ definen una representación del producto tensorial torcido $K[x] \otimes_\sigma K[y]$ en $L(K[x]^{\mathbb{N}_0})$ (ver[8], Teorema 1.10). Para ello notemos que los elementos de $L(K[x]^{\mathbb{N}_0})$ son matrices infinitas con entradas en $K[x]$ indexadas por $\mathbb{N}_0 \times \mathbb{N}_0$ tal que cada fila tiene solamente un número finito de entradas diferentes de cero.

Denotaremos con $Y$ a la matriz infinita

$$Y := \begin{pmatrix} 0 & 1 & 0 & 0 & \cdots \\ 0 & 0 & 1 & 0 & \cdots \\ 0 & 0 & 0 & 1 & \cdots \\ \vdots & \vdots & \vdots & \vdots & \ddots \end{pmatrix}$$

Es decir $Y_{i,j} = \delta_{i+1,j}$ para todo $i, j \in \mathbb{N}_0$.

Para cualquier matriz $D$ se cumple

a) $(Y^k)_{i,j} = \delta_{i+k,j}, k \geq 0$

b) $(Y^k D)_{r,j} = D_{r+k,j}, \quad k \geq 0$

c)
$$(DY^k)_{r,j} = \begin{cases} D_{r,j-k}, & k \leq j \\ 0, & k > j \geq 0 \end{cases}$$

Asociada a la aplicación de torcimiento $\sigma : K[y] \otimes K[x] \longrightarrow K[x] \otimes K[y]$, determinada por las aplicaciones $K$-lineales $\gamma_j^r : K[x] \longrightarrow K[x]$, para cada $a \in K[x]$ definimos la matriz $M(a) \in L(K[x]^{\mathbb{N}_0})$ mediante $(M(a))_{r,j} = \gamma_j^r(a)$. Es decir

$$M(a) := \begin{pmatrix} a & 0 & \cdots & 0 & \cdots \\ \gamma_0^1(a) & \gamma_1^1(a) & \cdots & \gamma_n^1(a) & \cdots \\ \vdots & \vdots & \ddots & \vdots & \cdots \\ \gamma_0^n(a) & \gamma_1^n(a) & \cdots & \gamma_n^n(a) & \cdots \\ \vdots & \vdots & \vdots & \vdots & \vdots \end{pmatrix}$$

Esta matriz cumple la condición de finitud, ya que

$$\sigma(y^r \otimes a) = \sum_{j \geq 0} \gamma_j^r(a) \otimes y^j \in K[x] \otimes K[y]$$

de manera que $\gamma_j^r(a) \neq 0$ sólo para un número finito de $j's$.

**Observación 1.3.2.** Las matrices $M(a)$ cumplen

a) $M(1) = Id$.

b) $M(ab) = M(a)M(b)$ para todo $a, b \in K[x]$.

En efecto:
a) $(M(1))_{i,j} = \gamma_j^i(1) = \delta_{i,j}$.
b) $(M(ab))_{i,j} = \gamma_j^i(ab) = \sum_{k \geq 0} \gamma_k^i(a) \gamma_j^k(b) = \sum_{k \geq 0} (M(a))_{i,k}(M(b))_{k,j} = (M(a)M(b))_{i,j}$.

**Proposición 1.3.3.** *Las fórmulas $\psi(a \otimes 1) = M(a)$, $\psi(1 \otimes y) = Y$ definen un morfismo inyectivo de álgebras, $\psi : K[x] \otimes_\sigma K[y] \longrightarrow L(K[x]^{\mathbb{N}_0})$(representación fiel).*

**Demostración.** Por la Observación 1.3.2 para demostrar que $\psi$ es una aplicación de álgebras basta demostrar que

$$\psi(\sigma(y \otimes a)) = \psi(y)\psi(a).$$

En efecto:

$$\psi(\sigma(y \otimes a)) = \psi\left(\sum_{u \geq 0} \gamma_u^1(a) \otimes y^u\right) = \sum_{u \geq 0} \psi(\gamma_u^1(a) \otimes y^u) = \sum_{u \geq 0} \psi(\gamma_u^1(a))\psi(y^u).$$

Por otro lado

$$\begin{aligned}
(\psi(y)\psi(a))_{i,j} &= (YM(a))_{i,j} = (M(a))_{i+1,j} = \gamma_j^{i+1}(a) \\
&= \sum_{u=0}^{j} (\gamma_{j-u}^i \circ \gamma_u^1)(a) = \sum_{u=0}^{j} (\gamma_{j-u}^i(\gamma_u^1(a)) \\
&= \sum_{u=0}^{j} (M(\gamma_u^1(a))_{i,j-u} = \sum_{u \geq 0} (M(\gamma_u^1(a))Y^u)_{i,j} \\
&= \sum_{u \geq 0} (\psi(\gamma_u^1(a))\psi(y^u))_{i,j}.
\end{aligned}$$

De este modo queda demostrado que $\psi$ es una aplicación de álgebras.
La inyectividad se sigue del hecho de que la composición de $\psi$ con la sobreyección sobre la primera fila da el isomorfismo lineal canónico

$$K[x] \otimes_\sigma K[y] \xrightarrow{\cong} \bigoplus_{i \in N_0} K_i, \quad \text{donde } K_i \cong K[x].$$

$\square$

**Proposición 1.3.4.** *Sea $M \in L(K[x]^{\mathbb{N}_0})$ tal que $M_{0,j} = x\delta_{0,j}$, y*

$$Y^k M = \sum_{j \geq 0} M_{k,j}(M) Y^j, \quad k \geq 0. \tag{1.3.4}$$

(*Notemos que la suma es finita y que si $M_{k,j} = a_0 + a_1 x + a_2 x^2 + \cdots + a_n x^n \in K[x]$, entonces $M_{k,j}(M) = a_0 Id + a_1 M + a_2 M^2 + \cdots + a_n M^n \in L(K[x]^{\mathbb{N}_0})$).*
*Entonces las aplicaciones $\gamma_j^k : K[x] \longrightarrow K[x]$ definidas por*

$$\gamma_j^k(x^i) = (M^i)_{k,j}$$

*determinan una aplicación de torcimiento.*

**Demostración.** Demostraremos que las aplicaciones $\gamma_j^k$ satisfacen
1) $\gamma_j^0 = \delta_{0,j} Id_{K[x]}$
2) $\gamma_j^k(1) = \delta_{k,j}$
3) Para todo $k, j$ tenemos

$$\gamma_j^k(ab) = \sum_{s \geq 0} \gamma_s^k(a) \gamma_j^s(b)$$

para todo $a, b \in K[x]$
4) Para todo $k, j$, y todo $l < k$

$$\gamma_j^k = \sum_{u=0}^{j} \gamma_u^l \circ \gamma_{j-u}^{k-l}$$

**Demostración de 1)**
Reemplazando $i = 1$ y $k = 0$ en la definición $\gamma_j^k(x^i) = (M^i)_{k,j}$
se tiene $\gamma_j^0(x) = M_{0,j} = x\delta_{0,j}$ lo cual se cumple por hipótesis.
Luego $\gamma_j^0 = \delta_{0,j} Id_{K[x]}$.
**Demostración de 2)**

$$\gamma_j^k(1) = \gamma_j^k(x^0) = (M^0)_{k,j} = Id_{k,j} = \delta_{k,j}.$$

**Demostración de 3)**
Sean $a = x^i, \quad b = x^l$.
Entonces

$$\gamma_j^k(ab) = \gamma_j^k(x^{i+l}) = (M^{i+l})_{k,j} = (M^i M^l)_{k,j} = \sum_{s \geq 0} (M^i)_{k,s} (M^l)_{s,j} = \sum_{s \geq 0} \gamma_s^k(x^i) \gamma_j^s(x^l).$$

**Demostración de 4)**
Haciendo $k = l + m$ en la sumatoria, demostraremos que

$$\gamma_j^{l+m}(x^i) = \sum_{u=0}^{j} (\gamma_u^l \circ \gamma_{j-u}^m)(x^i).$$

Usaremos inducción sobre $i$.

$\boxed{i=0}$

Demostraremos que
$$\gamma_j^{l+m}(1) = \sum_{u=0}^{j}(\gamma_u^l \circ \gamma_{j-u}^m)(1).$$

En efecto:

$\gamma_j^{l+m}(1) = \delta_{l+m,j}$ por 2).

Por otro lado
$$\sum_{u=0}^{j}(\gamma_u^l \circ \gamma_{j-u}^m)(1) = \sum_{u=0}^{j}\gamma_u^l(\gamma_{j-u}^m(1)).$$

Pero
$$\gamma_{j-u}^m(1) = \delta_{m,j-u} = \begin{cases} 1 & si \quad m = j-u \\ 0 & si \quad m \neq j-u \end{cases} = \begin{cases} 1 & si \quad u = j-m \\ 0 & si \quad u \neq j-m \end{cases}.$$

Luego
$$\sum_{u=0}^{j}(\gamma_u^l \circ \gamma_{j-u}^m)(1) = \gamma_{j-m}^l(1) = \delta_{l,j-m} = \delta_{l+m,j}.$$

$\boxed{i=1}$

Demostraremos que
$$\gamma_j^{l+m}(x) = \sum_{u=0}^{j}(\gamma_u^l \circ \gamma_{j-u}^m)(x).$$

En efecto:

$$\gamma_j^{l+m}(x) = M_{l+m,j} = (Y^m M)_{l,j} = \sum_{u \geq 0}(M_{m,u}(M)Y^u)_{l,j} = \sum_{u=0}^{j}(M_{m,u}(M))_{l,j-u}.$$

Supongamos ahora que $M_{r,s} = \sum_{p=0}^{t} a_p x^p$, entonces

$$M_{r,s}(M) = \sum_{p=0}^{t} a_p M^p.$$

Evaluando en la entrada $(k, j)$

$$\begin{aligned}(M_{r,s}(M))_{k,j} &= \sum_{p=0}^{t} a_p(M^p)_{k,j} = \sum_{p=0}^{t} a_p \gamma_j^k(x^p) \\ &= \gamma_j^k\left(\sum_{p=0}^{t} a_p x^p\right) = \gamma_j^k(M_{r,s}).\end{aligned}$$

Luego

$$\begin{aligned} \gamma_j^{l+m}(x) &= \sum_{u=0}^{j}(M_{m,u}(M))_{l,j-u} = \sum_{u=0}^{j}\gamma_{j-u}^{l}(M_{m,u}) \\ &= \sum_{u=0}^{j}\gamma_{j-u}^{l}(\gamma_u^m(x)) = \sum_{u=0}^{j}(\gamma_{j-u}^{l}\circ\gamma_u^m)(x) \\ &= \sum_{u=0}^{j}(\gamma_u^l\circ\gamma_{j-u}^m)(x). \end{aligned}$$

Supongamos que se cumple

$$\gamma_j^{l+m}(x^i) = \sum_{u=0}^{j}(\gamma_u^l\circ\gamma_{j-u}^m)(x^i)$$

para algún $i \geq 1$.

Demostraremos

$$\gamma_j^{l+m}(x^{i+1}) = \sum_{u=0}^{j}(\gamma_u^l\circ\gamma_{j-u}^m)(x^{i+1}).$$

En efecto:

$$\gamma_j^{l+m}(x^{i+1}) = \gamma_j^{l+m}(x^i x) = \sum_{r\geq 0}\gamma_r^{l+m}(x^i)\gamma_j^r(x) = \sum_{r\geq 0}\left(\sum_{s=0}^{r}\gamma_s^l(\gamma_{r-s}^m(x^i))\right)\gamma_j^r(x).$$

Haciendo $n = r - s$ se tiene

$$\begin{aligned} \gamma_j^{l+m}(x^{i+1}) &= \sum_{s\geq 0}\sum_{n\geq 0}\gamma_s^l(\gamma_n^m(x^i))\gamma_j^{s+n}(x) \\ &= \sum_{s\geq 0}\sum_{n\geq 0}\gamma_s^l(\gamma_n^m(x^i))\sum_{u=0}^{j}\gamma_u^s(\gamma_{j-u}^n(x)) \\ &= \sum_{s\geq 0}\sum_{n\geq 0}\sum_{u=0}^{j}\gamma_s^l(\gamma_n^m(x^i))\gamma_u^s(\gamma_{j-u}^n(x)) \\ &= \sum_{n\geq 0}\sum_{u=0}^{j}\sum_{s\geq 0}\gamma_s^l(\underbrace{\gamma_n^m(x^i)}_{a})\gamma_u^s(\underbrace{\gamma_{j-u}^n(x)}_{b}) = \sum_{n\geq 0}\sum_{u=0}^{j}\gamma_u^l(ab) \\ &= \sum_{u=0}^{j}\sum_{n\geq 0}\gamma_u^l(ab) = \sum_{u=0}^{j}\gamma_u^l\left(\sum_{n\geq 0}ab\right) = \sum_{u=0}^{j}\gamma_u^l\left(\sum_{n\geq 0}\gamma_n^m(x^i)\gamma_{j-u}^n(x)\right) \\ &= \sum_{u=0}^{j}(\gamma_u^l\circ\gamma_{j-u}^m)(x^{i+1}). \end{aligned}$$
$\square$

**Observación 1.3.5.** Se cumple el recíproco de la Proposición 1.3.4, es decir, si existe una aplicación de torcimiento con aplicaciones $\gamma_j^k$, entonces la matriz $M$ dada por $M_{k,j} := \gamma_j^k(x)$ cumple (1.3.4).

## 1.4. Planos torcidos graduados

Ahora supongamos que la potencial aplicación de torcimiento dada por

$$\sigma(y^k \otimes a) = \sum_{j \geq 0} \gamma_j^k(a) \otimes y^j$$

es graduada.
Haciendo $a = x^i$ se tiene

$$\sigma(y^k \otimes x^i) = \sum_{j \geq 0} \gamma_j^k(x^i) \otimes y^j.$$

Como $y^k \otimes x^i$ es de grado $k+i$ entonces $\gamma_j^k(x^i)$ debe ser de grado $k+i-j$.
Es decir

$$\gamma_j^k(x^i) = \alpha_{ijk} x^{k+i-j}, \quad \alpha_{ijk} \in K.$$

**Proposición 1.4.1.** *Sean $Y, M \in L(K^{\mathbb{N}_0})$ tales que $Y_{i,j} = \delta_{i+1,j}$, $M_{0,j} = \delta_{0,j}$, $M_{k,j} = 0$ para $j > k+1$, y*

$$Y^k M = \sum_{j=0}^{k+1} M_{k,j} M^{k+1-j} Y^j, \quad k \geq 0.$$

*Entonces las aplicaciones*

$$\gamma_j^k(x^i) = (M^i)_{k,j} x^{k+i-j}$$

*determinan una aplicación de torcimiento **graduada**.*

**Demostración.** La aplicación $\sigma$ determinada por los $\gamma_j^k$ evidentemente es graduada, así que, usando la Proposición 1.3.4, solamente necesitamos demostrar que la matriz $\widetilde{M}$ definida por $\widetilde{M}_{k,j} = M_{k,j} x^{k+1-j}$, es decir

$$\widetilde{M} = \begin{pmatrix} x & 0 & 0 & 0 & 0 & \cdots \\ M_{1,0} x^2 & M_{1,1} x & M_{1,2} & 0 & 0 & \cdots \\ M_{2,0} x^3 & M_{2,1} x^2 & M_{2,2} x & M_{2,3} & 0 & \cdots \\ \vdots & \vdots & \vdots & \vdots & \vdots & \vdots \end{pmatrix}$$

satisface

(1) $\widetilde{M}_{0,j} = x \delta_{0,j}$,

(2) $\gamma_j^k(x^i) = (\widetilde{M}^i)_{k,j}$,

(3) $Y^k \widetilde{M} = \sum_{j \geq 0} \widetilde{M}_{k,j}(\widetilde{M}) Y^j, \quad k \geq 0.$

**Demostración de (1)** :

$$\widetilde{M}_{0,j} = M_{0,j} x^{1-j} = \delta_{0,j} x^{1-j} = \begin{cases} x & si \quad j=0 \\ 0 & si \quad j \neq 0 \end{cases}.$$

Luego $\widetilde{M}_{0,j} = x \delta_{0,j}$.

**Demostración de (2)** :

Es suficiente demostrar que

$$(\widetilde{M}^i)_{k,j} = (M^i)_{k,j} x^{k+i-j}. \qquad (1.4.5)$$

Por inducción sobre $i$

$\boxed{i=0}$

$$(\widetilde{M}^0)_{k,j} = (Id)_{k,j} = \delta_{k,j}.$$

Por otro lado

$$(M^0)_{k,j} x^{k-j} = (Id)_{k,j} x^{k-j} = \delta_{k,j} x^{k-j} = \delta_{k,j}.$$

$\boxed{i=1}$

$$\widetilde{M}_{k,j} = M_{k,j} x^{k+1-j}$$

se cumple por definición de $\widetilde{M}$.

Supongamos que se cumple (1.4.5) para algún $i \geq 1$.

Demostraremos que se cumple para $i+1$.

En efecto:

$$\begin{aligned}
(\widetilde{M}^{i+1})_{k,j} &= (\widetilde{M}^i \widetilde{M})_{k,j} = \sum_{l \geq 0} (\widetilde{M}^i)_{k,l} (\widetilde{M})_{l,j} \\
&= \sum_{l \geq 0} (M^i)_{k,l} x^{k+i-l} M_{l,j} x^{l+1-j} \\
&= \sum_{l \geq 0} (M^i)_{k,l} M_{l,j} x^{k+i+1-j} = (M^{i+1})_{k,j} x^{k+i+1-j}.
\end{aligned}$$

**Demostración de (3)** :

Es suficiente demostrar la igualdad en la entrada $(r,s)$.

Tenemos

$$\begin{aligned}
(Y^k \widetilde{M})_{r,s} &= \widetilde{M}_{r+k,s} = M_{r+k,s} x^{r+k+1-s} \\
&= (Y^k M)_{r,s} x^{r+k+1-s} = \sum_{j=0}^{k+1} (M_{k,j} M^{k+1-j} Y^j)_{r,s} x^{r+k+1-s} \\
&= \sum_{j=0}^{k+1} M_{k,j} (M^{k+1-j})_{r,s-j} x^{r+k+1-s}.
\end{aligned}$$

Ahora $\widetilde{M}_{k,j} = M_{k,j} x^{k+1-j}$.

Luego $\widetilde{M}_{k,j}(\widetilde{M}) = M_{k,j} \widetilde{M}^{k+1-j}$, y tenemos

$$(\widetilde{M}_{k,j}(\widetilde{M}))_{r,s-j} = M_{k,j}(\widetilde{M}^{k+1-j})_{r,s-j} = M_{k,j}(M^{k+1-j})_{r,s-j} x^{r+k+1-s}.$$

Reemplazando en la última sumatoria se tiene

$$(Y^k \widetilde{M})_{r,s} = \sum_{j=0}^{k+1} (\widetilde{M}_{k,j}(\widetilde{M}))_{r,s-j} = \sum_{j=0}^{k+1} (\widetilde{M}_{k,j}(\widetilde{M}) Y^j)_{r,s} = \sum_{j \geq 0} (\widetilde{M}_{k,j}(\widetilde{M}) Y^j)_{r,s}.$$

La última igualdad se cumple por la condición $M_{k,j} = 0$ para $j > k+1$. $\square$

**Observación 1.4.2.** En el caso graduado las fórmulas $\psi(x \otimes 1) = M$ y $\psi(1 \otimes y) = Y$ definen una aplicación inyectiva de álgebras (representación fiel) $\psi : K[x] \otimes_\sigma K[y] \to L(K^{\mathbb{N}_0})$.

# Capítulo 2

# El caso cuadrático

## 2.1. Matriz asociada a una aplicación de torcimiento graduada

Una matriz que satisface las condiciones de la Proposición 1.4.1 determina una aplicación de torcimiento graduada, pero también es claro que una aplicación de torcimiento graduada determina una tal matriz $M$.(Notemos que $M_{k,j} = 0$ para $j > k+1$).
Para clasificar a las aplicaciones de torcimiento graduadas tenemos que clasificar a las matrices $M$.
Denotaremos $M_{1,0} := a$, $M_{1,1} := b$, $M_{1,2} := c$ para una tal matriz. En algunos casos los valores de $a, b$ y $c$ determinan completamente a la matriz $M$ (y por lo tanto a la aplicación de torcimiento).

### 2.1.1. Caso $a = 0$

Si $\boxed{a = 0}$ entonces para $k \geq 1$, las únicas entradas posiblemente distintas de cero de $M_{k,*}$ son $M_{k,k}$ y $M_{k,k+1}$.
Para $k = 1$ se cumple $YM = bMY + cY^2$.
Haciendo $x = M$, $y = Y$ en esta igualdad matricial, se obtiene $yx = bxy + cy^2$.
Usaremos un proceso de inducción sobre $k$.
Para $k = 2$ se tiene

$$\begin{aligned} y^2 x &= y(bxy + cy^2) = b(yx)y + cy^3 = b(bxy + cy^2)y + cy^3 \\ &= b^2 xy^2 + c(b+1)y^3. \end{aligned}$$

Supongamos que se cumple

$$y^k x = b^k x y^k + c(b^{k-1} + b^{k-2} + \cdots + b + 1)y^{k+1},$$

para un $k \geq 1$.
Entonces

$$\begin{aligned} y^{k+1}x = y(y^k x) &= y(b^k x y^k + c(b^{k-1} + b^{k-2} + \cdots + b + 1)y^{k+1}) \\ &= b^k(yx)y^k + c(b^{k-1} + b^{k-2} + \cdots + b + 1)y^{k+2} \\ &= b^k(bxy + cy^2)y^k + c(b^{k-1} + b^{k-2} + \cdots + b + 1)y^{k+2} \\ &= b^{k+1}xy^{k+1} + c(b^k + b^{k-1} + \cdots + b + 1)y^{k+2} \end{aligned}$$

De esta manera hemos demostrado que

$$Y^k M = b^k M Y^k + c(b^{k-1} + b^{k-2} + \cdots + b + 1)Y^{k+1}$$

para $k \geq 1$.

Pongamos ahora $b_n = M_{n,n}, \quad c_n = M_{n,n+1}$.
Entonces a partir de

$$Y^k M = \sum_{j=0}^{k+1} M_{k,j} M^{k+1-j} Y^j, \quad k \geq 0,$$

obtenemos para $k = 1$

$$YM = \sum_{j=0}^{2} M_{1,j} M^{2-j} Y^j = \underbrace{M_{1,0}}_{0} M^2 + \underbrace{M_{1,1}}_{b} MY + \underbrace{M_{1,2}}_{c} Y^2.$$

Es decir $\boxed{YM = bMY + cY^2}$.
Evaluando en la entrada $(n-1, n)$ se tiene por un lado

$$(YM)_{n-1,n} = M_{n,n} = b_n.$$

Por otro lado

$$b(MY)_{n-1,n} + c(Y^2)_{n-1,n} = bM_{n-1,n-1} + c\delta_{n+1,n} = b\,b_{n-1}.$$

Luego se obtiene la fórmula de recurrencia

$$\boxed{b_n = b\,b_{n-1}}, \quad n \geq 1.$$

Así $b_1 = b\,b_0 = b\,M_{0,0} = b$,
$b_2 = b\,b_1 = b\,b = b^2$.
En general se obtiene

$$b_n = b^n, \quad n \geq 0.$$

Evaluando en la entrada $(n-1, n+1)$ se tiene por un lado

$$(YM)_{n-1,n+1} = M_{n,n+1} = c_n.$$

Por otro lado

$$b(MY)_{n-1,n+1} + c(Y^2)_{n-1,n+1} = bM_{n-1,n} + c\delta_{n+1,n+1} = bc_{n-1} + c.$$

Luego se obtiene la fórmula de recurrencia

$$\boxed{c_n = bc_{n-1} + c}, \quad n \geq 1.$$

### 2.1.2. Caso $a = 0$, $b = 1$

Si $\boxed{b = 1}$ entonces $c_n = nc$, $n \geq 0$ ($c_0 = M_{0,1} = 0$).
La matriz correspondiente es

$$M = \begin{pmatrix} 1 & 0 & 0 & 0 & 0 & 0 & \cdots \\ 0 & 1 & c & 0 & 0 & 0 & \cdots \\ 0 & 0 & 1 & 2c & 0 & 0 & \cdots \\ 0 & 0 & 0 & 1 & 3c & 0 & \cdots \\ \vdots & \vdots & \vdots & \ddots & \ddots & \ddots & \ddots \end{pmatrix}.$$

La aplicación de torcimiento graduada está dada por

$$\tau(1 \otimes x) = x \otimes 1,$$
$$\tau(y \otimes x) = x \otimes y + c(1 \otimes y^2),$$
$$\tau(y^2 \otimes x) = x \otimes y^2 + 2c(1 \otimes y^3),$$
$$\tau(y^3 \otimes x) = x \otimes y^3 + 3c(1 \otimes y^4).$$

En general

$$\tau(y^n \otimes x) = x \otimes y^n + nc(1 \otimes y^{n+1}).$$

Usando la inclusión $K[x] \hookrightarrow K[x] \otimes K[y] \cong_K K[x,y]$ y similarmente $K[y] \hookrightarrow K[x,y]$, podemos escribir

$$yx = xy + cy^2,$$
$$y^2 x = xy^2 + 2cy^3,$$
$$y^3 x = xy^3 + 3cy^4.$$

En general

$$y^n x = xy^n + ncy^{n+1}.$$

Por la Proposición 1.4.1 las igualdades matriciales respectivas

$$YM = MY + cY^2,$$
$$Y^2M = MY^2 + 2cY^3,$$
$$Y^3M = MY^3 + 3cY^4,$$

en general

$$Y^n M = MY^n + ncY^{n+1},$$

garantizan que $\tau$ es una aplicación de torcimiento graduada. De hecho la primera igualdad implica todas las demás y se puede ver que efectivamente $M$ cumple $YM = MY + cY^2$. Es decir

$$\begin{pmatrix} 0 & 1 & c & 0 & 0 & \cdots \\ 0 & 0 & 1 & 2c & 0 & \cdots \\ 0 & 0 & 0 & 1 & 3c & \cdots \\ 0 & 0 & 0 & 0 & 1 & \cdots \\ \vdots & \vdots & \vdots & \vdots & \vdots & \ddots \end{pmatrix} = \begin{pmatrix} 0 & 1 & 0 & 0 & 0 & \cdots \\ 0 & 0 & 1 & c & 0 & \cdots \\ 0 & 0 & 0 & 1 & 2c & \cdots \\ 0 & 0 & 0 & 0 & 1 & \cdots \\ \vdots & \vdots & \vdots & \vdots & \vdots & \ddots \end{pmatrix} + c \begin{pmatrix} 0 & 0 & 1 & 0 & 0 & \cdots \\ 0 & 0 & 0 & 1 & 0 & \cdots \\ 0 & 0 & 0 & 0 & 1 & \cdots \\ 0 & 0 & 0 & 0 & 0 & \cdots \\ \vdots & \vdots & \vdots & \vdots & \vdots & \ddots \end{pmatrix}.$$

### 2.1.3. Caso $a = 0$, $b \neq 1$

Si $\boxed{b \neq 1}$ entonces

$$c_1 = bc_0 + c = c,$$
$$c_2 = bc_1 + c = bc + c = c(b+1),$$
$$c_3 = bc_2 + c = bc(b+1) + c = c(b^2 + b + 1),$$

en general

$$c_n = c(b^{n-1} + b^{n-2} + \cdots + b^2 + b + 1) = c\left(\frac{b^n - 1}{b - 1}\right), \quad n \geq 0.$$

La matriz correspondiente es

$$M = \begin{pmatrix} 1 & 0 & 0 & 0 & 0 & 0 & \cdots \\ 0 & b & c & 0 & 0 & 0 & \cdots \\ 0 & 0 & b^2 & c(b+1) & 0 & 0 & \cdots \\ 0 & 0 & 0 & b^3 & c(b^2+b+1) & 0 & \cdots \\ \vdots & \vdots & \vdots & \vdots & \vdots & \vdots & \ddots \end{pmatrix}.$$

La aplicación de torcimiento graduada está dada por

$$\tau(1 \otimes x) = x \otimes 1,$$
$$\tau(y \otimes x) = bx \otimes y + c(1 \otimes y^2),$$
$$\tau(y^2 \otimes x) = b^2 x \otimes y^2 + c(b+1)(1 \otimes y^3),$$
$$\tau(y^3 \otimes x) = b^3 x \otimes y^3 + c(b^2 + b + 1)(1 \otimes y^4),$$

en general
$$\tau(y^n \otimes x) = b^n x \otimes y^n + c\left(\frac{b^n-1}{b-1}\right)(1 \otimes y^{n+1}).$$

Usando la inclusión $K[x] \hookrightarrow K[x] \otimes K[y] \cong_K K[x,y]$ y similarmente $K[y] \hookrightarrow K[x,y]$ podemos escribir

$$yx = bxy + cy^2,$$
$$y^2 x = b^2 x y^2 + c(b+1)y^3,$$
$$y^3 x = b^3 x y^3 + c(b^2+b+1)y^4,$$

en general

$$y^n x = b^n x y^n + c\left(\frac{b^n-1}{b-1}\right)y^{n+1}.$$

Por la Proposición 1.4.1 las igualdades matriciales respectivas

$$YM = bMY + cY^2,$$
$$Y^2 M = b^2 M Y^2 + c(b+1)Y^3,$$
$$Y^3 M = b^3 M Y^3 + c(b^2+b+1)Y^4,$$

en general

$$Y^n M = b^n M Y^n + c\left(\frac{b^n-1}{b-1}\right)Y^{n+1},$$

garantizan que $\tau$ es una aplicación de torcimiento graduada. De hecho la primera igualdad implica todas las demás y se puede ver que efectivamente $M$ cumple $YM = bMY + cY^2$. Es decir

$$\begin{pmatrix} 0 & b & c & 0 & \cdots \\ 0 & 0 & b^2 & c(b+1) & \cdots \\ 0 & 0 & 0 & b^3 & \cdots \\ \vdots & \vdots & \vdots & \vdots & \ddots \end{pmatrix} = b \begin{pmatrix} 0 & 1 & 0 & 0 & \cdots \\ 0 & 0 & b & c & \cdots \\ 0 & 0 & 0 & b^2 & \cdots \\ \vdots & \vdots & \vdots & \vdots & \ddots \end{pmatrix} + c \begin{pmatrix} 0 & 0 & 1 & 0 & \cdots \\ 0 & 0 & 0 & 1 & \cdots \\ 0 & 0 & 0 & 0 & \cdots \\ \vdots & \vdots & \vdots & \vdots & \ddots \end{pmatrix}.$$

### 2.1.4. Caso $a = 1$

**Lema 2.1.1.** *Sean $A$ una $K$-álgebra asociativa, $k \geq 2$, $x \in A, y \in A$, $M_{i,j} \in K$ para $0 \leq j < k$, $0 \leq i \leq j+1$ tales que $M_{0,j} = \delta_{0,j}$, $M_{1,0} = 1$, y*

$$y^j x = \sum_{i=0}^{j+1} M_{j,i} x^{j+1-i} y^i.$$

*Pongamos $b = M_{1,1}$, $c = M_{1,2}$.*
*Entonces*

$$(1 - M_{k-1,k}) y^k x = \sum_{s=0}^{k+1} \overline{M}_{k,s} x^{k+1-s} y^s, \tag{2.1.1}$$

donde $\overline{M}_{k,0} = \sum_{i=0}^{k-1} M_{k-1,i} M_{i,0}$,

$\overline{M}_{k,k+1} = c + b M_{k-1,k}$,

$\overline{M}_{k,s} = b M_{k-1,s-1} + \sum_{i=s-1}^{k-1} M_{k-1,i} M_{i,s} \quad s = 1, 2, \cdots, k.$

**Demostración.** Tenemos $y^k x = y^{k-1}(yx)$.

Por hipótesis

$$y x = \sum_{i=0}^{2} M_{1,i} x^{2-i} y^i = M_{1,0} x^2 + M_{1,1} x y + M_{1,2} y^2 = x^2 + b x y + c y^2.$$

Reemplazando obtenemos

$$\begin{aligned} y^k x &= y^{k-1}(x^2 + b x y + c y^2) \\ &= y^{k-1} x^2 + b y^{k-1} x y + c y^{k+1} \\ &= (y^{k-1} x) x + b (y^{k-1} x) y + c y^{k+1}. \end{aligned}$$

Usando nuevamente la hipótesis para $j = k - 1$

$$\begin{aligned} y^k x &= (y^{k-1} x) x + b (y^{k-1} x) y + c y^{k+1} \\ &= \left( \sum_{i=0}^{k} M_{k-1,i} x^{k-i} y^i \right) x + b \left( \sum_{i=0}^{k} M_{k-1,i} x^{k-i} y^i \right) y + c y^{k+1} \\ &= \left( \sum_{i=0}^{k-1} M_{k-1,i} x^{k-i} y^i + M_{k-1,k} y^k \right) x + b \sum_{i=0}^{k} M_{k-1,i} x^{k-i} y^{i+1} + c y^{k+1} \\ &= \sum_{i=0}^{k-1} M_{k-1,i} x^{k-i} (y^i x) + M_{k-1,k} y^k x + b \sum_{i=0}^{k} M_{k-1,i} x^{k-i} y^{i+1} + c y^{k+1}. \end{aligned}$$

Luego

$$\begin{aligned} (1 - M_{k-1,k}) y^k x &= \sum_{i=0}^{k-1} M_{k-1,i} x^{k-i} (y^i x) + b \sum_{s=0}^{k} M_{k-1,s} x^{k-s} y^{s+1} + c y^{k+1} \\ &= \sum_{i=0}^{k-1} M_{k-1,i} x^{k-i} (y^i x) + b \sum_{s=0}^{k-1} M_{k-1,s} x^{k-s} y^{s+1} + b M_{k-1,k} y^{k+1} + c y^{k+1} \\ &= \sum_{i=0}^{k-1} M_{k-1,i} x^{k-i} (y^i x) + b \sum_{s=1}^{k} M_{k-1,s-1} x^{k-s+1} y^s + (c + b M_{k-1,k}) y^{k+1}. \end{aligned}$$

Usando nuevamente la hipótesis en la primera sumatoria, en la expresión

$$y^i x = \sum_{s=0}^{i+1} M_{i,s} x^{i+1-s} y^s \quad 0 \le i < k,$$

se obtiene

$$
\begin{aligned}
\sum_{i=0}^{k-1} M_{k-1,i} x^{k-i} (y^i x) &= \sum_{i=0}^{k-1} M_{k-1,i} x^{k-i} \left( \sum_{s=0}^{i+1} M_{i,s} x^{i+1-s} y^s \right) \\
&= \sum_{i=0}^{k-1} M_{k-1,i} x^{k-i} \left( M_{i,0} x^{i+1} + \sum_{s=1}^{i+1} M_{i,s} x^{i+1-s} y^s \right) \\
&= \sum_{i=0}^{k-1} M_{k-1,i} M_{i,0} x^{k+1} + \sum_{i=0}^{k-1} \sum_{s=1}^{i+1} M_{k-1,i} M_{i,s} x^{k+1-s} y^s \\
&= \sum_{i=0}^{k-1} M_{k-1,i} M_{i,0} x^{k+1} + \sum_{s=1}^{k} \sum_{i=s-1}^{k-1} M_{k-1,i} M_{i,s} x^{k+1-s} y^s.
\end{aligned}
$$

Reemplazando obtenemos

$$
\begin{aligned}
(1 - M_{k-1,k}) y^k x &= \sum_{i=0}^{k-1} M_{k-1,i} M_{i,0} x^{k+1} + \sum_{s=1}^{k} \sum_{i=s-1}^{k-1} M_{k-1,i} M_{i,s} x^{k+1-s} y^s \\
&\quad + b \sum_{s=1}^{k} M_{k-1,s-1} x^{k-s+1} y^s + (c + b M_{k-1,k}) y^{k+1}.
\end{aligned}
$$

Definiendo
$$\overline{M}_{k,0} = \sum_{i=0}^{k-1} M_{k-1,i} M_{i,0},$$
$$\overline{M}_{k,k+1} = c + b M_{k-1,k},$$
$$\overline{M}_{k,s} = b M_{k-1,s-1} + \sum_{i=s-1}^{k-1} M_{k-1,i} M_{i,s}, \quad s = 1, 2, \cdots, k,$$
y reemplazando en

$$
\begin{aligned}
(1 - M_{k-1,k}) y^k x &= \sum_{i=0}^{k-1} M_{k-1,i} M_{i,0} x^{k+1} + \sum_{s=1}^{k} \left( \sum_{i=s-1}^{k-1} M_{k-1,i} M_{i,s} + b M_{k-1,s-1} \right) x^{k+1-s} y^s \\
&\quad + (c + b M_{k-1,k}) y^{k+1},
\end{aligned}
$$

se obtiene

$$
(1 - M_{k-1,k}) y^k x = \overline{M}_{k,0} x^{k+1} + \sum_{s=1}^{k} \overline{M}_{k,s} x^{k+1-s} y^s + \overline{M}_{k,k+1} y^{k+1}.
$$

Finalmente escribimos esto como

$$
(1 - M_{k-1,k}) y^k x = \sum_{s=0}^{k+1} \overline{M}_{k,s} x^{k+1-s} y^s. \qquad \square
$$

**Observación 2.1.2.** Sea $M \in L(K^{N_0})$ tal que $M_{0j} = \delta_{0j}$, $M_{ji} = 0$ para todo $i > j+1$ y para algún $k > 1$,

$$Y^j M = \sum_{i=0}^{j+1} M_{ji} M^{j+1-i} Y^i, \quad \text{para } j < k. \qquad (2.1.2)$$

Entonces $x = M$, $y = Y$ satisfacen las hipótesis del lema, la igualdad (2.1.1) se lee

$$(1-M_{k-1,k})Y^k M = \sum_{s=0}^{k+1} \overline{M}_{k,s} M^{k+1-s} Y^s,$$

y si tomamos la entrada $(0, i)$, entonces el lado de la izquierda da

$$\left((1-M_{k-1,k})Y^k M\right)_{0,i} = (1-M_{k-1,k})M_{k,i},$$

y el lado de la derecha da $\overline{M}_{k,i}$, ya que

$$\left(M^{k+1-s} Y^s\right)_{0,i} = \left(M^{k+1-s}\right)_{0,i-s} = \delta_{i,s}.$$

Por lo tanto en una aplicación de torcimiento con $M_{k-1,k} \neq 1$, los coeficientes $M_{k,i}$ están determinados de manera única por los coeficientes $M_{j,i}$ con $j < k$.

### 2.1.5. Caso general

**Observación 2.1.3.** Dada una aplicación de torcimiento $\sigma$ por

$$\sigma(y \otimes x) = a(x^2 \otimes 1) + b(x \otimes y) + c(1 \otimes y^2). \tag{2.1.3}$$

Dado $\lambda \neq 0$, podemos sustituir $\sigma$ por la aplicación de torcimiento isomorfa

$$\sigma' = (f^{-1} \otimes Id) \circ \sigma \circ (Id \otimes f),$$

donde $f(x) = \lambda x$.
Entonces para $\sigma'$ tenemos $a' = \dfrac{a}{\lambda}$, $b' = b$ y $c' = c\lambda$. De esta manera no hay pérdida de generalidad si asumimos $a = 1$ en (2.1.3), cuando $a \neq 0$.
Verifiquemos los parámetros de $\sigma'$:

$$\begin{aligned}
\sigma'(y \otimes x) &= [(f^{-1} \otimes Id) \circ \sigma \circ (Id \otimes f)](y \otimes x) \\
&= [(f^{-1} \otimes Id) \circ \sigma](y \otimes \lambda x) \\
&= (f^{-1} \otimes Id)(\sigma(y \otimes \lambda x)) = (f^{-1} \otimes Id)(\lambda \sigma(y \otimes x)) \\
&= (f^{-1} \otimes Id)(a\lambda(x^2 \otimes 1) + \lambda b(x \otimes y) + \lambda c(1 \otimes y^2)) \\
&= f^{-1}(a\lambda x^2) \otimes 1 + f^{-1}(\lambda b x) \otimes y + f^{-1}(\lambda c)(1 \otimes y^2) \\
&= \frac{a}{\lambda} x^2 \otimes 1 + b(x \otimes y) + \lambda c(1 \otimes y^2).
\end{aligned}$$

Para el caso $a = 0$ y $c \neq 0$ podemos tomar $\lambda = 1/c$, y entonces $\sigma'(y \otimes x) = b x \otimes y + 1 \otimes y^2$. Por lo tanto en este caso podemos asumir que $c = 1$.

**Proposición 2.1.4.** *Supongamos que M determina una aplicación de torcimiento graduada.*
*Asumamos $a = 1$, y $M_{k,k+1} \neq 1$ para todo $k \geq 1$.*
*Entonces b y c determinan en forma única a la matriz M (y por lo tanto a la aplicación de torcimiento). Recordemos que $a = M_{1,0}$, $b = M_{1,1}$ y $c = M_{1,2}$.*

**Demostración.** Por la Observación 2.1.2 se tiene que en ese caso las entradas $M_{j,i}$ con $j < k$ determinan de manera única a las entradas $M_{k,i}$, luego, por inducción, las entradas $M_{1,j}$, es decir $M_{1,0} = 1$, $M_{1,1} = b$, $M_{1,2} = c$, determinan a toda la matriz $M$. □

Sin embargo no toda elección de $b, c$ es válida, como veremos en la siguiente sección en la Proposición 2.2.1

**Lema 2.1.5.** *Sea M una matriz que representa a una aplicación de torcimiento graduada y pongamos $b_n = M_{n,n}$, $c_n = M_{n,n+1}$,   $n \geq 0$.*
*Entonces*
$$c_{n+1}(1 - c_n) = b\, c_n + c, \quad n \geq 0. \tag{2.1.4}$$

*Además*

1) *Si $c \neq 1$ entonces $c_n \neq 1$,  para todo $n \geq 0$,*

2) *Si $c = 1$ entonces $b = -1$ y*

$$(1 - c_n)(1 - c_{n+1}) = 0, \quad n \geq 0.$$

**Demostración.** Sabemos que $M$ satisface la igualdad

$$YM = M^2 + bMY + cY^2.$$

Evaluando en la entrada $(n, n+2)$, obtenemos
$(YM)_{n,n+2} = M_{n+1,n+2} = c_{n+1}$,
$(MY)_{n,n+2} = M_{n,n+1} = c_n$,
$(Y^2)_{n,n+2} = \delta_{n+2,n+2} = 1$,
$(M^2)_{n,n+2} = \sum_{k \geq 0} M_{n,k} M_{k,n+2}$.
Ahora $k > n+1$ implica $M_{n,k} = 0$ y $k < n+1$ implica $M_{k,n+2} = 0$.
Luego la sumatoria se reduce al término correspondiente a $k = n+1$, es decir
$(M^2)_{n,n+2} = M_{n,n+1} M_{n+1,n+2} = c_n c_{n+1}$.
Reemplazando en

$$(YM)_{n,n+2} = (M^2)_{n,n+2} + b(MY)_{n,n+2} + c(Y^2)_{n,n+2}$$

se tiene:
$$c_{n+1} = c_n c_{n+1} + b c_n + c.$$

Es decir
$$c_{n+1}(1 - c_n) = b c_n + c, \quad n \geq 0.$$

**Demostración de 1)**

Supongamos que $c \neq 1$ y por el absurdo supongamos que $c_{n_0} = 1$ para algún $n_0 > 1$.

Notemos que $c_0 = 0, c_1 = c$.

Reemplazando en (2.1.4) se obtiene:
$$c_{n_0+1} \underbrace{(1 - c_{n_0})}_{0} = b \underbrace{c_{n_0}}_{1} + c$$

y por lo tanto $b = -c$.

Reemplazando nuevamente en (2.1.4) se tiene
$$c_{n+1}(1 - c_n) = -c c_n + c = c(1 - c_n).$$

De esto se obtiene
$$(1 - c_n)(c - c_{n+1}) = 0, \quad n \geq 0.$$

Si $n = 1$ entonces $(1 - c_1)(c - c_2) = 0$ y de allí $c_2 = c$ pues $1 - c_1 = 1 - c \neq 0$.

Si $n = 2$ entonces $(1 - c_2)(c - c_3) = 0$ y de allí $c_3 = c$ pues $1 - c_2 = 1 - c \neq 0$.

Por recurrencia se llega a $c_n = c \neq 1$ para todo $n > 1$.

Esto contradice lo que hemos supuesto por el absurdo y concluye la demostración de 1).

**Demostración de 2)**

Si $c = 1$ entonces de (2.1.4) con $n = 1$ se obtiene
$$c_2 \underbrace{(1 - c)}_{0} = b c + c = (b+1)c = b + 1.$$

Es decir $b = -1$.

De esta manera (2.1.4) se lee
$$c_{n+1}(1 - c_n) = -c_n + 1.$$

Es decir
$$(1 - c_n)(1 - c_{n+1}) = 0. \qquad \square$$

Si $c \neq 1$ entonces no todos los valores de $b$ y $c$ producen aplicaciones de torcimiento.

Por ejemplo si $c \neq 1$ y $1 - c = c(b+1)$ en (2.1.4) para $n = 1$ obtenemos $c_2(1 - c_1) = b c_1 + c$, es decir $c_2(1 - c) = b c + c$. Luego $c_2 = \dfrac{(b+1)c}{1-c} = 1$, lo cual es una contradicción con lo afirmado en el Lema 2.1.5.

**Nota**: Entre las hipótesis de la Proposicion 2.1.4 hemos pedido que $c_k := M_{k,k+1} \neq 1$ para todo $k$. Del item a) del lema anterior se sigue que para ello basta pedir que $c = c_1 \neq 1$.

## 2.2.   El caso $yx = x^2 + bxy + cy^2, \quad c \neq 1$

Sabemos por (2.1.4) que $c_{n+1}(1 - c_n) = bc_n + c$, para $n \geq 0$.
Si $c \neq 1$ obtenemos la fórmula de recurrencia

$$c_{n+1} = \frac{bc_n + c}{1 - c_n}, \text{ para } n \geq 0.$$

Las primeras fórmulas para $c_n$ son

$\boxed{n=0}$   $c_1 = \dfrac{bc_0 + c}{1 - c_0} = \dfrac{c}{1} = c,$

$\boxed{n=1}$   $c_2 = \dfrac{bc_1 + c}{1 - c_1} = \dfrac{bc + c}{1 - c} = c\left(\dfrac{1+b}{1-c}\right),$

$\boxed{n=2}$   $c_3 = \dfrac{bc_2 + c}{1 - c_2} = \dfrac{b[c\frac{1+b}{1-c}] + c}{1 - [c\frac{1+b}{1-c}]} = \dfrac{bc(1+b) + c(1-c)}{1 - c - c(1+b)} = \dfrac{c(b^2 + b + 1 - c)}{1 - 2c - bc},$

$\boxed{n=3}$   $c_4 = \dfrac{bc_3 + c}{1 - c_3} = \dfrac{b[\frac{c(b^2+b+1-c)}{1-2c-bc}] + c}{1 - [\frac{c(b^2+b+1-c)}{1-2c-bc}]} = \dfrac{bc(b^2 + b + 1 - c) + c - 2c^2 - bc^2}{1 - 2c - bc - c(b^2 + b + 1 - c)},$

$$c_4 = \frac{c(1+b)(1 + b^2 - 2c)}{1 - 3c - 2bc - cb^2 + c^2}.$$

En general podemos escribir $\quad c_n = c\dfrac{P_n}{Q_n}$, para $n \geq 1$.

Notemos que $P_1 = Q_1 = 1$ y $P_n, Q_n$ son polinomios en $b$ y $c$ para $n \geq 2$.
Reemplazando en (2.1.4) se obtiene

$$c\frac{P_{n+1}}{Q_{n+1}}\left(1 - c\frac{P_n}{Q_n}\right) = bc\frac{P_n}{Q_n} + c, \text{ para } n \geq 1.$$

Es decir
$$\frac{cP_{n+1}(Q_n - cP_n)}{Q_n Q_{n+1}} = \frac{bcP_n + cQ_n}{Q_n}, \text{ para } n \geq 1.$$

Como $Q_n \neq 0$ para todo $n \geq 1$ entonces

$$cP_{n+1}(Q_n - cP_n) = c(bP_n + Q_n)Q_{n+1}, \text{ para } n \geq 1.$$

Si ponemos $Q_1 = 1, P_1 = 1$ y para cada $n \geq 1$ definimos recursivamente

$$\begin{cases} P_{n+1} &= bP_n + Q_n \\ Q_{n+1} &= -cP_n + Q_n \end{cases}, \qquad (2.2.5)$$

entonces los $c_n$ definidos por $c_n = c\dfrac{P_n}{Q_n}$ satisfacen (2.1.4).
Dados $b$ y $c$, los valores de $P_n(b, c)$ y $Q_n(b, c)$ están bien definidos aún cuando algún $c_n = 1$, y esto último sucede si y sólo si $Q_{n+1} = Q_{n+1}(b, c) = 0$.

**Proposición 2.2.1.** *Sea $K$ un cuerpo y sean $b, c \in K$ con $c \neq 1$. Si existe una matriz infinita $M$ de modo que $M_{1,0} = 1$, $M_{1,1} = b$ y $M_{1,2} = c$, y tal que $M$ determina una aplicación de torcimiento y además cumple $M_{n,n+1} \neq 1$ para todo $n$, entonces $Q_n(b,c) \neq 0$ para todo $n \in \mathbb{N}$, donde los polinomios $P_n, Q_n \in K[b,c]$ están definidos por $P_1 = 1$, $Q_1 = 1$ y las reglas recursivas* (2.2.5).

**Demostración.** Por la Proposición 2.1.4 y la discusión anterior. □

Para demostrar el recíproco del Proposición 2.2.1, consideramos la siguiente valuación en el álgebra $L(K[x]^{\mathbb{N}_0})$.

**Definición 2.2.2.** *Sea $M = (m_{i,j})_{i,j \in \mathbb{N}_0}$, la valuación de $M$ está dada por*

$$w(M) := \inf\{i - j, m_{i,j} \neq 0\}.$$

*Además convenimos: $w(0) = +\infty$.*
*Al valor $w(M)$ se le llama el peso de $M$.*

Por ejemplo $w(E_{1,0}) = 1$ y $w(E_{0,1}) = -1 = w(E_{1,0} + E_{0,1})$.
Notemos que para alguna $M$ se podría tener $w(M) = -\infty$.
En nuestro trabajo esto no ocurre, de manera que en lo que sigue consideraremos que $w(M) \neq -\infty$.

**Proposición 2.2.3.** *Sean $M, N \in L(K[x]^{\mathbb{N}_0})$. Se cumple:*

1. $w(M + N) \geq \inf\{w(M), w(N)\}$.

2. *Si $w(M) \neq w(N)$, entonces $w(M + N) = \inf\{w(M), w(N)\}$.*

3. $w(MN) \geq w(M) + w(N)$.

**Demostración.**
**Demostración de 1.**
Supongamos que $M + N = 0$.
En este caso $w(M + N) = +\infty$ y se cumple trivialmente
$w(M + N) \geq \inf\{w(M), w(N)\}$.
Supongamos ahora que $M + N \neq 0$.
Entonces $w(M + N) = i_0 - j_0$ para algunos $i_0, j_0$ tales que $m_{i_0, j_0} + n_{i_0, j_0} \neq 0$.
Entonces $m_{i_0, j_0} \neq 0$ y por lo tanto

$$w(M) = \inf\{i - j, m_{i,j} \neq 0\} \leq i_0 - j_0 = w(M + N),$$

o sino $n_{i_0, j_0} \neq 0$ y en este caso

$$w(N) = \inf\{i - j, n_{i,j} \neq 0\} \leq i_0 - j_0 = w(M + N).$$

En ambos casos se tiene que ínf$\{w(M), w(N)\} \leq w(M + N)$.

**Demostración de 2.**

Si $w(M) \neq w(N)$ entonces $w(M) > w(N)$ ó $w(M) < w(N)$.

Supongamos que $w(M) > w(N)$.

Entonces $w(N) = i_0 - j_0$ para algunos $i_0, j_0$ tales que $n_{i_0, j_0} \neq 0$.

Pero $m_{i_0, j_0} = 0$, porque sino $m_{i_0, j_0} \neq 0$ implica

$$w(N) = i_0 - j_0 \geq \text{ínf}\{i - j, m_{i,j} \neq 0\} = w(M),$$

lo cual contradice lo que hemos supuesto.

Entonces

$$(M + N)_{i_0, j_0} = n_{i_0, j_0} + m_{i_0, j_0} = n_{i_0, j_0} \neq 0,$$

y entonces

$$w(N + M) = \text{ínf}\{i - j, (n_{i,j} + m_{i,j}) \neq 0\} \leq i_0 - j_0 = w(N) = \text{ínf}\{w(N), w(M)\}.$$

Pero por el ítem 1. sabemos que ínf$\{w(M), w(N)\} \leq w(M + N)$, lo cual concluye la demostración de 2. cuando $w(M) > w(N)$. Si $w(M) < w(N)$, la demostración es simétrica.

**Demostración de 3.**

Supongamos que $MN = 0$.

En este caso $w(MN) = +\infty$ y se cumple trivialmente
$w(MN) \geq w(M) + w(N)$.

Supongamos ahora que $MN \neq 0$.

Entonces $w(MN) = i_0 - j_0$ para algunos $i_0, j_0$ tales que $(MN)_{i_0, j_0} \neq 0$.

Pero como $M \in L(K[x]^{\mathbb{N}_0})$, existe $k_0$ tal que $M_{i_0, k} = 0$ para $k > k_0$, y entonces

$$(MN)_{i_0, j_0} = \sum_{k=0}^{k_0} m_{i_0, k} n_{k, j_0}.$$

Como la suma finita es no nula, existe un $k$ tal que $m_{i_0, k} n_{k, j_0} \neq 0$, y entonces $m_{i_0, k} \neq 0$, lo que implica

$$w(M) = \text{ínf}\{i - j, m_{i,j} \neq 0\} \leq i_0 - k$$

y también tenemos $n_{k, j_0} \neq 0$, lo que implica

$$w(N) = \text{ínf}\{i - j, m_{i,j} \neq 0\} \leq k - j_0.$$

Entonces

$$w(MN) = i_0 - j_0 = (i_0 - k) + (k - j_0) \geq w(M) + w(N),$$

lo cual concluye la demostración. $\square$

**Definición 2.2.4.** *Sea $M = (m_{i,j})_{i,j \in \mathbb{N}_0}$.*

1. *Decimos que $M$ es homogénea si $m_{i,j} = 0$ cuando $i - j \neq w(M)$.*

2. *La componente homogénea de peso $k$ de $M$ es la matriz, denotada por $M^{(k)}$, dada por*

$$(M^{(k)})_{i,j} = \begin{cases} M_{i,j} & si\ i-j = k \\ 0 & si\ i-j \neq k \end{cases}.$$

Por ejemplo, consideremos las matrices $Y$, $Z$ dadas por $Y_{i,j} = \delta_{i+1,j}$, $Z_{i,j} = \delta_{i,j+1}$. Entonces ambas son homogéneas con $w(Y) = -1$, $w(Z) = 1$.

Consideremos el subálgebra $\mathscr{R} \subset L(K[x]^{\mathbb{N}_0})$ que consiste de las matrices homogéneas de peso cero (matrices diagonales).

Si $k > 0$ entonces la matriz $N_k = (N_k)_{i,j} \in \mathscr{R} \subset L(K[x]^{\mathbb{N}_0})$ dada por

$$(N_k)_{i,j} = \begin{cases} M_{k+i,i}, & \text{si } i = j \\ 0, & \text{si } i \neq j \end{cases}$$

satisface

$$M^{(k)} = Z^k N_k. \tag{2.2.6}$$

En efecto

$$(Z^k N_k)_{i,j} = \sum_p (Z^k)_{i,p} (N_k)_{p,j} = \sum_p \delta_{i,p+k} (N_k)_{p,j}.$$

El único término que sobrevive en la sumatoria es aquel que se obtiene cuando $p + k = i$.

Entonces

$$(Z^k N_k)_{i,j} = (N_k)_{i-k,j} = \begin{cases} M_{i,j} & \text{si } i - k = j \\ 0 & \text{si } i - k \neq j \end{cases}.$$

Pero

$$(M^{(k)})_{i,j} = \begin{cases} M_{i,j} & \text{si } i - j = k \\ 0 & \text{si } i - j \neq k \end{cases}.$$

Como $i - k = j \Leftrightarrow i - j = k$, entonces queda demostrado lo deseado.

Por otro lado, si $k$ es negativo, entonces la matriz $N_k = (N_k)_{i,j} \in \mathscr{R} \subset L(K[x]^{\mathbb{N}_0})$ dada por

$$(N_k)_{i,j} = \begin{cases} M_{i,j-k}, & \text{si } i = j \\ 0, & \text{si } i \neq j \end{cases}$$

satisface
$$M^{(k)} = N_k Y^{-k}. \tag{2.2.7}$$

En efecto
$$(N_k Y^{-k})_{i,j} = (N_k)_{i,j+k} = \begin{cases} M_{i,j} & \text{si } i = j+k \\ 0 & \text{si } i \neq j+k \end{cases}.$$

Pero
$$(M^{(k)})_{i,j} = \begin{cases} M_{i,j} & \text{si } i-j = k \\ 0 & \text{si } i-j \neq k \end{cases}.$$

Como $i = j+k \Leftrightarrow i-j = k$, entonces queda demostrado lo deseado.

Finalmente definimos $N_0 = M^{(0)}$, la componente homogénea de $M$ de peso 0.

Se sigue que para $M$ con $w(M) > -\infty$ tenemos una descomposición

$$M = \sum_{j=w(M)}^{0} M^{(j)} + \sum_{k>0} M^{(k)}.$$

Es decir
$$M = \sum_{j=w(M)}^{0} N_j Y^{-j} + \sum_{k>0} Z^k N_k \tag{2.2.8}$$

para algunos $N_j, N_k \in \mathscr{R}$, donde la suma infinita converge en la topología $Z$-ádica. Notemos que si $w(M) > 0$, entonces la primera suma es vacía.

Dados un cuerpo $K$ y $b, c \in K$ sean $P_n(b,c)$ y $Q_n(b,c)$ definidas como en (2.2.5) con $P_1(b,c) = Q_1(b,c) = 1$.

**Proposición 2.2.5.** *Sea $K$ un cuerpo y sean $b, c \in K$ con $Q_n(b,c) \neq 0$ para todo $n \in \mathbb{N}$. Entonces existe una única matriz $M \in L(K[x]^{\mathbb{N}_0})$ con $w(M) \geq -1$, tal que $M_{0,*} = (1, 0, 0, \dots)$, $M_{1,0} = 1$, $M_{1,1} = b$, $M_{1,2} = c$,*

$$YM = M^2 + bMY + cY^2. \tag{2.2.9}$$

**Demostración.** De acuerdo a (2.2.8) podemos escribir

$$M = CY + B_0 + \sum_{j \geq 1} Z^j B_j,$$

donde $C = N_{-1}$, $B_0 = M^{(0)}$ y $B_j = N_j$.

Es decir si denotamos $M_{i,i+1} = c_i$ y $M_{i,i} = b_i$ entonces las matrices $C$ y $B_0$ están dadas respectivamente por

$$C = \begin{pmatrix} c_0 & 0 & 0 & 0 & \cdots \\ 0 & c_1 & 0 & 0 & \cdots \\ 0 & 0 & c_2 & 0 & \cdots \\ \vdots & \vdots & \vdots & \ddots & \vdots \end{pmatrix}, \qquad B_0 = \begin{pmatrix} b_0 & 0 & 0 & 0 & \cdots \\ 0 & b_1 & 0 & 0 & \cdots \\ 0 & 0 & b_2 & 0 & \cdots \\ \vdots & \vdots & \vdots & \ddots & \vdots \end{pmatrix}.$$

Por otro lado si $j \geq 1$, entonces

$$B_j = \begin{pmatrix} M_{j,0} & 0 & 0 & 0 & \cdots \\ 0 & M_{j+1,1} & 0 & 0 & \cdots \\ 0 & 0 & M_{j+2,2} & 0 & \cdots \\ \vdots & \vdots & \vdots & \ddots & \vdots \end{pmatrix}.$$

Luego la igualdad (2.2.9) es válida si y sólo si es válida en cada componente homogénea, es decir, si y sólo si

$$(YM)^{(k)} = (M^2 + bMY + cY^2)^{(k)} \qquad (2.2.10)$$

para todo $k \geq -2$.

En el resto de la demostración construiremos recursivamente $M^{(-1)}$, $M^{(0)}$, $M^{(1)}$, $M^{(2)}$,..., $M^{(j)}$, tales que (2.2.10) es válida para $k = -2, -1, 0, 1, \ldots, j-1$. Esto produce una construcción inductiva de la única $M$ talque (2.2.9) es válida.

Antes notemos que una matriz $A \in \mathcal{R}$, se puede escribir como $A = (a_0, a_1, a_2, \ldots)$. Definimos entonces $S : \mathcal{R} \longrightarrow \mathcal{R}$ el shift dado por

$$SA := YAZ = (a_1, a_2, \ldots)$$

y $T : \mathcal{R} \longrightarrow \mathcal{R}$ la inversa a derecha de $S$, es decir:

$$TA := ZAY = (0, a_0, a_1, a_2, \ldots).$$

Las matrices que hemos definido satisfacen las siguientes propiedades, para $C \in \mathcal{R}$:

(1) $YZ = 1$,

(2) $(SC)Y = YC$,

(3) $Z^j(S^j C) = CZ^j$, $\quad j \geq 1$,

(4) $Z^j(S^j(B_k)) = B_k Z^j$, $\quad j \geq 1$.

Notemos que

$$\begin{aligned} YM &= Y\left(CY + B_0 + \sum_{j \geq 1} Z^j B_j\right) \\ &= YCY + YB_0 + \sum_{j \geq 1} YZ^j B_j \\ &= YCY + YB_0 + \sum_{j \geq 0} YZ^{j+1} B_{j+1}. \end{aligned}$$

Luego

$$(YM)^{(-2)} = YCY = (SC)Y^2,$$
$$(YM)^{(-1)} = YB_0 = (SB_0)Y,$$
$$(YM)^{(j)} = YZ^{j+1}B_{j+1} = (YZ)Z^j B_{j+1} = Z^j B_{j+1}, \quad j \geq 0.$$

Notemos también que

$$\begin{aligned} MY &= \left(CY + B_0 + \sum_{j\geq 1} Z^j B_j\right)Y \\ &= CY^2 + B_0 Y + \sum_{j\geq 1} Z^j B_j Y \\ &= CY^2 + B_0 Y + \sum_{j\geq 0} Z^{j+1} B_{j+1} Y. \end{aligned}$$

Luego

$$(MY)^{(-2)} = CY^2,$$
$$(MY)^{(-1)} = B_0 Y,$$
$$(MY)^{(j)} = Z^{j+1} B_{j+1} Y = Z^j (ZB_{j+1} Y) = Z^j T B_{j+1}, \quad j \geq 0.$$

Finalmente

$$M^2 = \left(\sum_{k\geq -1} M^{(k)}\right)\left(\sum_{i\geq -1} M^{(i)}\right) = \sum_{k\geq -1}\sum_{i\geq -1} M^{(k)} M^{(i)},$$

donde $M^{(k)} M^{(i)}$ es una matriz homogénea de peso $k+i$.

Haciendo $j = k+i$ se tiene para cada $j$

$$(M^2)^{(j)} = \sum_{k=-1}^{j+1} M^{(k)} M^{(j-k)}.$$

Luego

$$(M^2)^{(-2)} = \sum_{k=-1}^{-1} M^{(k)} M^{(-2-k)} = M^{(-1)} M^{(-1)} = CYCY = C(YC)Y = C(SC)Y^2,$$

$$(M^2)^{(-1)} = \sum_{k=-1}^{0} M^{(k)} M^{(-1-k)} = M^{(-1)} M^{(0)} + M^{(0)} M^{(-1)} = CYB_0 + B_0 CY$$
$$= C(SB_0)Y + B_0 CY,$$

$$(M^2)^{(j)} = \sum_{k=-1}^{j+1} M^{(k)} M^{(j-k)} = M^{(-1)} M^{(j+1)} + M^{(j+1)} M^{(-1)} + \sum_{k=0}^{j} M^{(k)} M^{(j-k)}, \quad j \geq 0.$$

Es decir

$$(M^2)^{(j)} = CYZ^{j+1} B_{j+1} + Z^{j+1} B_{j+1} CY + \sum_{k=0}^{j} Z^k B_k Z^{j-k} B_{j-k}, \quad j \geq 0.$$

Pero
$$CYZ^{j+1}B_{j+1} = CZ^j B_{j+1} = Z^j(S^j C)B_{j+1},$$

por las propiedades (1) y (3),

$$Z^{j+1}B_{j+1}CY = Z^j[Z(B_{j+1}C)Y]$$
$$= Z^j T(B_{j+1}C),$$
$$\sum_{k=0}^{j} Z^k B_k Z^{j-k} B_{j-k} = \sum_{k=0}^{j} Z^k Z^{j-k}(S^{j-k}B_k)B_{j-k}$$
$$= \sum_{k=0}^{j} Z^j(S^{j-k}B_k)B_{j-k},$$

por la propiedad (4).

Para $j = -2$ la igualdad

$$(YM)^{(j)} = (M^2)^{(j)} + b(MY)^{(j)} + c(Y^2)^{(j)},$$

se lee

$$(YM)^{(-2)} = (M^2)^{(-2)} + b(MY)^{(-2)} + c(Y^2)^{(-2)},$$
$$(SC)Y^2 = C(SC)Y^2 + bCY^2 + c(Y^2).$$

Multiplicando a la derecha por $Z^2$ se obtiene

$$SC = C(SC) + bC + c\mathbf{1},$$

lo cual nos da la igualdad (2.1.4)

$$(1 - c_n)c_{n+1} = b c_n + c, \quad n \geq 0.$$

Por lo tanto la ecuación matricial (2.2.10) para $k = -2$ es equivalente al conjunto de ecuaciones (2.1.4) para todo $n$. Aquí es donde se usa la condición $Q_n(b, c) \neq 0$, porque sino no se pueden construir los $c_n$ satisfaciendo (2.1.4).

Para $j = -1$ la igualdad se lee

$$(YM)^{(-1)} = (M^2)^{(-1)} + b(MY)^{(-1)} + c(Y^2)^{(-1)},$$
$$(SB_0)Y = C(SB_0)Y + B_0 CY + b B_0 Y.$$

Dado que la multiplicación a la derecha por $Y$ es inyectiva se obtiene

$$SB_0 = C(SB_0) + B_0 C + b B_0.$$

Así tenemos una fórmula recursiva para $B_0$

$$(1 - C)SB_0 = B_0 C + b B_0.$$

Es decir
$$(1-c_n)(B_0)_{n+1} = (B_0)_n(c_n+b).$$

$\boxed{n=0}$

$$(1-c_0)(B_0)_1 = (B_0)_0(c_0+b),$$
$$(B_0)_1 = b(B_0)_0 = bM_{0,0} = b.$$

$\boxed{n=1}$

$$(1-c_1)(B_0)_2 = (B_0)_1(c_1+b),$$
$$(B_0)_2 = \frac{b(c+b)}{1-c}.$$

$\boxed{n=2}$

$$(1-c_2)(B_0)_3 = (B_0)_2(c_2+b),$$
$$(B_0)_3 = \frac{b(c+b)(c_2+b)}{(1-c)(1-c_2)}.$$

Dado que $(1-c_n) \neq 0$ para todo $n \geq 0$ y ya tenemos $(B_0)_0 = 1$, $(B_0)_1 = b$, esta fórmula determina una única $B_0$ tal que la igualdad (2.2.10) correspondiente al peso $k = -1$ es satisfecha.

Para $j \geq 0$ la igualdad (2.2.10) se lee

$$Z^j B_{j+1} = Z^j(S^j C)B_{j+1} + Z^j T(B_{j+1}C) + Z^j \sum_{k=0}^{j}(S^{j-k}B_k)B_{j-k} + bZ^j T(B_{j+1}).$$

Dado que la multiplicación a la izquierda por $Z^j$ es inyectiva se obtiene

$$\begin{aligned}(1-S^jC)B_{j+1} &= T(B_{j+1}C) + \sum_{k=0}^{j}(S^{j-k}B_k)B_{j-k} + bT(B_{j+1}) \\ &= T(B_{j+1})(TC+b1) + \sum_{k=0}^{j}(S^{j-k}B_k)B_{j-k}.\end{aligned}$$

$\boxed{j=0}$

$$\begin{aligned}(1-C)B_1 &= T(B_1)(TC+b1) + \sum_{k=0}^{0}(S^{-k}B_k)B_{-k}, \\ &= T(B_1)(TC+b1) + B_0 B_0.\end{aligned}$$

$\boxed{j=1}$

$$\begin{aligned}(1-SC)B_2 &= T(B_2)(TC+b1) + \sum_{k=0}^{1}(S^{1-k}B_k)B_{1-k}, \\ &= T(B_2)(TC+b1) + (SB_0)B_1 + B_1 B_0.\end{aligned}$$

Supongamos que hemos construido $C$ y $B_i$ para $i = 0, 1, \cdots, j$ tal que (2.2.10) es satisfecha para $k = -2, -1, \cdots, j-1$. Entonces ponemos

$$R = \sum_{i=0}^{j}(S^{j-i}B_i)B_{j-i},$$

y notemos que $R$ depende solamente de $B_i$ para $i = 0, 1, \cdots, j$ y obtenemos una fórmula recursiva

$$(1 - c_{n+j})(B_{j+1})_n = (B_{j+1})_{n-1}(c_{n-1} + b) + R_n,$$

la cual produce una única $B_{j+1}$ tal que (2.2.10) es satisfecha para $k = -2, -1, \cdots, j-1$.

Note que la fórmula es válida para $n = 0$ haciendo $(B_{j+1})_{-1} = c_{-1} = 0$.

Esto demuestra que existe una única matriz

$$M = CY + B_0 + \sum_{j \geq 1} Z^j B_j,$$

que satisface (2.2.9). $\square$

**Lema 2.2.6.** *Sea la primera fila de $M \in L(K[x]^{\mathbb{N}_0})$ dada por $M_{0,*} = E_0 = (1, 0, 0, 0, \ldots)$. Si $M$ satisface*

$$Y^k M = \sum_{i=0}^{k+1} a_i M^{k+1-i} Y^i,$$

*entonces $M_{k,j} = a_j$ para $j = 0, \ldots, k+1$.*

**Demostración.** Notemos que $(M^r)_{0,*} = E_0$ para todo $r$, y así $(M^r Y^i)_{0,j} = M_{0,j-i} = \delta_{ij}$. Por lo tanto

$$M_{k,j} = (Y^k M)_{0,j} = \sum_{i=0}^{k+1} a_i (M^{k+1-i} Y^i)_{0,j} = \sum_{i=0}^{k+1} a_i \delta_{ij} = a_j,$$

como se quería. $\square$

**Teorema 2.2.7.** *Sea $K$ un cuerpo y sean $b, c \in K$ con $c \neq 1$.*
*Supongamos que $Q_n(b, c) \neq 0$ para todo $n \in \mathbb{N}$. Entonces existe una única matriz $M$, con $M_{1,0} = 1, M_{1,1} = b$ y $M_{1,2} = c$, que determina una (única) aplicación de torcimiento via Proposición 1.4.1. Entonces $b$ y $c$ determinan una única matriz $M$, con $M_{1,0} = 1$, $M_{1,1} = b$ y $M_{1,2} = c$, y que a su vez determina una única aplicación de torcimiento via Proposición 1.4.1.*

**Demostración.** Usaremos la Proposición 1.4.1. Para esto demostraremos que la matriz $M$ construída en la Proposición 2.2.5 satisface la igualdad

$$Y^k M = \sum_{j=0}^{k+1} M_{k,j} M^{k+1-j} Y^j, \quad k \geq 0.$$

Para $k=0$ esto es evidente, y por el Lema 2.2.6 obtenemos $M_{1,*} = E_0 + b E_1 + c E_2$, y por lo tanto, por (2.2.9), la igualdad es válida para $k=1$.

Asumamos, por inducción, que la igualdad es válida para todo $k < k_0$. Entonces por el Lema 2.1.2 y el hecho de que $M_{k_0, k_0+1} \neq 1$ se sigue que

$$Y^{k_0} M = \sum_{s=0}^{k_0+1} \frac{\overline{M}_{k_0,s}}{1 - M_{k_0,k_0+1}} M^{k_0+1-s} Y^s.$$

Entonces por el Lema 2.2.6 tenemos $\dfrac{\overline{M}_{k_0,s}}{1 - M_{k_0,k_0+1}} = M_{k_0,s}$ para $s = 0, 1, \cdots, k_0 + 1$, lo cual produce la igualdad deseada para $k = k_0$ y se completa el paso inductivo.

Finalmente la Proposición 1.4.1 asegura que $M$ representa una aplicación de torcimiento que además es única por la Proposición 2.1.4. □

## 2.3. Raíces de $Q_n$

En vista del Teorema 2.2.7 queremos analizar los polinomios $Q_n$ y sus raíces. En particular estamos interesados en responder la siguiente pregunta: Dado un par $(b,c) \in K^2$, ¿existe un $n \in \mathbb{N}$ tal que $Q_n(b,c) = 0$?. Si la respuesta es no entonces $(b,c)$ define una única aplicación de torcimiento vía el teorema anterior. Si la respuesta es si, y además $(b,c) \neq (-1,1)$, entonces no existe aplicación de torcimiento para tal $(b,c)$.

Para cada par $(b,c)$ fijo podemos escribir las fórmulas recursivas (2.2.5) en forma matricial

$$\begin{pmatrix} P_{n+1} \\ Q_{n+1} \end{pmatrix} = \begin{pmatrix} b & 1 \\ -c & 1 \end{pmatrix} \begin{pmatrix} P_n \\ Q_n \end{pmatrix} \qquad n \geq 1.$$

Haciendo $P_0 = 0$, $Q_0 = 1$, la fórmula anterior también vale para $n \geq 0$.

$\boxed{n=0}$

$$\begin{pmatrix} P_1 \\ Q_1 \end{pmatrix} = \begin{pmatrix} b & 1 \\ -c & 1 \end{pmatrix} \begin{pmatrix} 0 \\ 1 \end{pmatrix} = \begin{pmatrix} 1 \\ 1 \end{pmatrix}.$$

$\boxed{n=1}$

$$\begin{pmatrix} P_2 \\ Q_2 \end{pmatrix} = \begin{pmatrix} b & 1 \\ -c & 1 \end{pmatrix} \begin{pmatrix} P_1 \\ Q_1 \end{pmatrix} = \begin{pmatrix} b & 1 \\ -c & 1 \end{pmatrix}^2 \begin{pmatrix} 0 \\ 1 \end{pmatrix}.$$

En general
$$\begin{pmatrix} P_n \\ Q_n \end{pmatrix} = \begin{pmatrix} b & 1 \\ -c & 1 \end{pmatrix}^{n-1} \begin{pmatrix} P_1 \\ Q_1 \end{pmatrix} = \begin{pmatrix} b & 1 \\ -c & 1 \end{pmatrix}^n \begin{pmatrix} 0 \\ 1 \end{pmatrix}, \text{ para } n \geq 0.$$

Si los autovalores $\lambda_1, \lambda_2$ de $D = \begin{pmatrix} b & 1 \\ -c & 1 \end{pmatrix}$ son diferentes, entonces existe una matriz inversible $T$ tal que
$$D = T \begin{pmatrix} \lambda_1 & 0 \\ 0 & \lambda_2 \end{pmatrix} T^{-1}.$$

Asumamos que $T = \begin{pmatrix} u_1 & u_2 \\ v_1 & v_2 \end{pmatrix}$ donde $\begin{pmatrix} u_1 \\ v_1 \end{pmatrix}, \begin{pmatrix} u_2 \\ v_2 \end{pmatrix}$ son los autovectores correspondientes a $\lambda_1, \lambda_2$ respectivamente.

Luego
$$\begin{pmatrix} 0 \\ 1 \end{pmatrix} = \alpha_1 \begin{pmatrix} u_1 \\ v_1 \end{pmatrix} + \alpha_2 \begin{pmatrix} u_2 \\ v_2 \end{pmatrix}.$$

Por lo tanto
$$D^n \begin{pmatrix} 0 \\ 1 \end{pmatrix} = \alpha_1 D^n \begin{pmatrix} u_1 \\ v_1 \end{pmatrix} + \alpha_2 D^n \begin{pmatrix} u_2 \\ v_2 \end{pmatrix} = \alpha_1 \lambda_1^n \begin{pmatrix} u_1 \\ v_1 \end{pmatrix} + \alpha_2 \lambda_2^n \begin{pmatrix} u_2 \\ v_2 \end{pmatrix}.$$

Por lo tanto
$$\begin{pmatrix} P_n \\ Q_n \end{pmatrix} = \begin{pmatrix} \alpha_1 u_1 \lambda_1^n + \alpha_2 u_2 \lambda_2^n \\ \alpha_1 v_1 \lambda_1^n + \alpha_2 v_2 \lambda_2^n \end{pmatrix}.$$

Es decir $P_n = \alpha_1 u_1 \lambda_1^n + \alpha_2 u_2 \lambda_2^n, \quad Q_n = \alpha_1 v_1 \lambda_1^n + \alpha_2 v_2 \lambda_2^n$.

Haciendo $r_1 = \alpha_1 v_1, \quad r_2 = \alpha_2 v_2$ obtenemos
$$Q_n = r_1 \lambda_1^n + r_2 \lambda_2^n.$$

En particular
$$Q_0 = r_1 + r_2 = 1 \quad (\alpha),$$
$$Q_1 = r_1 \lambda_1 + r_2 \lambda_2 = 1 \quad (\beta).$$

De $(\alpha)$ se obtiene $r_1 = 1 - r_2$,

reemplazando en $(\beta)$ se obtiene $\boxed{r_2 = \dfrac{1-\lambda_1}{\lambda_2 - \lambda_1}}$.

Similarmente, de $(\alpha)$ se obtiene $r_2 = 1 - r_1$,

reemplazando en $(\beta)$ se obtiene $\boxed{r_1 = \dfrac{1-\lambda_2}{\lambda_1 - \lambda_2}}$.

Luego para todo $n \geq 0$
$$Q_n = 0 \Leftrightarrow r_1 \lambda_1^n + r_2 \lambda_2^n = 0 \Leftrightarrow r_1 \left(\frac{\lambda_1}{\lambda_2}\right)^n + r_2 = 0 \Leftrightarrow \boxed{\left(\frac{\lambda_1}{\lambda_2}\right)^n = -\frac{r_2}{r_1}} \qquad (2.3.11)$$

siempre que $\lambda_2 \neq 0, r_1 \neq 0$.

Esta condición es más fácil de verificar que el número infinito de evaluaciones $Q_n(b,c)$. Por ejemplo si $K \subset \mathbb{C}$ y $r_1 = |r_1|e^{i\theta_1}$, $r_2 = |r_2|e^{i\theta_2}$ con $|r_1| \neq |r_2|$, entonces

$$\frac{r_2}{r_1} = \frac{|r_2|}{|r_1|} e^{i(\theta_1-\theta_2)}.$$

Luego para verificar que se se cumple la igualdad $Q_n = 0$ basta hallar, usando logaritmos reales, $n \geq 0$ que satisface

$$\left|\frac{\lambda_1}{\lambda_2}\right|^n = \left|\frac{r_2}{r_1}\right| \neq 1.$$

Calculemos los autovalores de $D$.

$$\begin{vmatrix} b-\lambda & 1 \\ -c & 1-\lambda \end{vmatrix} = 0,$$

$$(b-\lambda)(1-\lambda) + c = 0,$$

$$\lambda^2 - (b+1)\lambda + (b+c) = 0,$$

$$\lambda_{1,2} = \frac{1}{2}(b+1 \pm \sqrt{(b-1)^2 - 4c}).$$

El argumento de arriba nos demuestra el siguiente lema:

**Lema 2.3.1.** *Si ninguno de los valores $\lambda_1, \lambda_2, \lambda_1-1, \lambda_2-1, \lambda_1-\lambda_2$ se anula, es decir, si*

$$0 \notin \{\lambda_1, \lambda_2, \lambda_1-1, \lambda_2-1, \lambda_1-\lambda_2\}, \tag{2.3.12}$$

*entonces $Q_n = 0$ si y solo si $\left(\dfrac{\lambda_1}{\lambda_2}\right)^n = -\dfrac{r_2}{r_1}$ si y solo si $\left(\dfrac{\lambda_1}{\lambda_2}\right)^n = \dfrac{1-\lambda_1}{1-\lambda_2}$.*

Notemos que si $b,c \in \mathbb{R}$ entonces $4c < (b-1)^2$ implica $\dfrac{\lambda_1}{\lambda_2} \neq 1$. Si además $b \neq -1$, entonces $\dfrac{|\lambda_1|}{|\lambda_2|} \neq 1$ y así se puede determinar de manera efectiva si se cumple $Q_n = 0$ para algún $n$.

## 2.4. Casos excepcionales

**Observación 2.4.1.** Dados un cuerpo $K$ y $b,c \in K$ sean $P_n(b,c)$ y $Q_n(b,c)$ definidas como en (2.2.5) con $P_1(b,c) = Q_1(b,c) = 1$.

Para verificar que una matriz infinita $M$ es la de torcimiento asociada a $(b,c)$ es suficiente verificar que $Q_n(b,c) \neq 0$ para todo $n$, $M_{0,*} = (1,0,0,\dots)$, $M_{1,0} = 1$, $M_{1,1} = b$, $M_{n,n+1} = c \dfrac{P_n(b,c)}{Q_n(b,c)}$, y que

$$YM = M^2 + bMY + cY^2.$$

En efecto, por el Teorema 2.2.7 existe una única matriz $M_1$ asociada al torcimiento, que necesariamente cumple esta fórmula. Por la Proposición 2.2.5 se tiene que $M_1 = M$.

Ahora damos cuenta detallada de cada caso excepcional en (2.3.12).

**Caso 1.** $\boxed{\lambda_1 = 0 \vee \lambda_2 = 0}$

En este caso se tiene (para $\lambda_i = 0$)

$$0 = \begin{vmatrix} b - \lambda_i & 1 \\ -c & 1 - \lambda_i \end{vmatrix} = \begin{vmatrix} b & 1 \\ -c & 1 \end{vmatrix} = b + c$$

Es decir $b = -c$ y por lo tanto se tiene

$$\boxed{YM = M^2 - cMY + cY^2, \quad c \neq 1}. \tag{2.4.13}$$

Se ve que en este caso $Q_n(-c, c) = (1-c)^{n-1} \neq 0$ para $n > 0$, y entonces por la Observación 2.4.1, para encontrar la matriz asociada al torcimiento, es suficiente encontrar una matriz que cumpla (2.4.13).

Afirmamos que

$$M = M_0 - cM_0Y + cY$$

donde $(M_0)_{i,j} = \delta_{0,j}$, cumple (2.4.13).

Es decir,

$$M = \begin{pmatrix} 1 & 0 & 0 & 0 & 0 & \cdots \\ 1 & -c & c & 0 & 0 & \cdots \\ 1 & -c & 0 & c & 0 & \cdots \\ 1 & -c & 0 & 0 & c & \cdots \\ \vdots & \vdots & \vdots & \vdots & \vdots & \ddots \end{pmatrix}$$

Para esto vemos que $M_0^2 = M_0$ y $YM_0 = M_0$, lo cual implica que

$$(M - Y)M_0 = (M_0 - cM_0Y + cY)M_0 - YM_0 = 0$$

Pero $M - cY = M_0(Id - cY)$, y entonces

$$0 = (M - Y)M_0(Id - cY) = (M - Y)(M - cY) = M^2 - YM - cMY + cY^2,$$

lo cual dice que $M$ cumple (2.4.13).

**Caso 2.** $\boxed{\lambda_1 = 1 \vee \lambda_2 = 1}$

En este caso se tiene (para $\lambda_i = 1$)

$$0 = \begin{vmatrix} b - \lambda_i & 1 \\ -c & 1 - \lambda_i \end{vmatrix} = \begin{vmatrix} b-1 & 1 \\ -c & 0 \end{vmatrix} = c.$$

Es decir $c = 0$ y $c_n = 0$ para todo $n \geq 1$.

Notemos que los autovalores de $\begin{pmatrix} b & 1 \\ 0 & 1 \end{pmatrix}$ son $\{b, 1\}$.

Entonces este caso se da cuando $c = 0$ y $b \in K$, y siempre se cumple que $Q_n(b,c) = 1 \neq 0$.

Por lo tanto se tiene $\boxed{YM = M^2 + bMY}$.

Afirmamos que en este caso la matriz $M$ definida por:

$$M_{k,j} = \begin{cases} \delta_{0,j} & \text{si} \quad k = 0 \\ b^j \left( \dfrac{\prod_{i=0}^{k-1} A_i}{\prod_{i=0}^{j-1} A_i} \right) & \text{si} \quad 0 \leq j \leq k \\ 0 & \text{si} \quad j > k \geq 1 \end{cases}$$

donde $A_i := 1 + b + b^2 + \cdots + b^i$ para $i \geq 0$ con $\prod_{i=0}^{-1} A_i = 1$, es la matriz asociada al torcimiento.

Notemos que

$$M = \begin{pmatrix} 1 & 0 & 0 & 0 & 0 & \cdots \\ 1 & b & 0 & 0 & 0 & \cdots \\ 1+b & b(1+b) & b^2 & 0 & 0 & \cdots \\ (1+b)(1+b+b^2) & b(1+b)(1+b+b^2) & b^2(1+b+b^2) & b^3 & 0 & \cdots \\ \vdots & \vdots & \vdots & \vdots & \vdots & \ddots \end{pmatrix}.$$

Por la Observación 2.4.1 sólo tenemos que demostrar que $M$ cumple

$$YM = M^2 + bMY.$$

Es suficiente demostrar

$$(YM)_{k,j} = (M^2)_{k,j} + b(MY)_{k,j}. \tag{2.4.14}$$

Demostraremos primero (2.4.14) para $j = 0$.

$$(YM)_{k,0} = M_{k+1,0} = b^0 \left( \frac{\prod_{i=0}^{k} A_i}{\prod_{i=0}^{-1} A_i} \right) = \prod_{i=0}^{k} A_i.$$

Demostraremos que el segundo miembro es igual a esta entrada.

$$\begin{aligned}
(M^2)_{k,0} &= M_{k,*} \cdot M_{*,0} \\
&= (M_{k,0}, M_{k,1}, \cdots, M_{k,k}, M_{k,k+1}, 0, \cdots) \cdot (M_{0,0}, M_{1,0}, \cdots, M_{k,0}, M_{k+1,0}, \cdots) \\
&= \left(b^0\left(\frac{\prod_{i=0}^{k-1} A_i}{\prod_{i=0}^{-1} A_i}\right), b\left(\frac{\prod_{i=0}^{k-1} A_i}{\prod_{i=0}^{0} A_i}\right), \cdots, b^k\left(\frac{\prod_{i=0}^{k-1} A_i}{\prod_{i=0}^{k-1} A_i}\right), 0, \cdots\right) \\
&\quad \cdot \left(1, b^0\left(\frac{\prod_{i=0}^{0} A_i}{\prod_{i=0}^{-1} A_i}\right), \cdots, b^0\left(\frac{\prod_{i=0}^{k-1} A_i}{\prod_{i=0}^{-1} A_i}\right), b^0\left(\frac{\prod_{i=0}^{k} A_i}{\prod_{i=0}^{-1} A_i}\right), \cdots\right) \\
\\
&= b^0\left(\frac{\prod_{i=0}^{k-1} A_i}{\prod_{i=0}^{-1} A_i}\right) + b\left(\frac{\prod_{i=0}^{k-1} A_i}{\prod_{i=0}^{0} A_i}\right) b^0\left(\frac{\prod_{i=0}^{0} A_i}{\prod_{i=0}^{-1} A_i}\right) + \cdots \\
&\quad + b^k\left(\frac{\prod_{i=0}^{k-1} A_i}{\prod_{i=0}^{k-1} A_i}\right) b^0\left(\frac{\prod_{i=0}^{k-1} A_i}{\prod_{i=0}^{-1} A_i}\right) \\
&= \left(\frac{\prod_{i=0}^{k-1} A_i}{\prod_{i=0}^{-1} A_i}\right)(1 + b + \cdots + b^k) = \prod_{i=0}^{k} A_i.
\end{aligned}$$

Ahora como $(MY)_{k,0} = M_{k,-1} = 0$, reemplazando se tiene

$$(M^2)_{k,0} + b(MY)_{k,0} = \prod_{i=0}^{k} A_i.$$

Ahora demostraremos (2.4.14) para $1 \leq j \leq k$:

$$(YM)_{k,j} = M_{k+1,j} = b^j\left(\frac{\prod_{i=0}^{k} A_i}{\prod_{i=0}^{j-1} A_i}\right).$$

Demostraremos que el segundo miembro es igual a esta entrada.

$$\begin{aligned}
(M^2)_{k,j} &= M_{k,*} \cdot M_{*,j} = (M_{k,0}, M_{k,1}, \cdots, M_{k,j-1}, M_{k,j}, \cdots, M_{k,k}, 0, \cdots) \\
&\quad \cdot (\underbrace{M_{0,j}}_{0}, \underbrace{M_{1,j}}_{0}, \cdots, \underbrace{M_{j-1,j}}_{0}, M_{j,j}, \cdots, M_{k,j}, M_{k+1,j}, \cdots) \\
&= M_{k,j} M_{j,j} + M_{k,j+1} M_{j+1,j} + \cdots + M_{k,k} M_{k,j} \\
&= \left(b^j\left(\frac{\prod_{i=0}^{k-1} A_i}{\prod_{i=0}^{j-1} A_i}\right) b^j\left(\frac{\prod_{i=0}^{j-1} A_i}{\prod_{i=0}^{j-1} A_i}\right) + b^{j+1}\left(\frac{\prod_{i=0}^{k-1} A_i}{\prod_{i=0}^{j} A_i}\right) b^j\left(\frac{\prod_{i=0}^{j} A_i}{\prod_{i=0}^{j-1} A_i}\right) + \cdots\right. \\
&\quad \left. + b^k\left(\frac{\prod_{i=0}^{k-1} A_i}{\prod_{i=0}^{k-1} A_i}\right) b^j\left(\frac{\prod_{i=0}^{k-1} A_i}{\prod_{i=0}^{j-1} A_i}\right)\right) \\
&= b^j(b^j + b^{j+1} + \cdots + b^k)\left(\frac{\prod_{i=0}^{k-1} A_i}{\prod_{i=0}^{j-1} A_i}\right).
\end{aligned}$$

Ahora $(MY)_{k,j} = M_{k,j-1} = b^{j-1} A_{j-1} \left( \dfrac{\prod_{i=0}^{k-1} A_i}{\prod_{i=0}^{j-1} A_i} \right)$.

Reemplazando se tiene

$$\begin{aligned}
(M^2)_{k,j} + b(MY)_{k,j} &= [b^j(b^j + b^{j+1} + \cdots + b^k) + b^j A_{j-1}] \left( \dfrac{\prod_{i=0}^{k-1} A_i}{\prod_{i=0}^{j-1} A_i} \right) \\
&= [b^j(\underbrace{1 + b + b^2 + b^{j-1}}_{A_{j-1}} + b^j + b^{j+1} + \cdots + b^k] \left( \dfrac{\prod_{i=0}^{k-1} A_i}{\prod_{i=0}^{j-1} A_i} \right) \\
&= b^j \left( \dfrac{\prod_{i=0}^{k} A_i}{\prod_{i=0}^{j-1} A_i} \right).
\end{aligned}$$

Demostraremos (2.4.14) para $j = k+1$:

$$(YM)_{k,k+1} = M_{k+1,k+1} = b^{k+1} \left( \dfrac{\prod_{i=0}^{k} A_i}{\prod_{i=0}^{k} A_i} \right) = b^{k+1}.$$

$$(M^2)_{k,k+1} = M_{k,*} \cdot M_{*,k+1} = (M_{k,0}, M_{k,1}, \cdots, M_{k,k}) \cdot (\underbrace{M_{0,k+1}}_{0}, \underbrace{M_{1,k+1}}_{0}, \cdots, \underbrace{M_{k,k+1}}_{0}) = 0$$

$$(MY)_{k,k+1} = M_{k,k} = b^k \left( \dfrac{\prod_{i=0}^{k-1} A_i}{\prod_{i=0}^{k-1} A_i} \right) = b^k.$$

Reemplazando se tiene

$$(M^2)_{k,k+1} + b(MY)_{k,k+1} = b^{k+1}.$$

Finalmente demostraremos (2.4.14) para $j > k+1$:

$$(YM)_{k,j} = M_{k+1,j} = 0$$

por definición de $M$.

$$\begin{aligned}
(M^2)_{k,j} &= M_{k,*} \cdot M_{*,j} \\
&= (M_{k,0}, M_{k,1}, \cdots, M_{k,k}) \cdot (\underbrace{M_{0,j}}_{0}, \underbrace{M_{1,j}}_{0}, \cdots, \underbrace{M_{k,j}}_{0}) = 0.
\end{aligned}$$

Ahora $(MY)_{k,j} = M_{k,j-1} = 0$ pues $j - 1 > k (\Leftrightarrow j > k+1)$.

Reemplazando se tiene

$$(M^2)_{k,j} + b(MY)_{k,j} = 0.$$

**Caso 3.** $\boxed{\lambda_1 = \lambda_2}$

En este caso asumimos que la característica de $K$ es 0 y aquí tenemos

$$b + 1 + \sqrt{(b-1)^2 - 4c} = b + 1 - \sqrt{(b-1)^2 - 4c} \Leftrightarrow \sqrt{(b-1)^2 - 4c} = 0.$$

Es decir $\boxed{c = \dfrac{(b-1)^2}{4}}$.

Además $\lambda = \lambda_1 = \lambda_2 = \dfrac{b+1}{2} \neq 0$, pues si $\lambda = 0$ entonces $b = -1$ y por lo tanto $c = 1$, lo cual es una contradicción.

Reemplazando $b = 2\lambda - 1$ y $c = (\lambda - 1)^2$ en las relaciones recursivas

$$\begin{cases} P_{n+1} &= b P_n + Q_n \\ Q_{n+1} &= -c P_n + Q_n \end{cases}$$

obtenemos por recurrencia:

$$P_{n+1} = (n+1)\lambda^n, \quad Q_{n+1} = (n+1-n\lambda)\lambda^n. \qquad (2.4.15)$$

En efecto, para $n = 0$ la fórmula (2.4.15) se cumple trivialmente.

Para $n = 1$

$$P_2 = b P_1 + Q_1 = b + 1 = 2\lambda, \qquad Q_2 = -c P_1 + Q_1 = -c + 1 = -(\lambda-1)^2 + 1 = (2-\lambda)\lambda.$$

Supongamos que (2.4.15) se cumple para algún $n$, demostraremos que se cumple para $n+1$.

En efecto

$$\begin{aligned}
P_{n+2} &= b P_{n+1} + Q_{n+1} = b(n+1)\lambda^n + (n+1-n\lambda)\lambda^n \\
&= (2\lambda - 1)(n+1)\lambda^n + (n+1-n\lambda)\lambda^n \\
&= (2\lambda n + 2\lambda - n - 1 + n + 1 - n\lambda)\lambda^n = (n+2)\lambda^{n+1} \\
Q_{n+2} &= -c P_{n+1} + Q_{n+1} = -(\lambda - 1)^2 (n+1)\lambda^n + (n+1-n\lambda)\lambda^n \\
&= [-(\lambda^2 - 2\lambda + 1)(n+1) + n + 1 - n\lambda]\lambda^n = [-\lambda^2 n - \lambda^2 + 2\lambda n + 2\lambda - n\lambda]\lambda^n \\
&= [-\lambda n - \lambda + 2n + 2 - n]\lambda^{n+1} = [n + 2 - (n+1)\lambda]\lambda^{n+1}.
\end{aligned}$$

Así hemos demostrado

$$Q_{n+1} = (n+1-n\lambda)\lambda^n, \text{ para } n \geq 1.$$

En consecuencia $Q_{n+1} = 0 \Leftrightarrow \lambda = \dfrac{n+1}{n} \Leftrightarrow b = 2\lambda - 1 = \dfrac{2n+2}{n} - 1 = \dfrac{n+2}{n}$.

Es decir

$$Q_{n+1} = 0 \Leftrightarrow b = 1 + \dfrac{2}{n}, \text{ para } n \geq 1.$$

Desarrollaremos el caso particular $\boxed{b = 0, c = \dfrac{1}{4}}$.

Afirmamos que en este caso $M$ definida por

$$M_{k,j} = \begin{cases} \dfrac{2^k}{k+1} & \text{si} & j = 0 \\ \dfrac{k}{2k+2} & \text{si} & j = k+1 \\ 0 & & \text{en otros casos} \end{cases}$$

para $k \geq 1$ y $M_{0,j} = \delta_{0,j}$, es la matriz asociada al torcimiento.

Notemos que

$$M = \begin{pmatrix} 1 & 0 & 0 & 0 & 0 & \cdots & 0 & 0 & 0 & 0 & \cdots \\ 1 & 0 & \frac{1}{4} & 0 & 0 & \cdots & 0 & 0 & 0 & 0 & \cdots \\ \frac{4}{3} & 0 & 0 & \frac{1}{3} & 0 & \cdots & 0 & 0 & 0 & 0 & \cdots \\ 2 & 0 & 0 & 0 & \frac{3}{8} & \cdots & 0 & 0 & 0 & 0 & \cdots \\ \vdots & \vdots & \vdots & \vdots & \vdots & \vdots & \vdots & \vdots & \vdots & \vdots \\ \frac{2^n}{n+1} & 0 & 0 & 0 & 0 & \cdots & 0 & \frac{n}{2(n+1)} & 0 & 0 & \cdots \\ \vdots & \vdots & \vdots & \vdots & \vdots & \vdots & \vdots & \vdots & \ddots & \vdots & \vdots \end{pmatrix}$$

Por la Observación 2.4.1 sólo tenemos que demostrar que $M$ cumple

$$YM = M^2 + \frac{1}{4} Y^2. \tag{2.4.16}$$

Para demostrar (2.4.16), definimos las matrices $A$, $SA$, $B$ y $SB$, con

$$A_{ij} = \frac{2^i}{i+1} \delta_{0j}, \quad (SA)_{ij} = A_{i+1,j}, \quad B_{ij} = \frac{i}{2(i+1)} \delta_{ij}, \quad \text{y} \quad (SB)_{ij} = B_{i+1,j+1}.$$

Se verifica directamente que $M = A + BY$, $AB = 0$, $YB = (SB)Y$ y $A^2 = A$.

Entonces, como $YA = SA$, se tiene

$$M^2 = A + B(SA) + B(SB)Y^2.$$

Como además se calcula que

$$A + B(SA) = SA = YA \quad \text{y} \quad B(SB) = SB - \frac{1}{4} Id,$$

obtenemos que

$$YM = YA + (SB)Y^2 = (A + B(SA)) + (B(SB) + \frac{1}{4} Id)Y^2 = M^2 + \frac{1}{4} Y^2,$$

que es lo que queríamos demostrar.

**Observación 2.4.2.** Existen otros casos interesantes como por ejemplo cuando

$$\boxed{\lambda_1 = -\lambda_2, \quad \lambda_2 \neq 0}.$$

En este caso se tiene

$$b + 1 + \sqrt{(b-1)^2 - 4c} = -(b+1) + \sqrt{(b-1)^2 - 4c}.$$

Es decir $\boxed{b = -1}$.

Luego estamos en el caso

$$\boxed{YM = M^2 - MY + cY^2}, \quad c \neq 1.$$

Afirmamos que en este caso $M$ dada por

$$M_{ij} = \begin{cases} \delta_{ij} & \text{si } i \text{ es par,} \\ 1 & \text{si } i \text{ es impar y } j = i-1, \\ -1 & \text{si } i \text{ es impar y } j = i, \\ c & \text{si } i \text{ es impar y } j = i+1, \\ 0 & \text{en otros casos} \end{cases}$$

es la matriz asociada al torcimiento.

Notemos que

$$M = \begin{pmatrix} 1 & 0 & 0 & 0 & 0 & 0 & \cdots \\ 1 & -1 & c & 0 & 0 & 0 & \cdots \\ 0 & 0 & 1 & 0 & 0 & 0 & \cdots \\ 0 & 0 & 1 & -1 & c & 0 & \cdots \\ 0 & 0 & 0 & 0 & 1 & 0 & \cdots \\ \vdots & \vdots & \vdots & \vdots & \vdots & \vdots & \ddots \end{pmatrix}$$

Por la Observación 2.4.1 sólo tenemos que demostrar que $M$ cumple

$$YM = M^2 - MY + cY^2, \qquad (2.4.17)$$

o equivalentemente, que

$$(M-Y)^2 = (1-c)Y^2.$$

Pero $M - Y$ está formada por bloques $A$ de $2 \times 2$ en la diagonal y bloques $C$ de $2 \times 2$ en la primera diagonal de arriba con

$$A = \begin{pmatrix} 1 & -1 \\ 1 & -1 \end{pmatrix} \quad \text{y} \quad C = \begin{pmatrix} 0 & 0 \\ c-1 & 0 \end{pmatrix}.$$

Además se cumple

$$A^2 = 0, \quad AC + CA = (1-c)Id \quad \text{y} \quad C^2 = 0.$$

Usando esto se verifica que $(M-Y)^2 = (1-c)Y^2$.

Para demostrar que esta matriz representa efectivamente una aplicación de torcimiento basta demostrar que no es posible hallar $n \geq 0$ que satisfaga

$$\left(\frac{\lambda_1}{\lambda_2}\right)^n = \frac{1-\lambda_1}{1-\lambda_2}.$$

Supongamos por el absurdo que existe tal $n$. Entonces

$$(-1)^n = \frac{1+\lambda_2}{1-\lambda_2}.$$

Si $n$ es par se tiene $1 + \lambda_2 = 1 - \lambda_2$ lo cual conduce a $\lambda_2 = 0$, una contradicción.
Si $n$ es impar se tiene $1 + \lambda_2 = -1 + \lambda_2$ lo cual conduce a $1 = -1$, un absurdo.

## 2.5. Comparación con la clasificación de Conner-Goetz

Combinando la Proposición 2.2.1 y el Teorema 2.2.7 obtenemos que $b$ y $c$ determinan una (necesariamente única) aplicación de torcimiento vía la Proposición 2.1.4 si y sólo si $Q_n(b,c) \neq 0$ para todo $n \in \mathbb{N}$.

Esta condición es la misma que la condición usada en [6, Theorem 3.4]. Para verificar esto, notamos primero que los polinomios $f_n(a,b)$ usados por [6] satisfacen $f_n(a,b) = Q_{n+1}(b,a)$.

En efecto, dado que $Q_1(b,a) = f_0(a,b) = 1$, $P_1(b,a) = e_0(a,b) = 1$, y las fórmulas recursivas son las mismas, es decir,

$$\begin{pmatrix} P_n(b,a) \\ Q_n(b,a) \end{pmatrix} = \begin{pmatrix} b & 1 \\ -a & 1 \end{pmatrix} \begin{pmatrix} P_{n-1}(b,a) \\ Q_{n-1}(b,a) \end{pmatrix} \quad \text{y} \quad \begin{pmatrix} e_n(a,b) \\ f_n(a,b) \end{pmatrix} = \begin{pmatrix} b & 1 \\ -a & 1 \end{pmatrix} \begin{pmatrix} e_{n-1}(a,b) \\ f_{n-1}(a,b) \end{pmatrix},$$

concluímos $e_n(a,b) = P_{n+1}(b,a)$ y $f_n(a,b) = Q_{n+1}(b,a)$, como se quería.

Esto significa que nuestros resultados coinciden con los de [T. 3.4, Conner-Goetz]: Si $a \neq 0$ y $ac \neq 1$, entonces

$$K[x] \otimes_\sigma K[y] \cong C(a,b,c) = K\langle x,y \rangle / <-yx + ax^2 + bxy + cy^2>,$$

si y solo si $Q_n(b,ac) \neq 0$ para todo $n \in \mathbb{N}$.

A continuación daremos una clasificación de las álgebras $C(a,b,c)$ con $ac \neq 1$, salvo isomorfismo de álgebras graduadas, siguiendo [Sección 3.2, Conner-Goetz]. Por isomorfismos de productos torcidos (ver Observación 2.1.4) es suficiente considerar los casos $C(0,b,1)$, $C(0,b,0)$ y $C(1,b,c)$ con $c \neq 1$. Se conoce (ver por ejemplo [Sección 3.2, Conner-Goetz]) que en ese caso los productos torcidos obtenidos son isomorfos (como álgebras graduadas) a una de las siguientes álgebras:

i) $K_q[x,y] = K\langle x,y \rangle / <yx - qxy>$, $q \neq 0$, el plano cuántico,

ii) $J(x,y) = K\langle x,y \rangle / <yx - xy - y^2>$, el plano de Jordan,

iii) $K\langle x,y \rangle / <yx>$.

Los casos $C(0,b,1)$ y $C(0,b,0)$ están contenidos en la subsección 2.1.1.

Ahora veremos todos estos casos.

$\boxed{(1)\ C(0,b,0) \cong K_b[x,y] \text{ si } b \neq 0}$.

En este caso se tiene $C(0,b,0) = K\langle x,y \rangle / \langle yx - bxy \rangle = K_b[x,y]$.

$\boxed{(2)\ C(0,1,1) \cong J(x,y)}$.

En este caso se tiene $C(0,1,1) = K\langle x,y \rangle / \langle yx - xy - y^2 \rangle = J(x,y)$.

$\boxed{(3)\ C(0,b,1) \cong K_b[x,y] \text{ si } b \neq 1}$.

Para demostrar el isomorfismo usamos el cambio de variables $w = x + y, \quad z = (1-b)y$ en la relación

$$zw - bwz - z^2 = 0$$

y obtenemos

$$(1-b)y(x+y)-b(x+y)(1-b)y-(1-b)^2 y^2 = 0,$$
$$(1-b)(yx - bxy) = 0,$$
$$yx - bxy = 0$$

pues $b \neq 1$.

Luego $K\langle w, z\rangle/\langle zw - bwz - z^2\rangle \cong K\langle x, y\rangle/\langle yx - bxy\rangle$

Es decir $C(0, b, 1) \cong K\langle x, y\rangle/\langle yx - bxy\rangle = K_b[x, y]$.

Continuamos con el caso $C(1, b, c)$, $c \neq 1$ de la sección 2.4.

$\boxed{(4)\ C(1,-c,c) \cong K\langle x, y\rangle/\langle yx\rangle}$.

Para demostrar el isomorfismo primero notemos que

$$0 = yx - x^2 + cxy - cy^2 = (y-x)x + c(x-y)y = (y-x)(x-cy).$$

Haciendo el cambio de variables $z = y - x$ y $w = x - cy$ se tiene $zw = 0$. Como $c \neq 1$, se tiene $K\langle x, y\rangle = K\langle z, w\rangle$ y por lo tanto $C(1, -c, c) \cong K\langle z, w\rangle/\langle zw\rangle$.

Este caso corresponde al caso excepcional 1: $\lambda_1 = 0$ o $\lambda_2 = 0$ para $c \neq 1$.

$\boxed{(5)\ C(1,-1,c) \cong K_{-1}[x, y] = K\langle x, y\rangle/\langle yx + xy\rangle}$.

Para demostrar el isomorfismo usamos el cambio de variables

$$w = x + (-1 + \sqrt{1-c})y, \quad z = \frac{x}{2} + \left(\frac{-1 - \sqrt{1-c}}{2}\right)y$$

en $0 = wz + zw$. Obtenemos

$$\begin{aligned}
0 &= wz + zw \\
&= \frac{x^2}{2} + \left(\frac{-1-\sqrt{1-c}}{2}\right)xy + \left(\frac{-1+\sqrt{1-c}}{2}\right)yx + \frac{c}{2}y^2 \\
&+ \frac{x^2}{2} + \left(\frac{-1+\sqrt{1-c}}{2}\right)xy + \left(\frac{-1-\sqrt{1-c}}{2}\right)yx + \frac{c}{2}y^2 \\
&= x^2 - xy - yx + cy^2.
\end{aligned}$$

Es decir

$$0 = wz + zw = -yx + x^2 - xy + cy^2.$$

Por lo tanto: $C(1, -1, c) \cong K\langle w, z\rangle/\langle wz + zw\rangle \cong K_{-1}[w, z]$.

Este caso corresponde al caso $\lambda_1 = -\lambda_2$, $\lambda_2 \neq 0$, de la Observación 2.4.2.

$\boxed{(6)\ C(1, b, c) \cong J(x, y), \text{ si } b \neq -1 \text{ y } 4c - (b-1)^2 = 0}$.

Para demostrar el isomorfismo usamos el cambio de variables

$$w = -\frac{i(b+1)}{2}y, \quad z = ix + \frac{i(b-1)}{2}y$$

en $0 = zw - wz - z^2$, donde $i \in K$ es tal que $i^2 = -1$. Obtenemos

$$0 = -\left(ix + \frac{i(b-1)y}{2}\right)\left(\frac{i(b+1)y}{2}\right) + \left(\frac{i(b+1)y}{2}\right)\left(ix + \frac{i(b-1)y}{2}\right) - (ix + \frac{i(b-1)}{2}y)^2,$$

$$0 = \frac{b+1}{2}xy + \frac{b^2-1}{4}y^2 - \frac{b+1}{2}yx - \frac{b^2-1}{4}y^2 + x^2 + \frac{b-1}{2}xy + \frac{b-1}{2}yx + \frac{(b-1)}{4}y^2.$$

Es decir

$$0 = zw - wz - z^2 = -yx + x^2 + bxy + cy^2.$$

Por lo tanto $C(1,b,c) \cong K\langle w,z \rangle / \langle zw - wz - z^2 \rangle \cong J(w,z)$.

Este caso corresponde al caso excepcional 3: $\lambda_1 = \lambda_2$ para $c = \frac{(b-1)^2}{4}$, $b \neq 1 + \frac{2}{n}$.

$\boxed{(7)\ C(1,b,c) \cong K_q[x,y], \text{ si } b \neq -1, b \neq -c \text{ y } 4c - (b-1)^2 \neq 0}$,

donde $q$ satisface $(c+b)q^2 + (2c - b^2 - 1)q + c + b = 0$.

Afirmamos que $q \in K - \{0, 1, -1\}$.

Para ver esto supongamos que $q \in \{0, 1, -1\}$.

Si $q = 0$ entonces reemplazando en la relación que define a $q$ se obtiene $c + b = 0 (\Rightarrow \Leftarrow)$.

Si $q = 1$ se obtiene $c + b + 2c - b^2 - 1 + c + b = 0$, es decir $4c - (b-1)^2 = 0 (\Rightarrow \Leftarrow)$.

Si $q = -1$, entonces $c + b - (2c - b^2 - 1) + c + b = 0$, es decir $b = -1 (\Rightarrow \Leftarrow)$.

Para demostrar el isomorfismo usamos el cambio de variables

$$w = \frac{x}{B} + \frac{b-q}{AB}y, \quad z = x + \frac{qb-1}{A}y$$

donde $A = 1 + q$, $B = 1 - q$.

Reemplazando en $0 = zw - qwz$ obtenemos

$$0 = \left(x + \frac{qb-1}{A}y\right)\left(\frac{x}{B} + \frac{b-q}{AB}y\right) - q\left(\frac{x}{B} + \frac{b-q}{AB}y\right)\left(x + \frac{qb-1}{A}y\right),$$

$$0 = (ABx + B(qb-1)y)(Ax + (b-q)y) - q(Ax + (b-q)y)(ABx + B(qb-1)y),$$

$$0 = (1-q)A^2Bx^2 + AB(1-q^2)bxy - (1-q^2)AByx + (1-q)B(qb-1)(b-q)y^2.$$

Dividiendo entre $A^2B^2$:

$$0 = x^2 + bxy - yx + \frac{qb^2 - q^2b - b + q}{(1+q)^2}y^2.$$

Pero por hipótesis $c = \frac{qb^2 - q^2b - b + q}{(1+q)^2}$.

Por lo tanto:

$$C(1,b,c) \cong K\langle w,z \rangle / \langle zw - qwz \rangle \cong K_q[w,z].$$

$\boxed{\text{Subcaso } c = 0, b \in K - \{-1, 0, 1\}}$.

Se tiene entonces $C(1,b,0) \cong K_q[x,y]$ si $b \in K - \{-1, 0, 1\}$.

Este subcaso corresponde también al caso excepcional 2: $\lambda_1 = 1$ ó $\lambda_2 = 1$ para $c = 0$, $b \in K$. En este subcaso $q$ satisface $bq^2 + (-b^2-1)q + b = 0$, es decir
$$q = \frac{b^2+1 \pm \sqrt{(b^2+1)^2 - 4b^2}}{2b}.$$
De donde $q = b$ ó $q = \frac{1}{b}$. En cualquier caso se tiene $C(1,b,0) \cong K_q[x,y]$.

$\boxed{\text{Subcaso } c \neq 0, b \neq -1, b \neq -c \text{ y } 4c - (b-1)^2 \neq 0}$.

Se tiene entonces $C(1,b,c) \cong K_q[x,y]$ si $b \neq -c$ y $4c - (b-1)^2 \neq 0$.

Este subcaso corresponde al caso en el cual los valores propios $\lambda_1, \lambda_2$ son diferentes tales que
$$0 \notin \{\lambda_1, \lambda_2, 1-\lambda_1, 1-\lambda_2, \lambda_1 - \lambda_2\}.$$

En efecto, al calcular dichos autovalores se obtiene la fórmula
$$\lambda_{1,2} = \frac{1}{2}(b+1 \pm \sqrt{(b-1)^2 - 4c}).$$

Las condiciones que corresponden a este subcaso garantizan que $\lambda_1, \lambda_2$ sean diferentes.

# Capítulo 3

# El caso $yx = x^2 - xy + y^2$

## 3.1. El caso $yx = x^2 - xy + y^2$

En esta sección asumiremos que $\sigma$ es una aplicación de torcimiento y que $Y$, $M$ son como en la Proposición 1.4.1. Como antes escribiremos $M_{1,*} = (1, b, c, 0, ...)$ y supondremos que $(b, c) = (-1, 1)$, el cual, por el Lema 2.1.5 (2), es el único caso no cubierto por el Teorema 2.2.7.

Esto significa que estamos tratando con la regla de conmutación

$$yx = x^2 - xy + y^2.$$

Por lo tanto tenemos $YM = M^2 - MY + Y^2$ lo cual implica $\widetilde{M}^2 = 0$, donde $\widetilde{M} = M - Y$.

La matriz $\widetilde{M} = M - Y = \psi(x - y) \in L(K^{\mathbb{N}_0})$ juega un papel central en la clasificación de todas las aplicaciones de torcimiento con $(b, c) = (-1, 1)$.

Notemos que

$$\widetilde{M}_{0,j} = \begin{cases} 1, & \text{si } j = 0 \\ -1, & \text{si } j = 1 \\ 0, & \text{si } j > 1 \end{cases}$$

y que $\widetilde{M}_{1,*} = \widetilde{M}_{0,*}$.

**Observación 3.1.1.** Sea $\widetilde{M} \in L(K^{\mathbb{N}_0})$ tal que $\widetilde{M}_{0,j} = \delta_{0,j} - \delta_{1,j}$ y $\widetilde{M}_{k,j} = 0$ para $j > k + 1$. Entonces la matriz $M := \widetilde{M} + Y$ determina una aplicación de torcimiento vía la Proposición 1.4.1 si y sólo si para todo $k \in \mathbb{N}$

$$Y^k \widetilde{M} = \sum_{j=0}^{k+1} \widetilde{M}_{k,j} M^{k+1-j} Y^j. \tag{3.1.1}$$

En efecto supongamos que

$$Y^k M = \sum_{j=0}^{k+1} M_{k,j} M^{k+1-j} Y^j, \quad k \geq 0.$$

Entonces

$$Y^k \widetilde{M} = Y^k(M-Y) = Y^k M - Y^{k+1} = \sum_{j=0}^{k} M_{k,j} M^{k+1-j} Y^j + (M_{k,k+1}-1) Y^{k+1}.$$

Pero $\widetilde{M}_{k,j} = M_{k,j}$, si $0 \leq j \leq k$ y $\widetilde{M}_{k,k+1} = M_{k,k+1} - 1$.
Luego

$$Y^k \widetilde{M} = \sum_{j=0}^{k+1} \widetilde{M}_{k,j} M^{k+1-j} Y^j$$

para todo $k \geq 0$.
En forma análoga se demuestra el recíproco.

**Lema 3.1.2.** *Sea $d \neq 0$. Entonces para $k \geq 1$ se cumple*

(1) $\widetilde{M}_{k,*} = \widetilde{M}_{k-1,*}$ *si y sólo si* $\widetilde{M} Y^{k-1} \widetilde{M} = 0$.

(2) $\widetilde{M}_{k,*} = d \widetilde{M}_{k-1,*}$ *si y sólo si* $\widetilde{M} Y^{k-1} \widetilde{M} = \dfrac{1-d}{d} Y^k \widetilde{M}$.

(3) $\widetilde{M}_{k,*} = d \widetilde{M}_{k-2,*}$ *si y sólo si* $\widetilde{M} Y^{k-1} \widetilde{M} + Y \widetilde{M} Y^{k-2} \widetilde{M} = \dfrac{1-d}{d} Y^k \widetilde{M}$.

**Demostración.** Sólo demostraremos (2) y (3) ya que (1) se sigue de (2) con $d = 1$.
Dado que

$$Y^{k-1} \widetilde{M} = \sum_{j=0}^{k} \widetilde{M}_{k-1,j} M^{k-j} Y^j \quad \text{y} \quad Y^k \widetilde{M} = \sum_{j=0}^{k+1} \widetilde{M}_{k,j} M^{k+1-j} Y^j,$$

evidentemente $\widetilde{M}_{k,*} = d \widetilde{M}_{k-1,*}$ implica $Y^k \widetilde{M} = dM Y^{k-1} \widetilde{M} = d(\widetilde{M} + Y) Y^{k-1} \widetilde{M}$, lo cual es válido si y sólo si $(1-d) Y^k \widetilde{M} = d \widetilde{M} Y^{k-1} \widetilde{M}$. Por otro lado, si $Y^k \widetilde{M} = dM Y^{k-1} \widetilde{M}$, entonces la primera fila de la igualdad matricial

$$\sum_{j=0}^{k+1} \widetilde{M}_{k,j} M^{k+1-j} Y^j = \sum_{j=0}^{k} d \widetilde{M}_{k-1,j} M^{k+1-j} Y^j$$

es

$$\sum_{j=0}^{k+1} \widetilde{M}_{k,j} E_j = \sum_{j=0}^{k} d \widetilde{M}_{k-1,j} E_j,$$

donde $E_j$ es el vector infinito con $(E_j)_i := \delta_{i,j}$, lo cual produce $\widetilde{M}_{k,*} = d \widetilde{M}_{k-1,*}$. De esta manera queda demostrado (2).

Similarmente se demuestra (3) notando que $\widetilde{M}_{k,*} = d \widetilde{M}_{k-2,*}$ si y sólo si $Y^k \widetilde{M} = dM^2 Y^{k-2} \widetilde{M} = d(\widetilde{M} Y + Y \widetilde{M} + Y^2) Y^{k-2} \widetilde{M}$ si y sólo si $(1-d) Y^k \widetilde{M} = d \widetilde{M} Y^{k-1} \widetilde{M} + d Y \widetilde{M} Y^{k-2} \widetilde{M}$. □

**Lema 3.1.3.** *Sea* $n \in \mathbb{N}$ *con* $n \geq 2$. *Entonces* $\widetilde{M}_{k,*} = \widetilde{M}_{0,*}$ *para* $1 \leq k \leq n-1$ *si y sólo si* $\widetilde{M} Y^k \widetilde{M} = 0$ *para* $0 \leq k \leq n-2$. *Además en este caso:*

a) $M^{k+1} = Y^{k+1} + \sum_{j=0}^{k} Y^j \widetilde{M} Y^{k-j}$   *para* $0 \leq k \leq n-1$.

b) $M^{n+1} = Y^{n+1} + \widetilde{M} Y^{n-1} \widetilde{M} + \sum_{j=0}^{n} Y^j \widetilde{M} Y^{n-j}$.

**Demostración.** Por un lado $\widetilde{M}_{k,*} = \widetilde{M}_{k-1,*}$ para $1 \leq k \leq n-1$ si y solo si $\widetilde{M}_{k,*} = \widetilde{M}_{0,*}$ para $1 \leq k \leq n-1$. Por otro lado, cambiando $k$ por $k-1$, se tiene que $\widetilde{M} Y^{k-1} \widetilde{M} = 0$ para $0 \leq k-1 \leq n-2$ si y solo si $\widetilde{M} Y^k \widetilde{M} = 0$ para $0 \leq k \leq n-2$, lo cual demuestra la primera parte del lema.

<u>Demostración de a)</u>

Usaremos inducción matemática sobre $k$.

Para $k = 0$: $M = \widetilde{M} + Y = Y + \sum_{j=0}^{0} Y^j \widetilde{M} Y^{-j}$.

Para $k = 1$: $M^2 = (\widetilde{M} + Y)^2 = Y^2 + \widetilde{M} Y + Y \widetilde{M} = Y^2 + \sum_{j=0}^{1} Y^j \widetilde{M} Y^{1-j}$.

Supongamos que se cumple

$$M^k = Y^k + \sum_{j=0}^{k-1} Y^j \widetilde{M} Y^{k-1-j} \quad \text{para algún } 0 \leq k \leq n-1.$$

Entonces

$$\begin{aligned} M^{k+1} &= (\widetilde{M} + Y) M^k = (\widetilde{M} + Y) \left( Y^k + \sum_{j=0}^{k-1} Y^j \widetilde{M} Y^{k-1-j} \right) \\ &= \widetilde{M} Y^k + \sum_{j=0}^{k-1} \widetilde{M} Y^j \widetilde{M} Y^{k-1-j} + Y^{k+1} + \sum_{j=0}^{k-1} Y^{j+1} \widetilde{M} Y^{k-1-j}. \end{aligned}$$

Ahora la primera sumatoria se anula pues $\widetilde{M} Y^j \widetilde{M} = 0$ ya que $0 \leq j \leq k-1 \leq n-2$ implica $0 \leq j \leq n-2$.

Por lo tanto

$$M^{k+1} = \widetilde{M} Y^k + Y^{k+1} + \sum_{j=1}^{k} Y^j \widetilde{M} Y^{k-j} = Y^{k+1} + \sum_{j=0}^{k} Y^j \widetilde{M} Y^{k-j}.$$

De esta manera queda demostrado el paso inductivo y así se ha demostrado *a)*.

<u>Demostración de b)</u>

$$\begin{aligned} M^{n+1} &= (\widetilde{M} + Y) M^n = (\widetilde{M} + Y) \left( Y^n + \sum_{j=0}^{n-1} Y^j \widetilde{M} Y^{n-1-j} \right) \\ &= \widetilde{M} Y^n + \sum_{j=0}^{n-1} \widetilde{M} Y^j \widetilde{M} Y^{n-1-j} + Y^{n+1} + \sum_{j=0}^{n-1} Y^{j+1} \widetilde{M} Y^{n-1-j}. \end{aligned}$$

Ahora la primera sumatoria es igual a $\widetilde{M} Y^{n-1}\widetilde{M}$ pues $\widetilde{M} Y^j \widetilde{M} = 0$ para todo $0 \leq j \leq n-2$.

Luego

$$M^{n+1} = \widetilde{M} Y^n + \widetilde{M} Y^{n-1}\widetilde{M} + Y^{n+1} + \sum_{j=1}^{n} Y^j \widetilde{M} Y^{n-j}.$$

Es decir

$$M^{n+1} = Y^{n+1} + \widetilde{M} Y^{n-1}\widetilde{M} + \sum_{j=0}^{n} Y^j \widetilde{M} Y^{n-j}. \qquad \square$$

**Proposición 3.1.4.** *Existe una aplicación de torcimiento tal que $\widetilde{M}_{k,*} = \widetilde{M}_{0,*}$ para todo $k \geq 1$.*

**Demostración.** Sea la matriz $\widetilde{M}$ tal que $\widetilde{M}_{k,0} = 1$, $\widetilde{M}_{k,1} = -1$ y $\widetilde{M}_{k,j} = 0$ para $j > 1$.
Para demostrar que $M = \widetilde{M} + Y$ representa a una aplicación de torcimiento usaremos la Observación 3.1.1.
Evidentemente $\widetilde{M}$ cumple $\widetilde{M}_{0,j} = \delta_{0,j} - \delta_{1,j}$ para todo $j \geq 0$ y $\widetilde{M}_{k,j} = 0$ para $j > k + 1$. Por la Observación 3.1.1 es suficiente demostrar que:

$$Y^k \widetilde{M} = \sum_{j=0}^{k+1} \widetilde{M}_{k,j} M^{k+1-j} Y^j$$

para todo $k \geq 1$.
En nuestro caso debemos demostrar que

$$Y^k \widetilde{M} = M^{k+1} - M^k Y = M^k \widetilde{M}.$$

Pero $Y^k \widetilde{M} = \widetilde{M}$ para todo $k \geq 1$, pues $\widetilde{M}_{k,*} = \widetilde{M}_{0,*}$ para todo $k \geq 1$.
Entonces es suficiente demostrar que

$$M^k \widetilde{M} = \widetilde{M} \quad \text{para todo } k \geq 1.$$

Para $k = 1$ se cumple $M\widetilde{M} = (\widetilde{M} + Y)\widetilde{M} = \widetilde{M}^2 + Y\widetilde{M} = Y\widetilde{M} = \widetilde{M}$.
Por inducción sobre k, supongamos que se cumple $M^k \widetilde{M} = \widetilde{M}$ para algún $k \geq 1$.
Demostraremos que:

$$M^{k+1}\widetilde{M} = \widetilde{M}.$$

En efecto:

$$M^{k+1}\widetilde{M} = M(M^k \widetilde{M}) = M\widetilde{M} = Y\widetilde{M} = \widetilde{M}. \qquad \square$$

**Proposición 3.1.5.** *Sea $n \in \mathbb{N}$ con $n \geq 2$ y supongamos que $\widetilde{M}_{k,*} = \widetilde{M}_{0,*}$ para $1 < k < n$ y que $\widetilde{M}_{n,*} \neq \widetilde{M}_{0,*}$.*

*Entonces* $\widetilde{M}_{n,j} = 0$ *para* $1 < j < n$.

*Además, si definimos* $m_i := \widetilde{M}_{n,i}$ *entonces*

$$(1-m_0)Y^n\widetilde{M} = m_0\widetilde{M}Y^{n-1}\widetilde{M} + (m_0+m_1)\sum_{j=0}^{n-1} Y^j\widetilde{M}Y^{n-j} \qquad (3.1.2)$$
$$+ m_n\widetilde{M}Y^n + (m_0+m_1+m_n+m_{n+1})Y^{n+1}$$

$$y \quad (m_{n+1}+1)(m_0+m_1) + m_n + m_{n+1} = 0. \qquad (3.1.3)$$

**Demostración.** Las hipótesis de la proposición se traducen en la matriz $\widetilde{M}$ de la siguiente manera:

$$\widetilde{M} = \begin{pmatrix} 1 & -1 & 0 & \cdots & 0 & 0 & 0 & \cdots \\ 1 & -1 & 0 & \cdots & 0 & 0 & 0 & \cdots \\ \vdots & \vdots & \vdots & \vdots & \vdots & \vdots & \vdots & \vdots \\ 1 & -1 & 0 & \cdots & 0 & 0 & 0 & \cdots \\ m_0 & m_1 & m_2 & \cdots & m_n & m_{n+1} & 0 & \cdots \\ \vdots & \vdots & \vdots & \vdots & \vdots & \vdots & \vdots & \vdots \end{pmatrix} \leftarrow n-1.$$

Primero demostraremos que $m_j = \widetilde{M}_{n,j} = 0$ para $2 \leq j \leq n-1$.

Caso 1: $\boxed{m_{n+1} = 0}$

En este caso demostraremos que $\widetilde{M}_{n,j} = 0$ para $2 \leq j \leq n$.

En efecto:

$$\begin{aligned} 0 &= (\widetilde{M}^2)_{n,n} = \widetilde{M}_{n,*}\widetilde{M}_{*,n} \\ &= (\widetilde{M}_{n,0}, \widetilde{M}_{n,1}, \cdots, \widetilde{M}_{n,n}, \widetilde{M}_{n,n+1}, 0, \cdots) \cdot (\widetilde{M}_{0,n}, \widetilde{M}_{1,n}, \cdots, \widetilde{M}_{n,n}, \widetilde{M}_{n+1,n}, \cdots) \\ &= (m_0, m_1, \cdots, m_n, m_{n+1}, 0, \cdots) \cdot (\widetilde{M}_{0,n}, \widetilde{M}_{1,n}, \cdots, \widetilde{M}_{n,n}, \widetilde{M}_{n+1,n}, \cdots). \end{aligned}$$

Pero $\widetilde{M}_{k,*} = (1, -1, 0, \cdots)$ para $0 \leq k \leq n-1$.
Es decir $\widetilde{M}_{k,j} = 0$ para $j \geq 2$. Como $n \geq 2$, entonces
$\widetilde{M}_{k,n} = 0$ para $0 \leq k \leq n-1$.

De esta manera

$$\begin{aligned} 0 &= (m_0, m_1, \cdots, m_{n-1}, m_n, m_{n+1}, 0, \cdots) \cdot (0, 0, \cdots, 0, m_n, \widetilde{M}_{n+1,n}, \widetilde{M}_{n+2,n}, \cdots) \\ &= (m_n)^2 + m_{n+1}\widetilde{M}_{n+1,n} \\ &= m_n^2. \end{aligned}$$

Por lo tanto $m_n = 0$.

Asumiendo que $m_{n+1} = 0$ hemos demostrado que en la fila $\widetilde{M}_{n,*}$ la entrada $m_n =$

$\widetilde{M}_{n,n} = 0$. Falta demostrar que en esta misma fila, la entrada $m_j = \widetilde{M}_{n,j}$ cumple $m_j = 0$ para $2 \leq j \leq n-1$.

Supongamos ahora que $2 < k \leq n$ y que $m_j = 0$ para $j \geq k$.

Es decir $\widetilde{M}_{n,*} = (m_0, m_1, \cdots, m_{k-1}, m_k, \cdots, m_n, m_{n+1}, 0, \cdots) = (m_0, m_1, \cdots, m_{k-1}, 0, 0, 0, \cdots)$.

Demostraremos que $m_j = 0$ para $j \geq k-1$.

Es suficiente demostrar que $m_{k-1} = 0$.

(Notemos que $k > 2 \Rightarrow j \geq k-1 > 1 \Rightarrow j \geq 2$).

Por hipótesis $\widetilde{M}_{k,*} = \widetilde{M}_{0,*}$ para $1 \leq k \leq n-1$, luego usando el Lema 3.1.3 se tiene que $\widetilde{M} Y^k \widetilde{M} = 0$ para $0 \leq k \leq n-2$.

Ahora $1 < k-1 \leq n-1 \Leftrightarrow -n+1 \leq -(k-1) < -1 \Leftrightarrow 1 \leq n-k+1 < n-1$.

Entonces haciendo $k' = n-k+1$ obtenemos $1 \leq k' \leq n-2$.

Luego $\widetilde{M} Y^{k'} \widetilde{M} = 0$, es decir $\widetilde{M} Y^{n-k+1} \widetilde{M} = 0$.

Evaluando esta igualdad en la entrada $(n, k-1)$ obtenemos:

$$0 = \left(\widetilde{M} Y^{n-k+1} \widetilde{M}\right)_{n,k-1} = \sum_{j=0}^{\infty} \widetilde{M}_{n,j} \left(Y^{n-k+1} \widetilde{M}\right)_{j,k-1} = \sum_{j=0}^{\infty} \widetilde{M}_{n,j} \widetilde{M}_{j+n-k+1,k-1}.$$

Por hipótesis inductiva, para $j > k-1$ se tiene $\widetilde{M}_{n,j} = 0$. Para $j < k-1$ se tiene $j+n-k+1 < n$ y entonces $\widetilde{M}_{j+n-k+1,k-1} = 0$. Por lo tanto

$$0 = \sum_{j=k-1}^{k-1} \widetilde{M}_{n,j} \widetilde{M}_{j+n-k+1,k-1} = \widetilde{M}_{n,k-1} \widetilde{M}_{n,k-1} = m_{k-1}^2,$$

y obtenemos que $m_{k-1} = 0$. Así, inductivamente, hemos demostrado que $\widetilde{M}_{n,j} = m_j = 0$ para $2 \leq j \leq n$.

A partir de $(\widetilde{M}^2)_{n,0} = 0$ obtenemos

$$\begin{aligned}
0 &= (\widetilde{M}^2)_{n,0} = \widetilde{M}_{n,*} \widetilde{M}_{*,n} \\
&= (\widetilde{M}_{n,0}, \widetilde{M}_{n,1}, \cdots, \widetilde{M}_{n,n}, \widetilde{M}_{n,n+1}, 0, \cdots) \cdot (\widetilde{M}_{0,0}, \widetilde{M}_{1,0}, \cdots, \widetilde{M}_{n-1,0}, \widetilde{M}_{n,0}, \widetilde{M}_{n+1,0}, \cdots) \\
&= (m_0, m_1, 0, \cdots, 0, 0, 0, \cdots) \cdot (1, 1, 1, \cdots, 1, \widetilde{M}_{n,0}, \widetilde{M}_{n+1,0}, \cdots) \\
&= m_0 + m_1.
\end{aligned}$$

En resumen, en el caso $m_{n+1} = 0$ hemos demostrado que $\widetilde{M}_{n,j} = m_j = 0$ para $2 \leq j \leq n$ y $m_0 + m_1 = 0$.

Es decir $\widetilde{M}_{n,*} = (m_0, m_1, 0, \cdots) = (m_0, -m_0, \cdots) = m_0(1, -1, 0, \cdots) = m_0 \widetilde{M}_{0,*}$.

Por el Lema 3.1.2 (2) con $d = m_0$ se tiene:

$$(1 - m_0) Y^n \widetilde{M} = m_0 \widetilde{M} Y^{n-1} \widetilde{M}$$

y queda demostrado (3.1.2) para este caso.

La afirmación (3.1.3) es inmediata pues $m_0 + m_1 = 0$ y $m_n = m_{n+1} = 0$.

Caso 2: $\boxed{m_{n+1} \neq 0}$

Sabemos que

$$Y^n \widetilde{M} = \sum_{k=0}^{n+1} m_k M^{n+1-k} Y^k = m_0 M^{n+1} + \sum_{k=1}^{n+1} m_k M^{n+1-k} Y^k. \quad (3.1.4)$$

Pero $1 \leq k \leq n-1 \Leftrightarrow 1-n \leq -k \leq -1 \Leftrightarrow 1 \leq n-k \leq n-1$.

Por el Lema 3.1.3 a) para $k' = n-k$

$$M^{n-k+1} Y^k = \left( Y^{n-k+1} + \sum_{j=0}^{n-k} Y^j \widetilde{M} Y^{n-k-j} \right) Y^k = Y^{n+1} + \sum_{j=0}^{n-k} Y^j \widetilde{M} Y^{n-j}.$$

También se cumple para $k = n$ pues en este caso

$$MY^n = (\widetilde{M} + Y)Y^n = \widetilde{M} Y^n + Y^{n+1}.$$

Es decir

$$M^{n-k+1} Y^k = Y^{n+1} + \sum_{j=0}^{n-k} Y^j \widetilde{M} Y^{n-j}, \quad 1 \leq k \leq n. \quad (3.1.5)$$

Usando el Lema 3.1.3 b) y reemplazando (3.1.5) en (3.1.4) se obtiene

$$\begin{aligned} Y^n \widetilde{M} &= \sum_{k=0}^{n+1} m_k M^{n+1-k} Y^k = m_0 M^{n+1} + \sum_{k=1}^{n+1} m_k M^{n+1-k} Y^k \\ &= m_0 \left( Y^{n+1} + \widetilde{M} Y^{n-1} \widetilde{M} + \sum_{j=0}^{n} Y^j \widetilde{M} Y^{n-j} \right) + \sum_{k=1}^{n+1} m_k \left( Y^{n+1} + \sum_{j=0}^{n-k} Y^j \widetilde{M} Y^{n-j} \right) \\ &= m_0 Y^{n+1} + m_0 \widetilde{M} Y^{n-1} \widetilde{M} + m_0 \sum_{j=0}^{n} Y^j \widetilde{M} Y^{n-j} \\ &\quad + \sum_{k=1}^{n+1} m_k Y^{n+1} + \sum_{k=1}^{n} \sum_{j=0}^{n-k} m_k Y^j \widetilde{M} Y^{n-j}. \end{aligned}$$

Pero $m_0 \sum_{j=0}^{n} Y^j \widetilde{M} Y^{n-j} = \sum_{k=0}^{0} \sum_{j=0}^{n-k} m_k Y^j \widetilde{M} Y^{n-j}$.

Luego

$$Y^n \widetilde{M} = m_0 Y^{n+1} + m_0 \widetilde{M} Y^{n-1} \widetilde{M} + \sum_{k=0}^{n} \sum_{j=0}^{n-k} m_k Y^j \widetilde{M} Y^{n-j} + \left( \sum_{k=1}^{n+1} m_k \right) Y^{n+1}.$$

Cambiando el orden de los índices en la doble sumatoria se obtiene:

$$\begin{aligned} Y^n \widetilde{M} &= m_0 Y^{n+1} + m_0 \widetilde{M} Y^{n-1} \widetilde{M} + \sum_{j=0}^{n} \sum_{k=0}^{n-j} m_k Y^j \widetilde{M} Y^{n-j} + \left( \sum_{k=1}^{n+1} m_k \right) Y^{n+1} \\ &= m_0 \widetilde{M} Y^{n-1} \widetilde{M} + \sum_{j=0}^{n} \left( \sum_{k=0}^{n-j} m_k \right) Y^j \widetilde{M} Y^{n-j} + \left( \sum_{k=0}^{n+1} m_k \right) Y^{n+1}. \end{aligned} \quad (3.1.6)$$

Evaluaremos (3.1.6) en las entradas $(i, i+n+1)$ para $i = 1, 2, \cdots, n-1$.
Primero demostraremos que $\left(\widetilde{M} Y^{n-1} \widetilde{M}\right)_{i,i+n+1} = 0$ para $i = 1, 2, \cdots, n-1$.
En efecto

$$\begin{aligned}
\left(\widetilde{M} Y^{n-1} \widetilde{M}\right)_{i,i+n+1} &= \widetilde{M}_{i,*} \left(Y^{n-1} \widetilde{M}\right)_{*,i+n+1} \\
&= (1, -1, 0, 0, \cdots) \cdot (\widetilde{M}_{n-1,i+n+1}, \widetilde{M}_{n,i+n+1}, \cdots) \\
&= \widetilde{M}_{n-1,i+n+1} - \widetilde{M}_{n,i+n+1}.
\end{aligned}$$

Como $\widetilde{M}_{k,j} = 0$ para $j > k+1$ tenemos:
$\widetilde{M}_{n,i+n+1} = 0$ pues $i+n+1 > n+1$,
$\widetilde{M}_{n-1,i+n+1} = 0$ pues $i+n+1 > n+1 > n = (n-1)+1$.
Esto demuestra que $\left(\widetilde{M} Y^{n-1} \widetilde{M}\right)_{i,i+n+1} = 0$ para $i = 1, \ldots, n-1$.
Por otro lado es evidente que $(Y^{n+1})_{i,i+n+1} = 1$.
Para los términos restantes hay que calcular $\left(Y^j \widetilde{M} Y^{n-j}\right)_{i,i+n+1}$ para $j = 0, \ldots, n$.
Se tiene $\left(Y^j \widetilde{M} Y^{n-j}\right)_{i,i+n+1} = \left(\widetilde{M} Y^{n-j}\right)_{i+j,i+n+1} = \widetilde{M}_{i+j,i+j+1}$.
Ahora afirmamos que

$$\widetilde{M}_{i+j,i+j+1} = m_{n+1} \delta_{j,n-i}. \tag{3.1.7}$$

En efecto se tiene $1 \leq i \leq n-1$ y $0 \leq j \leq n$ (ver sumatoria en (3.1.6)).
Luego $1 \leq i+j \leq 2n-1$. Denotando $k = i+j$ se tiene $1 \leq k \leq 2n-1$ y $\delta_{j,n-i} = \delta_{k,n}$.
Entonces para $1 \leq k < n$ se cumple $\delta_{k,n} = 0$ y $\widetilde{M}_{i+j,i+j+1} = \widetilde{M}_{k,k+1} = 0$ pues $2 \leq k+1 < n+1$.
Para $k = n = i+j$ $(j = n-i)$ se cumple $\widetilde{M}_{n,n+1} = m_{n+1} \delta_{j,j}$.
De esta manera queda por demostrar $\widetilde{M}_{k,k+1} = m_{n+1} \delta_{j,n-i}$ para $n+1 \leq k \leq 2n-1$.
Es suficiente demostrar

$$\widetilde{M}_{n+k,n+k+1} = 0 \quad \text{para } 1 \leq k \leq n-1. \tag{3.1.8}$$

Pero

$$\begin{aligned}
0 &= (\widetilde{M} Y^{k-1} \widetilde{M})_{n,n+k+1} = \widetilde{M}_{n,*} \widetilde{M}_{*+k-1,n+k+1} \\
&= (m_0, m_1, \cdots, m_n, m_{n+1}) \cdot (\widetilde{M}_{k-1,n+k+1}, \widetilde{M}_{k,n+k+1}, \cdots, \widetilde{M}_{n+k-1,n+k+1}, \widetilde{M}_{n+k,n+k+1}) \\
&= (m_0, m_1, \cdots, m_n, m_{n+1}) \cdot (0, 0, \cdots, 0, \widetilde{M}_{n+k,n+k+1}) \\
&= m_{n+1} \widetilde{M}_{n+k,n+k+1}.
\end{aligned}$$

Como $m_{n+1} \neq 0$, entonces $\widetilde{M}_{n+k,n+k+1} = 0$, y así queda demostrado (3.1.8) y por ende (3.1.7).
Finalmente por (3.1.7)

$$(Y^n \widetilde{M})_{i,i+n+1} = \widetilde{M}_{i+n,i+n+1} = 0, \quad i = 1, 2, \cdots, n-1.$$

Reuniendo las entradas en $(i, i+n+1)$ en todos los términos de (3.1.6) obtenemos

$$(Y^n \widetilde{M})_{i,i+n+1} = m_0(\widetilde{M} Y^{n-1} \widetilde{M})_{i,i+n+1} + \sum_{j=0}^{n}\left(\sum_{k=0}^{n-j} m_k\right)(Y^j \widetilde{M} Y^{n-j})_{i,i+n+1}$$

$$+ \left(\sum_{k=0}^{n+1} m_k\right)(Y^{n+1})_{i,i+n+1}$$

$$0 = m_0(0) + \sum_{j=0}^{n}\left(\sum_{k=0}^{n-j} m_k\right) m_{n+1} \delta_{j,n-i} + \left(\sum_{k=0}^{n+1} m_k\right)(Y^{n+1})_{i,i+n+1}$$

$$0 = m_{n+1} \sum_{k=0}^{i} m_k + \sum_{k=0}^{n+1} m_k, \quad i = 1, 2, \cdots, n-1. \tag{3.1.9}$$

Para $i = 1, 2$ se cumple
$$0 = \sum_{k=0}^{n+1} m_k + m_{n+1}(m_0 + m_1),$$
$$0 = \sum_{k=0}^{n+1} m_k + m_{n+1}(m_0 + m_1 + m_2).$$
Luego $m_{n+1} m_2 = 0$ y por lo tanto $m_2 = 0$.
En general para $i = j-1, i = j$ se cumple
$$0 = \sum_{k=0}^{n+1} m_k + m_{n+1}(m_0 + m_1 + \cdots + m_{j-1}),$$
$$0 = \sum_{k=0}^{n+1} m_k + m_{n+1}(m_0 + m_1 + \cdots + m_{j-1} + m_j).$$
Luego $m_{n+1} m_j = 0$ y por lo tanto $m_j = 0$ para $1 < j < n$.
Además para $i = 1$ se obtiene

$$0 = \sum_{k=0}^{n+1} m_k + m_{n+1} \sum_{k=0}^{1} m_k.$$

Pero $\displaystyle\sum_{k=2}^{n-1} m_k = 0$ por lo anterior. Luego

$$0 = m_0 + m_1 + m_n + m_{n+1} + m_{n+1}(m_0 + m_1),$$
$$0 = (m_{n+1} + 1)(m_0 + m_1) + m_n + m_{n+1},$$

lo cual es (3.1.3).

Finalmente, usando el hecho de que $m_j = 0$ para $2 \leq j \leq n-1$ en la igualdad (3.1.6) se obtiene (3.1.2).

En efecto

$$
\begin{aligned}
Y^n \widetilde{M} &= m_0 \widetilde{M} Y^{n-1} \widetilde{M} + \sum_{j=0}^{n} \left( \sum_{k=0}^{n-j} m_k \right) Y^j \widetilde{M} Y^{n-j} + \left( \sum_{k=0}^{n+1} m_k \right) Y^{n+1} \\
&= m_0 \widetilde{M} Y^{n-1} \widetilde{M} + \sum_{j=0}^{n-1} \left( \sum_{k=0}^{n-j} m_k \right) Y^j \widetilde{M} Y^{n-j} + m_0 Y^n \widetilde{M} \\
&\quad + (m_0 + m_1 + m_n + m_{n+1}) Y^{n+1} \\
&= m_0 \widetilde{M} Y^{n-1} \widetilde{M} + \sum_{j=0}^{n-1} \left( \sum_{k=0}^{1} m_k \right) Y^j \widetilde{M} Y^{n-j} + m_n \widetilde{M} Y^n + m_0 Y^n \widetilde{M} \\
&\quad + (m_0 + m_1 + m_n + m_{n+1}) Y^{n+1}.
\end{aligned}
$$

Luego

$$
\begin{aligned}
(1 - m_0) Y^n \widetilde{M} &= m_0 \widetilde{M} Y^{n-1} \widetilde{M} + (m_0 + m_1) \sum_{j=0}^{n-1} Y^j \widetilde{M} Y^{n-j} + m_n \widetilde{M} Y^n \\
&\quad + (m_0 + m_1 + m_n + m_{n+1}) Y^{n+1}. \qquad \square
\end{aligned}
$$

**Proposición 3.1.6.** *Sea $n \in \mathbb{N}$ con $n \geq 2$ y supongamos que $\widetilde{M}_{k,*} = \widetilde{M}_{0,*}$ para $1 < k < n$ y que $\widetilde{M}_{n,*} \neq \widetilde{M}_{0,*}$. Renombramos las únicas posibles entradas diferentes de cero de la siguiente manera:*

$$a := m_{n+1} = \widetilde{M}_{n,n+1}, \quad b := m_n = \widetilde{M}_{n,n}, \quad c := m_1 = \widetilde{M}_{n,1}, \quad d := m_0 = \widetilde{M}_{n,0}.$$

*Entonces se cumple exactamente uno de los dos siguientes casos:*

*(1) $(d, c, b, a) = (d, -d, -a, a)$ con $(d, a) \neq (1, 0)$,*

*(2) $(d, c, b, a) = (d, -1, 0, a)$ con $a \neq 0$ y $d(a+1) = 1$.*

**Demostración.** Consideremos las siguientes afirmaciones:

(i) Si $a = 0$ entonces $b = 0$ y $c = -d$,

(ii) Si $a \neq 0$ y $d = 1$ entonces $b = -a$ y $c = -1$,

(iii) Si $a \neq 0$, $d \neq 1$ y $b \neq 0$ entonces $b = -a$ y $c = -d$,

(iv) Si $a \neq 0$, $d \neq 1$ y $b = 0$ entonces $c = -1$ y $d(a+1) = 1$.

Notemos que en el item (i) tenemos $(d, a) \neq (1, 0)$ pues $\widetilde{M}_{n,*} \neq \widetilde{M}_{0,*}$.
Por un lado (i), (ii), (iii) implican el caso (1) y por otro lado (iv) implica el caso (2).
Como los ítems (i),(ii), (iii), (iv) cubren todos los casos posibles, es suficiente probar estos ítems para mostrar que necesariamente se da uno de los casos (1) o (2).
Demostración de (i)

Asumimos $a = \widetilde{M}_{n,n+1} = 0$.

Entonces

$$\begin{aligned}
0 &= (\widetilde{M}^2)_{n,n} = \widetilde{M}_{n,*}\widetilde{M}_{*,n} \\
&= (\widetilde{M}_{n,0}, \widetilde{M}_{n,1}, \cdots, \widetilde{M}_{n,n}, \widetilde{M}_{n,n+1}) \cdot (\widetilde{M}_{0,n}, \widetilde{M}_{1,n}, \cdots, \widetilde{M}_{n,n}, \widetilde{M}_{n+1,n}) \\
&= (d, c, 0, \cdots, 0, b, 0) \cdot (0, 0, \cdots, 0, b, \widetilde{M}_{n+1,n}) \\
&= b^2.
\end{aligned}$$

Por lo tanto $b = 0$.

Por otro lado

$$\begin{aligned}
0 &= (\widetilde{M}^2)_{n,0} = \widetilde{M}_{n,*}\widetilde{M}_{*,0} \\
&= (\widetilde{M}_{n,0}, \widetilde{M}_{n,1}, \cdots, \widetilde{M}_{n,n}, \widetilde{M}_{n,n+1}) \cdot (\widetilde{M}_{0,0}, \widetilde{M}_{1,0}, \cdots, \widetilde{M}_{n,0}, \widetilde{M}_{n+1,0}) \\
&= (d, c, 0, \cdots, 0, 0, 0) \cdot (1, 1, \cdots, 1, \widetilde{M}_{n,0}, \widetilde{M}_{n+1,0}) \\
&= d + c.
\end{aligned}$$

Por lo tanto $c = -d$.

<u>Demostración de (ii)</u>

Asumimos $a \neq 0$, $\quad d = 1$.

Usando (3.1.2)

$$(1 - m_0)Y^n\widetilde{M} = m_0\widetilde{M}Y^{n-1}\widetilde{M} + (m_0 + m_1)\sum_{j=0}^{n-1} Y^j\widetilde{M}Y^{n-j} + m_n\widetilde{M}Y^n$$

$$+ (m_0 + m_1 + m_n + m_{n+1})Y^{n+1}.$$

Con la nueva notación

$$(1-d)Y^n\widetilde{M} = d\widetilde{M}Y^{n-1}\widetilde{M} + (c+d)\sum_{j=0}^{n-1} Y^j\widetilde{M}Y^{n-j} + b\widetilde{M}Y^n + (d+c+b+a)Y^{n+1}.$$

Reemplazando $d = 1$

$$0 = \widetilde{M}Y^{n-1}\widetilde{M} + (c+1)\sum_{j=0}^{n-1} Y^j\widetilde{M}Y^{n-j} + b\widetilde{M}Y^n + (1+c+b+a)Y^{n+1}.$$

Evaluamos esta igualdad en la entrada $(1, 2)$ y obtenemos

$$0 = (\widetilde{M}Y^{n-1}\widetilde{M})_{1,2} + (c+1)\sum_{j=0}^{n-1}(Y^j\widetilde{M}Y^{n-j})_{1,2} + b(\widetilde{M}Y^n)_{1,2} + (1+c+b+a)(Y^{n+1})_{1,2}.$$

Ahora

$$(\widetilde{M} Y^{n-1}\widetilde{M})_{1,2} = \widetilde{M}_{1,*}(Y^{n-1}\widetilde{M})_{*,2} = \widetilde{M}_{1,*}\widetilde{M}_{*+n-1,2}$$
$$= (\widetilde{M}_{1,0}, \widetilde{M}_{1,1}, 0, \cdots) \cdot (\widetilde{M}_{n-1,2}, \widetilde{M}_{n,2}, \widetilde{M}_{n+1,2}, \cdots)$$
$$= \widetilde{M}_{n-1,2} - \widetilde{M}_{n,2},$$
$$(Y^j \widetilde{M} Y^{n-j})_{1,2} = (\widetilde{M} Y^{n-j})_{j+1,2} = \widetilde{M}_{j+1,2-n+j}.$$

Pero cuando $j < n-2$ se tiene $2-n+j < 0$ y en este caso $\widetilde{M}_{j+1,2-n+j} = 0$.
Luego

$$\sum_{j=0}^{n-1}(Y^j \widetilde{M} Y^{n-j})_{1,2} = \sum_{j=0}^{n-1} \widetilde{M}_{j+1,2-n+j} = \widetilde{M}_{n-1,0} + \widetilde{M}_{n,1},$$
$$(\widetilde{M} Y^n)_{1,2} = \widetilde{M}_{1,2-n},$$
$$(Y^{n+1})_{1,2} = \delta_{n+2,2} = 0.$$

Por lo tanto $\qquad 0 = \widetilde{M}_{n-1,2} - \widetilde{M}_{n,2} + (c+1)(\widetilde{M}_{n-1,0} + \widetilde{M}_{n,1}) + b\widetilde{M}_{1,2-n}.$
Si $n = 2$ entonces $\qquad 0 = \widetilde{M}_{1,2} - \widetilde{M}_{2,2} + (c+1)(\widetilde{M}_{1,0} + \widetilde{M}_{2,1}) + b\widetilde{M}_{1,0}.$
Es decir

$$0 = -b + (c+1)(c+1) + b.$$

Si $n > 2$ entonces $\widetilde{M}_{n-1,2} = \widetilde{M}_{n,2} = \widetilde{M}_{1,2-n} = 0$, y obtenemos directamente que $(c+1)^2 = 0$. Por lo tanto siempre se cumple que $c = -1$.
Reemplazando $d = 1$ y $c = -1$ en (3.1.2) se obtiene

$$0 = \widetilde{M} Y^{n-1}\widetilde{M} + b\widetilde{M} Y^n + (b+a)Y^{n+1}.$$

Multiplicando a la izquierda por $\widetilde{M}$ se obtiene $0 = (b+a)\widetilde{M} Y^{n+1}$.
Como $\widetilde{M} Y^{n+1} \neq 0$, se tiene $b + a = 0$ y de allí $b = -a$.
Demostración de (iii) y (iv)
Asumiremos que $d = m_0 \neq 1$ y demostraremos que

$$\widetilde{M}_{n+1,0} = d, \quad \widetilde{M}_{n+1,1} = -d. \qquad (3.1.10)$$

Para esto evaluaremos (3.1.2) en la entrada $(1,0)$:

$$(1-m_0)(Y^n \widetilde{M})_{1,0} = m_0(\widetilde{M} Y^{n-1}\widetilde{M})_{1,0} + (m_0 + m_1)\sum_{j=0}^{n-1}(Y^j \widetilde{M} Y^{n-j})_{1,0}$$
$$+ m_n(\widetilde{M} Y^n)_{1,0} + (m_0 + m_1 + m_n + m_{n+1})(Y^{n+1})_{1,0}.$$

Ahora

$$(Y^n \widetilde{M})_{1,0} = \widetilde{M}_{n+1,0},$$
$$(\widetilde{M} Y^{n-1} \widetilde{M})_{1,0} = \widetilde{M}_{1,*} \cdot (Y^{n-1} \widetilde{M})_{*,0}$$
$$= (\widetilde{M}_{1,0}, \widetilde{M}_{1,1}, 0, 0, \cdots) \cdot (\widetilde{M}_{n-1,0}, \widetilde{M}_{n,0}, \widetilde{M}_{n+1,0}, \cdots)$$
$$= \widetilde{M}_{n-1,0} - \widetilde{M}_{n,0} = 1 - d$$
$$(Y^j \widetilde{M} Y^{n-j})_{1,0} = (Y^j)_{1,*} \cdot (\widetilde{M} Y^{n-j})_{*,0} = Id_{j+1,*} \cdot \widetilde{M}_{*,-n+j}.$$

Por el ítem c) de la página 16 se tiene $(DY^k)_{1,0} = 0$ para todo $k > 0$, luego $\sum_{j=0}^{n-1}(Y^j \widetilde{M} Y^{n-j})_{1,0} = 0$ y $(\widetilde{M} Y^n)_{1,0} = 0$. Reemplazando obtenemos

$$(1-d)\widetilde{M}_{n+1,0} = d(1-d) + d(0) + 0 = d(1-d).$$

Como $d \neq 1$ entonces $\widetilde{M}_{n+1,0} = d$.

De la misma manera hallaremos $\widetilde{M}_{n+1,1}$ evaluando (3.1.2) en la entrada $(1,1)$:

$$(Y^n \widetilde{M})_{1,1} = \widetilde{M}_{n+1,1}$$
$$(\widetilde{M} Y^{n-1} \widetilde{M})_{1,1} = \widetilde{M}_{1,*} \cdot (Y^{n-1} \widetilde{M})_{*,1}$$
$$= (\widetilde{M}_{1,0}, \widetilde{M}_{1,1}, 0, 0, \cdots) \cdot (\widetilde{M}_{n-1,1}, \widetilde{M}_{n,1}, \widetilde{M}_{n+1,1}, \cdots)$$
$$= \widetilde{M}_{n-1,1} - \widetilde{M}_{n,1} = -1 - c.$$

Por el ítem c) de la página 16 se tiene $(DY^k)_{1,1} = 0$ para todo $k > 1$, luego

$$\sum_{j=0}^{n-1}(Y^j \widetilde{M} Y^{n-j})_{1,1} = (Y^{n-1} \widetilde{M} Y)_{1,1} = \widetilde{M}_{n,0} = d$$

y $(\widetilde{M} Y^n)_{1,1} = 0$.
$$(Y^{n+1})_{1,1} = \delta_{n+2,1} = 0.$$

Reemplazando obtenemos

$$(1-d)\widetilde{M}_{n+1,1} = d(-1-c) + d(c+d) = d(-1-c+c+d) = -d(1-d).$$

Como $d \neq 1$, se tiene $\widetilde{M}_{n+1,1} = -d$.

Esto concluye la demostración de (3.1.10).

Ahora, usando (3.1.10) en $(\widetilde{M}^2)_{n,0} = (\widetilde{M}^2)_{n,1} = 0$:

$$\begin{aligned}
0 &= (\widetilde{M}^2)_{n,0} = \widetilde{M}_{n,*} \cdot \widetilde{M}_{*,0} \\
&= (m_0, m_1, 0, \cdots, 0, m_n, m_{n+1}, 0, \cdots) \cdot (\widetilde{M}_{0,0}, \widehat{M}_{1,0}, \widetilde{M}_{2,0}, \cdots, \widetilde{M}_{n,0}, \widetilde{M}_{n+1,0}, \cdots) \\
&= (d, c, 0, \cdots, 0, b, a, 0, \cdots) \cdot (1, 1, \widetilde{M}_{2,0}, \cdots, d, d, \cdots).
\end{aligned}$$

Luego
$$0 = d + c + bd + ad. \tag{3.1.11}$$

Análogamente

$$\begin{aligned} 0 &= (\widetilde{M}^2)_{n,1} = \widetilde{M}_{n,*} \cdot \widetilde{M}_{*,1} \\ &= (d,c,0,\cdots,0,b,a,0,\cdots) \cdot (\widetilde{M}_{0,1}, \widehat{M}_{1,1}, \widetilde{M}_{2,1}, \cdots, \widetilde{M}_{n,1}, \widetilde{M}_{n+1,1}, \cdots) \\ &= (d,c,0,\cdots,0,b,a,0,\cdots) \cdot (-1,-1,\widetilde{M}_{2,1}, \cdots, c, -d, \cdots). \end{aligned}$$

Luego
$$0 = -d - c + bc - ad. \tag{3.1.12}$$

De (3.1.3) se obtiene
$$(a+1)(d+c) + b + a = 0. \tag{3.1.13}$$

Sumando miembro a miembro (3.1.11) y (3.1.12):
$d + c + bd + ad - d - c + bc - ad = 0$,
$(c+d)b = 0$.
Si $b \neq 0$ entonces $c = -d$ y en este caso $b = -a$ por (3.1.13) lo cual demuestra el ítem (iii).
Si $b = 0$, de (3.1.11) obtenemos:

$$d + c + ad = 0. \tag{3.1.14}$$

De (3.1.13) obtenemos

$$0 = (a+1)(d+c) + b + a = ad + c + d + ac + b + a = ac + a = a(c+1).$$

Como $a \neq 0$, obtenemos $c = -1$.
Finalmente, (3.1.14) implica $d - 1 + ad = 0$, es decir $d(a+1) = 1$, lo cual concluye la demostración de (iv). □

## 3.2. La familia $A(n,d,a)$

En esta sección describiremos el caso (1) de la Proposición 3.1.6. Demostraremos que la aplicación de torcimiento depende solamente de $n, d, a$. Obtenemos la familia de productos tensoriales $A(n,d,a)$ parametrizada por $n \in \mathbb{N}$, $n \geq 2$ y $(a,d) \in K^2$ tales que para la familia infinita de polinomios $R_k$ (ver Definición 3.2.5) tenemos $R_k(a,d) \neq 0$.

**Proposición 3.2.1.** *Sean A una K-álgebra asociativa, $a, d \in K$ con $(d,a) \neq (1,0)$, $\widetilde{M}, Y \in A$ tales que*

$$d\widetilde{M} Y^{n-1} \widetilde{M} = e Y^n \widetilde{M} + a\widetilde{M} Y^n, \tag{3.2.15}$$

*donde $e := 1-d$.*

*Entonces para todo $k \geq 1$ tenemos:*

$$d_k \widetilde{M} Y^{kn-1} \widetilde{M} = e_k Y^{kn} \widetilde{M} - a_k \widetilde{M} Y^{kn}, \qquad (3.2.16)$$

*donde $e_k := e^k$, $\quad a_k := (-a)^k$, $\quad d_k := d \sum_{j=0}^{k-1} e_j a_{k-1-j}$*

**Demostración.** Si $d = 0$ entonces $e = 1 - d = 1$ y reemplazando en (3.2.15) se tiene $Y^n \widetilde{M} = -a \widetilde{M} Y^n$.

Luego

$$\begin{aligned} Y^{2n} \widetilde{M} &= Y^n (Y^n \widetilde{M}) = Y^n (-a \widetilde{M} Y^n) = -a Y^n \widetilde{M} Y^n = -a(-a \widetilde{M} Y^n) Y^n \\ &= (-a)^2 \widetilde{M} Y^{2n}, \\ Y^{3n} \widetilde{M} &= Y^n (Y^{2n} \widetilde{M}) = Y^n ((-a)^2 \widetilde{M} Y^{2n}) = (-a)^2 Y^n \widetilde{M} Y^{2n} = (-a)^2 (-a \widetilde{M} Y^n) Y^{2n} \\ &= (-a)^3 \widetilde{M} Y^{3n}. \end{aligned}$$

En general se cumple

$$Y^{kn} \widetilde{M} = (-a)^k \widetilde{M} Y^{kn}.$$

Esto es precisamente (3.2.16) en este caso.

Ahora supongamos que $d \neq 0$.

Por hipótesis sabemos que $d_k = d \sum_{j=0}^{k-1} e_j a_{k-1-j}$.

Luego

$$d_{k+1} = d \sum_{j=0}^{k} e_j a_{k-j} = -da \sum_{j=0}^{k-1} e_j a_{k-1-j} + d e_k = -a d_k + d e_k.$$

Por lo tanto

$$d_{k+1} = -a d_k + d e_k. \qquad (3.2.17)$$

Demostraremos (3.2.16) por inducción sobre $k$.

Para $k = 1$ la igualdad (3.2.15) es justamente (3.2.16).

Ahora supongamos que (3.2.16) es válida para algún $k \geq 1$, es decir

$$d_k \widetilde{M} Y^{kn-1} \widetilde{M} = e_k Y^{kn} \widetilde{M} - a_k \widetilde{M} Y^{kn}.$$

Multiplicando esta igualdad a la izquierda por $\widetilde{M} Y^{n-1}$:

$$d_k \widetilde{M} Y^{n-1} \widetilde{M} Y^{kn-1} \widetilde{M} = e_k \widetilde{M} Y^{(k+1)n-1} \widetilde{M} - a_k \widetilde{M} Y^{n-1} \widetilde{M} Y^{kn}.$$

Usando (3.2.15) obtenemos:

$$d_k\left(\frac{e}{d}Y^n\widetilde{M}+\frac{a}{d}\widetilde{M}Y^n\right)Y^{kn-1}\widetilde{M}=e_k\widetilde{M}Y^{(k+1)n-1}\widetilde{M}-a_k\left(\frac{e}{d}Y^n\widetilde{M}+\frac{a}{d}\widetilde{M}Y^n\right)Y^{kn},$$

$$\frac{d_k e}{d}Y^n\widetilde{M}Y^{kn-1}\widetilde{M}+\frac{d_k a}{d}\widetilde{M}Y^{(k+1)n-1}\widetilde{M}=e_k\widetilde{M}Y^{(k+1)n-1}\widetilde{M}$$
$$-\frac{a_k e}{d}Y^n\widetilde{M}Y^{kn}-\frac{a_k a}{d}\widetilde{M}Y^{(k+1)n},$$

$$\left(e_k-\frac{d_k a}{d}\right)\widetilde{M}Y^{(k+1)n-1}\widetilde{M}=\frac{d_k e}{d}Y^n\widetilde{M}Y^{kn-1}\widetilde{M}+\frac{a_k e}{d}Y^n\widetilde{M}Y^{kn}$$
$$+\frac{a_k a}{d}\widetilde{M}Y^{(k+1)n}.$$

Multiplicando por $d\neq 0$ obtenemos

$$(de_k-d_k a)\widetilde{M}Y^{(k+1)n-1}\widetilde{M}=d_k e\,Y^n\widetilde{M}Y^{kn-1}\widetilde{M}+a_k e\,Y^n\widetilde{M}Y^{kn}+a_k a\widetilde{M}Y^{(k+1)n}.$$

Usando la hipótesis inductiva en el primer término del segundo miembro:

$$(de_k-d_k a)\widetilde{M}Y^{(k+1)n-1}\widetilde{M} = eY^n\left(e_k Y^{kn}\widetilde{M}-a_k\widetilde{M}Y^{kn}\right)+a_k e\,Y^n\widetilde{M}Y^{kn}$$
$$+a_k a\widetilde{M}Y^{(k+1)n}.$$

Simplificando:

$$(de_k-d_k a)\widetilde{M}Y^{(k+1)n-1}\widetilde{M}=e\,e_k Y^{(k+1)n}\widetilde{M}+a_k a\widetilde{M}Y^{(k+1)n}.$$

Por (3.2.17) se tiene $d_{k+1}=e_k d-ad_k$.
Además $e\,e_k=e_{k+1}$ y $aa_k=a(-a)^k=-(-a)^{k+1}=-a_{k+1}$.
Reemplazando en la última igualdad obtenemos

$$d_{k+1}\widetilde{M}Y^{(k+1)n-1}\widetilde{M}=e_{k+1}Y^{(k+1)n}\widetilde{M}-a_{k+1}\widetilde{M}Y^{(k+1)n}.$$

Esto completa el paso de inducción y por lo tanto la demostración. $\square$

**Observación 3.2.2.** Notemos que

$$d_k=\begin{cases} d\dfrac{e^k-(-a)^k}{e+a} & \text{si } e\neq -a \\ kde^{k-1} & \text{si } e=-a \end{cases}.$$

En efecto:
$\boxed{e\neq -a}$
$$\sum_{j=0}^{k-1}e_j a_{k-1-j}=\sum_{j=0}^{k-1}e_j(-a)^{k-1-j}=(-a)^{k-1}\sum_{j=0}^{k-1}\left(-\frac{e}{a}\right)^j=(-a)^{k-1}\frac{1-(-\frac{e}{a})^k}{1+\frac{e}{a}}=\frac{e^k-(-a)^k}{e+a}.$$
$\boxed{e=-a}$
$$\sum_{j=0}^{k-1}(-a)^j a_{k-1-j}=(-a)^{k-1}\sum_{j=0}^{k-1}1=k(-a)^{k-1}=ke^{k-1}.$$

**Teorema 3.2.3.** *Sea $\sigma$ una aplicación de torcimiento y supongamos que $Y$, $M$ son como en la Proposición 1.4.1 y asumamos que $\widetilde{M} = M - Y$.*

*Sea $n \in \mathbb{N}$ con $n \geq 2$, tomemos $a, d \in K$ con $(d, a) \neq (1, 0)$ y supongamos que $\widetilde{M}_{k,*} = \widetilde{M}_{0,*}$ para $k < n$ y que*

$$d\widetilde{M} Y^{n-1}\widetilde{M} = e Y^n \widetilde{M} + a\widetilde{M} Y^n, \qquad (3.2.18)$$

*con $e = 1 - d$.*

*Definamos*

$$d_k := d \sum_{j=0}^{k-1} e_j a_{k-1-j}, \quad e_k := e^k, \quad a_k := (-a)^k.$$

*Entonces:*

*(1)*
$$d_k + e_k \neq 0 \text{ para todo } k \in \mathbb{N}, \qquad (3.2.19)$$

*(2)* $\widetilde{M}_{kn+j,0} = -\widetilde{M}_{kn+j,1} = \prod_{i=1}^{k} \dfrac{d_i}{d_i + e_i}$ *para $k \geq 1$, $0 \leq j \leq n-1$.*

*(3)* $\widetilde{M}_{kn+j,rn} = -\widetilde{M}_{kn+j,rn+1} = \dfrac{a_r}{e_r + d_r} \prod_{i=r+1}^{k} \dfrac{d_i}{d_i + e_i}$ *para $1 \leq r \leq k$, $0 \leq j \leq n-1$, y las otras entradas son ceros.*

**Demostración.** Primero demostraremos que

$$\widetilde{M}_{kn,*} = \widetilde{M}_{kn+j,*} \qquad (3.2.20)$$

para $k \geq 1$, $0 < j < n$.

En efecto, como (3.2.18) se cumple, entonces por la Proposición 3.2.1 se cumple (3.2.16):

$$d_k \widetilde{M} Y^{kn-1} \widetilde{M} = e_k Y^{kn} \widetilde{M} - a_k \widetilde{M} Y^{kn}.$$

Sea $i \in \{0, 1, \cdots, n-2\}$.

Si $e \neq 0$ multiplicamos a la izquierda por $\widetilde{M} Y^i$ y obtenemos:

$$d_k \widetilde{M} Y^i \widetilde{M} Y^{kn-1} \widetilde{M} = e_k \widetilde{M} Y^{kn+i} \widetilde{M} - a_k \widetilde{M} Y^i \widetilde{M} Y^{kn}.$$

Ahora, como $\widetilde{M}_{k,*} = \widetilde{M}_{0,*}$ para $k < n$ entonces por el Lema 3.1.3 $\widetilde{M} Y^i \widetilde{M} = 0$ para $0 \leq i \leq n-2$. Reemplazando:

$$0 = e_k \widetilde{M} Y^{kn+i} \widetilde{M}.$$

Pero $e_k = e^k \neq 0$, luego

$$\widetilde{M} Y^{kn+i} \widetilde{M} = 0, \quad 0 \leq i \leq n-2.$$

Similarmente, si $a \neq 0$, entonces multiplicando (3.2.16) a la derecha por $Y^i \widetilde{M}$ obtenemos

$$d_k \widetilde{M} Y^{kn-1} \widetilde{M} Y^i \widetilde{M} = e_k Y^{kn} \widetilde{M} Y^i \widetilde{M} - a_k \widetilde{M} Y^{kn+i} \widetilde{M}.$$

Ahora, como $\widetilde{M}_{k,*} = \widetilde{M}_{0,*}$ para $k < n$ entonces por el Lema 3.1.3 $\widetilde{M} Y^i \widetilde{M} = 0$ para $0 \leq i \leq n-2$. Reemplazando:

$$0 = -a_k \widetilde{M} Y^{kn+i} \widetilde{M}.$$

Pero $a_k \neq 0$ pues $a \neq 0$, luego

$$\widetilde{M} Y^{kn+i} \widetilde{M} = 0, \quad 0 \leq i \leq n-2.$$

Dado que $(e, a) \neq (0, 0)$, pues $(d, a) \neq (1, 0)$, obtenemos

$$\widetilde{M} Y^{kn+i} \widetilde{M} = 0, \quad 0 \leq i \leq n-2.$$

Por el Lema 3.1.2 aplicado a $\widetilde{M}_{kn,*}$ se tiene (3.2.20).
En efecto:
Para $i = 0$: $\widetilde{M} Y^{kn} \widetilde{M} = 0 \Leftrightarrow \widetilde{M}_{kn,*} = \widetilde{M}_{kn+1,*}$
Para $i = 1$: $\widetilde{M} Y^{kn+1} \widetilde{M} = 0 \Leftrightarrow \widetilde{M}_{kn+1,*} = \widetilde{M}_{kn+2,*}$.
Repitiendo esta acción sucesivamente se llega a:
Para $i = n-2$: $\widetilde{M} Y^{kn+n-2} \widetilde{M} = 0 \Leftrightarrow \widetilde{M}_{kn+n-2,*} = \widetilde{M}_{kn+n-1,*}$.
Por lo tanto $\widetilde{M}_{kn+j,*} = \widetilde{M}_{kn,*}$ para $0 \leq j \leq n-1$.
Por la Observación 3.1.1 sabemos que

$$Y^{kn} \widetilde{M} = \sum_{i=0}^{kn+1} \widetilde{M}_{kn,i} M^{kn+1-i} Y^i.$$

Queremos demostrar que las únicas posibles entradas diferentes de cero de $\widetilde{M}_{kn,*}$ son

$$\widetilde{M}_{kn,0}, \widetilde{M}_{kn,1}, \widetilde{M}_{kn,n}, \widetilde{M}_{kn,n+1}, \widetilde{M}_{kn,2n}, \widetilde{M}_{kn,2n+1}, \cdots, \widetilde{M}_{kn,rn}, \widetilde{M}_{kn,rn+1}, \cdots,$$
$$\widetilde{M}_{kn,kn}, \widetilde{M}_{kn,kn+1}.$$

Además queremos demostrar que $\widetilde{M}_{kn,rn} = -\widetilde{M}_{kn,rn+1}$ para $0 \leq r \leq k$.
Si llamamos $C_{k,r} = \widetilde{M}_{kn,rn}$, para $0 \leq r \leq k$, entonces debemos demostrar que

$$\widetilde{M}_{kn,*} = (C_{k,0}, -C_{k,0}, 0, \cdots, 0, C_{k,1}, -C_{k,1}, 0, \cdots, 0, C_{k,r}, -C_{k,r}, 0, \cdots, 0, C_{k,k}, -C_{k,k})$$

ó equivalentemente

$$\begin{aligned} Y^{kn} \widetilde{M} &= C_{k,0} M^{kn+1} - C_{k,0} M^{kn} Y + C_{k,1} M^{(k-1)n+1} Y^n - C_{k,1} M^{(k-1)n} Y^{n+1} + \cdots \\ &\quad + C_{k,r} M^{(k-r)n+1} Y^{rn} - C_{k,r} M^{(k-r)n} Y^{rn+1} + \cdots + C_{k,k} M Y^{kn} - C_{k,k} Y^{kn+1} \\ &= C_{k,0} M^{kn}(M - Y) + C_{k,1} M^{(k-1)n}(M - Y) Y^n + \cdots + C_{k,r} M^{(k-r)n}(M - Y) Y^{rn} \\ &\quad + \cdots + C_{k,k}(M - Y) Y^{kn}. \end{aligned}$$

Entonces para demostrar el teorema es suficiente demostrar (3.2.19) y que

$$Y^{kn}\widetilde{M} = \sum_{r=0}^{k} C_{k,r} M^{(k-r)n} \widetilde{M} Y^{rn}, \qquad (3.2.21)$$

donde $C_{k,0} = \prod_{i=1}^{k} \dfrac{d_i}{d_i + e_i}$, $\quad C_{k,r} = \dfrac{a_r}{e_r + d_r} \prod_{i=r+1}^{k} \dfrac{d_i}{d_i + e_i}$ para $1 \leq r \leq k$.

Primero asumiremos que $d, a \neq 0$ y demostraremos (3.2.19) y (3.2.21) por inducción sobre $k$.

Para $k = 1$ la igualdad (3.2.21) que debemos demostrar es

$$Y^n \widetilde{M} = \sum_{r=0}^{1} C_{1,r} M^{(1-r)n} \widetilde{M} Y^{rn} = C_{1,0} M^n \widetilde{M} + C_{1,1} \widetilde{M} Y^n.$$

Demostraremos

$$Y^n \widetilde{M} = d M^n \widetilde{M} - a \widetilde{M} Y^n,$$

pues

$$d = \prod_{i=1}^{1} \dfrac{d_i}{d_i + e_i} = C_{1,0},$$

$$-a = \dfrac{a_1}{d_1 + e_1} \prod_{i=2}^{1} \dfrac{d_i}{d_i + e_i} = C_{1,1}.$$

Para esto vemos que $\widetilde{M}_{n-1,*} = \widetilde{M}_{0,*} = (1, -1, 0, \cdots)$ y entonces $Y^{n-1} \widetilde{M} = M^{n-1} \widetilde{M}$.
Multiplicando a la izquierda por $M$ obtenemos $M Y^{n-1} \widetilde{M} = M^n \widetilde{M}$.
Luego $M^n \widetilde{M} = M Y^{n-1} \widetilde{M} = (\widetilde{M} + Y) Y^{n-1} \widetilde{M}$,
$M^n \widetilde{M} = \widetilde{M} Y^{n-1} \widetilde{M} + Y^n \widetilde{M}$.
Usando (3.2.18) obtenemos:

$$M^n \widetilde{M} = \dfrac{e}{d} Y^n \widetilde{M} + \dfrac{a}{d} \widetilde{M} Y^n + Y^n \widetilde{M},$$

de manera que

$$\dfrac{e+d}{d} Y^n \widetilde{M} = M^n \widetilde{M} - \dfrac{a}{d} \widetilde{M} Y^n,$$

y como $e + d = 1 \neq 0$ se obtiene

$$Y^n \widetilde{M} = d M^n \widetilde{M} - a \widetilde{M} Y^n,$$

que es lo planteado.

Ahora fijemos $k$ y supongamos que (3.2.21) se cumple para $k-1$ y que $e_r + d_r \neq 0$ para $r < k$.

Es decir se cumple

$$Y^{(k-1)n} \widetilde{M} = \sum_{r=0}^{k-1} C_{k-1,r} M^{(k-1-r)n} \widetilde{M} Y^{rn},$$

y $e_r + d_r \neq 0$, para $r < k$.

Demostraremos

$$Y^{kn}\widetilde{M} = \sum_{r=0}^{k} C_{k,r} M^{(k-r)n} \widetilde{M} Y^{rn}$$

y $e_k + d_k \neq 0$.

Por (3.2.16)

$$d_k \widetilde{M} Y^{kn-1} \widetilde{M} = e_k Y^{kn} \widetilde{M} - a_k \widetilde{M} Y^{kn}.$$

$$\begin{aligned} e_k Y^{kn} \widetilde{M} &= d_k \widetilde{M} Y^{kn-1} \widetilde{M} + a_k \widetilde{M} Y^{kn} \\ &= d_k (M - Y) Y^{kn-1} \widetilde{M} + a_k \widetilde{M} Y^{kn} \\ &= d_k M Y^{kn-1} \widetilde{M} - d_k Y^{kn} \widetilde{M} + a_k \widetilde{M} Y^{kn}. \end{aligned}$$

Luego

$$(e_k + d_k) Y^{kn} \widetilde{M} = d_k M Y^{kn-1} \widetilde{M} + a_k \widetilde{M} Y^{kn}. \tag{3.2.22}$$

Sabemos por (3.2.20) que $\widetilde{M}_{(k-1)n,*} = \widetilde{M}_{kn-1,*}$ pues $kn-1 = (k-1)n + j$, con $j = n-1$. Es decir los coeficientes de $Y^{(k-1)n}\widetilde{M}$, los cuales son las entradas de $\widetilde{M}$ en la fila $(k-1)n$, son iguales a los coeficientes de $Y^{kn-1}\widetilde{M}$, los cuales son las entradas de $\widetilde{M}$ en la fila $kn-1$. Por lo tanto, por la fórmula (3.1.1) se obtiene :

$$Y^{kn-1}\widetilde{M} = M^{n-1} Y^{(k-1)n} \widetilde{M},$$

pues $n - 1 + (k-1)n = n - 1 + kn - n = kn - 1$.

Usando la hipótesis inductiva se tiene que

$$Y^{kn-1}\widetilde{M} = M^{n-1} \sum_{r=0}^{k-1} C_{k-1,r} M^{(k-1-r)n} \widetilde{M} Y^{rn} = \sum_{r=0}^{k-1} C_{k-1,r} M^{nk-1-rn} \widetilde{M} Y^{rn}.$$

Reemplazando en (3.2.22) :

$$(e_k + d_k) Y^{kn} \widetilde{M} = d_k \sum_{r=0}^{k-1} C_{k-1,r} M^{(k-r)n} \widetilde{M} Y^{rn} + a_k \widetilde{M} Y^{kn}. \tag{3.2.23}$$

Evaluemos la igualdad (3.2.23) en la entrada $(0, kn)$ y obtendremos $e_k + d_k \neq 0$.

En efecto: para $p \geq 1$

$$\begin{aligned} (M^p)_{0,j} &= (MM^{p-1})_{0,j} = M_{0,*}(M^{p-1})_{*,j} \\ &= (1,0,\cdots).((M^{p-1})_{0,j},(M^{p-1})_{1,j},\cdots) = (M^{p-1})_{0,j} = (MM^{p-2})_{0,j} \\ &= M_{0,*}(M^{p-2})_{*,j} = (1,0,\cdots).((M^{p-2})_{0,j},(M^{p-2})_{1,j},\cdots) \\ &= (M^{p-2})_{0,j} = \cdots = M_{0,j} = \delta_{0,j}. \end{aligned}$$

Luego haciendo $p = n(k - r)$

$(0 \leq r < k < n \Rightarrow k - r > 0 \Rightarrow k - r \geq 1 \Rightarrow p = n(k-r) \geq n \geq 2)$

se cumple $(M^{n(k-r)})_{0,*} = (1,0,0,\cdots)$.

Por lo tanto

$$\begin{aligned}(M^{n(k-r)}\widetilde{M}Y^{rn})_{0,kn} &= (M^{n(k-r)})_{0,*}(\widetilde{M}Y^{nr})_{*,kn} = (M^{n(k-r)})_{0,*}(\widetilde{M})_{*,(k-r)n} \\ &= (1,0,0,\cdots).(\widetilde{M}_{0,(k-r)n},\widetilde{M}_{1,(k-r)n},\widetilde{M}_{2,(k-r)n},\cdots) \\ &= \widetilde{M}_{0,(k-r)n} = 0,\end{aligned}$$

pues $(k-r)n \geq n \geq 2$.

Entonces se tiene $\sum_{r=0}^{k-1} C_{k-1,r}(M^{(k-r)n}\widetilde{M}Y^{rn})_{0,kn} = 0$.

Por lo tanto en (3.2.23)

$$(e_k + d_k)(Y^{kn}\widetilde{M})_{0,kn} = a_k(\widetilde{M}Y^{kn})_{0,kn} = a_k\widetilde{M}_{0,0} = a_k.$$

Se sigue que $e_k + d_k \neq 0$.

La igualdad (3.2.23) también implica que (3.2.21) es válido para $k$ haciendo

$$C_{k,k} = \frac{a_k}{d_k + e_k} \quad C_{k,r} = \frac{d_k C_{k-1,r}}{d_k + e_k} \quad \text{para } r < k.$$

En efecto:

De (3.2.23) se tiene

$$\begin{aligned}Y^{kn}\widetilde{M} &= \frac{d_k}{e_k + d_k}\sum_{r=0}^{k-1} C_{k-1,r} M^{(k-r)n}\widetilde{M}Y^{rn} + \frac{a_k}{e_k+d_k}\widetilde{M}Y^{kn} \\ &= \sum_{r=0}^{k-1} \underbrace{\frac{d_k C_{k-1,r}}{e_k + d_k}}_{C_{k,r}} M^{(k-r)n}\widetilde{M}Y^{rn} + \underbrace{\frac{a_k}{e_k+d_k}}_{C_{k,k}}\widetilde{M}Y^{kn} \\ &= \sum_{r=0}^{k} C_{k,r} M^{(k-r)n}\widetilde{M}Y^{rn}.\end{aligned}$$

Esto completa el paso de inducción y concluye la demostración en el caso $a,d \neq 0$.

Caso $\boxed{d=0}$

En este caso para todo $k$ tenemos

$$d_k := d\sum_{j=0}^{k-1} e_j a_{k-1-j} = 0, \quad e_k := e^k = 1, \quad d_k + e_k = 1 \neq 0.$$

De (3.2.18) se obtiene $Y^n\widetilde{M} = -a\widetilde{M}Y^n$.

Además

$$\begin{aligned}C_{k,0} &= \prod_{i=1}^{k}\frac{d_i}{d_i+e_i} = 0, \\ C_{k,r} &= \frac{a_r}{e_r + d_r}\prod_{i=r+1}^{k}\frac{d_i}{d_i+e_i} = a_r\prod_{i=r+1}^{k}\frac{d_i}{1} = a_r(0) = 0, \ (r \leq k-1), \\ C_{k,k} &= a_k(1) = a_k.\end{aligned}$$

En (3.2.21) es suficiente demostrar $Y^{kn}\widetilde{M} = C_{k,k}\widetilde{M}Y^{kn}$.

Esto es equivalente a demostrar $Y^{kn}\widetilde{M} = a_k \widetilde{M} Y^{kn}$, lo cual se sigue de la Proposición 3.2.1.

Finalmente tenemos el

Caso $\boxed{a=0}$

En este caso $d_k = d\dfrac{e^k-(-a)^k}{e+a} = d\dfrac{e^k}{e} \Leftrightarrow \dfrac{e_k}{d_k} = \dfrac{1-d}{d}$.

Por (3.2.16)

$$d_k \widetilde{M} Y^{kn-1}\widetilde{M} = e_k Y^{kn}\widetilde{M},$$

$$\widetilde{M} Y^{kn-1}\widetilde{M} = \frac{e_k}{d_k} Y^{kn}\widetilde{M},$$

$$\widetilde{M} Y^{kn-1}\widetilde{M} = \frac{1-d}{d} Y^{kn}\widetilde{M}.$$

Por el Lema 3.1.2 (2)

$$\widetilde{M}_{kn,*} = d\widetilde{M}_{kn-1,*}.$$

Por (3.2.20)

$$\widetilde{M}_{kn,*} = d\widetilde{M}_{kn-1,*} = d\widetilde{M}_{(k-1)n,*},$$

y así (3.2.21) se cumple, dado que

$$C_{k,0} = d^k, \quad C_{k,r} = 0, \ r > 0.$$

En efecto, $d_k = de^{k-1}$, $d_k + e_k = de^{k-1} + e^k = e^{k-1}(d+e) = e^{k-1}$, $\dfrac{d_k}{e_k + d_k} = d$.

Entonces

$$C_{k,0} = \prod_{i=1}^{k} \frac{d_i}{d_i+e_i} = \prod_{i=1}^{k} d = d^k.$$

También tenemos

$$C_{k,j} = \underbrace{\frac{a_j}{e_j+d_j}}_{0} \prod_{i=j+1}^{k} \frac{d_i}{d_i+e_i} = 0,$$

lo cual concluye la demostración ya que $d_k + e_k = e^{k-1} \neq 0$.

(Notemos que si $e=0$ entonces $d = 1-e = 1$ lo cual contradice el supuesto $(d,a) \neq (1,0)$). $\square$

**Ejemplo 3.2.4.** Ilustraremos el Teorema 3.2.3 para el caso $n=3$ y para $k=1,2$.

Para $\boxed{k=1}$

1) $\widetilde{M}_{3,0} = \widetilde{M}_{4,0} = \widetilde{M}_{5,0}$,

$$\widetilde{M}_{3,0} = \prod_{i=1}^{1} \frac{d_i}{d_i+e_i} = d.$$

2) En este caso $r = 1$, luego se tiene
$\widetilde{M}_{3,3} = \widetilde{M}_{4,3} = \widetilde{M}_{5,3}$,
$$\widetilde{M}_{3,3} = \frac{a_1}{e_1 + d_1} \prod_{i=2}^{1} \frac{d_i}{d_i + e_i} = -a.$$

Entonces
$$\widetilde{M}_{3,*} = \widetilde{M}_{4,*} = \widetilde{M}_{5,*} = (d, -d, 0, -a, a, 0, \cdots).$$

Para $\boxed{k = 2}$
1) $\widetilde{M}_{6,0} = \widetilde{M}_{7,0} = \widetilde{M}_{8,0}$,
$$\widetilde{M}_{6,0} = \prod_{i=1}^{2} \frac{d_i}{d_i + e_i} = d\left(\frac{d_2}{d_2 + e_2}\right).$$

Asumiendo que $e \neq -a$, por la Observación 3.2.2 se tiene

$$d_2 = d\left(\frac{e^2 - a^2}{e + a}\right) = d(e - a), \, d_2 + e_2 = d(e - a) + e^2 = e^2 + de - da = e - da = 1 - d - da.$$

Reemplazando
$$\widetilde{M}_{6,0} = \frac{d^2(1 - d - a)}{1 - d - da}.$$

2) En este caso $r = 1, 2$.
Para $\boxed{r = 1}$
$\widetilde{M}_{6,3} = \widetilde{M}_{7,3} = \widetilde{M}_{8,3}$,
$$\widetilde{M}_{6,3} = \frac{a_1}{e_1 + d_1} \prod_{i=2}^{2} \frac{d_i}{d_i + e_i} = -a \frac{d_2}{d_2 + e_2} = \frac{-ad(1 - d - a)}{1 - d - da}.$$

Para $\boxed{r = 2}$
$\widetilde{M}_{6,6} = \widetilde{M}_{7,6} = \widetilde{M}_{8,6}$,
$$\widetilde{M}_{6,6} = \frac{a_2}{e_2 + d_2} \prod_{i=3}^{2} \frac{d_i}{d_i + e_i} = \frac{a^2}{1 - d - da}.$$

Entonces
$$\begin{aligned}
\widetilde{M}_{6,*} &= \widetilde{M}_{7,*} = \widetilde{M}_{8,*} \\
&= \left(\frac{d^2(1 - d - a)}{1 - d - da}, -\frac{d^2(1 - d - a)}{1 - d - da}, 0, \frac{-ad(1 - d - a)}{1 - d - da}, \right. \\
&\quad \left. \frac{ad(1 - d - a)}{1 - d - da}, 0, \frac{a^2}{1 - d - da}, \frac{-a^2}{1 - d - da}, 0, \cdots \right).
\end{aligned}$$

**Definición 3.2.5.** *Para $a, d \in K$ y $k \in \mathbb{N}$ definimos el polinomio*

$$R_k(a, d) = (1 - d)^k + d \sum_{j=0}^{k-1} (1 - d)^j (-a)^{k-1-j}.$$

Notemos que $R_k(a,d) = e_k + d_k$, donde $e_k, d_k$ están definidos en el Teorema 3.2.3.

Tenemos $R_1(a,d) = e_1 + d_1 = 1 - d + d \sum_{j=0}^{0} (1-d)^j (-a)^{-j} = 1 - d + d(1) = 1 \neq 0$.

$$R_2(a,d) = e_2 + d_2 = (1-d)^2 + d \sum_{j=0}^{1} (1-d)^j (-a)^{1-j}$$
$$= (1-d)^2 + d(-a) + d(1-d) = 1 - d - ad = 1 - (a+1)d \neq 0.$$

Luego $R_2(a,d) \neq 0 \Leftrightarrow d(a+1) \neq 1$.

En particular $(a,d) \neq (0,1)$ si requerimos que $R_2(a,d) \neq 0$.

**Lema 3.2.6.** *Supongamos que $k \geq 1$, que $\widetilde{M}$ es una matriz infinita con $\widetilde{M}_{i,j} = 0$ para $j > i+1$ y que cumple (3.1.1) para $k-1$. Si se cumple $\widetilde{M} Y^{k-1} \widetilde{M} = 0$ y $\widetilde{M}_{k,*} = \widetilde{M}_{k-1,*}$, entonces $\widetilde{M}$ cumple (3.1.1) para $k$.*

**Demostración.** En efecto,

$$\begin{aligned} Y^k \widetilde{M} &= (M - \widetilde{M}) Y^{k-1} \widetilde{M} = M Y^{k-1} \widetilde{M} - \underbrace{\widetilde{M} Y^{k-1} \widetilde{M}}_{0} \\ &= M \sum_{j=0}^{k} \widetilde{M}_{k-1,j} M^{k-j} Y^j = \sum_{j=0}^{k} \widetilde{M}_{k-1,j} M^{k+1-j} Y^j \\ &= \sum_{j=0}^{k} \widetilde{M}_{k,j} M^{k+1-j} Y^j + \underbrace{\widetilde{M}_{k,k+1} Y^{k+1}}_{0} \\ &= \sum_{j=0}^{k+1} \widetilde{M}_{k,j} M^{k+1-j} Y^j. \end{aligned}$$

$\square$

**Corolario 3.2.7.** *Sean $a, d \in K$ tales que $R_k(a,d) \neq 0$ para todo $k$. Entonces las fórmulas del Teorema 3.2.3 definen una matriz $M$ que determina una aplicación de torcimiento vía la Proposición 1.4.1.*

**Demostración.** Observemos que

$$M_{i,j} = \begin{cases} \widetilde{M}_{i,j} & \text{si } i \geq j \\ \widetilde{M}_{i,j} + 1 & \text{si } i + 1 = j \end{cases}.$$

Entonces $M_{0,*} = (1, 0, 0, \cdots)$ pues $\widetilde{M}_{0,*} = (1, -1, 0, \cdots)$.

Por la Proposición 1.4.1 tenemos que demostrar que $M$ cumple

$$Y^k M = \sum_{j=0}^{k+1} M_{k,j} M^{k+1-j} Y^j \tag{3.2.24}$$

para todo $k \geq 0$.

Es evidente que se cumple para $k = 0$.

Consideremos (3.2.24) para algún $1 \leq k \leq n-1$.

Sabemos por hipótesis que $\widetilde{M}_{k,*} = \widetilde{M}_{0,*} = (1, -1, 0, \cdots)$.

Luego $M_{k,*} = (1, -1, 0, \cdots, 0, \underbrace{1}_{k+1}, 0, \cdots)$.

Es decir
$$Y^k M = M^{k+1} - M^k Y + Y^{k+1},$$

o, equivalentemente, $Y^k \widetilde{M} = M^k \widetilde{M}$.

Por el Lema 3.2.6, es suficiente demostrar que $\widetilde{M} Y^{k-1} \widetilde{M} = 0$ para $1 \leq k \leq n-1$.

Demostraremos algo más general:
$$\widetilde{M} Y^{kn+j} \widetilde{M} = 0 \text{ para } k \geq 0 \text{ y } 0 \leq j \leq n-2. \qquad (3.2.25)$$

En efecto por las fórmulas dadas en el Teorema 3.2.3 se cumple que las únicas entradas diferentes de cero en $\widetilde{M}_{k,*}$ son de la forma $\widetilde{M}_{k,rn}$ ó $\widetilde{M}_{k,rn+1}$ para algún $r \geq 0$ con $\widetilde{M}_{k,rn+1} = -\widetilde{M}_{k,rn}$.

Además se cumple $\widetilde{M}_{kn+j,l} = \widetilde{M}_{kn+j+1,l}$ para $k \geq 1$, $0 \leq j \leq n-2$.

Haciendo $k' = r + k \geq 1$ se tiene que $\widetilde{M}_{k'n+j,l} = \widetilde{M}_{k'n+j+1,l}$.

Es decir $\widetilde{M}_{(k+r)n+j,l} = \widetilde{M}_{(k+r)n+j+1,l}$, lo cual es equivalente a
$(Y^{kn+j}\widetilde{M})_{rn,l} = (Y^{kn+j}\widetilde{M})_{rn+1,l}$.

Calculamos

$$\begin{aligned}(\widetilde{M} Y^{kn+j} \widetilde{M})_{t,l} &= \sum_s \widetilde{M}_{t,s}(Y^{kn+j}\widetilde{M})_{s,l} \\ &= \sum_r \left(\widetilde{M}_{t,rn}(Y^{kn+j}\widetilde{M})_{rn,l} + \widetilde{M}_{t,rn+1}(Y^{kn+j}\widetilde{M})_{rn+1,l}\right) = 0,\end{aligned}$$

y por lo tanto
$$\widetilde{M} Y^{kn+j} \widetilde{M} = 0, \quad j = 0, 1, \cdots, n-2.$$

Resta demostrar (3.2.24) para $k \geq n$ y para ello expresamos
$$k = rn + j, \quad r \geq 1, \quad j = 0, 1, \cdots, n-1.$$

Es suficiente demostrar el caso $k = rn$, $r \geq 1$.

En efecto, aplicando el Lema 3.2.6 con $k = rn+1$ se tiene que en ese caso $\widetilde{M}$ cumple (3.1.1) para $rn+1$.

Aplicando sucesivamente el Lema 3.2.6 para $k = rn+2, rn+3, \cdots, rn+(n-1)$, se tiene que $\widetilde{M}$ cumple (3.1.1) para $k = rn+2, rn+3, \cdots, rn+(n-1)$, como se quería.

Ahora demostremos (3.2.24) para $k = rn$, $r \geq 1$.

Para ello demostraremos primero que si se cumple (3.2.15)
$$d\widetilde{M} Y^{n-1} \widetilde{M} = e Y^n \widetilde{M} + a\widetilde{M} Y^n,$$

entonces se cumple (3.2.24) para todo $kn$ con $k \geq 1$, que se puede escribir como

$$Y^{kn}\widetilde{M} = \sum_{r=0}^{k} C_{k,r} M^{(k-r)n} \widetilde{M} Y^{rn}. \tag{3.2.26}$$

Probaremos esta fórmula primero para $k=1$. Notemos que (3.2.24) para $n-1$ es

$$Y^{n-1}\widetilde{M} = M^{n-1}\widetilde{M}.$$

De esto y de (3.2.15) obtenemos

$$dMY^{n-1}\widetilde{M} - dY^n\widetilde{M} = eY^n\widetilde{M} + a\widetilde{M}Y^n,$$
$$Y^n\widetilde{M} = (d+e)Y^n\widetilde{M} = dM^n\widetilde{M} - a\widetilde{M}Y^n,$$

que es (3.2.26) para $k=1$, pues $C_{1,0} = d$ y $C_{1,1} = a_1 = -a$.

Ahora fijemos $k$ y asumamos por hipótesis inductiva que (3.2.26) se cumple para $k-1$. Por el argumento anterior esto implica que (3.2.26) se cumple para $(k-1)n + n-1 = kn-1$, y esto lo podemos escribir de la siguiente forma:

$$Y^{kn-1}\widetilde{M} = \sum_{r=0}^{k-1} C_{k-1,r} M^{(k-r)n-1} \widetilde{M} Y^{rn}.$$

Como asumimos que se cumple (3.2.15), por la Proposición 3.2.1 se cumple (3.2.16):

$$d_k \widetilde{M} Y^{kn-1}\widetilde{M} = e_k Y^{kn}\widetilde{M} - a_k \widetilde{M} Y^{kn}.$$

Entonces

$$d_k M Y^{kn-1}\widetilde{M} - d_k Y^{kn}\widetilde{M} = e_k Y^{kn}\widetilde{M} - a_k \widetilde{M} Y^{kn},$$
$$(d_k + e_k) Y^{kn}\widetilde{M} = d_k \left( \sum_{r=0}^{k-1} C_{k-1,r} M^{(k-r)n} \widetilde{M} Y^{rn} \right) + a_k \widetilde{M} Y^{kn},$$

que es (3.2.26) para $k$, pues $C_{k,r} = \frac{d_k}{d_k + e_k} C_{k-1,r}$ para $r < k$ y $C_{k,k} = \frac{a_k}{d_k + e_k}$. Esto concluye el paso inductivo y demuestra que se cumple (3.2.26) para todo $kn$ con $k \geq 1$.

Finalmente, sólo resta probar que $\widetilde{M}$ cumple (3.2.15).

Demostraremos (3.2.15) en cada entrada $(l, kn+j)$.

Como las columnas $\widetilde{M}_{*, kn+j}$ se anulan para $j = 2, 3, \cdots, n-1$ y $\widetilde{M}_{*,nk} = -\widetilde{M}_{*,nk+1}$, es suficiente demostrar (3.2.15) en las entradas $(l, kn)$.

En efecto, (3.2.15) en $(l, kn+j)$ es

$$d(\widetilde{M} Y^{n-1} \widetilde{M})_{l,kn+j} = e(Y^n \widetilde{M})_{l,kn+j} + a(\widetilde{M} Y^n)_{l,kn+j}$$
$$\Leftrightarrow d \sum_r (\widetilde{M} Y^{n-1})_{l,r} \widetilde{M}_{r,kn+j} = e \widetilde{M}_{l+n,kn+j} + a \widetilde{M}_{l,(k-1)n+j}.$$

Si $j = 2, 3, \cdots, n-1$, ambos lados se anulan y si $j = 1$ esto es equivalente, vía $\widetilde{M}_{*,kn} = -\widetilde{M}_{*,kn+1}$, a

$$-d\sum_r (\widetilde{M} Y^{n-1})_{l,r} \widetilde{M}_{r,kn} = -e\widetilde{M}_{l+n,kn} - a\widetilde{M}_{l,(k-1)n}$$

lo cual es equivalente a (3.2.15) en $(l, kn)$.

Así que tenemos que demostrar

$$d(\widetilde{M} Y^{n-1}\widetilde{M})_{l,kn} = e(Y^n \widetilde{M})_{l,kn} + a(\widetilde{M} Y^n)_{l,kn}$$

para $l \geq 0$, $k \geq 0$.

Pero $(\widetilde{M} Y^{n-1}\widetilde{M})_{l,kn} = \widetilde{M}_{l,*}(Y^{n-1}\widetilde{M})_{*,kn} = \widetilde{M}_{l,*}\widetilde{M}_{*+n-1,kn}$.

$(Y^n \widetilde{M})_{l,kn} = \widetilde{M}_{l+n,kn}$,

$(\widetilde{M} Y^n)_{l,kn} = \widetilde{M}_{l,(k-1)n}$,

y $\widetilde{M}_{l,*} = \widetilde{M}_{rn,*}$ para $l = rn + j$ con $j = 0, 1, \cdots, n-1$.

Por lo tanto es suficiente demostrar

$$d(\widetilde{M}_{rn,*}\widetilde{M}_{*+n-1,kn}) = e\widetilde{M}_{(r+1)n,nk} + a\widetilde{M}_{rn,(k-1)n}, \text{ para todo } r, k \geq 0. \qquad (3.2.27)$$

Notemos que $\widetilde{M}_{i,j} = 0$ si $i < 0$, o si $j < 0$.

Por definición

$$\widetilde{M}_{rn,*} = \sum_{i=0}^{r} C_{r,i}(E_{in} - E_{in+1}), \qquad (3.2.28)$$

donde $E_j$ es el vector infinito con $(E_j)_i := \delta_{i,j}$.

Dado que
$(Y^{n-1}\widetilde{M})_{*,kn} = \widetilde{M}_{*+n-1,kn} = (\widetilde{M}_{n-1,kn}, \widetilde{M}_{n,kn}, \widetilde{M}_{n+1,kn}, \cdots)$,
tenemos

$$\begin{aligned}
\widetilde{M}_{rn,*}\widetilde{M}_{*+n-1,kn} &= \sum_{i=0}^{r} C_{r,i}(E_{in} - E_{in+1}) \cdot (\widetilde{M}_{n-1,kn}, \widetilde{M}_{n,kn}, \widetilde{M}_{n+1,kn}, \cdots) \\
&= \sum_{i=0}^{r} C_{r,i}(E_{in} \cdot (\widetilde{M}_{n-1,kn}, \widetilde{M}_{n,kn}, \widetilde{M}_{n+1,kn}, \cdots) \\
&\quad - E_{in+1} \cdot (\widetilde{M}_{n-1,kn}, \widetilde{M}_{n,kn}, \widetilde{M}_{n+1,kn}, \cdots)) \\
&= C_{r,0}(\widetilde{M}_{n-1,kn} - \widetilde{M}_{n,kn}) + C_{r,1}(\widetilde{M}_{2n-1,kn} - \widetilde{M}_{2n,kn}) \\
&\quad + \cdots + C_{r,r}(\widetilde{M}_{rn+n-1,kn} - \widetilde{M}_{(r+1)n,kn}) \\
&= \sum_{i=0}^{r} C_{r,i}(\widetilde{M}_{in+n-1,kn} - \widetilde{M}_{(i+1)n,kn}) \\
&= \sum_{i=0}^{r} C_{r,i}(C_{i,k} - C_{i+1,k}).
\end{aligned}$$

Además $e\widetilde{M}_{(r+1)n,kn} = e\,C_{r+1,k}$, $a\widetilde{M}_{rn,(k-1)n} = a\,C_{r,k-1}$, (notemos que $C_{r,-1} = \widetilde{M}_{rn,-n} = 0$), de manera que (3.2.27) se lee:

$$d\sum_{i=0}^{r} C_{r,i}(C_{i,k} - C_{i+1,k}) = e\,C_{r+1,k} + a\,C_{r,k-1}. \tag{3.2.29}$$

Para demostrar (3.2.29) usaremos

$$e\,d_r + d\,a_r = d_{r+1}, \tag{3.2.30}$$

lo cual se sigue directamente de las definiciones de $a_r$, $e_r$ y $d_r$.

En efecto, $d_r = d\sum_{j=0}^{r-1}(1-d)^j(-a)^{r-1-j}$,

$e\,d_r = d\sum_{j=0}^{r-1}(1-d)^{j+1}(-a)^{r-1-j}$,

$d\,a_r = d(-a)^r$.

Luego

$$e\,d_r + d\,a_r = d\left[\sum_{j=0}^{r-1}(1-d)^{j+1}(-a)^{r-1-j} + (-a)^r\right] = d\left[\sum_{j=1}^{r}(1-d)^j(-a)^{r-j} + (-a)^r\right]$$

$$= d\sum_{j=0}^{r}(1-d)^j(-a)^{r-j} = d_{r+1}.$$

Primero demostraremos (3.2.29) en el caso $\boxed{k > r+1}$.

$C_{i,k} = 0$, pues $k > r+1 \Rightarrow k > i+1 > i$,

$C_{i+1,k} = 0$, pues $k > r+1 \Rightarrow k > i+1$,

$C_{r+1,k} = 0$, pues $k > r+1$,

$C_{r,k-1} = 0$, pues $k > r+1 \Rightarrow k-1 > r$.

Por lo tanto se cumple (3.2.29) pues los dos términos de la igualdad se anulan.

Ahora demostraremos el caso $\boxed{k = r+1}$.

En (3.2.29) todos los términos de la forma $C_{i,r+1}$ se anulan, pues $r+1 > i$. Asimismo de los términos de la forma $C_{i+1,r+1}$ solo sobrevive el término con $i = r$ pues si $i < r$, entonces $r+1 > i+1$ y $C_{i+1,r+1} = 0$.

Luego se tiene

$$-d\,C_{r,r}\,C_{r+1,r+1} = e\,C_{r+1,r+1} + a\,C_{r,r}. \tag{3.2.31}$$

Notemos que $C_{r,r} = \dfrac{a_r}{d_r + e_r}\prod_{i=r+1}^{r}\dfrac{d_i}{d_i + e_i} = \dfrac{a_r}{d_r + e_r}$.

Es decir $C_{r,r}(d_r + e_r) = a_r$ (notemos que $e_0 = a_0 = 1$ y $d_0 = 0$).

Luego (3.2.29) es equivalente a

$$-d\,\frac{a_r}{d_r + e_r}\,\frac{a_{r+1}}{d_{r+1} + e_{r+1}} = e\,\frac{a_{r+1}}{d_{r+1} + e_{r+1}} + \frac{a\,a_r}{d_r + e_r},$$

$-d a_r a_{r+1} = e a_{r+1}(e_r + d_r) + a a_r(e_{r+1} + d_{r+1})$.

Ahora dividimos entre $a_{r+1} = -a a_r$ y así tenemos que demostrar que

$-d a_r = e(e_r + d_r) - (e_{r+1} + d_{r+1})$.

Pero esto se sigue directamente de (3.2.30) usando $e_{r+1} = e e_r$.

En efecto

$d_{r+1} = e d_r + d a_r$,

$e_{r+1} = e e_r$.

Luego $d_{r+1} + e_{r+1} = e(d_r + e_r) + d a_r$,

$-d a_r = e(d_r + e_r) - (d_{r+1} + e_{r+1})$.

Por lo tanto se ha demostrado el caso $k = r+1$.

Finalmente demostramos (3.2.29) en el caso $\boxed{k \leq r}$.

Usaremos el hecho de que para $i \leq r$ tenemos

$$C_{r+1,i} = \frac{a_i}{d_i + e_i} \prod_{j=i+1}^{r+1} \frac{d_j}{d_j + e_j} = \frac{a_i}{d_i + e_i} \prod_{j=i+1}^{r} \frac{d_j}{d_j + e_j} \frac{d_{r+1}}{d_{r+1} + e_{r+1}} = C_{r,i} \frac{d_{r+1}}{d_{r+1} + e_{r+1}}.$$

Es decir

$$C_{r+1,i} = C_{r,i} \frac{d_{r+1}}{d_{r+1} + e_{r+1}}. \tag{3.2.32}$$

Demostraremos (3.2.29) por inducción sobre $r$ (asumiendo $k \leq r$ y usando que (3.2.29) se cumple para $k = r+1$).

Para $r = 0 = k$ esto significa $\quad d C_{0,0}(C_{0,0} - C_{1,0}) = e C_{1,0} + a C_{0,-1}$.

Pero $C_{0,0} = 1$, $C_{1,0} = d = 1 - e$.

Luego en el primer miembro $(1-e)(1-1+e) = e(1-e)$ y en el segundo miembro $e(1-e) + a(0) = e(1-e)$ y por lo tanto (3.2.29) se cumple para $r = k = 0$.

Supongamos que (3.2.29) es válido para algún $r - 1 \geq 0$.

Multiplicando (3.2.30) por $e_{r+1} = e e_r$ obtenemos

$$e_{r+1}(e d_r + d a_r) = e_{r+1} d_{r+1},$$
$$e d_r e_{r+1} + d a_r e_{r+1} = e e_r d_{r+1}.$$

Sumando $e d_r d_{r+1}$ a ambos miembros se obtiene

$$e d_r e_{r+1} + d a_r e_{r+1} + e d_r d_{r+1} = e e_r d_{r+1} + e d_r d_{r+1},$$
$$e d_r (e_{r+1} + d_{r+1}) + d a_r e_{r+1} = e d_{r+1}(d_r + e_r).$$

Dividiendo esta igualdad entre $d_r + e_r$ obtenemos

$$\frac{d_r}{d_r + e_r} e(e_{r+1} + d_{r+1}) + d e_{r+1} \frac{a_r}{d_r + e_r} = e d_{r+1},$$
$$\frac{d_r}{d_r + e_r} e(e_{r+1} + d_{r+1}) + d e_{r+1} C_{r,r} = e d_{r+1},$$

donde hemos usado $a_r = C_{r,r}(e_r + d_r)$.

Ahora dividimos entre $e_{r+1} + d_{r+1}$ y asi podemos escribir el resultado como
$$\frac{e d_r}{d_r + e_r} + d C_{r,r}\left(1 - \frac{d_{r+1}}{e_{r+1} + d_{r+1}}\right) = \frac{e d_{r+1}}{e_{r+1} + d_{r+1}}.$$

Luego multiplicamos por $C_{r,k}$
$$\frac{e d_r}{d_r + e_r} C_{r,k} + d C_{r,r} C_{r,k}\left(1 - \frac{d_{r+1}}{e_{r+1} + d_{r+1}}\right) = \frac{e d_{r+1}}{e_{r+1} + d_{r+1}} C_{r,k},$$

y dado que $C_{r+1,k} = C_{r,k} \dfrac{d_{r+1}}{e_{r+1} + d_{r+1}}$ se obtiene:

$$\frac{e d_r}{d_r + e_r} C_{r,k} + d C_{r,r}(C_{r,k} - C_{r+1,k}) = e C_{r+1,k} \tag{3.2.33}$$

Afirmamos que

$$e C_{r,k} = d \sum_{i=0}^{r-1} C_{r-1,i}(C_{i,k} - C_{i+1,k}) - a C_{r-1,k-1}. \tag{3.2.34}$$

En efecto, si $k < r$, esto se sigue de la hipótesis inductiva, y si $k = r$, entonces (3.2.31) nos da la igualdad (3.2.34). Ahora las igualdades (3.2.34) y (3.2.32) implican

$$\begin{aligned}
\frac{d_r}{e_r + d_r} e C_{r,k} &= \frac{d_r}{e_r + d_r} d \sum_{i=0}^{r-1} C_{r-1,i}(C_{i,k} - C_{i+1,k}) - a C_{r,k-1} \\
&= d \sum_{i=0}^{r-1} C_{r,i}(C_{i,k} - C_{i+1,k}) - a C_{r,k-1}.
\end{aligned}$$

Insertando este valor en (3.2.33) nos da

$$d \sum_{i=0}^{r-1} C_{r,i}(C_{i,k} - C_{i+1,k}) - a C_{r,k-1} + d C_{r,r}(C_{r,k} - C_{r+1,k}) = e C_{r+1,k},$$

lo cual es (3.2.29) para $r$. Esto completa el paso inductivo, demuestra (3.2.29) y concluye la demostración. □

## 3.3. Raíces de $R_k$

En vista del Corolario 3.2.7, queremos analizar los polinomios $R_k$ y sus raíces. En particular estamos interesados en la siguiente pregunta.

Dado un par $(a, d)$, ¿existe un $k \in \mathbb{N}$ tal que $R_k(a, d) = 0$? Si la respuesta es no, entonces el par $(a, d)$ define una única aplicación de torcimiento vía el Teorema 3.2.3. De lo contrario no existe una aplicación de torcimiento que satisfaga el item (1) de la Proposición 3.1.6 para dicho par $(a, d)$.

Fijemos $(a, d)$.

Si $\boxed{-a = 1 - d = e}$ entonces $d_k = k d e^{k-1}$.

De esta manera $R_k(a,d) = e_k + d_k = e_k + kd\,e^{k-1} = e^{k-1}(kd+e)$.

Luego $R_k(a,d) = 0 \Leftrightarrow kd + e = 0 \Leftrightarrow d(k-1) + 1 = 0 \Leftrightarrow d = -\dfrac{1}{k-1}$.

Además en este caso $-a = 1 - d = 1 + \dfrac{1}{k-1} = \dfrac{k}{k-1}$, es decir $a = -\dfrac{k}{k-1}$.

Por lo tanto, si $a$ no es de la forma $\dfrac{k}{1-k}$ para $k \in N$, entonces existe una aplicación de torcimiento que satisface el item (1) de la Proposición 3.1.6 para $(a,d)$, con $d = 1 + a$. Ahora supongamos $\boxed{-a \neq 1-d = e}$.

Entonces
$$R_k(a,d) = e_k + d_k = e^k + d\frac{e^k - (-a)^k}{e+a} = \frac{e^{k+1} + ae^k + de^k - d(-a)^k}{e+a},$$

de manera que
$$\begin{aligned}R_k(a,d) &= 0 \Leftrightarrow e^{k+1} + ae^k + de^k - d(-a)^k = 0 \\ &\Leftrightarrow e^k(e+d) + ae^k - d(-a)^k = 0 \\ &\Leftrightarrow (a+1)e^k - (1-e)(-a)^k = 0.\end{aligned}$$

En este caso $e \neq 0$ y $e \neq 1$ pues si $e = 0$ entonces $a = 0$ y si $e = 1$ entonces $a = -1$, lo cual contradice $-a \neq e$.

De esta manera
$$R_k(a,d) = 0 \Leftrightarrow \frac{1+a}{1-e} = \frac{(-a)^k}{e^k}. \tag{3.3.35}$$

Esta condición es más fácil de manejar que la condición original.

Asumiendo que $K \subset \mathbb{C}$, si (3.3.35) es satisfecha y $\left|\dfrac{a}{e}\right| \neq 1$, entonces

$k\log\left|\dfrac{a}{e}\right| = \log\left|\dfrac{1+a}{1-e}\right|$.

Además si $\left|\dfrac{a}{e}\right| = 1$ y se satisface (3.3.35) entonces necesariamente $\left|\dfrac{1+a}{1-e}\right| = 1$.

Luego $\left|\dfrac{1+a}{1-e}\right| = 1 \Leftrightarrow |1+a| = |1-e| \Leftrightarrow (1+a)(1+\overline{a}) = (1-e)(1-\overline{e})$.

Pero $|a| = |e| \Leftrightarrow a\overline{a} = e\overline{e}$,

luego $\left|\dfrac{1+a}{1-e}\right| = 1 \Leftrightarrow a + \overline{a} = -(e + \overline{e}) \Leftrightarrow \mathrm{Re}(a) = -\mathrm{Re}(e)$.

Es decir $-a = \overline{e}$ ó $-a = e$.

El caso $-a = e$ no es posible por la suposición $-a \neq e$.

Entonces $-a = \overline{e}$.

Reemplazando en (3.3.35):
$$R_k(a,d) = 0 \Leftrightarrow \frac{1+a}{1-e} = \frac{(-a)^k}{e^k} \Leftrightarrow \frac{1-\overline{e}}{1-e} = \frac{\overline{e}^k}{e^k} \Leftrightarrow \frac{1-\overline{e}}{1-e} = \frac{\overline{e}^k e^k}{e^{2k}} \Leftrightarrow \frac{1-\overline{e}}{1-e} = \frac{|e|^{2k}}{e^{2k}}.$$

Haciendo $u = \dfrac{e}{|e|}$, número complejo unitario, se tiene entonces

$R_k(a,d) = 0 \Leftrightarrow \dfrac{1-e}{1-\overline{e}} = u^{2k}$.

A continuación veremos algunos casos particulares de la familia $A(n,d,a)$.

**Ejemplo 3.3.1.** $\boxed{d=0,\ a=0}$

$$\widetilde{M} = \begin{pmatrix} 1 & -1 & 0 & 0 & 0 & 0 & 0 & \cdots \\ 1 & -1 & 0 & 0 & 0 & 0 & 0 & \cdots \\ \vdots & \vdots & \vdots & \vdots & \vdots & \vdots & \vdots & \vdots \\ 1 & -1 & 0 & 0 & 0 & 0 & 0 & \cdots \\ 0 & 0 & 0 & 0 & 0 & 0 & 0 & \cdots \\ \vdots & \vdots & \vdots & \vdots & \vdots & \vdots & \vdots & \vdots \end{pmatrix} \leftarrow n-1$$

Se tiene $Y^n \widetilde{M} = 0$.
Se sigue que $Y^k \widetilde{M} = 0$ para todo $k \geq n$.
Es decir
$$y^k x = y^{k+1}, \text{ para } k \geq n.$$

**Ejemplo 3.3.2.** $\boxed{d=0,\ a \neq 0}$

$\boxed{n=2}$
En este caso se verifica que la matriz $\widetilde{M}$ está dada por

$$\widetilde{M} = \begin{pmatrix} 1 & -1 & 0 & 0 & 0 & 0 & 0 & 0 & 0 & \cdots \\ 1 & -1 & 0 & 0 & 0 & 0 & 0 & 0 & 0 & \cdots \\ 0 & 0 & -a & a & 0 & 0 & 0 & 0 & 0 & \cdots \\ 0 & 0 & -a & a & 0 & 0 & 0 & 0 & 0 & \cdots \\ 0 & 0 & 0 & 0 & a^2 & -a^2 & 0 & 0 & 0 & \cdots \\ 0 & 0 & 0 & 0 & a^2 & -a^2 & 0 & 0 & 0 & \cdots \\ 0 & 0 & 0 & 0 & 0 & 0 & -a^3 & a^3 & 0 & \cdots \\ 0 & 0 & 0 & 0 & 0 & 0 & -a^3 & a^3 & 0 & \cdots \\ \vdots & \vdots & \vdots & \vdots & \vdots & \vdots & \vdots & \vdots & \vdots & \vdots \end{pmatrix}$$

En general para $\boxed{k \geq 1}$ por la Proposición 3.2.1 tenemos
$Y^{2k} \widetilde{M} = (-a)^k \widetilde{M} Y^{2k}$, y por el Teorema 3.2.3 $Y^{2k+1} \widetilde{M} = (-a)^k M \widetilde{M} Y^{2k}$.

$\boxed{n=3}$

En este caso se verifica que la matriz $\widetilde{M}$ está dada por

$$\widetilde{M} = \begin{pmatrix} 1 & -1 & 0 & 0 & 0 & 0 & 0 & 0 & 0 & 0 & 0 & 0 & \cdots \\ 1 & -1 & 0 & 0 & 0 & 0 & 0 & 0 & 0 & 0 & 0 & 0 & \cdots \\ 1 & -1 & 0 & 0 & 0 & 0 & 0 & 0 & 0 & 0 & 0 & 0 & \cdots \\ 0 & 0 & 0 & -a & a & 0 & 0 & 0 & 0 & 0 & 0 & 0 & \cdots \\ 0 & 0 & 0 & -a & a & 0 & 0 & 0 & 0 & 0 & 0 & 0 & \cdots \\ 0 & 0 & 0 & -a & a & 0 & 0 & 0 & 0 & 0 & 0 & 0 & \cdots \\ 0 & 0 & 0 & 0 & 0 & 0 & a^2 & -a^2 & 0 & 0 & 0 & 0 & \cdots \\ 0 & 0 & 0 & 0 & 0 & 0 & a^2 & -a^2 & 0 & 0 & 0 & 0 & \cdots \\ 0 & 0 & 0 & 0 & 0 & 0 & a^2 & -a^2 & 0 & 0 & 0 & 0 & \cdots \\ 0 & 0 & 0 & 0 & 0 & 0 & 0 & 0 & 0 & -a^3 & a^3 & 0 & \cdots \\ 0 & 0 & 0 & 0 & 0 & 0 & 0 & 0 & 0 & -a^3 & a^3 & 0 & \cdots \\ 0 & 0 & 0 & 0 & 0 & 0 & 0 & 0 & 0 & -a^3 & a^3 & 0 & \cdots \\ \vdots & \vdots & \vdots & \vdots & \vdots & \vdots & \vdots & \vdots & \vdots & \vdots & \vdots & \vdots \end{pmatrix}$$

**Ejemplo 3.3.3.** $\boxed{d=1,\ a\neq 0}$

$\boxed{n=2}$

En este caso se verifica que la matriz $\widetilde{M}$ está dada por

$$\widetilde{M} = \begin{pmatrix} 1 & -1 & 0 & 0 & 0 & 0 & 0 & 0 & 0 & \cdots \\ 1 & -1 & 0 & 0 & 0 & 0 & 0 & 0 & 0 & \cdots \\ 1 & -1 & -a & a & 0 & 0 & 0 & 0 & 0 & \cdots \\ 1 & -1 & -a & a & 0 & 0 & 0 & 0 & 0 & \cdots \\ 1 & -1 & -a & a & -a & a & 0 & 0 & 0 & \cdots \\ 1 & -1 & -a & a & -a & a & 0 & 0 & 0 & \cdots \\ 1 & -1 & -a & a & -a & a & -a & a & 0 & \cdots \\ 1 & -1 & -a & a & -a & a & -a & a & 0 & \cdots \\ \vdots & \vdots & \vdots & \vdots & \vdots & \vdots & \vdots & \vdots & \vdots \end{pmatrix}$$

$\boxed{n=4}$

En este caso se verifica que la matriz $\widetilde{M}$ está dada por

$$\widetilde{M} = \begin{pmatrix}
1 & -1 & 0 & 0 & 0 & 0 & 0 & 0 & 0 & 0 & 0 & 0 & 0 & 0 & \cdots \\
1 & -1 & 0 & 0 & 0 & 0 & 0 & 0 & 0 & 0 & 0 & 0 & 0 & 0 & \cdots \\
1 & -1 & 0 & 0 & 0 & 0 & 0 & 0 & 0 & 0 & 0 & 0 & 0 & 0 & \cdots \\
1 & -1 & 0 & 0 & 0 & 0 & 0 & 0 & 0 & 0 & 0 & 0 & 0 & 0 & \cdots \\
1 & -1 & 0 & 0 & -a & a & 0 & 0 & 0 & 0 & 0 & 0 & 0 & 0 & \cdots \\
1 & -1 & 0 & 0 & -a & a & 0 & 0 & 0 & 0 & 0 & 0 & 0 & 0 & \cdots \\
1 & -1 & 0 & 0 & -a & a & 0 & 0 & 0 & 0 & 0 & 0 & 0 & 0 & \cdots \\
1 & -1 & 0 & 0 & -a & a & 0 & 0 & 0 & 0 & 0 & 0 & 0 & 0 & \cdots \\
1 & -1 & 0 & 0 & -a & a & 0 & 0 & -a & a & 0 & 0 & 0 & 0 & \cdots \\
1 & -1 & 0 & 0 & -a & a & 0 & 0 & -a & a & 0 & 0 & 0 & 0 & \cdots \\
1 & -1 & 0 & 0 & -a & a & 0 & 0 & -a & a & 0 & 0 & 0 & 0 & \cdots \\
1 & -1 & 0 & 0 & -a & a & 0 & 0 & -a & a & 0 & 0 & 0 & 0 & \cdots \\
1 & -1 & 0 & 0 & -a & a & 0 & 0 & -a & a & 0 & 0 & -a & a & \cdots \\
\vdots & \vdots & \vdots & \vdots & \vdots & \vdots & \vdots & \vdots & \vdots & \vdots & \vdots & \vdots & \vdots & \vdots & \vdots
\end{pmatrix}$$

**Ejemplo 3.3.4.** $\boxed{d \neq 1,\ a=0}$

$\boxed{n=3}$

En este caso se verifica que la matriz $\widetilde{M}$ está dada por

$$\widetilde{M} = \begin{pmatrix}
1 & -1 & 0 & 0 & 0 & 0 & 0 & \cdots \\
1 & -1 & 0 & 0 & 0 & 0 & 0 & \cdots \\
1 & -1 & 0 & 0 & 0 & 0 & 0 & \cdots \\
d & -d & 0 & 0 & 0 & 0 & 0 & \cdots \\
d & -d & 0 & 0 & 0 & 0 & 0 & \cdots \\
d & -d & 0 & 0 & 0 & 0 & 0 & \cdots \\
d^2 & -d^2 & 0 & 0 & 0 & 0 & 0 & \cdots \\
d^2 & -d^2 & 0 & 0 & 0 & 0 & 0 & \cdots \\
d^2 & -d^2 & 0 & 0 & 0 & 0 & 0 & \cdots \\
d^3 & -d^3 & 0 & 0 & 0 & 0 & 0 & \cdots \\
d^3 & -d^3 & 0 & 0 & 0 & 0 & 0 & \cdots \\
d^3 & -d^3 & 0 & 0 & 0 & 0 & 0 & \cdots \\
\vdots & \vdots & \vdots & \vdots & \vdots & \vdots & \vdots & \vdots
\end{pmatrix}$$

# Capítulo 4

# El caso $y^n x = dx^{n+1} - x^n y + (a+1)y^{n+1}$

En este capítulo describiremos el caso (2) de la Proposición 3.1.6. Por lo tanto $\sigma$ es una aplicación de torcimiento y $Y, M$ son como en la Proposición 1.4.1 y $\widetilde{M} = M - Y$. Además $n \in \mathbb{N}$ con $n \geq 2$, $\widetilde{M}_{k,*} = \widetilde{M}_{0,*}$ para $1 < k < n$. En particular $\widetilde{M}^2 = 0$ y $\widetilde{M}_{n,*} = dE_0 - E_1 + aE_{n+1}$ para $d, a \in K^*$ con $d(a+1) = 1$. Esto implica que estamos tratando con el caso $y^n x = dx^{n+1} - x^n y + (a+1)y^{n+1}$, con $d, a \in K^\times$ y $d(a+1) = 1$. La correspondiente igualdad matricial es

$$Y^n \widetilde{M} = dM^{n+1} - M^n Y + aY^{n+1}.$$

Usaremos esta igualdad y similares igualdades matriciales para calcular inductivamente las diferentes posibilidades para las filas de $M$ (ó equivalentemente de $\widetilde{M}$).

## 4.1. Fórmulas técnicas

En esta sección estableceremos algunas fórmulas técnicas. En particular en el Lema 4.1.4 encontraremos una fórmula para las potencias de $M$. Por un lado estos resultados técnicos nos ayudarán a calcular todas las posibilidades para las primeras 3n+2 filas de $\widetilde{M}$ en la sección 4.2, y por otro lado, en la sección 4.3, ellos nos permitirán describir una cierta familia de aplicaciones de torcimiento llamada $B(n, L)$.

**Lema 4.1.1.** *Sea $A \in L(K^\mathbb{N})$ una matriz infinita tal que $A_{i,j} = 0$ para $j > i+1$. Entonces*

$$(A^{n+1})_{j,j+n+1} = \prod_{k=0}^{n} A_{j+k, j+k+1}.$$

**Demostración.** Se tiene $(A^2)_{j,j+2} = \sum_{i_1} A_{j,i_1} A_{i_1,j+2}$,

$$(A^3)_{j,j+3} = \sum_{i_2}(A^2)_{j,i_2} A_{i_2,j+3} = \sum_{i_2}\left(\sum_{i_1} A_{j,i_1} A_{i_1,i_2}\right) A_{i_2,j+3} = \sum_{i_1,i_2} A_{j,i_1} A_{i_1,i_2} A_{i_2,j+3}.$$

En general
$$(A^{n+1})_{j,j+n+1} = \sum_{i_1,\cdots,i_n} A_{j,i_1} A_{i_1,i_2} \cdots A_{i_{n-1},i_n} A_{i_n,j+n+1}.$$
Pongamos $i_0 = j$, $i_{n+1} = j + n + 1$, entonces
$$\sum_{k=0}^{n}(i_{k+1} - i_k) = i_{n+1} - i_0 = n + 1.$$
Luego $\sum_{k=0}^{n}(i_{k+1} - i_k - 1) = \sum_{k=0}^{n}(i_{k+1} - i_k) - (n + 1) = 0.$
Por hipótesis $A_{i_k, i_{k+1}} = 0$ para $i_{k+1} > i_k + 1$, es decir para $i_{k+1} - i_k - 1 > 0$.
Luego los términos $A_{j,i_1} A_{i_1,i_2} \cdots A_{i_{n-1},i_n} A_{i_n,j+n+1}$ son iguales a cero
para $i_{k+1} - i_k - 1 > 0$.
De manera que si $A_{j,i_1} A_{i_1,i_2} \cdots A_{i_{n-1},i_n} A_{i_n,j+n+1} \neq 0$ entonces $A_{i_k,i_{k+1}} \neq 0$ para todo $k$
con $i_{k+1} - i_k - 1 \leq 0$, pero como $\sum_{k=0}^{n}(i_{k+1} - i_k - 1) = 0$ entonces $i_{k+1} - i_k - 1 = 0$.
Luego $i_{k+1} = i_k + 1$, para todo $k$.
Es decir $i_1 = i_0 + 1 = j + 1$, $i_2 = i_1 + 1 = j + 2, \cdots, \boxed{i_k = j + k}$ para todo $k$.
Luego el único sumando en $(A^{n+1})_{j,j+n+1}$ que puede ser diferente de cero es
$$\prod_{k=0}^{n} A_{j+k,j+k+1}. \qquad \square$$

**Proposición 4.1.2.** *Sean $\widetilde{M}$, $a$ y $d$ como al inicio del capítulo.*
*Entonces*

$$\widetilde{M} Y^{n-1} \widetilde{M} = a \sum_{j=0}^{n} Y^j \widetilde{M} Y^{n-j} - a^2 Y^{n+1}, \qquad (4.1.1)$$

$$(\widetilde{M} - aY)^{n+1} = 0. \qquad (4.1.2)$$

*Además, sea $m_i := \widetilde{M}_{i,i+1}$ y sea $L = (L_1 = n, L_2, L_3, \cdots)$ una sucesión creciente de enteros tales que $m_{L_i} \neq 0$ y $m_i = 0$ si $L_k < i < L_{k+1}$ para algún $k$.*
*Entonces:*

*(1)* $\prod_{k=1}^{n+1}(m_{j+k} - a) = 0$ *para todo $j \geq 0$.*

*(2) Si $m_j \neq 0$ entonces $m_{j+k} = 0$ para $k = 1, 2, \cdots, n - 1$.*

*(3) Tenemos $m_{L_i} = a$ y $(L_{i+1} - L_i) \in \{n, n+1\}$ para todo $i \geq 1$.*

**Demostración.** De $Y^n \widetilde{M} = dM^{n+1} - M^n Y + aY^{n+1}$ y por el Lema 3.1.3 obtenemos

$$Y^n \widetilde{M} = d\left(Y^{n+1} + \widetilde{M} Y^{n-1} \widetilde{M} + \sum_{j=0}^{n} Y^j \widetilde{M} Y^{n-j}\right) - \left(Y^n + \sum_{j=0}^{n-1} Y^j \widetilde{M} Y^{n-1-j}\right) Y + aY^{n+1}$$

$$= dY^{n+1} + d\widetilde{M} Y^{n-1} \widetilde{M} + d\left(\sum_{j=0}^{n-1} Y^j \widetilde{M} Y^{n-j} + Y^n \widetilde{M}\right) - Y^{n+1}$$

$$- \sum_{j=0}^{n-1} Y^j \widetilde{M} Y^{n-j} + aY^{n+1}$$

$$= (d-1) \sum_{j=0}^{n-1} Y^j \widetilde{M} Y^{n-j} + dY^n \widetilde{M} + d\widetilde{M} Y^{n-1} \widetilde{M} + (d-1+a)Y^{n+1}.$$

De manera que

$$0 = (d-1) \sum_{j=0}^{n-1} Y^j \widetilde{M} Y^{n-j} + (d-1)Y^n \widetilde{M} + d\widetilde{M} Y^{n-1} \widetilde{M} + (d-1+a)Y^{n+1},$$

$$0 = -(1-d) \sum_{j=0}^{n} Y^j \widetilde{M} Y^{n-j} + d\widetilde{M} Y^{n-1} \widetilde{M} + (d-1+a)Y^{n+1}.$$

Por lo que asumimos al inicio del capítulo: $d(a+1) = 1 \Leftrightarrow da + d = 1 \Leftrightarrow \boxed{1-d = ad}$.
Entonces

$$0 = d\widetilde{M} Y^{n-1} \widetilde{M} - ad \sum_{j=0}^{n} Y^j \widetilde{M} Y^{n-j} + (a-ad)Y^{n+1}.$$

Usando $a - ad = a^2 d$ y dividiendo entre $d \neq 0$ obtenemos (4.1.1)

$$0 = \widetilde{M} Y^{n-1} \widetilde{M} - a \sum_{j=0}^{n} Y^j \widetilde{M} Y^{n-j} + a^2 Y^{n+1}.$$

Multiplicando la igualdad anterior por $a^{n-1} \neq 0$ obtenemos

$$0 = \widetilde{M}(aY)^{n-1} \widetilde{M} - \sum_{j=0}^{n} (aY)^j \widetilde{M}(aY)^{n-j} + (aY)^{n+1}.$$

El lado derecho de esta igualdad es el desarrollo de $(aY - \widetilde{M})^{n+1}$ pues en este desarrollo sólo aparecen los términos $(aY)^{n+1}$, $-(aY)^j \widetilde{M}(aY)^{n-j}$, $j = 0, \cdots, n$ y términos de la forma $\widetilde{M}(aY)^j \widetilde{M}$, $j = 1, \cdots, n-1$ pero sabemos que $\widetilde{M}(aY)^j \widetilde{M} = 0$ para $j = 1, \cdots, n-2$.
Así queda demostrado (4.1.2).
Demostración de (1)

$$(\widetilde{M} - aY)_{i,i+1} = \widetilde{M}_{i,i+1} - aY_{i,i+1} = m_i - a.$$

Por el Lema 4.1.1

$$[(\widetilde{M}-aY)^{n+1}]_{j,j+n+1} = \prod_{k=0}^{n}(\widetilde{M}-aY)_{j+k,j+k+1} = \prod_{k=0}^{n}(m_{j+k}-a) = \prod_{k=1}^{n+1}(m_{j+k-1}-a) = 0$$

para todo $j \geq 0$.

Haciendo $\overline{j} = j-1$ y, teniendo en cuenta que la igualdad anterior es también válida para $j \geq 1 > 0$, entonces se tiene que para $\overline{j} \geq 0$ se cumple

$$\prod_{k=1}^{n+1}(m_{\overline{j}+k}-a) = 0.$$

Renombrando $\overline{j}$ como $j$ queda demostrado (1).

Demostración de (2)

Por el Lema 3.1.3 tenemos $\widetilde{M}Y^k\widetilde{M} = 0$ para $k = 0, 1, \cdots, n-2$.

Entonces para $1 \leq k \leq n-1$

$$\begin{aligned}
0 &= (\widetilde{M}Y^{k-1}\widetilde{M})_{j,j+k+1} = \widetilde{M}_{j,*}\widetilde{M}_{*+k-1,j+k+1} \\
&= (\widetilde{M}_{j,0}, \widetilde{M}_{j,1}, \cdots, \widetilde{M}_{j,j}, \widetilde{M}_{j,j+1}, 0, \cdots) \\
&\quad \cdot (\underbrace{\widetilde{M}_{k-1,j+k+1}}_{0}, \underbrace{\widetilde{M}_{k,j+k+1}}_{0}, \cdots, \underbrace{\widetilde{M}_{j+k-1,j+k+1}}_{0}, \widetilde{M}_{j+k,j+k+1}, \cdots) \\
&= \widetilde{M}_{j,j+1}\widetilde{M}_{j+k,j+k+1} = m_j m_{j+k}.
\end{aligned}$$

Si $m_j \neq 0$ entonces $m_{k+j} = 0$.

De esta manera el item (2) es verdadero.

Demostración de (3)

Por el item (2) si $m_j \neq 0$ entonces $m_{j+1} = m_{j+2} = \cdots = m_{j+n-1} = 0$.

Por hipótesis sabemos que

$m_{L_i} \neq 0$ y $m_i = 0$ si $L_k < i < L_{k+1}$ para algún $k$.

Luego se tiene

$$\underbrace{m_{L_i}}_{\neq 0}, \underbrace{m_{L_i+1}}_{=0}, \underbrace{m_{L_i+2}}_{=0}, \cdots, \underbrace{m_{L_i+n-1}}_{=0}, \cdots, \underbrace{m_{L_{i+1}}}_{\neq 0}.$$

Es decir entre $L_i$ y $L_{i+1}$ existe un mínimo de $n-1$ términos nulos.

Luego $L_{i+1} - L_i \geq (n-1) + 1 = n$.

Ahora supongamos, por contradicción, que $L_{i+1} - L_i > n+1$.

Entonces

$$\prod_{k=1}^{n+1}(m_{L_i+k}-a) = (m_{L_i+1}-a)(m_{L_i+2}-a)\cdots(m_{L_i+n+1}-a) = (-a)^{n+1} \neq 0 \quad (\Rightarrow\Leftarrow).$$

Por lo tanto $L_{i+1} - L_i \in \{n, n+1\}$ para todo $i \geq 1$.

Finalmente supongamos, por contradicción, que $m_{L_i} \neq a$ para algún $i \geq 1$, entonces

$$\prod_{k=1}^{n+1}(m_{L_i-2+k}-a) = (m_{L_i-1}-a)(m_{L_i}-a)\cdots(m_{L_i+n-1}-a) = (-a)^n(m_{L_i}-a) \neq 0 \quad (\Rightarrow\Leftarrow).$$

Por lo tanto $m_{L_i} = a$ para todo $i \geq 1$. $\square$

**Proposición 4.1.3.** *Supongamos que para algún $k \geq 1$ existe una sucesión creciente $(L_1, L_2, \cdots, L_k)$ tal que para $i \leq L_k + n - 1$ tenemos*

$$\widetilde{M}_{i,*} = \begin{cases} d^j \widetilde{M}_{0,*}, & si \quad L_j < i < L_{j+1} \text{ con } j \in \{1, \cdots, k-1\}, \\ d^j E_0 - d^{j-1} E_1 + a E_{L_j+1}, & si \quad i = L_j \text{ con } j \in \{1, \cdots, k\}, \\ d^k \widetilde{M}_{0,*}, & si \quad L_k < i < L_k + n. \end{cases}$$

*Pongamos $t_i = \widetilde{M}_{L_k+n,i}$. Entonces*

(1) $(\widetilde{M}_{L_k+n+1,0}, \widetilde{M}_{L_k+n+1,1}, \widetilde{M}_{L_k+n+1,2}) = (d^{k+1}, -t_0, -t_1 - d^k)$,

(2) $\widetilde{M}_{L_k+n+1,j} = -t_{j-1}$ *para* $3 \leq j \leq L_k + n + 1$,

(3) $t_0 \in \{d^k, d^{k+1}\}$,

(4) $(\widetilde{M}_{L_k+n, L_k+n+1}, \widetilde{M}_{L_k+n+1, L_k+n+2}) \in \{(a,0), (0,a)\}$.

**Demostración.** Calculemos las entradas $(L_k + 1, j)$ de (4.1.1)

$$(\widetilde{M} Y^{n-1} \widetilde{M})_{L_k+1,j} = a \sum_{i=0}^{n} (Y^i \widetilde{M} Y^{n-i})_{L_k+1,j} - a^2 (Y^{n+1})_{L_k+1,j}$$

y obtenemos

$$(\widetilde{M} Y^{n-1} \widetilde{M})_{L_k+1,j} = \widetilde{M}_{L_k+1,*} \cdot \widetilde{M}_{*+n-1,j} = d^k \widetilde{M}_{0,*} \cdot \widetilde{M}_{*+n-1,j}$$
$$= d^k (E_0 - E_1) \cdot \widetilde{M}_{*+n-1,j} = d^k (\widetilde{M}_{n-1,j} - \widetilde{M}_{n,j}),$$
$$(Y^i \widetilde{M} Y^{n-i})_{L_k+1,j} = (\widetilde{M} Y^{n-i})_{L_k+i+1,j} = \widetilde{M}_{L_k+i+1,j-n+i},$$
$$(Y^{n+1})_{L_k+1,j} = \delta_{L_k+n+2,j}.$$

Luego

$$d^k (\widetilde{M}_{n-1,j} - \widetilde{M}_{n,j}) = a \sum_{i=0}^{n} \widetilde{M}_{L_k+i+1,j-n+i} - a^2 \delta_{L_k+n+2,j}, \tag{4.1.3}$$

donde $\widetilde{M}_{i,j} = 0$ si $j < 0$.

Para $j = 0$

$$d^k (\widetilde{M}_{n-1,0} - \widetilde{M}_{n,0}) = a \sum_{i=0}^{n} \widetilde{M}_{L_k+i+1,-n+i} - a^2 \delta_{L_k+n+2,0}.$$

El único término que sobrevive en la sumatoria es cuando $i = n$ pues si $i < n$ entonces $i - n < 0$.

Entonces se obtiene: $d^k (1 - d) = a \widetilde{M}_{L_k+n+1,0}$.

Pero $1 - d = ad$, luego $\widetilde{M}_{L_k+n+1,0} = d^{k+1}$.

Cuando calculamos (4.1.3) para $j = 1$ los dos únicos términos que sobreviven son los que corresponden a $i = n$, $i = n-1$.
Por lo tanto se obtiene

$$d^k(\widetilde{M}_{n-1,1} - \widetilde{M}_{n,1}) = a\sum_{i=0}^{n} \widetilde{M}_{L_k+i+1,1-n+i} - a^2\delta_{L_k+n+2,1},$$

$$d^k(-1+1) = a(\widetilde{M}_{L_k+n+1,1} + \widetilde{M}_{L_k+n,0}),$$

$$0 = a(\widetilde{M}_{L_k+n+1,1} + t_0).$$

Como $a \neq 0$ se sigue que $\widetilde{M}_{L_k+n+1,1} = -t_0$.
Similarmente, para $j = 2$ solamente sobreviven los tres términos correspondientes a $i = n, n-1, n-2$ y así obtenemos

$$d^k(\widetilde{M}_{n-1,2} - \widetilde{M}_{n,2}) = a(\widetilde{M}_{L_k+n+1,2} + \widetilde{M}_{L_k+n,1} + \widetilde{M}_{L_k+n-1,0}),$$

$$d^k(0-0) = a(\widetilde{M}_{L_k+n+1,2} + \widetilde{M}_{L_k+n,1} + \widetilde{M}_{L_k+n-1,0}),$$

$$0 = a(\widetilde{M}_{L_k+n+1,2} + \widetilde{M}_{L_k+n,1} + \widetilde{M}_{L_k+n-1,0}),$$

lo cual da
$\widetilde{M}_{L_k+n+1,2} = -\widetilde{M}_{L_k+n,1} - \widetilde{M}_{L_k+n-1,0} = -t_1 - d^k$.
(Notemos que $L_k < L_k + n - 1 < L_k + n$ y entonces por definición $\widetilde{M}_{L_k+n-1,*} = d^k\widetilde{M}_{0,*}$).
Esto demuestra el item (1).

Para demostrar el item (2) observemos que en la suma de la igualdad

$$d^k(\widetilde{M}_{n-1,j} - \widetilde{M}_{n,j}) = a\sum_{i=0}^{n} \widetilde{M}_{L_k+i+1,j-n+i} - a^2\delta_{L_k+n+2,j},$$

si $j - n + i < 0$ entonces $\widetilde{M}_{L_k+i+1,j-n+i} = 0$. Luego los términos correspondientes a $i < n - j$ se eliminan y la suma queda

$$\sum_{i=n-j}^{n} \widetilde{M}_{L_k+i+1,j-n+i} = \sum_{i=0}^{j} \widetilde{M}_{L_k+i+n-j+1,i}$$
$$= \widetilde{M}_{L_k+n-j+1,0} + \widetilde{M}_{L_k+n-j+2,1} + \widetilde{M}_{L_k+n-j+3,2}$$
$$+ \cdots + \widetilde{M}_{L_k+i+n-j+1,i} + \cdots + \widetilde{M}_{L_k+n,j-1} + \widetilde{M}_{L_k+1+n,j}.$$

Caso $\boxed{3 \leq j \leq n}$
En este caso se tiene $-n \leq -j \leq -3$.
Luego $L_k + i + 1 \leq L_k + i + n - j + 1 \leq L_k + i + n - 2$.
Supongamos ahora que $0 \leq i \leq j - 2$.
Entonces $-j \leq i - j \leq -2$.

Por lo tanto $L_k < L_k + n + 1 - j \leq L_k + i + n - j + 1 \leq L_k + n - 1$.

Usando la hipótesis se sigue que $\widetilde{M}_{L_k+i+n-j+1,i} = d^k \widetilde{M}_{0,i}$, si $0 \leq i \leq j-2$.

Luego los únicos términos que sobreviven en la suma, para $j = 3, 4, \cdots, n$, son

$$\widetilde{M}_{L_k+n-j+1,0} = d^k, \widetilde{M}_{L_k+n-j+2,1} = -d^k, \widetilde{M}_{L_k+n,j-1}, \widetilde{M}_{L_k+1+n,j}.$$

De esta manera la suma se reduce a

$$\widetilde{M}_{L_k+n,j-1} + \widetilde{M}_{L_k+1+n,j}.$$

Por otro lado, es claro que $\widetilde{M}_{n-1,j} = \widetilde{M}_{n,j} = 0$ para $j = 3, 4, \cdots, n$.

Por lo tanto de (4.1.3) obtenemos

$0 = a(\widetilde{M}_{L_k+n,j-1} + \widetilde{M}_{L_k+1+n,j})$.

Es decir $\widetilde{M}_{L_k+1+n,j} = -\widetilde{M}_{L_k+n,j-1} = -t_{j-1}$.

Caso $\boxed{j = n+1}$

Consideremos ahora la igualdad

$$d^k(\widetilde{M}_{n-1,j} - \widetilde{M}_{n,j}) = a \sum_{i=0}^{n} \widetilde{M}_{L_k+i+1,j-n+i} - a^2 \delta_{L_k+n+2,j},$$

para $j = n+1$.

A la izquierda tenemos

$$d^k(\widetilde{M}_{n-1,j} - \widetilde{M}_{n,j}) = d^k(\widetilde{M}_{n-1,n+1} - \widetilde{M}_{n,n+1}) = -d^k \widetilde{M}_{n,n+1} = -d^k a.$$

La sumatoria de la derecha, en este caso, es igual a

$\sum_{i=0}^{n} \widetilde{M}_{L_k+i+1,i+1} = \widetilde{M}_{L_k+1,1} + \widetilde{M}_{L_k+2,2} + \widetilde{M}_{L_k+3,3} + \cdots + \widetilde{M}_{L_k+n,n} + \widetilde{M}_{L_k+n+1,n+1}$.

Pero $L_k < i \leq L_k + n - 1 \Rightarrow \widetilde{M}_{i,*} = d^k \widetilde{M}_{0,*}$.

Luego $\widetilde{M}_{L_k+2,2} = \widetilde{M}_{L_k+3,3} = \cdots = \widetilde{M}_{L_k+n-1,n-1} = 0$, es decir los únicos posibles términos diferentes de cero corresponden a $i = 0, i = n-1, i = n$.

Esto produce

$$-d^k a = a(\widetilde{M}_{L_k+1,1} + \widetilde{M}_{L_k+n,n} + \widetilde{M}_{L_k+n+1,n+1}),$$

y dado que $\widetilde{M}_{L_k+1,1} = -d^k$ obtenemos

$$-d^k a = -d^k a + a(\widetilde{M}_{L_k+n,n} + \widetilde{M}_{L_k+n+1,n+1}).$$

Luego

$$\widetilde{M}_{L_k+n+1,n+1} = -\widetilde{M}_{L_k+n,n} = -t_n,$$

como se quería.

Caso $\boxed{n+1 < j < L_k + n + 2}$

Supongamos que $0 \le i \le n-2$

Entonces

$L_k + 1 \le L_k + i + 1 \le L_k + n - 1 \Rightarrow \widetilde{M}_{L_k+i+1,*} = d^k \widetilde{M}_{0,*} \Rightarrow \widetilde{M}_{L_k+i+1,j-n+i} = 0$ ya que $i \ge 0$ implica $1 < j - n \le j - n + i$.

En este caso los dos únicos términos que sobreviven en la suma $\sum_{i=0}^{n} \widetilde{M}_{L_k+i+1,j-n+i}$ son $\widetilde{M}_{L_k+n+1,j}$ y $\widetilde{M}_{L_k+n,j-1}$, y así obtenemos

$$0 = d^k(\underbrace{\widetilde{M}_{n-1,j}}_{0} - \underbrace{\widetilde{M}_{n,j}}_{0}) = (\widetilde{M}_{L_k+n+1,j} + \widetilde{M}_{L_k+n,j-1}).$$

Por lo tanto $\widetilde{M}_{L_k+n+1,j} = -\widetilde{M}_{L_k+n,j-1} = -t_{j-1}$, lo cual concluye la demostración del item (2).

Para demostrar el item (3), consideremos la entrada $(1,0)$ de la igualdad

$$Y^{L_k+n} \widetilde{M} = \sum_{i=0}^{L_k+n+1} \widetilde{M}_{L_k+n,i} M^{L_k+n+1-i} Y^i,$$

y tenemos:

$$(Y^{L_k+n} \widetilde{M})_{1,0} = \sum_{i=0}^{L_k+n+1} t_i (M^{L_k+n+1-i} Y^i)_{1,0}.$$

Por un lado

$$(Y^{L_k+n} \widetilde{M})_{1,0} = \widetilde{M}_{L_k+n+1,0} = d^{k+1}. \tag{4.1.4}$$

Por otro lado, como $(M^{L_k+n+1-i} Y^i)_{1,0} = (M^{L_k+n+1-i})_{1,-i}$, el único término diferente de cero en la suma es el que corresponde a $i = 0$. Es decir

$$\sum_{i=0}^{L_k+n+1} t_i (M^{L_k+n+1-i} Y^i)_{1,0} = t_0 (M^{L_k+n+1})_{1,0}.$$

Reemplazando en la igualdad se obtiene

$$\begin{aligned}
d^{k+1} &= t_0 (M^{L_k+n+1})_{1,0} = t_0 (M^{L_k+n}(\widetilde{M}+Y))_{1,0} \\
&= t_0 (M^{L_k+n} \widetilde{M} + M^{L_k+n} Y)_{1,0} = t_0 (M^{L_k+n} \widetilde{M})_{1,0} + t_0 (M^{L_k+n})_{1,-1} \\
&= t_0 (M^{L_k+n} \widetilde{M})_{1,0}.
\end{aligned}$$

Ahora en

$$Y^{L_k+n-1} \widetilde{M} = \sum_{i=0}^{L_k+n} \widetilde{M}_{L_k+n-1,i} M^{L_k+n-i} Y^i,$$

se tiene que $\widetilde{M}_{L_k+n-1,*} = d^k \widetilde{M}_{0,*}$, luego

$$\begin{aligned}
Y^{L_k+n-1} \widetilde{M} &= d^k \widetilde{M}_{0,0} M^{L_k+n} + d^k \widetilde{M}_{0,1} M^{L_k+n-1} Y \\
&= d^k M^{L_k+n} - d^k M^{L_k+n-1} Y = d^k M^{L_k+n-1}(M-Y) \\
&= d^k M^{L_k+n-1} \widetilde{M}.
\end{aligned}$$

Por lo tanto

$$d^{k+1} = t_0(M^{L_k+n}\widetilde{M})_{1,0} = t_0((\widetilde{M}+Y)M^{L_k+n-1}\widetilde{M})_{1,0} = \frac{t_0}{d^k}((\widetilde{M}+Y)Y^{L_k+n-1}\widetilde{M})_{1,0}$$
$$= \frac{t_0}{d^k}(\widetilde{M}Y^{L_k+n-1}\widetilde{M} + Y^{L_k+n}\widetilde{M})_{1,0}.$$

Multiplicando por $d^k$ y usando (4.1.4) $(Y^{L_k+n}\widetilde{M})_{1,0} = \widetilde{M}_{L_k+n+1,0} = d^{k+1}$ se obtiene

$$d^{2k+1} = t_0(\widetilde{M}Y^{L_k+n-1}\widetilde{M})_{1,0} + t_0 d^{k+1},$$
$$(d^k - t_0)d^{k+1} = t_0(\widetilde{M}Y^{L_k+n-1}\widetilde{M})_{1,0} = t_0\widetilde{M}_{1,*}\widetilde{M}_{*+L_k+n-1,0}$$
$$= t_0(1,-1)\cdot(\widetilde{M}_{L_k+n-1,0}, \widetilde{M}_{L_k+n,0}) = t_0(\widetilde{M}_{L_k+n-1,0} - \widetilde{M}_{L_k+n,0})$$
$$= t_0(d^k - t_0).$$

Por lo tanto

$$(d^{k+1} - t_0)(d^k - t_0) = 0,$$

de lo cual se sigue el item (3): $t_0 \in \{d^k, d^{k+1}\}$.

Finalmente por la Proposición 4.1.2 (3):

$$m_{L_i} := \widetilde{M}_{L_i, L_i+1} = a, \quad (L_{i+1} - L_i) \in \{n, n+1\}, \text{para todo } i \geq 1.$$

Caso $\boxed{L_{i+1} - L_i = n}$
$L_{k+1} = L_k + n.$
Entonces $\widetilde{M}_{L_k+n, L_k+n+1} = \widetilde{M}_{L_{k+1}, L_{k+1}+1} = m_{L_{k+1}} = a,$
$\widetilde{M}_{L_k+n+1, L_k+n+2} = \widetilde{M}_{L_{k+1}+1, L_{k+1}+2} = m_{L_{k+1}+1} = 0.$
Caso $\boxed{L_{i+1} - L_i = n+1}$
$L_{k+1} = L_k + n + 1.$
Entonces
$\widetilde{M}_{L_k+n, L_k+n+1} = \widetilde{M}_{L_{k+1}-1, L_{k+1}} = m_{L_{k+1}-1} = 0,$
$\widetilde{M}_{L_k+n+1, L_k+n+2} = \widetilde{M}_{L_{k+1}, L_{k+1}+1} = m_{L_{k+1}} = a.$ $\square$

**Lema 4.1.4.** *Sean $\widetilde{M}$, $a$ y $d$ como antes. Entonces*

*(1) $\widetilde{M}_{j,*} = d\widetilde{M}_{0,*}$ para $j = n+1, \cdots, 2n-1$,*

*(2) para $n \leq j < 2n$ tenemos*

$$dM^{j+1} = \sum_{i=0}^{j} Y^i \widetilde{M} Y^{j-i} + (1-a)Y^{j+1}, \tag{4.1.5}$$

*en particular*

$$dM^{2n} = \sum_{i=0}^{2n-1} Y^i \widetilde{M} Y^{2n-1-i} + (1-a)Y^{2n}, \tag{4.1.6}$$

*(3)*
$$dM^{2n+1} = \sum_{i=0}^{2n} Y^i \widetilde{M} Y^{2n-i} + (1-a)Y^{2n+1} + \widetilde{M} Y^{2n-1} \widetilde{M}. \tag{4.1.7}$$

**Demostración.** Sabemos que $\widetilde{M}_{i,*} = \widetilde{M}_{0,*}$ para $1 \leq i \leq n-1$.

Luego $\widetilde{M} Y^{i-1} \widetilde{M} = 0$ para $0 \leq i-1 \leq n-2$.

Evaluando en la entrada $(n,j)$ se obtiene

$0 = (\widetilde{M} Y^{i-1} \widetilde{M})_{n,j} = \widetilde{M}_{n,*} \cdot \widetilde{M}_{*+i-1,j} = (d,-1,0,\cdots,a) \cdot (\widetilde{M}_{i-1,j}, \widetilde{M}_{i,j}, \cdots, \widetilde{M}_{i+n-1,j}, \widetilde{M}_{n+i,j}).$

Si $\boxed{j=0}$

$0 = (\widetilde{M} Y^{i-1} \widetilde{M})_{n,0} = (d,-1,0,\cdots,a) \cdot (\underbrace{\widetilde{M}_{i-1,0}}_{1}, \underbrace{\widetilde{M}_{i,0}}_{1}, \widetilde{M}_{i+1,0}, \cdots, \widetilde{M}_{i+n,0}),$

$0 = d - 1 + a\widetilde{M}_{n+i,0}.$

Luego $\widetilde{M}_{n+i,0} = \dfrac{1-d}{a} = \dfrac{da}{a} = d$.

Si $\boxed{j=1}$

$0 = (\widetilde{M} Y^{i-1} \widetilde{M})_{n,1} = (d,-1,0,\cdots,a) \cdot (\underbrace{\widetilde{M}_{i-1,1}}_{-1}, \underbrace{\widetilde{M}_{i,1}}_{-1}, \widetilde{M}_{i+1,1}, \cdots, \widetilde{M}_{i+n-1,1}, \widetilde{M}_{i+n,1}),$

$0 = -d + 1 + a\widetilde{M}_{n+i,1}.$

Luego $\widetilde{M}_{n+i,1} = -\dfrac{1-d}{a} = -\dfrac{da}{a} = -d$.

Si $\boxed{j>1}$

$0 = (\widetilde{M} Y^{i-1} \widetilde{M})_{n,j} = (d,-1,0,\cdots,a) \cdot (\underbrace{\widetilde{M}_{i-1,j}}_{0}, \underbrace{\widetilde{M}_{i,j}}_{0}, \widetilde{M}_{i+1,j}, \cdots, \widetilde{M}_{i+n,j}),$

$0 = a\widetilde{M}_{n+i,j}.$

Luego $\widetilde{M}_{n+i,j} = 0$.

Así hemos demostrado que para $n+1 \leq k = n+i \leq 2n-1$ se cumple $\widetilde{M}_{k,*} = d\widetilde{M}_{0,*}$, con lo cual queda demostrado (1).

Ahora multiplicamos $Y^{n-1}\widetilde{M} = M^{n-1}\widetilde{M}$ por $M = \widetilde{M} + Y$ y obtenemos

$$(\widetilde{M} + Y) Y^{n-1} \widetilde{M} = M^n \widetilde{M},$$
$$\widetilde{M} Y^{n-1} \widetilde{M} + Y^n \widetilde{M} = M^n \widetilde{M}.$$

De esto obtenemos

$\widetilde{M} Y^{n-1} \widetilde{M} = M^n \widetilde{M} - Y^n \widetilde{M} = M^n \widetilde{M} - (dM^{n+1} - M^n Y + aY^{n+1}) = M^{n+1} - dM^{n+1} - aY^{n+1}.$

Es decir

$$\widetilde{M} Y^{n-1} \widetilde{M} = (1-d)M^{n+1} - aY^{n+1}.$$

Por el Lema 3.1.3 b)

$$M^{n+1} = \sum_{j=0}^{n} Y^j \widetilde{M} Y^{n-j} + Y^{n+1} + \widetilde{M} Y^{n-1} \widetilde{M}.$$

Reemplazando lo hallado anteriormente se tiene que

$$M^{n+1} = \sum_{j=0}^{n} Y^j \widetilde{M} Y^{n-j} + Y^{n+1} + (1-d)M^{n+1} - aY^{n+1}.$$

Luego

$$dM^{n+1} = \sum_{j=0}^{n} Y^j \widetilde{M} Y^{n-j} + (1-a)Y^{n+1},$$

lo cual es (4.1.5) para $j = n$.

Por inducción supongamos que (4.1.5) se cumple para algún $j$ con $n \leq j < 2n-1$.
Es decir para dicho $j$ se cumple

$$dM^{j+1} = \sum_{i=0}^{j} Y^i \widetilde{M} Y^{j-i} + (1-a)Y^{j+1}.$$

Multiplicando a la derecha por $M = Y + \widetilde{M}$ se obtiene

$$\begin{aligned} dM^{j+2} &= \sum_{i=0}^{j} Y^i \widetilde{M} Y^{j-i}(Y + \widetilde{M}) + (1-a)Y^{j+1}(Y + \widetilde{M}) \\ &= \sum_{i=0}^{j} Y^i \widetilde{M} Y^{j+1-i} + \sum_{i=0}^{j} Y^i \widetilde{M} Y^{j-i}\widetilde{M} + (1-a)Y^{j+2} + (1-a)Y^{j+1}\widetilde{M} \\ &= \sum_{i=0}^{j+1} Y^i \widetilde{M} Y^{j+1-i} - aY^{j+1}\widetilde{M} + \sum_{i=0}^{j} Y^i \widetilde{M} Y^{j-i}\widetilde{M} + (1-a)Y^{j+2}. \end{aligned}$$

En la segunda sumatoria $0 \leq i \leq j \Rightarrow -j \leq -i \leq 0 \Rightarrow 0 \leq j-i \leq j < 2n-1$.
Haciendo $k = j-i$ se tiene que $0 \leq k \leq 2n-2$.
Sabemos, por el Lema 3.1.3, que para $0 \leq k \leq n-2$ se cumple $\widetilde{M} Y^k \widetilde{M} = 0$.
Por la parte (1) $\widetilde{M}_{k,*} = d\widetilde{M}_{0,*}$ para $n+1 \leq k \leq 2n-1$.
Luego por el Lema 3.1.2 (1) $\widetilde{M} Y^k \widetilde{M} = 0$ para $n+1 \leq k \leq 2n-2$.
Por lo tanto en la segunda sumatoria sólo sobreviven los términos correspondientes a $k = n-1$ y $k = n$.
Entonces reemplazando se obtiene

$$\begin{aligned} dM^{j+2} &= \sum_{i=0}^{j+1} Y^i \widetilde{M} Y^{j+1-i} - aY^{j+1}\widetilde{M} + \sum_{i=0}^{j} Y^i \widetilde{M} Y^{j-i}\widetilde{M} + (1-a)Y^{j+2} \\ &= \sum_{i=0}^{j+1} Y^i \widetilde{M} Y^{j+1-i} - aY^{j+1}\widetilde{M} \\ &\quad + Y^{j-n+1}\widetilde{M} Y^{n-1}\widetilde{M} + Y^{j-n}\widetilde{M} Y^n \widetilde{M} + (1-a)Y^{j+2}. \end{aligned}$$

Por el Lema 3.1.2 (3), como $\widetilde{M}_{n+1,*} = d\widetilde{M}_{n-1,*}$ entonces

$$\widetilde{M} Y^n \widetilde{M} + Y\widetilde{M} Y^{n-1}\widetilde{M} = \frac{1-d}{d} Y^{n+1}\widetilde{M} = aY^{n+1}\widetilde{M}.$$

Reemplazando

$$dM^{j+2} = \sum_{i=0}^{j+1} Y^i \widetilde{M} Y^{j+1-i} - aY^{j+1}\widetilde{M} + Y^{j-n}\underbrace{\left[Y\widetilde{M}Y^{n-1}\widetilde{M} + \widetilde{M}Y^n\widetilde{M}\right]}_{aY^{n+1}\widetilde{M}} + (1-a)Y^{j+2},$$

$$dM^{j+2} = \sum_{i=0}^{j+1} Y^i \widetilde{M} Y^{j+1-i} + (1-a)Y^{j+2},$$

lo cual es (4.1.5) para $j+1$.

Esto demuestra (2).

Finalmente (4.1.5) para $j = 2n-1$ se lee

$$dM^{2n} = \sum_{i=0}^{2n-1} Y^i \widetilde{M} Y^{2n-1-i} + (1-a)Y^{2n}.$$

Multiplicando a la derecha por $M = Y + \widetilde{M}$ se obtiene

$$\begin{aligned}
dM^{2n+1} &= \sum_{i=0}^{2n-1} Y^i \widetilde{M} Y^{2n-1-i}(Y+\widetilde{M}) + (1-a)Y^{2n}(Y+\widetilde{M}) \\
&= \sum_{i=0}^{2n-1} Y^i \widetilde{M} Y^{2n-i} + \sum_{i=0}^{2n-1} Y^i \widetilde{M} Y^{2n-1-i}\widetilde{M} + (1-a)Y^{2n+1} + (1-a)Y^{2n}\widetilde{M} \\
&= \sum_{i=0}^{2n} Y^i \widetilde{M} Y^{2n-i} - aY^{2n}\widetilde{M} + \sum_{i=0}^{2n-1} Y^i \widetilde{M} Y^{2n-1-i}\widetilde{M} + (1-a)Y^{2n+1}.
\end{aligned}$$

Utilizando el mismo argumento en la demostración de (2) para $k = 2n-1-i$ se tiene que $0 \leq k \leq 2n-1$ y por lo tanto en la segunda sumatoria sobreviven los términos correspondientes a $k = n-1$, $k = n$ y $k = 2n-1$.

Entonces reemplazando se obtiene

$$\begin{aligned}
dM^{2n+1} &= \sum_{i=0}^{2n} Y^i \widetilde{M} Y^{2n-i} - aY^{2n}\widetilde{M} + \sum_{i=0}^{2n-1} Y^i \widetilde{M} Y^{2n-1-i}\widetilde{M} + (1-a)Y^{2n+1} \\
&= \sum_{i=0}^{2n} Y^i \widetilde{M} Y^{2n-i} - aY^{2n}\widetilde{M} \\
&\quad + Y^n \widetilde{M} Y^{n-1}\widetilde{M} + Y^{n-1}\widetilde{M}Y^n\widetilde{M} + \widetilde{M}Y^{2n-1}\widetilde{M} + (1-a)Y^{2n+1}.
\end{aligned}$$

Por el Lema 3.1.2 (3), como $\widetilde{M}_{n+1,*} = d\widetilde{M}_{n-1,*}$ entonces

$$\widetilde{M}Y^n\widetilde{M} + Y\widetilde{M}Y^{n-1}\widetilde{M} = \frac{1-d}{d}Y^{n+1}\widetilde{M} = aY^{n+1}\widetilde{M}.$$

Reemplazando

$$\begin{aligned}
dM^{2n+1} &= \sum_{i=0}^{2n} Y^i \widetilde{M} Y^{2n-i} - a Y^{2n}\widetilde{M} + Y^{n-1}\underbrace{\left[Y\widetilde{M}Y^{n-1}\widetilde{M} + \widetilde{M}Y^n\widetilde{M}\right]}_{a Y^{n+1}\widetilde{M}}\\
&\quad + (1-a)Y^{2n+1} + \widetilde{M}Y^{2n-1}\widetilde{M}\\
&= \sum_{i=0}^{2n} Y^i \widetilde{M} Y^{2n-i} + (1-a)Y^{2n+1} + \widetilde{M}Y^{2n-1}\widetilde{M}.
\end{aligned}$$

Queda demostrado (3) y con ello el lema. □

## 4.2. Cálculo de $\widetilde{M}_{j,*}$ para $j \leq 3n+2$

En esta sección continuaremos describiendo el caso (2) de la Proposición 3.1.6. Por lo tanto $\sigma$ es una aplicación de torcimiento, $Y, M$ son como en la Proposición 1.4.1 y $\widetilde{M} = M - Y$. Además existe $n \in \mathbb{N}$ con $n \geq 2$, tal que $\widetilde{M}_{k,*} = \widetilde{M}_{0,*}$ para $1 < k < n$ y $\widetilde{M}_{n,*} = dE_0 - E_1 + aE_{n+1}$ para algunos $a, d \in K^{\times}$ tales que $d(a+1) = 1$.
Calcularemos las diferentes posibilidades para las aplicaciones de torcimiento resultantes. Las primeras $2n-1$ filas están completamente determinadas, pero hay cuatro posibilidades para la $2n$-ésima fila. En cada uno de los cuatro casos las filas están determinadas hasta la fila $3n-1$. Mostraremos como proceder para determinar las filas $3n, 3n+1, 3n+2$, en cada uno de los cuatro casos y obtendremos de nuevo cuatro casos (de manera que tenemos 16 casos). Como el número de posibilidades crece, el sistema de ecuaciones se hace más y más complicado, de manera que nuestros métodos no proveen una clasificación completa. Sin embargo, en la sección 4.3 describiremos una familia de aplicaciones de torcimiento, tal que cuatro de los 16 casos mencionados antes coinciden en las primeras 3n+2 filas con miembros de esta familia.
Si hacemos $t_i = \widetilde{M}_{2n,i}$, la matriz $\widetilde{M}$ hasta la fila $2n+1$ se ve como sigue:

$$\begin{array}{c}
0 \to \\
\\
\\
n \to \\
\\
\\
\\
2n \to \\
2n+1 \to
\end{array}
\left(\begin{array}{cccccccccc}
1 & -1 & 0 & \cdots & 0 & 0 & \cdots & 0 & 0 & 0\\
\vdots & \vdots & \vdots & \cdots & \vdots & \vdots & \cdots & \vdots & \vdots & \vdots\\
1 & -1 & 0 & \cdots & 0 & 0 & \cdots & 0 & 0 & 0\\
\boldsymbol{d} & \boldsymbol{-1} & \boldsymbol{0} & \cdots & \boldsymbol{a} & \boldsymbol{0} & \cdots & \boldsymbol{0} & \boldsymbol{0} & \boldsymbol{0}\\
d & -d & 0 & \cdots & 0 & 0 & \cdots & 0 & 0 & 0\\
\vdots & \vdots & \vdots & \cdots & \vdots & \vdots & \cdots & \vdots & \vdots & \vdots\\
d & -d & 0 & \cdots & 0 & 0 & \cdots & 0 & 0 & 0\\
t_0 & t_1 & t_2 & \cdots & t_{n+1} & t_{n+2} & \cdots & t_{2n} & t_{2n+1} & 0\\
d^2 & -t_0 & -d-t_1 & -t_2 & \cdots & -t_{n+1} & \cdots & -t_{2n-1} & -t_{2n} & a-t_{2n+1}
\end{array}\right)$$

con $t_0 \in \{d, d^2\}$ y $t_{2n+1} \in \{0, a\}$.

En efecto el Lema 4.1.4(1) produce las filas $n+1, \cdots, 2n-1$.

Por la Proposición 4.1.3 para $k = 1$ se tiene $L_1 = n$ y para $i \leq L_1 + n - 1 = 2n - 1$.

se cumple:

$$\boxed{i = L_1 = n} \Rightarrow \widetilde{M}_{n,*} = (d, -1, 0, \cdots, \overset{n+1}{\widehat{a}}, 0, \cdots),$$

$$\boxed{L_1 = n < i \leq 2n-1} \Rightarrow \widetilde{M}_{i,*} = d\widetilde{M}_{0,*} = (d, -d, 0, 0, \cdots).$$

Luego se sigue que si $t_i = \widetilde{M}_{2n,i}$, entonces:

(1) $(\widetilde{M}_{2n+1,0}, \widetilde{M}_{2n+1,1}, \widetilde{M}_{2n+1,2}) = (d^2, -t_0, -t_1 - d)$,

(2) $\widetilde{M}_{2n+1,j} = -t_{j-1}$ para $3 \leq j \leq 2n+1$,

(3) $t_0 \in \{d, d^2\}$,

(4) $(\widetilde{M}_{2n,2n+1}, \widetilde{M}_{2n+1,2n+2}) \in \{(a,0), (0,a)\}$ y en cualquier caso $\widetilde{M}_{2n+1,2n+2} = a - t_{2n+1}$.

Las dos elecciones para $t_0$ y para $t_{2n+1}$ generan cuatro casos diferentes para $\widetilde{M}_{2n,*}$, los cuales describiremos en las dos proposiciones siguientes, cuyas demostraciones requieren sus lemas respectivos.

**Lema 4.2.1.** *Si $t_{2n+1} = \widetilde{M}_{2n,2n+1} = 0$ entonces*

*(1)*
$$t_k = 0 \text{ para } k \notin \{0, 1, n+1\}, \tag{4.2.8}$$

*(2)*
$$t_{n+1}(t_0 - d^2) = 0. \tag{4.2.9}$$

**Demostración.** Para demostrar (4.2.8) debemos demostrar que $t_k = 0$ para $n+2 \leq k \leq 2n$ (Caso 1) y que $t_k = 0$ para $2 \leq k \leq n$ (Caso 2).

$\boxed{\text{Caso 1}: n+2 \leq k \leq 2n}$.

Demostraremos inductivamente que $t_k = 0$.

En efecto, para $k = 2n$ se tiene

$$\begin{aligned}
0 &= (\widetilde{M}^2)_{2n,2n} = \widetilde{M}_{2n,*} \cdot \widetilde{M}_{*,2n} \\
&= (t_0, t_1, \cdots, t_n, t_{n+1}, \cdots, t_{2n}) \cdot (\widetilde{M}_{0,2n}, \widetilde{M}_{1,2n}, \cdots, \widetilde{M}_{n,2n}, \widetilde{M}_{n+1,2n}, \cdots, \widetilde{M}_{2n,2n}).
\end{aligned}$$

Como $\widetilde{M}_{k,j} = 0$ para $k < j - 1$, se tiene

$$\widetilde{M}_{0,2n} = \widetilde{M}_{1,2n} = \cdots = \widetilde{M}_{2n-2,2n} = 0.$$

Además $\widetilde{M}_{2n-1,*} = d\widetilde{M}_{0,*} = (d, -d, 0, \ldots)$ y entonces $\widetilde{M}_{2n-1,2n} = 0$ pues $2n > 1$.

Por lo tanto

$0 = t_{2n} \cdot \widetilde{M}_{2n,2n} = t_{2n}^2$, y entonces $t_{2n} = 0$.

Supongamos que se cumple $t_{2n-1} = 0$, $t_{2n-2} = 0$, ..., $t_{2n-(h-1)} = 0$, para algún $h$ con $1 \leq h \leq n-2$. Demostraremos que se cumple $t_{2n-h} = 0$. Notemos que $1 \leq h \leq n-2$ si y solo si $n+2 \leq 2n-h \leq 2n-1$.

En efecto $\widetilde{M} Y^h \widetilde{M} = 0$ para $1 \leq h \leq n-2$.

Luego

$$\begin{aligned} 0 &= (\widetilde{M} Y^h \widetilde{M})_{2n, 2n-h} = \widetilde{M}_{2n,*} \cdot \widetilde{M}_{*+h, 2n-h} \\ &= (t_0, t_1, \cdots, t_n, \cdots, t_{2n-h}) \cdot (\widetilde{M}_{h, 2n-h}, \widetilde{M}_{h+1, 2n-h}, \cdots, \widetilde{M}_{n+h, 2n-h}, \cdots, \widetilde{M}_{2n, 2n-h}). \end{aligned}$$

Ahora $\widetilde{M}_{j,*} = d\widetilde{M}_{0,*} = (d, -d, 0, \ldots)$, para $n+1 \leq j \leq 2n-1$.

Luego
$\widetilde{M}_{h, 2n-h} = \widetilde{M}_{h+1, 2n-h} = \cdots = \widetilde{M}_{n+h, 2n-h} = \cdots = \widetilde{M}_{2n-2, 2n-h} = \widetilde{M}_{2n-1, 2n-h} = 0$,
pues $n+1 \leq h+n \leq 2n-2$ y $2n-h \geq n+2 > 2$.

Por lo tanto $\quad 0 = t_{2n-h} \cdot \widetilde{M}_{2n, 2n-h} = t_{2n-h}^2$, es decir $\quad t_{2n-h} = 0$, lo cual concluye el Caso 1.

$\boxed{\text{Caso 2}: 2 \leq k \leq n}$.

Para esto sabemos que $\widetilde{M} Y^i \widetilde{M} = 0$ para $0 \leq i \leq n-2$.

Luego

$$\begin{aligned} 0 &= (\widetilde{M} Y^i \widetilde{M})_{2n, n+1} = \widetilde{M}_{2n,*} \cdot \widetilde{M}_{*+i, n+1} \\ &= (t_0, t_1, \cdots, t_j, \cdots, t_n, t_{n+1}) \cdot (\widetilde{M}_{i, n+1}, \widetilde{M}_{i+1, n+1}, \cdots, \widetilde{M}_{i+j, n+1}, \cdots, \widetilde{M}_{i+n, n+1}, \widetilde{M}_{i+n+1, n+1}). \end{aligned}$$

Ahora si $j < n+2$: $(Y^i \widetilde{M})_{j, n+1} = \widetilde{M}_{i+j, n+1} = \begin{cases} a & \text{si } i+j = n \\ 0 & \text{si } i+j \neq n. \end{cases}$

Además $i + j = n \Leftrightarrow j = n - i$.

Por lo tanto $0 = a\, t_{n-i}$.

Es decir $t_{n-i} = 0$.

Así hemos demostrado que $t_k = 0$ para $2 \leq k \leq n$. De esta manera concluye la demostración de (4.2.8).

Ahora demostraremos que $t_{n+1}(t_0 - d^2) = 0$.

Para esto evaluamos la igualdad matricial

$$Y^{2n} \widetilde{M} = t_0 M^{2n+1} + t_1 M^{2n} Y + t_{n+1} M^n Y^{n+1} \tag{4.2.10}$$

la cual es válida por (3.1.1), en la entrada $(1, n+1)$.

Primero calculamos

$$(Y^{2n} \widetilde{M})_{1, n+1} = \widetilde{M}_{2n+1, n+1} = -t_n = 0 \text{ (ya demostrado).} \tag{4.2.11}$$

Ahora demostramos que

$$(M^{2n+1})_{1, n+1} = -\frac{t_{n+1}}{d}. \tag{4.2.12}$$

En efecto por el Lema 4.1.4(3) tenemos que

$$dM^{2n+1} = \sum_{i=0}^{2n} Y^i \widetilde{M} Y^{2n-i} + (1-a)Y^{2n+1} + \widetilde{M} Y^{2n-1} \widetilde{M}.$$

Pero $(Y^{2n+1})_{1,n+1} = \delta_{2n+2,n+1} = 0$,

$$(\widetilde{M} Y^{2n-1} \widetilde{M})_{1,n+1} = \widetilde{M}_{1,*} \cdot \widetilde{M}_{*+2n-1,n+1} = (1,-1) \cdot (\underbrace{\widetilde{M}_{2n-1,n+1}}_{0}, \widetilde{M}_{2n,n+1}) = -t_{n+1}.$$

Reemplazando obtenemos

$$d(M^{2n+1})_{1,n+1} = \sum_{i=0}^{2n} (Y^i \widetilde{M} Y^{2n-i})_{1,n+1} - t_{n+1}.$$

Para demostrar (4.2.12) es suficiente demostrar que $\sum_{i=0}^{2n} (Y^i \widetilde{M} Y^{2n-i})_{1,n+1} = 0$, lo cual se sigue fácilmente del hecho de que para $0 \leq i \leq 2n$ tenemos

$$(Y^i \widetilde{M} Y^{2n-i})_{1,n+1} = [Y^i (\widetilde{M} Y^{2n-i})]_{1,n+1} = \widetilde{M}_{i+1,i-n+1}$$
$$= \begin{cases} \widetilde{M}_{n,0} = d, & \text{si} \quad i = n-1, \\ \widetilde{M}_{n+1,1} = -d, & \text{si} \quad i = n, \\ \widetilde{M}_{2n,n} = 0, & \text{si} \quad i = 2n-1 \\ \widetilde{M}_{2n+1,n+1} = 0, & \text{si} \quad i = 2n, \\ 0 & \text{en otros casos.} \end{cases} \quad (4.2.13)$$

Notemos que si $0 \leq i \leq n-2$ entonces $1-n \leq i-n+1 \leq -1$,
si $n+1 \leq i \leq 2n-2$ entonces $n+2 \leq i+1 \leq 2n-1$ y $2 \leq i-n+1 \leq n-1$.
En la sumatoria sólo sobreviven, por (4.2.13), $\widetilde{M}_{n,0} = d$ y $\widetilde{M}_{n+1,1} = -d$.
Luego la suma total es igual a cero y queda demostrado (4.2.12).
Ahora demostraremos que
$$(M^{2n} Y)_{1,n+1} = 0. \quad (4.2.14)$$

Para esto notemos que por el Lema 4.1.4 (2):

$$dM^{2n} = \sum_{i=0}^{2n-1} Y^i \widetilde{M} Y^{2n-1-i} + (1-a)Y^{2n}.$$

Multiplicando a la derecha por $Y$

$$dM^{2n} Y = \sum_{i=0}^{2n-1} Y^i \widetilde{M} Y^{2n-i} + (1-a)Y^{2n+1}.$$

Evaluemos esta igualdad matricial en la entrada $(1, n+1)$

$$d(M^{2n}Y)_{1,n+1} = \sum_{i=0}^{2n-1}(Y^i \widetilde{M} Y^{2n-i})_{1,n+1} + (1-a)(Y^{2n+1})_{1,n+1}.$$

Pero $(Y^{2n+1})_{1,n+1} = \delta_{2n+2,n+1} = 0$ y por (4.2.13) $\sum_{i=0}^{2n-1}(Y^i \widetilde{M} Y^{2n-i})_{1,n+1} = 0$.

Finalmente calculamos $(M^n Y^{n+1})_{1,n+1}$.

Usando el Lema 3.1.3 a) para $k = n-1$

$$M^n = \sum_{i=0}^{n-1} Y^i \widetilde{M} Y^{n-1-i} + Y^n.$$

Multiplicando a la derecha por $Y^{n+1}$ se obtiene

$$M^n Y^{n+1} = \sum_{i=0}^{n-1} Y^i \widetilde{M} Y^{2n-i} + Y^{2n+1}.$$

En la entrada $(1, n+1)$ obtenemos:
$(Y^{2n+1})_{1,n+1} = \delta_{2n+2,n+1} = 0$,
$(Y^i \widetilde{M} Y^{2n-i})_{1,n+1} = \widetilde{M}_{i+1, i-n+1}$.
En este caso $0 \leq i \leq n-1$ y por (4.2.13) en la suma sólo sobrevive $\widetilde{M}_{n,0} = d$.
Por lo tanto

$$(M^n Y^{n+1})_{1,n+1} = d. \tag{4.2.15}$$

Insertando los valores obtenidos en (4.2.11), (4.2.12), (4.2.14) y (4.2.15) en la igualdad (4.2.10) obtenemos:
$0 = t_0 \left( -\dfrac{t_{n+1}}{d} \right) + t_1(0) + t_{n+1} d,$
$0 = t_{n+1} \left( d - \dfrac{t_0}{d} \right)$
lo cual es (4.2.9) : $t_{n+1}(t_0 - d^2) = 0$. $\qquad\square$

**Proposición 4.2.2.** *Supongamos que* $t_{2n+1} = \widetilde{M}_{2n,2n+1} = 0$.

(1) *Si* $t_0 = d$ *entonces*

$$\widetilde{M}_{2n,*} = dE_0 - dE_1, \quad \widetilde{M}_{2n+1,*} = d^2 E_0 - dE_1 + aE_{2n+2}.$$

(2) *Si* $t_0 = d^2$ *entonces*

$$\widetilde{M}_{2n,*} = d^2 E_0 - dE_1 + adE_{n+1}, \quad \widetilde{M}_{2n+1,*} = d^2 E_0 - d^2 E_1 - adE_{n+2} + aE_{2n+2}.$$

**Demostración.** Tenemos $t_{2n+1} = 0$ y así $\widetilde{M}_{2n+1,2n+2} = a - t_{2n+1} = a$.

Si $t_0 = d$, de (4.2.9) obtenemos $t_{n+1} = 0$.

Entonces

$$\begin{aligned} 0 &= (\widetilde{M}^2)_{2n,0} = \widetilde{M}_{2n,*} \cdot \widetilde{M}_{*,0} \\ &= (t_0, t_1, \cdots, t_{n+1}, \cdots, t_{2n}, t_{2n+1}) \cdot (\widetilde{M}_{0,0}, \widetilde{M}_{1,0}, \cdots, \widetilde{M}_{n+1,0}, \cdots, \widetilde{M}_{2n,0}, \widetilde{M}_{2n+1,0}). \end{aligned}$$

Pero hemos demostrado que $t_k = 0$ para $k \notin \{0, 1, n+1\}$ y acabamos de demostrar que $t_{n+1} = 0$, luego obtenemos

$$0 = (t_0, t_1) \cdot (\widetilde{M}_{0,0}, \widetilde{M}_{1,0}) = (t_0, t_1) \cdot (1, 1) = t_0 + t_1,$$

lo cual produce $t_1 = -d$.

De esta manera hemos demostrado (1):

$$\widetilde{M}_{2n,*} = d E_0 - d E_1, \quad \widetilde{M}_{2n+1,*} = d^2 E_0 - d E_1 + a E_{2n+2}.$$

La última igualdad se sigue de la Proposición 4.1.3.

En efecto

$$(\widetilde{M}_{2n+1,0}, \widetilde{M}_{2n+1,1}, \widetilde{M}_{2n+1,2}) = (d^2, -d, 0),$$
$$\widetilde{M}_{2n+1,j} = -t_{j-1} = 0, \quad 3 \leq j \leq 2n+1,$$
$$\widetilde{M}_{2n+1,2n+2} = a - t_{2n+1} = a.$$

Por otro lado si $t_0 = d^2$, entonces demostraremos que

$$t_{n+1} = ad - d - t_1 \tag{4.2.16}$$

.

Para esto evaluaremos la igualdad matricial (4.2.10) en $(1, 2n+2)$.

En este caso (4.2.10) se lee:

$$Y^{2n} \widetilde{M} = d^2 M^{2n+1} + t_1 M^{2n} Y + t_{n+1} M^n Y^{n+1}. \tag{4.2.17}$$

Primero calculamos

$$(Y^{2n} \widetilde{M})_{1,2n+2} = \widetilde{M}_{2n+1,2n+2} = a. \tag{4.2.18}$$

Ahora demostraremos que

$$d^2 (M^{2n+1})_{1,2n+2} = 1. \tag{4.2.19}$$

Para esto notemos que por el Lema 4.1.4 (3)

$$dM^{2n+1} = \sum_{i=0}^{2n} Y^i \widetilde{M} Y^{2n-i} + (1-a)Y^{2n+1} + \widetilde{M} Y^{2n-1} \widetilde{M}.$$

Se tiene $(Y^{2n+1})_{1,2n+2} = \delta_{2n+2,2n+2} = 1$,
$(\widetilde{M} Y^{2n-1} \widetilde{M})_{1,2n+2} = \widetilde{M}_{1,*} \cdot \widetilde{M}_{*+2n-1,2n+2} = (\widetilde{M}_{1,0}, \widetilde{M}_{1,1}) \cdot (\widetilde{M}_{2n-1,2n+2}, \widetilde{M}_{2n,2n+2})$
$= (1,-1) \cdot (0,0) = 0$ pues $2n+2 > 2n+1 > 2n$.
De manera que para demostrar (4.2.19) es suficiente demostrar que

$$\sum_{i=0}^{2n} (Y^i \widetilde{M} Y^{2n-i})_{1,2n+2} = 2a, \tag{4.2.20}$$

dado que entonces $d^2(M^{2n+1})_{1,2n+2} = d(2a + (1-a)) = d(a+1) = 1$.

Pero (4.2.20) se sigue rápidamente del hecho de que para $i = 0, 1, \cdots, 2n$ tenemos

$$(Y^i \widetilde{M} Y^{2n-i})_{1,2n+2} = (\widetilde{M} Y^{2n-i})_{i+1,2n+2} = \widetilde{M}_{i+1,i+2} = \begin{cases} a, & \text{si } i \in \{n-1, 2n\}, \\ 0, & \text{de lo contrario.} \end{cases} \tag{4.2.21}$$

Esto establece (4.2.19).

A continuación demostraremos que

$$d(M^{2n} Y)_{1,2n+2} = 1. \tag{4.2.22}$$

Para esto notemos que por el Lema 4.1.4 (2) se tiene

$$dM^{2n} = \sum_{i=0}^{2n-1} Y^i \widetilde{M} Y^{2n-1-i} + (1-a)Y^{2n}.$$

Multiplicando a la derecha por $Y$ obtenemos

$$dM^{2n} Y = \sum_{i=0}^{2n-1} Y^i \widetilde{M} Y^{2n-i} + (1-a)Y^{2n+1}.$$

Luego

$$d(M^{2n} Y)_{1,2n+2} = \sum_{i=0}^{2n-1} (Y^i \widetilde{M} Y^{2n-i})_{1,2n+2} + (1-a)(Y^{2n+1})_{1,2n+2}.$$

Pero $(Y^{2n+1})_{1,2n+2} = \delta_{2n+2,2n+2} = 1$.

Por (4.2.21) $(Y^i \widetilde{M} Y^{2n-i})_{1,2n+2} = \widetilde{M}_{i+1,i+2} = \begin{cases} a, & \text{si } i \in \{n-1, 2n\}, \\ 0, & \text{de lo contrario.} \end{cases}$

Por lo tanto

$$d(M^{2n} Y)_{1,2n+2} = a + (1-a) = 1.$$

Finalmente usando el Lema 3.1.3 a) para $k = n-1$ se tiene que

$$M^n = \sum_{i=0}^{n-1} Y^i \widetilde{M} Y^{n-1-i} + Y^n.$$

Multiplicando a la derecha por $Y^{n+1}$ se obtiene

$$M^n Y^{n+1} = \sum_{i=0}^{n-1} Y^i \widetilde{M} Y^{2n-i} + Y^{2n+1}.$$

Evaluando en $(1, 2n+2)$ y usando (4.2.21)

$$(M^n Y^{n+1})_{1,2n+2} = \sum_{i=0}^{n-1} (Y^i \widetilde{M} Y^{2n-i})_{1,2n+2} + (Y^{2n+1})_{1,2n+2} = a+1. \qquad (4.2.23)$$

Insertando los valores obtenidos en (4.2.18),(4.2.19), (4.2.22) y (4.2.23) en la igualdad (4.2.17) obtenemos

$$a = 1 + \frac{t_1}{d} + t_{n+1}(a+1).$$

Luego, usando $a + 1 = \frac{1}{d}$ se sigue (4.2.16). En efecto :

$$a = 1 + \frac{t_1}{d} + \frac{t_{n+1}}{d},$$
$$da = d + t_1 + t_{n+1},$$
$$\boxed{t_{n+1} = da - d - t_1}.$$

Ahora

$$0 = (\widetilde{M}^2)_{2n,0} = \widetilde{M}_{2n,*} \cdot \widetilde{M}_{*,0} = (t_0, t_1, 0, \cdots, 0, t_{n+1}) \cdot (\widetilde{M}_{0,0}, \widetilde{M}_{1,0}, \widetilde{M}_{2,0}, \cdots, \widetilde{M}_{n,0} \widetilde{M}_{n+1,0}),$$
$$0 = t_0 \widetilde{M}_{0,0} + t_1 \widetilde{M}_{1,0} + t_{n+1} \widetilde{M}_{n+1,0},$$
$$0 = t_0 + t_1 + d\, t_{n+1},$$

lo cual junto con (4.2.16) producen

$$0 = t_0 + t_1 + d(da - d - t_1) = t_0 + t_1 + d^2 a - d^2 - d t_1 = \underbrace{(1-d)}_{da} t_1 + d^2 a = ad^2 + adt_1$$

y de esta manera

$$t_1 = \frac{-ad^2}{ad} = -d,$$

y

$$t_{n+1} = da - d - t_1 = da - d + d = da.$$

Esto da lugar a

$$\widetilde{M}_{2n,*} = d^2 E_0 - d E_1 + ad E_{n+1}, \quad \widetilde{M}_{2n+1,*} = d^2 E_0 - d^2 E_1 - ad E_{n+2} + a E_{2n+2},$$

que es lo que se quería y así concluye la demostración. □

En el caso de la Proposición 4.2.2 (1), la matriz $\widetilde{M}$ hasta la fila $2n+1$ se ve como sigue:

$$\begin{array}{c} 0 \to \\ \\ \\ n \to \\ \\ \\ \\ \\ 2n \to \\ 2n+1 \to \end{array} \left( \begin{array}{ccccccccc} 1 & -1 & 0 & \cdots & 0 & 0 & \cdots & 0 & 0 & 0 \\ \vdots & \vdots & \vdots & \cdots & \vdots & \vdots & \cdots & \vdots & \vdots & \vdots \\ 1 & -1 & 0 & \cdots & 0 & 0 & \cdots & 0 & 0 & 0 \\ \boldsymbol{d} & \boldsymbol{-1} & \boldsymbol{0} & \cdots & \boldsymbol{a} & \boldsymbol{0} & \cdots & \boldsymbol{0} & \boldsymbol{0} & \boldsymbol{0} \\ d & -d & 0 & \cdots & 0 & 0 & \cdots & 0 & 0 & 0 \\ \vdots & \vdots & \vdots & \cdots & \vdots & \vdots & \cdots & \vdots & \vdots & \vdots \\ d & -d & 0 & \cdots & 0 & 0 & \cdots & 0 & 0 & 0 \\ \boldsymbol{d} & \boldsymbol{-d} & \boldsymbol{0} & \cdots & \boldsymbol{0} & \boldsymbol{0} & \cdots & \boldsymbol{0} & \boldsymbol{0} & \boldsymbol{0} \\ \boldsymbol{d^2} & \boldsymbol{-d} & \boldsymbol{0} & \cdots & \boldsymbol{0} & \boldsymbol{0} & \cdots & \boldsymbol{0} & \boldsymbol{0} & \boldsymbol{a} \end{array} \right)$$

En el caso de la Proposición 4.2.2 (2), la matriz $\widetilde{M}$ hasta la fila $2n+1$ se ve como sigue:

$$\begin{array}{c} 0 \to \\ \\ \\ n \to \\ \\ \\ \\ \\ 2n \to \\ 2n+1 \to \end{array} \left( \begin{array}{ccccccccc} 1 & -1 & 0 & \cdots & 0 & 0 & \cdots & 0 & 0 & 0 \\ \vdots & \vdots & \vdots & \cdots & \vdots & \vdots & \cdots & \vdots & \vdots & \vdots \\ 1 & -1 & 0 & \cdots & 0 & 0 & \cdots & 0 & 0 & 0 \\ \boldsymbol{d} & \boldsymbol{-1} & \boldsymbol{0} & \cdots & \boldsymbol{a} & \boldsymbol{0} & \cdots & \boldsymbol{0} & \boldsymbol{0} & \boldsymbol{0} \\ d & -d & 0 & \cdots & 0 & 0 & \cdots & 0 & 0 & 0 \\ \vdots & \vdots & \vdots & \cdots & \vdots & \vdots & \cdots & \vdots & \vdots & \vdots \\ d & -d & 0 & \cdots & 0 & 0 & \cdots & 0 & 0 & 0 \\ \boldsymbol{d^2} & \boldsymbol{-d} & \boldsymbol{0} & \cdots & \boldsymbol{ad} & \boldsymbol{0} & \cdots & \boldsymbol{0} & \boldsymbol{0} & \boldsymbol{0} \\ \boldsymbol{d^2} & \boldsymbol{-d^2} & \boldsymbol{0} & \cdots & \boldsymbol{0} & \boldsymbol{-ad} & \cdots & \boldsymbol{0} & \boldsymbol{0} & \boldsymbol{a} \end{array} \right)$$

**Lema 4.2.3.** *Supongamos que $t_{2n+1} = \widetilde{M}_{2n,2n+1} = a$ y hagamos*

$$\widehat{t}_j = \begin{cases} t_j & \text{si} & j \leq n, \\ d\, t_j & \text{si} & n < j \leq 2n+1. \end{cases},$$

$$S_j = \sum_{i=0}^{j} \widehat{t}_i.$$

*Entonces*

*(i)* $(t_0 - d)t_2 = (t_1 + d)^2$,

*(ii)* $(t_0 - d)t_k = (d + t_1)t_{k-1}$, *para $k = 3, \cdots, n$*,

*(iii)* $(t_0 - d^2)t_{n+1} = (d + t_1)t_n$,

*(iv)* $(t_0 - d)t_k = (d + t_1)t_{k-1}$, *para $k = n+2, \cdots, 2n$*,

*(v)* $(d + t_1)t_{2n} = t_0 a + S_{2n}$,

*(vi)* $0 = t_k + d\, t_{k+n}$, para $k = 2, \cdots, n-1$.

**Demostración.** Por la Observación 3.1.1 para $k \geq 0$ se tiene

$$Y^k \widetilde{M} = \sum_{j=0}^{k+1} \widetilde{M}_{k,j} M^{k+1-j} Y^j.$$

Para $\boxed{k = 2n}$

$$Y^{2n} \widetilde{M} = \sum_{j=0}^{2n+1} \widetilde{M}_{2n,j} M^{2n+1-j} Y^j = \sum_{j=0}^{n} t_j M^{2n+1-j} Y^j + \sum_{j=n+1}^{2n+1} t_j M^{2n+1-j} Y^j,$$

$$Y^{2n} \widetilde{M} = t_0 M^{2n+1} + \sum_{j=1}^{n} t_j M^{2n+1-j} Y^j + \sum_{j=n+1}^{2n} t_j M^{2n+1-j} Y^j + t_{2n+1} Y^{2n+1}. \quad (4.2.24)$$

Por el Lema 4.1.4(3)

$$dM^{2n+1} = \widetilde{M} Y^{2n-1} \widetilde{M} + \sum_{i=0}^{2n} Y^i \widetilde{M} Y^{2n-i} + (1-a) Y^{2n+1}.$$

En la primera sumatoria de (4.2.24) $1 \leq j \leq n$, entonces $n \leq 2n - j \leq 2n - 1 < 2n$ y en este caso el Lema 4.1.4(2) para $2n - j$ afirma que

$$dM^{2n+1-j} = \sum_{i=0}^{2n-j} Y^i \widetilde{M} Y^{2n-j-i} + (1-a) Y^{2n-j+1}.$$

Multiplicando a la derecha por $Y^j$ se obtiene

$$dM^{2n+1-j} Y^j = \sum_{i=0}^{2n-j} Y^i \widetilde{M} Y^{2n-i} + (1-a) Y^{2n+1}.$$

En la segunda sumatoria de (4.2.24) $n+1 \leq j \leq 2n$ entonces $0 \leq 2n-j \leq n-1$ y en este caso el Lema 3.1.3 a) para $2n - j$ afirma que

$$M^{2n+1-j} = \sum_{i=0}^{2n-j} Y^i \widetilde{M} Y^{2n-j-i} + Y^{2n-j+1}.$$

Multiplicando a la derecha por $Y^j$

$$M^{2n+1-j} Y^j = \sum_{i=0}^{2n-j} Y^i \widetilde{M} Y^{2n-i} + Y^{2n+1}.$$

Así (4.2.24) multiplicado por $d$ se lee

$$dY^{2n}\widetilde{M} = t_0 dM^{2n+1} + \sum_{j=1}^{n} t_j dM^{2n+1-j}Y^j + \sum_{j=n+1}^{2n} t_j dM^{2n+1-j}Y^j + t_{2n+1}dY^{2n+1}$$

$$= t_0\left(\widetilde{M}Y^{2n-1}\widetilde{M} + \sum_{i=0}^{2n} Y^i\widetilde{M}Y^{2n-i} + (1-a)Y^{2n+1}\right)$$

$$+ \sum_{j=1}^{n} t_j dM^{2n+1-j}Y^j + \sum_{j=n+1}^{2n} t_j dM^{2n+1-j}Y^j + t_{2n+1}dY^{2n+1}$$

$$= t_0\widetilde{M}Y^{2n-1}\widetilde{M} + t_0\sum_{i=0}^{2n} Y^i\widetilde{M}Y^{2n-i} + t_0(1-a)Y^{2n+1} + t_{2n+1}dY^{2n+1}$$

$$+ \sum_{j=1}^{n} t_j\left(\sum_{i=0}^{2n-j} Y^i\widetilde{M}Y^{2n-i} + (1-a)Y^{2n+1}\right)$$

$$+ \sum_{j=n+1}^{2n} t_j d\left(\sum_{i=0}^{2n-j} Y^i\widetilde{M}Y^{2n-i} + Y^{2n+1}\right)$$

$$= t_0\widetilde{M}Y^{2n-1}\widetilde{M} + \sum_{j=0}^{n} t_j\left(\sum_{i=0}^{2n-j} Y^i\widetilde{M}Y^{2n-i} + (1-a)Y^{2n+1}\right)$$

$$+ \sum_{j=n+1}^{2n} t_j d\left(\sum_{i=0}^{2n-j} Y^i\widetilde{M}Y^{2n-i} + Y^{2n+1}\right) + t_{2n+1}dY^{2n+1}.$$

Si escribimos

$$\widehat{t}_j = \begin{cases} t_j & \text{si} \quad j \leq n, \\ dt_j & \text{si} \quad n < j \leq 2n+1. \end{cases}$$

y $S_j = \sum_{i=0}^{j} \widehat{t}_i$ entonces

$$dY^{2n}\widetilde{M} = t_0\widetilde{M}Y^{2n-1}\widetilde{M} + \sum_{j=0}^{n} t_j \sum_{i=0}^{2n-j} Y^i\widetilde{M}Y^{2n-i} + \sum_{j=0}^{n} t_j(1-a)Y^{2n+1}$$

$$+ \sum_{j=n+1}^{2n} t_j d\sum_{i=0}^{2n-j} Y^i\widetilde{M}Y^{2n-i} + \sum_{j=n+1}^{2n} t_j dY^{2n+1} + t_{2n+1}dY^{2n+1}$$

$$= t_0\widetilde{M}Y^{2n-1}\widetilde{M} + \sum_{j=0}^{2n} \widehat{t}_j \sum_{i=0}^{2n-j} Y^i\widetilde{M}Y^{2n-i} + \left(\sum_{j=0}^{n} \widehat{t}_j\right)(1-a)Y^{2n+1}$$

$$+ \left(\sum_{j=n+1}^{2n+1} \widehat{t}_j\right)Y^{2n+1},$$

esto se lee

$$dY^{2n}\widetilde{M} = t_0\widetilde{M}Y^{2n-1}\widetilde{M} + \sum_{j=0}^{2n} \widehat{t}_j \sum_{i=0}^{2n-j} Y^i\widetilde{M}Y^{2n-i} + (S_{2n+1} - aS_n)Y^{2n+1}. \qquad (4.2.25)$$

Pero
$$\sum_{j=0}^{2n} \widehat{t}_j \sum_{i=0}^{2n-j} Y^i \widetilde{M} Y^{2n-i} = \sum_{i=0}^{2n} \sum_{j=0}^{2n-i} \widehat{t}_j Y^i \widetilde{M} Y^{2n-i}.$$

Haciendo $u = 2n - i$ en la segunda sumatoria se tiene

$$\sum_{j=0}^{2n} \widehat{t}_j \sum_{i=0}^{2n-j} Y^i \widetilde{M} Y^{2n-i} = \sum_{u=0}^{2n} \sum_{j=0}^{u} \widehat{t}_j Y^{2n-u} \widetilde{M} Y^u = \widehat{t}_0 Y^{2n} \widetilde{M} + \sum_{u=1}^{2n} \sum_{j=0}^{u} \widehat{t}_j Y^{2n-u} \widetilde{M} Y^u.$$

Así de (4.2.25) obtenemos

$$(d - t_0) Y^{2n} \widetilde{M} = t_0 \widetilde{M} Y^{2n-1} \widetilde{M} + \sum_{i=1}^{2n} \sum_{j=0}^{i} \widehat{t}_j Y^{2n-i} \widetilde{M} Y^i + (S_{2n+1} - a S_n) Y^{2n+1}. \quad (4.2.26)$$

Evaluaremos ahora la igualdad matricial (4.2.26) en la entrada $(1, k)$ para $k = 2, 3, \cdots, 2n+1$.

Para $k = 2$ tenemos

$$(Y^{2n} \widetilde{M})_{1,2} = \widetilde{M}_{2n+1,2} = -t_1 - d,$$

$$(\widetilde{M} Y^{2n-1} \widetilde{M})_{1,2} = \widetilde{M}_{1,*} \cdot \widetilde{M}_{*+2n-1,2} = (1, -1) \cdot (\underbrace{\widetilde{M}_{2n-1,2}}_{0}, \underbrace{\widetilde{M}_{2n,2}}_{t_2}) = -t_2.$$

Calculamos $(Y^{2n-i} \widetilde{M} Y^i)_{1,2} = \widetilde{M}_{2n-i+1, 2-i}$.

Para $i = 1$: $(Y^{2n-i} \widetilde{M} Y^i)_{1,2} = (Y^{2n-1} \widetilde{M} Y)_{1,2} = \widetilde{M}_{2n,1} = t_1$.

Para $i = 2$: $(Y^{2n-i} \widetilde{M} Y^i)_{1,2} = (Y^{2n-2} \widetilde{M} Y^2)_{1,2} = \widetilde{M}_{2n-1,0} = d$.

Para $3 \leq i \leq 2n$: $(Y^{2n-i} \widetilde{M} Y^i)_{1,2} = \widetilde{M}_{2n-i+1, 2-i} = 0$, pues $2 - 2n \leq 2 - i \leq -1 < 0$.

Además $(Y^{2n+1})_{1,2} = \delta_{2n+2, 2} = 0$.

Reemplazando en (4.2.26) en la entrada $(1,2)$ se tiene

$$(d - t_0)(-t_1 - d) = t_0(-t_2) + S_1 t_1 + S_2 d.$$

Pero $S_1 = t_0 + t_1$ y $S_2 = t_0 + t_1 + t_2$.

Luego

$$-d t_1 - d^2 + t_o t_1 + t_0 d = -t_0 t_2 + (t_0 + t_1) t_1 + (t_0 + t_1 + t_2) d,$$
$$-2 t_1 d - d^2 = -t_0 t_2 + t_1^2 + t_2 d = -(t_0 - d) t_2 + t_1^2,$$
$$(t_0 - d) t_2 = t_1^2 + 2 t_1 d + d^2.$$

Es decir

$$(t_0 - d) t_2 = (t_1 + d)^2.$$

Queda demostrado $(i)$.

Para $k = 3, 4, \cdots, 2n+1$ tenemos
$(Y^{2n+1})_{1,k} = \delta_{2n+2, k} = 0$.

$(Y^{2n}\widetilde{M})_{1,k} = \widetilde{M}_{2n+1,k} = -t_{k-1}$   (Proposición 4.1.3 (2)).

Además

$$(\widetilde{M}Y^{2n-1}\widetilde{M})_{1,k} = \widetilde{M}_{1,*} \cdot (Y^{2n-1}\widetilde{M})_{*,k} = \widetilde{M}_{1,*} \cdot \widetilde{M}_{*+2n-1,k} = (1,-1) \cdot (\underbrace{\widetilde{M}_{2n-1,k}}_{0}, \underbrace{\widetilde{M}_{2n,k}}_{t_k}) = -t_k.$$

$$(Y^{2n-i}\widetilde{M}Y^i)_{1,k} = \widetilde{M}_{2n+1-i,k-i} = \begin{cases} \widetilde{M}_{2n,k-1} = t_{k-1} & \text{si} & i = 1, \\ \widetilde{M}_{2n+2-k,1} & \text{si} & i = k-1, \\ \widetilde{M}_{2n+1-k,0} & \text{si} & i = k, \\ 0 & \text{en otros casos.} \end{cases}$$

En efecto si $2 \leq i \leq k-2$ entonces $\boxed{2 \leq k-i \leq k-2}$.
Por otro lado $2 \leq i \leq k-2$ implica $\boxed{0 \leq 2n+1-i \leq 2n-1}$.
Supongamos que $2n+1-i = n$ y $k-i = n+1$ entonces $i = n+1$ y $k = 2n+2$,
lo cual contradice el hecho de que $3 \leq k \leq 2n+1$.
Por lo tanto en los demás casos $\widetilde{M}_{2n+1-i,k-i} = 0$.
Supongamos que $3 \leq k \leq n+1$. Entonces se cumple $-n-1 \leq -k \leq -3$, $n \leq 2n+1-k \leq 2n-2$ y $n+1 \leq 2n+2-k \leq 2n-1$.
Luego $\widetilde{M}_{2n+1-k,0} = d$ y $\widetilde{M}_{2n+2-k,1} = -d$.
De esta manera para $3 \leq k \leq n+1$ la igualdad (4.2.26) se lee

$$(d-t_0)(Y^{2n}\widetilde{M})_{1,k} = t_0(\widetilde{M}Y^{2n-1}\widetilde{M})_{1,k} + \sum_{i=1}^{2n}\sum_{j=0}^{i}(\widehat{t}_j Y^{2n-i}\widetilde{M}Y^i)_{1,k} + (S_{2n+1} - aS_n)(Y^{2n+1})_{1,k},$$

$$(d-t_0)(-t_{k-1}) = -t_0 t_k + \sum_{j=0}^{1}\widehat{t}_j \widetilde{M}_{2n,k-1} + \sum_{j=0}^{k-1}\widehat{t}_j \widetilde{M}_{2n-k+2,1} + \sum_{j=0}^{k}\widehat{t}_j \widetilde{M}_{2n-k+1,0},$$

$$-(d-t_0)t_{k-1} = t_0(-t_k) + t_{k-1}\sum_{j=0}^{1}\widehat{t}_j + (-d)\sum_{j=0}^{k-1}\widehat{t}_j + d\sum_{j=0}^{k}\widehat{t}_j,$$

$$-(d-t_0)t_{k-1} = -t_0 t_k + t_0 t_{k-1} + t_1 t_{k-1} - dS_{k-1} + dS_k = -t_0 t_k + t_0 t_{k-1} + t_1 t_{k-1} + d\widehat{t}_k$$

$$-dt_{k-1} = -t_0 t_k + t_1 t_{k-1} + d\widehat{t}_k.$$

Si $3 \leq k \leq n$ se obtiene

$$-dt_{k-1} = -t_0 t_k + t_1 t_{k-1} + dt_k.$$

Es decir

$$(t_0 - d)t_k = (d + t_1)t_{k-1}.$$

Queda demostrado $(ii)$.
Si $k = n+1$ se obtiene

$$-dt_n = -t_0 t_{n+1} + t_1 t_n + d^2 t_{n+1}.$$

Es decir
$$(t_0 - d^2)t_{n+1} = (d + t_1)t_n,$$

con lo cual queda demostrado $(iii)$.

Por otro lado para $n + 2 \leq k \leq 2n$ tenemos $-2n \leq -k \leq -n-2$, lo cual implica $1 \leq 2n+1-k \leq n-1$ y $2 \leq 2n+2-k \leq n$.

Por lo tanto $\widetilde{M}_{2n+1-k,0} = 1$ y $\widetilde{M}_{2n+2-k,1} = -1$.

En este caso

$$(Y^{2n-i}\widetilde{M}Y^i)_{1,k} = \widetilde{M}_{2n+1-i,k-i} = \begin{cases} \widetilde{M}_{2n,k-1} = t_{k-1} & \text{si} & i = 1, \\ \widetilde{M}_{2n+2-k,1} = -1 & \text{si} & i = k-1, \\ \widetilde{M}_{2n+1-k,0} = 1 & \text{si} & i = k, \\ 0 & \text{en otros casos.} \end{cases}$$

y la igualdad (4.2.26) se lee

$$(d - t_0)(-t_{k-1}) = -t_0 t_k + S_1 t_{k-1} + S_{k-1}(-1) + S_k,$$

$$-(d - t_0)t_{k-1} = -t_0 t_k + t_0 t_{k-1} + t_1 t_{k-1} - S_{k-1} + S_k = -t_0 t_k + t_0 t_{k-1} + t_1 t_{k-1} + \widehat{t}_k,$$

$$-d\, t_{k-1} = -t_0 t_k + t_1 t_{k-1} + \widehat{t}_k.$$

Si $n + 2 \leq k \leq 2n$ se obtiene

$$-d\, t_{k-1} = -t_0 t_k + t_1 t_{k-1} + d\, t_k.$$

Es decir

$$(t_0 - d)t_k = (d + t_1)t_{k-1},$$

y por lo tanto se cumple $(iv)$.

Para $k = 2n + 1$ se tiene que

$(Y^{2n+1})_{1,k} = \delta_{2n+2,2n+1} = 0$.

$(Y^{2n}\widetilde{M})_{1,k} = \widetilde{M}_{2n+1,2n+1} = -t_{2n}$  (Proposición 4.1.3 (2)).

Además $(\widetilde{M}Y^{2n-1}\widetilde{M})_{1,k} = -t_{2n+1}$.

$$(Y^{2n-i}\widetilde{M}Y^i)_{1,k} = \widetilde{M}_{2n+1-i,2n+1-i} = \begin{cases} \widetilde{M}_{2n,2n} = t_{2n} & \text{si} & i = 1, \\ \widetilde{M}_{1,1} = -1 & \text{si} & i = 2n, \\ \widetilde{M}_{0,0} = 1 & \text{si} & i = 2n+1, \\ 0 & \text{en otros casos.} \end{cases}$$

En (4.2.26)

$$(d - t_0)(-t_{2n}) = -t_0 t_{2n+1} + S_1 t_{2n} + S_{2n}(-1),$$

$$-d\, t_{2n} = -t_0 \underbrace{t_{2n+1}}_{a} + t_1 t_{2n} - S_{2n},$$

$$(d + t_1)t_{2n} = t_0 a + S_{2n}.$$

Queda demostrado $(v)$.

Ahora evaluamos (4.2.26) en $(k, 2n+1+k)$ para $k = 1, 2, \cdots, n-1$.

Obtendremos varias veces entradas de la forma $\widetilde{M}_{j,j+1}$.

Por la Proposición 4.1.2(2) sabemos que si $m_j \neq 0$ entonces $m_{j+k} = 0$ para $k = 1, 2, \cdots, n-1$. Luego tendremos

$$m_j = \widetilde{M}_{j,j+1} = \begin{cases} a, & \text{si } j \in \{n, 2n\}, \\ 0, & \text{si } 1 \leq j \leq 3n-1, \; j \notin \{n, 2n\}. \end{cases}$$

En consecuencia

$$(Y^{2n}\widetilde{M})_{k,2n+1+k} = \widetilde{M}_{2n+k,2n+1+k} = 0,$$

pues para $k = 1, 2, \cdots, n-1$ se cumple $2n+1 \leq 2n+k \leq 3n-1$.

También tenemos

$$(Y^{2n-i}\widetilde{M}Y^i)_{k,2n+1+k} = \widetilde{M}_{2n+k-i,2n+1+k-i} = \begin{cases} a, & \text{si } i \in \{k, n+k\}, \\ 0, & \text{de lo contrario.} \end{cases}$$

Observemos que $j = 2n+k-i = n \Leftrightarrow i = n+k$ y $j = 2n+k-i = 2n \Leftrightarrow i = k$.

Dado que

$$(\widetilde{M}Y^{2n-1}\widetilde{M})_{k,2n+1+k} = \widetilde{M}_{k,*} \cdot (Y^{2n-1}\widetilde{M})_{*,2n+1+k} = (1,-1) \cdot (\widetilde{M}_{2n-1,2n+1+k}, \widetilde{M}_{2n,2n+1+k}) = 0,$$

pues $2n+1+k > 2n = (2n-1)+1$ y $2n+1+k > 2n+1$,

y como $(Y^{2n+1})_{k,2n+1+k} = \delta_{2n+1+k,2n+1+k} = 1$, entonces en (4.2.26) obtenemos

$$0 = t_0(0) + (S_k + S_{k+n})a + S_{2n+1} - aS_n.$$

Es decir

$$0 = (S_k + S_{k+n})a + S_{2n+1} - aS_n, \quad k = 1, 2, \cdots, n-1.$$

Para $k = 1$ se tiene

$$0 = (S_1 + S_{1+n})a + S_{2n+1} - aS_n.$$

Para $k = 2, \cdots, n-1$ y para $k-1$ se tiene

$$0 = (S_{k-1} + S_{k-1+n})a + S_{2n+1} - aS_n.$$

Restando esto a la ecuación correspondiente a $k$ se tiene

$$0 = a t_k + \underbrace{(S_{k+n} - S_{k-1+n})}_{\widehat{t}_{k+n}} a.$$

Es decir

$$0 = a t_k + a d t_{k+n}.$$

Esto demuestra $(vi)$ y concluye la demostración del lema. $\square$

**Proposición 4.2.4.** *Supongamos que* $t_{2n+1} = \widetilde{M}_{2n,2n+1} = a$.

(1) *Si $t_0 = d^2$ entonces*

$$\widetilde{M}_{2n,*} = d^2 E_0 - d E_1 + a E_{2n+1}, \quad \widetilde{M}_{2n+1,*} = d^2 E_0 - d^2 E_1.$$

(2) *Si $t_0 = d$ entonces*

$$\widetilde{M}_{2n,*} = d E_0 - d E_1 - a d E_n + a E_{2n+1}, \quad \widetilde{M}_{2n+1,*} = d^2 E_0 - d E_1 + a d E_{n+1}.$$

**Demostración.** Supongamos que $t_0 = d^2$.
Demostraremos primero que

$$t_1 = -d. \tag{4.2.27}$$

En efecto, si asumimos por contradicción que $(t_1 + d) \neq 0$, entonces por $(iii)$ del Lema 4.2.3 se obtiene $0 = (d + t_1)t_n$ y por lo tanto $t_n = 0$.

Luego de $(t_0 - d)t_k = (d + t_1)t_{k-1}$ (Lema 4.2.3 $(ii)$) obtenemos sucesivamente $t_{n-1} = t_{n-2} = \cdots = t_2 = 0$.

Finalmente $(t_0 - d)t_2 = (t_1 + d)^2$ (Lema 4.2.3 $(i)$) produce $0 = (t_1 + d)$.

Esto es una contradicción, y por lo tanto (4.2.27) queda demostrado.

Ahora las igualdades $(i), (ii)$ y $(iv)$ del Lema 4.2.3 producen $t_k = 0$ para $k = 2, 3, \cdots, 2n$, con $k \neq n+1$.

Pero $(d + t_1)t_{2n} = t_0 a + S_{2n}$ (Lema 4.2.3 $(v)$) se lee $0 = t_0 a + S_{2n}$ y así

$$0 = d^2 a + t_0 + t_1 + \cdots + t_n + d t_{n+1} + d t_{n+2} + \cdots + d t_{2n}.$$

Es decir

$$0 = d^2 a + t_0 + t_1 + d t_{n+1} = d^2 a + d^2 - d + d t_{n+1} = \underbrace{d^2(a+1) - d}_{0} + d t_{n+1} = d t_{n+1}$$

Por lo tanto $t_{n+1} = 0$.

De esta manera queda demostrado el caso (1), es decir

$$\widetilde{M}_{2n,*} = d^2 E_0 - d E_1 + a E_{2n+1}, \quad \widetilde{M}_{2n+1,*} = d^2 E_0 - d^2 E_1.$$

Por otro lado si $t_0 = d$ entonces por Lema 4.2.3 $(i)$ en $(t_0 - d)t_2 = (t_1 + d)^2$ se obtiene $t_1 = -d$ y en $(t_0 - d^2)t_{n+1} = (d + t_1)t_n$ (Lema 4.2.3$(iii)$) se obtiene $t_{n+1} = 0$.

Sumando las igualdades $0 = t_k + d t_{k+n}$ (Lema 4.2.3 $(vi)$) para $k = 2, 3, \cdots, n-1$, obtenemos

$$\begin{aligned}
S_{2n} &= t_0 + t_1 + \cdots + t_n + \widehat{t}_{n+1} + \cdots + \widehat{t}_{2n}, \\
&= t_0 + t_1 + \cdots + t_n + d t_{n+1} + \cdots + d t_{2n}, \\
&= (t_0 + t_1) + (t_2 + d t_{n+2}) + (t_3 + d t_{n+3}) + \cdots + (t_{n-1} + d t_{2n-1}) + t_n + d t_{2n}.
\end{aligned}$$

Es decir $S_{2n} = t_n + d\, t_{2n}$.

Pero por $(v)$ del Lema 4.2.3, $(d + t_1)t_{2n} = t_0 a + S_{2n}$, obtenemos $d a + S_{2n} = 0$.

Luego $d a + t_n + d\, t_{2n} = 0$.

Dividiendo entre $ad$:
$$1 + \frac{t_n}{ad} + \frac{t_{2n}}{a} = 0.$$

Demostraremos que
$$f := \frac{t_{2n}}{a} = -\frac{t_n}{ad} - 1 = 0. \qquad (4.2.28)$$

Asumamos por el absurdo que $f \neq 0$.

Entonces evaluamos para $n + 2 \leq k \leq 2n + 1$:

$$\begin{aligned} 0 &= (\widetilde{M}^2)_{2n,k} = \widetilde{M}_{2n,*} \cdot \widetilde{M}_{*,k} \\ &= (t_0, t_1, \cdots, t_k, \cdots, t_{2n}, t_{2n+1}) \cdot (\widetilde{M}_{0,k}, \widetilde{M}_{1,k}, \cdots, \widetilde{M}_{k,k}, \cdots, \widetilde{M}_{2n,k}, \widetilde{M}_{2n+1,k}). \end{aligned}$$

Pero $\widetilde{M}_{j,k} = 0$ para $0 \leq j \leq 2n - 1$, $j \neq n + 1$, pues $2 < n + 2 \leq k \leq 2n + 1$ y además se demostró que $t_{n+1} = 0$.

Luego
$$0 = \underbrace{t_{2n}}_{fa} \underbrace{\widetilde{M}_{2n,k}}_{t_k} + \underbrace{t_{2n+1}}_{a} \underbrace{\widetilde{M}_{2n+1,k}}_{-t_{k-1}}.$$

Es decir
$$0 = a(f\, t_k - t_{k-1}), \quad n + 2 \leq k \leq 2n + 1.$$

Por lo tanto $t_{k-1} = f\, t_k$ para $n + 2 \leq k \leq 2n + 1$.

Se sigue que $\quad t_{n+1} = f\, t_{n+2} = f^2\, t_{n+3} = \cdots = f^j\, t_{n+j+1}$ para $1 \leq j \leq n$.

En particular para $j = n$ se obtiene $t_{n+1} = f^n t_{2n+1} = f^n a \neq 0$, lo cual es una contradicción a $t_{n+1} = 0$, y queda demostrado (4.2.28).

De (4.2.28) deducimos $t_n = -ad$ y de
$$0 = (\widetilde{M}^2)_{2n,k} = \widetilde{M}_{2n,*} \cdot \widetilde{M}_{*,k} = a(f\, t_k - t_{k-1}) = -a\, t_{k-1},$$

obtenemos $t_{k-1} = 0$ para $n + 2 \leq k \leq 2n + 1$. Es decir $t_k = 0$ para $n + 1 \leq k \leq 2n$.

Por otro lado por Lema 4.2.3 $(vi)$: $t_k + d\, t_{k+n} = 0$ para $k = 2, \cdots, n - 1$.

Luego $t_k = -d\, t_{k+n} = 0$ para $k = 2, \cdots, n - 1$, pues $k + n = n + 2, \cdots, 2n - 1$.

Esto da el segundo de los casos deseados

$$\widetilde{M}_{2n,*} = d E_0 - d E_1 - ad E_n + a E_{2n+1} \qquad \widetilde{M}_{2n+1,*} = d^2 E_0 - d E_1 + ad E_{n+1},$$

y concluye la demostración. $\qquad \square$

En el caso de la Proposición 4.2.4 (1), la matriz $\widetilde{M}$ hasta la fila $2n+1$ se ve como sigue:

$$
\begin{array}{r}
0 \to \\
\\
\\
n \to \\
\\
\\
\\
\\
2n \to \\
2n+1 \to
\end{array}
\left(
\begin{array}{cccccccccc}
1 & -1 & 0 & \cdots & 0 & 0 & \cdots & 0 & 0 & 0 \\
\vdots & \vdots & \vdots & \cdots & \vdots & \vdots & \cdots & \vdots & \vdots & \vdots \\
1 & -1 & 0 & \cdots & 0 & 0 & \cdots & 0 & 0 & 0 \\
\boldsymbol{d} & \boldsymbol{-1} & \boldsymbol{0} & \cdots & \boldsymbol{a} & \boldsymbol{0} & \cdots & \boldsymbol{0} & \boldsymbol{0} & \boldsymbol{0} \\
d & -d & 0 & \cdots & 0 & 0 & \cdots & 0 & 0 & 0 \\
\vdots & \vdots & \vdots & \cdots & \vdots & \vdots & \cdots & \vdots & \vdots & \vdots \\
d & -d & 0 & \cdots & 0 & 0 & \cdots & 0 & 0 & 0 \\
\boldsymbol{d^2} & \boldsymbol{-d} & \boldsymbol{0} & \cdots & \boldsymbol{0} & \boldsymbol{0} & \cdots & \boldsymbol{0} & \boldsymbol{a} & \boldsymbol{0} \\
\boldsymbol{d^2} & \boldsymbol{-d^2} & \boldsymbol{0} & \cdots & \boldsymbol{0} & \boldsymbol{0} & \cdots & \boldsymbol{0} & \boldsymbol{0} & \boldsymbol{0}
\end{array}
\right)
$$

En el caso de la Proposición 4.2.4 (2), la matriz $\widetilde{M}$ hasta la fila $2n+1$ se ve como sigue:

$$
\begin{array}{r}
0 \to \\
\\
\\
n \to \\
\\
\\
\\
\\
2n \to \\
2n+1 \to
\end{array}
\left(
\begin{array}{cccccccccc}
1 & -1 & 0 & \cdots & 0 & 0 & \cdots & 0 & 0 & 0 \\
\vdots & \vdots & \vdots & \cdots & \vdots & \vdots & \cdots & \vdots & \vdots & \vdots \\
1 & -1 & 0 & \cdots & 0 & 0 & \cdots & 0 & 0 & 0 \\
\boldsymbol{d} & \boldsymbol{-1} & \boldsymbol{0} & \cdots & \boldsymbol{a} & \boldsymbol{0} & \cdots & \boldsymbol{0} & \boldsymbol{0} & \boldsymbol{0} \\
d & -d & 0 & \cdots & 0 & 0 & \cdots & 0 & 0 & 0 \\
\vdots & \vdots & \vdots & \cdots & \vdots & \vdots & \cdots & \vdots & \vdots & \vdots \\
d & -d & 0 & \cdots & 0 & 0 & \cdots & 0 & 0 & 0 \\
\boldsymbol{d} & \boldsymbol{-d} & \boldsymbol{0} & \cdots & \boldsymbol{-ad} & \boldsymbol{0} & \cdots & \boldsymbol{0} & \boldsymbol{a} & \boldsymbol{0} \\
\boldsymbol{d^2} & \boldsymbol{-d} & \boldsymbol{0} & \cdots & \boldsymbol{0} & \boldsymbol{ad} & \cdots & \boldsymbol{0} & \boldsymbol{0} & \boldsymbol{0}
\end{array}
\right)
$$

En la matriz del primer caso de la Proposición 4.2.4 vemos que $\widetilde{M}_{j,*} = \widetilde{M}_{j-1,*}$ para $n+2 \leq j \leq 2n-1$.

Por el Lema 3.1.2 se sigue que $\widetilde{M} Y^j \widetilde{M} = 0$ para $n+1 \leq j \leq 2n-2$.

Luego

$$0 = (\widetilde{M} Y^j \widetilde{M})_{n,i} = \widetilde{M}_{n,*} \cdot \widetilde{M}_{*+j,i} = (d, -1, 0, \cdots, 0, a) \cdot (\widetilde{M}_{j,i}, \widetilde{M}_{j+1,i}, \cdots, \widetilde{M}_{j+n+1,i}).$$

Es decir

$$0 = d\widetilde{M}_{j,i} - \widetilde{M}_{j+1,i} + a\widetilde{M}_{j+n+1,i}.$$

Caso $\boxed{i=0}$

Se tiene $\widetilde{M}_{j,0} = d$, $\widetilde{M}_{j+1,0} = d$ y se obtiene $\widetilde{M}_{j+n+1,0} = \dfrac{d - d^2}{a} = d^2$.

Caso $\boxed{i=1}$

Se tiene $\widetilde{M}_{j,1} = -d$, $\widetilde{M}_{j+1,1} = -d$ y se obtiene $\widetilde{M}_{j+n+1,1} = \dfrac{d^2 - d}{a} = -d^2$.

Caso $\boxed{i \geq 2}$

En este caso como $n+1 \leq j \leq 2n-2$ se cumple $\widetilde{M}_{j,i} = \widetilde{M}_{j+1,i} = 0$.

Luego se tiene que $0 = a\widetilde{M}_{j+n+1,i}$ y como $a \neq 0$ entonces $\widetilde{M}_{j+n+1,i} = 0$.
Por lo tanto se cumple que

$$\widetilde{M}_{k,*} = (d^2, -d^2, 0, 0, \cdots) \text{ para } 2n+2 \leq k \leq 3n-1.$$

De esta manera la matriz $\widetilde{M}$ hasta la fila $3n-1$ se ve como sigue:

$$\begin{array}{r} 0 \to \\ \\ \\ n \to \\ \\ \\ \\ 2n \to \\ 2n+1 \to \\ \\ \\ 3n-1 \to \end{array} \left( \begin{array}{ccccccccc} 1 & -1 & 0 & \cdots & 0 & 0 & \cdots & 0 & 0 & 0 \\ \vdots & \vdots & \vdots & \cdots & \vdots & \vdots & \cdots & \vdots & \vdots & \vdots \\ 1 & -1 & 0 & \cdots & 0 & 0 & \cdots & 0 & 0 & 0 \\ \boldsymbol{d} & \boldsymbol{-1} & \boldsymbol{0} & \cdots & \boldsymbol{a} & \boldsymbol{0} & \cdots & \boldsymbol{0} & \boldsymbol{0} & \boldsymbol{0} \\ d & -d & 0 & \cdots & 0 & 0 & \cdots & 0 & 0 & 0 \\ \vdots & \vdots & \vdots & \cdots & \vdots & \vdots & \cdots & \vdots & \vdots & \vdots \\ d & -d & 0 & \cdots & 0 & 0 & \cdots & 0 & 0 & 0 \\ \boldsymbol{d^2} & \boldsymbol{-d} & \boldsymbol{0} & \cdots & \boldsymbol{0} & \boldsymbol{0} & \cdots & \boldsymbol{0} & \boldsymbol{a} & \boldsymbol{0} \\ \boldsymbol{d^2} & \boldsymbol{-d^2} & \boldsymbol{0} & \cdots & \boldsymbol{0} & \boldsymbol{0} & \cdots & \boldsymbol{0} & \boldsymbol{0} & \boldsymbol{0} \\ d^2 & -d^2 & 0 & \cdots & 0 & 0 & \cdots & 0 & 0 & 0 \\ \vdots & \vdots & \vdots & \cdots & \vdots & \vdots & \cdots & \vdots & \vdots \\ d^2 & -d^2 & 0 & \cdots & 0 & 0 & \cdots & 0 & 0 & 0 \end{array} \right)$$

Ahora usamos la Proposición 4.1.3 para $k = 2$ donde $L_1 = n$, $L_2 = 2n$.
Se cumple para $i \leq 3n-1$ y $j = 1, 2$ que

$$\widetilde{M}_{i,*} = \begin{cases} dE_0 - E_1 + aE_{n+1}, & \text{si} \quad i = n, \\ d\widetilde{M}_{0,*}, & \text{si} \quad n < i < 2n, \\ d^2 E_0 - dE_1 + aE_{2n+1}, & \text{si} \quad i = 2n, \\ d^2 \widetilde{M}_{0,*}, & \text{si} \quad 2n < i \leq 3n-1. \end{cases}$$

Luego, haciendo $t_i = \widetilde{M}_{3n,i}$, se cumple

(1) $(\widetilde{M}_{3n+1,0}, \widetilde{M}_{3n+1,1}, \widetilde{M}_{3n+1,2}) = (d^3, -t_0, -t_1 - d^2)$,

(2) $\widetilde{M}_{3n+1,j} = -t_{j-1}$ para $3 \leq j \leq 3n+1$,

(3) $t_0 \in \{d^2, d^3\}$,

(4) $(\widetilde{M}_{3n,3n+1}, \widetilde{M}_{3n+1,3n+2}) \in \{(a,0), (0,a)\}$.

En el primer caso de la Proposición 4.2.4 la matriz $\widetilde{M}$ hasta la fila $3n+2$ se ve como sigue:

$$
\begin{array}{c}
0 \to \\
\\
n \to \\
\\
\\
\\
\\
2n \to \\
2n+1 \to \\
\\
\\
\\
3n \to \\
3n+1 \to \\
3n+2 \to
\end{array}
\left(\begin{array}{ccccccccccc}
1 & -1 & 0 & \cdots & 0 & 0 & \cdots & 0 & \cdots & 0 & 0 \\
\vdots & \vdots & \vdots & \cdots & \vdots & \vdots & \cdots & \vdots & \cdots & \vdots & \vdots \\
1 & -1 & 0 & \cdots & 0 & 0 & \cdots & 0 & \cdots & 0 & 0 \\
\boldsymbol{d} & \boldsymbol{-1} & \boldsymbol{0} & \cdots & \boldsymbol{a} & \boldsymbol{0} & \cdots & \boldsymbol{0} & \cdots & \boldsymbol{0} & \boldsymbol{0} \\
d & -d & 0 & \cdots & 0 & 0 & \cdots & 0 & \cdots & 0 & 0 \\
\vdots & \vdots & \vdots & \cdots & \vdots & \vdots & \cdots & \vdots & \cdots & \vdots & \vdots \\
d & -d & 0 & \cdots & 0 & 0 & \cdots & 0 & \cdots & 0 & 0 \\
\boldsymbol{d^2} & \boldsymbol{-d} & \boldsymbol{0} & \cdots & \boldsymbol{0} & \boldsymbol{0} & \cdots & \boldsymbol{a} & \cdots & \boldsymbol{0} & \boldsymbol{0} \\
\boldsymbol{d^2} & \boldsymbol{-d^2} & \boldsymbol{0} & \cdots & \boldsymbol{0} & \boldsymbol{0} & \cdots & \boldsymbol{0} & \cdots & \boldsymbol{0} & \boldsymbol{0} \\
d^2 & -d^2 & 0 & \cdots & 0 & 0 & \cdots & 0 & \cdots & 0 & 0 \\
\vdots & \vdots & \vdots & \cdots & \vdots & \vdots & \cdots & \vdots & \cdots & \vdots & \vdots \\
d^2 & -d^2 & 0 & \cdots & 0 & 0 & \cdots & 0 & \cdots & 0 & 0 \\
\boldsymbol{t_0} & \boldsymbol{t_1} & \boldsymbol{t_2} & \cdots & \boldsymbol{t_{n+1}} & \boldsymbol{t_{n+2}} & \cdots & \boldsymbol{t_{2n+1}} & \cdots & \boldsymbol{t_{3n+1}} & \boldsymbol{0} \\
\boldsymbol{d^3} & \boldsymbol{-t_0} & \boldsymbol{u} & \cdots & \boldsymbol{-t_n} & \boldsymbol{-t_{n+1}} & \cdots & \boldsymbol{-t_{2n}} & \cdots & \boldsymbol{-t_{3n}} & \boldsymbol{a-t_{3n+1}} \\
\boldsymbol{d^3} & \boldsymbol{-d^3} & \boldsymbol{0} & \cdots & \boldsymbol{0} & \boldsymbol{0} & \cdots & \boldsymbol{0} & \cdots & \boldsymbol{0} & \boldsymbol{0}
\end{array}\right)
$$

donde $\mathbf{u} = -\boldsymbol{d^2} - \boldsymbol{t_1}$.

Para la fila $3n+2$ observemos que $\widetilde{M}_{2n+1,*} = \widetilde{M}_{2n+2,*}$.

Luego, por el Lema 3.1.2, se sigue que $\widetilde{M} Y^{2n+1} \widetilde{M} = 0$.

Evaluando en la entrada $(n, i)$ se tiene

$$0 = (\widetilde{M} Y^{2n+1} \widetilde{M})_{n,i} = \widetilde{M}_{n,*} \cdot \widetilde{M}_{*+2n+1,i} = (d, -1, 0, \cdots, 0, a) \cdot (\widetilde{M}_{2n+1,i}, \widetilde{M}_{2n+2,i}, \cdots, \widetilde{M}_{3n+2,i}).$$

Es decir
$$0 = d\widetilde{M}_{2n+1,i} - \widetilde{M}_{2n+2,i} + a\widetilde{M}_{3n+2,i}.$$

Caso $\boxed{i = 0}$

Se tiene $\widetilde{M}_{2n+1,0} = d^2$, $\widetilde{M}_{2n+2,0} = d^2$. Luego $\widetilde{M}_{3n+2,0} = \dfrac{d^2 - d^3}{a} = d^3$.

Caso $\boxed{i = 1}$

Se tiene $\widetilde{M}_{2n+1,1} = -d^2$, $\widetilde{M}_{2n+2,1} = -d^2$. Luego $\widetilde{M}_{3n+2,1} = \dfrac{d^3 - d^2}{a} = -d^3$.

Caso $\boxed{i \geq 2}$.

En este caso se cumple $\widetilde{M}_{2n+1,i} = \widetilde{M}_{2n+2,i} = 0$.

Luego se tiene que $0 = a\widetilde{M}_{3n+2,i}$ y como $a \neq 0$ entonces $\widetilde{M}_{3n+2,i} = 0$.

Por lo tanto se cumple
$$\widetilde{M}_{3n+2,*} = (d^3, -d^3, 0, 0, \cdots),$$

lo cual concluye la justificación de la forma de la matriz hasta la fila $3n+2$, en el primer caso de la Proposición 4.2.4.

Para encontrar la forma de la matriz en el primer caso de la Proposición 4.2.2 vemos que $\widetilde{M}_{j,*} = \widetilde{M}_{j-1,*}$ para $n+2 \leq j \leq 2n$.

Por el Lema 4.1.1 se sigue que $\widetilde{M} Y^j \widetilde{M} = 0$ para $n+1 \leq j \leq 2n-1$.

Luego

$$0 = (\widetilde{M} Y^j \widetilde{M})_{n,i} = \widetilde{M}_{n,*} \cdot \widetilde{M}_{*+j,i} = (d, -1, 0, \cdots, 0, a) \cdot (\widetilde{M}_{j,i}, \widetilde{M}_{j+1,i}, \cdots, \widetilde{M}_{j+n+1,i}).$$

Es decir

$$0 = d\widetilde{M}_{j,i} - \widetilde{M}_{j+1,i} + a\widetilde{M}_{j+n+1,i}.$$

Caso $\boxed{i=0}$

Se tiene $\widetilde{M}_{j,0} = d$, $\widetilde{M}_{j+1,0} = d$ y se obtiene $\widetilde{M}_{j+n+1,0} = \dfrac{d-d^2}{a} = d^2$.

Caso $\boxed{i=1}$

Se tiene $\widetilde{M}_{j,1} = -d$, $\widetilde{M}_{j+1,1} = -d$ y se obtiene $\widetilde{M}_{j+n+1,1} = \dfrac{d^2-d}{a} = -d^2$.

Caso $\boxed{i \geq 2}$

En este caso como $n+1 \leq j \leq 2n-1$ se cumple $\widetilde{M}_{j,i} = \widetilde{M}_{j+1,i} = 0$.

Luego se tiene que $0 = a\widetilde{M}_{j+n+1,i}$ y como $a \neq 0$ entonces $\widetilde{M}_{j+n+1,i} = 0$.

Por lo tanto se cumple

$$\widetilde{M}_{k,*} = (d^2, -d^2, 0, 0, \cdots) \text{ para } 2n+2 \leq k \leq 3n.$$

De esta manera la matriz $\widetilde{M}$ hasta la fila $3n$ se ve como sigue:

$$\begin{array}{r} 0 \to \\ \\ \\ n \to \\ \\ \\ \\ \\ 2n \to \\ 2n+1 \to \\ \\ \\ 3n \to \end{array} \left( \begin{array}{cccccccccc} 1 & -1 & 0 & \cdots & 0 & 0 & \cdots & 0 & 0 \\ \vdots & \vdots & \vdots & \cdots & \vdots & \vdots & \cdots & \vdots & \vdots \\ 1 & -1 & 0 & \cdots & 0 & 0 & \cdots & 0 & 0 \\ \boldsymbol{d} & \boldsymbol{-1} & \boldsymbol{0} & \cdots & \boldsymbol{a} & \boldsymbol{0} & \cdots & \boldsymbol{0} & \boldsymbol{0} \\ d & -d & 0 & \cdots & 0 & 0 & \cdots & 0 & 0 \\ \vdots & \vdots & \vdots & \cdots & \vdots & \vdots & \cdots & \vdots & \vdots \\ d & -d & 0 & \cdots & 0 & 0 & \cdots & 0 & 0 \\ \boldsymbol{d} & \boldsymbol{-d} & \boldsymbol{0} & \cdots & \boldsymbol{0} & \boldsymbol{0} & \cdots & \boldsymbol{0} & \boldsymbol{0} \\ \boldsymbol{d^2} & \boldsymbol{-d} & \boldsymbol{0} & \cdots & \boldsymbol{0} & \boldsymbol{0} & \cdots & \boldsymbol{0} & \boldsymbol{a} \\ d^2 & -d^2 & 0 & \cdots & 0 & 0 & \cdots & 0 & 0 \\ \vdots & \vdots & \vdots & \cdots & \vdots & \vdots & \cdots & \vdots \\ d^2 & -d^2 & 0 & \cdots & 0 & 0 & \cdots & 0 & 0 \end{array} \right)$$

Ahora usamos la Proposición 4.1.3 para $k=2$ donde $L_1 = n$, $L_2 = 2n+1$.

Se cumple para $i \leq 3n$ y $j = 1, 2$ que

$$\widetilde{M}_{i,*} = \begin{cases} dE_0 - E_1 + aE_{n+1} & \text{si} \quad i = n \\ d\widetilde{M}_{0,*} & \text{si} \quad n < i < 2n+1 \\ d^2 E_0 - dE_1 + aE_{2n+2} & \text{si} \quad i = 2n+1 \\ d^2 \widetilde{M}_{0,*} & \text{si} \quad 2n+1 < i \leq 3n \end{cases}$$

Luego, haciendo $t_i = \widetilde{M}_{3n+1,i}$, se cumple

(1) $(\widetilde{M}_{3n+2,0}, \widetilde{M}_{3n+2,1}, \widetilde{M}_{3n+2,2}) = (d^3, -t_0, -t_1 - d^2)$,

(2) $\widetilde{M}_{3n+2,j} = -t_{j-1}$ para $3 \leq j \leq 3n+2$,

(3) $t_0 \in \{d^2, d^3\}$,

(4) $(\widetilde{M}_{3n+1,3n+2}, \widetilde{M}_{3n+2,3n+3}) \in \{(a,0),(0,a)\}$.

Por lo tanto, en el primer caso de la Proposición 4.2.2 la matriz $\widetilde{M}$ hasta la fila $3n+2$ se ve como sigue:

$$
\begin{array}{r}
0 \to \\
\\
\\
n \to \\
\\
\\
\\
2n \to \\
2n+1 \to \\
\\
\\
\\
3n+1 \to \\
3n+2 \to
\end{array}
\left(
\begin{array}{cccccccccccc}
1 & -1 & 0 & \cdots & 0 & 0 & \cdots & 0 & \cdots & 0 & 0 \\
\vdots & \vdots & \vdots & \cdots & \vdots & \vdots & \cdots & \vdots & \cdots & \vdots & \vdots \\
1 & -1 & 0 & \cdots & 0 & 0 & \cdots & 0 & \cdots & 0 & 0 \\
\boldsymbol{d} & \boldsymbol{-1} & \boldsymbol{0} & \cdots & \boldsymbol{a} & \boldsymbol{0} & \cdots & \boldsymbol{0} & \cdots & 0 & 0 \\
d & -d & 0 & \cdots & 0 & 0 & \cdots & 0 & \cdots & \boldsymbol{0} & \boldsymbol{0} \\
\vdots & \vdots & \vdots & \cdots & \vdots & \vdots & \cdots & \vdots & \cdots & \vdots & \vdots \\
d & -d & 0 & \cdots & 0 & 0 & \cdots & 0 & \cdots & 0 & 0 \\
\boldsymbol{d} & \boldsymbol{-d} & \boldsymbol{0} & \cdots & \boldsymbol{0} & \boldsymbol{0} & \cdots & \boldsymbol{0} & \cdots & \boldsymbol{0} & \boldsymbol{0} \\
\boldsymbol{d^2} & \boldsymbol{-d} & \boldsymbol{0} & \cdots & \boldsymbol{0} & \boldsymbol{0} & \cdots & \boldsymbol{a} & \cdots & \boldsymbol{0} & \boldsymbol{0} \\
d^2 & -d^2 & 0 & \cdots & 0 & 0 & \cdots & 0 & \cdots & 0 & 0 \\
\vdots & \vdots & \vdots & \cdots & \vdots & \vdots & \cdots & \vdots & \cdots & \vdots & \vdots \\
d^2 & -d^2 & 0 & \cdots & 0 & 0 & \cdots & 0 & \cdots & 0 & 0 \\
\boldsymbol{t_0} & \boldsymbol{t_1} & \boldsymbol{t_2} & \cdots & \boldsymbol{t_{n+1}} & \boldsymbol{t_{n+2}} & \cdots & \boldsymbol{t_{2n+2}} & \cdots & \boldsymbol{t_{3n+2}} & 0 \\
\boldsymbol{d^3} & \boldsymbol{-t_0} & \boldsymbol{u} & \cdots & \boldsymbol{-t_n} & \boldsymbol{-t_{n+1}} & \cdots & \boldsymbol{-t_{2n+1}} & \cdots & \boldsymbol{-t_{3n+1}} & \boldsymbol{a - t_{3n+2}}
\end{array}
\right)
$$

donde $\mathbf{u} = -\boldsymbol{d}^2 - \boldsymbol{t}_1$.

En la matriz del segundo caso de la Proposición 4.2.2 vemos que $\widetilde{M}_{j,*} = \widetilde{M}_{j-1,*}$ para $n+2 \leq j \leq 2n-1$.
Por el Lema 4.1.1 se sigue que $\widetilde{M} Y^j \widetilde{M} = 0$ para $n+1 \leq j \leq 2n-2$.
Luego

$$0 = (\widetilde{M} Y^j \widetilde{M})_{n,i} = \widetilde{M}_{n,*} \cdot \widetilde{M}_{*+j,i} = (d, -1, 0, \cdots, 0, a) \cdot (\widetilde{M}_{j,i}, \widetilde{M}_{j+1,i}, \cdots, \widetilde{M}_{j+n+1,i}).$$

Es decir
$$0 = d\widetilde{M}_{j,i} - \widetilde{M}_{j+1,i} + a\widetilde{M}_{j+n+1,i}.$$

Caso $\boxed{i = 0}$

Se tiene $\widetilde{M}_{j,0} = d$, $\widetilde{M}_{j+1,0} = d$ y se obtiene $\widetilde{M}_{j+n+1,0} = \dfrac{d - d^2}{a} = d^2$.

Caso $\boxed{i=1}$

Se tiene $\widetilde{M}_{j,1} = -d$, $\widetilde{M}_{j+1,1} = -d$ y se obtiene $\widetilde{M}_{j+n+1,1} = \dfrac{d^2-d}{a} = -d^2$.

Caso $\boxed{i \geq 2}$

En este caso como $n+1 \leq j \leq 2n-2$ se cumple $\widetilde{M}_{j,i} = \widetilde{M}_{j+1,i} = 0$.
Luego se tiene que $0 = a\widetilde{M}_{j+n+1,i}$ y como $a \neq 0$ entonces $\widetilde{M}_{j+n+1,i} = 0$.
Por lo tanto se cumple

$$\widetilde{M}_{k,*} = (d^2, -d^2, 0, 0, \cdots) \text{ para } 2n+2 \leq k \leq 3n-1.$$

De esta manera la matriz $\widetilde{M}$ hasta la fila $3n-1$ se ve como sigue:

$$\begin{array}{r} 0 \to \\ \\ \\ n \to \\ \\ \\ \\ \\ 2n \to \\ 2n+1 \to \\ \\ \\ 3n-1 \to \end{array} \left( \begin{array}{cccccccccc} 1 & -1 & 0 & \cdots & 0 & 0 & \cdots & 0 & 0 \\ \vdots & \vdots & \vdots & \vdots & \vdots & \vdots & \vdots & \vdots & \vdots \\ 1 & -1 & 0 & \cdots & 0 & 0 & \cdots & 0 & 0 \\ \boldsymbol{d} & \boldsymbol{-1} & \boldsymbol{0} & \cdots & \boldsymbol{a} & \boldsymbol{0} & \cdots & \boldsymbol{0} & \boldsymbol{0} \\ d & -d & 0 & \cdots & 0 & 0 & \cdots & 0 & 0 \\ \vdots & \vdots & \vdots & \cdots & \vdots & \vdots & \vdots & \vdots & \vdots \\ d & -d & 0 & \cdots & 0 & 0 & \cdots & 0 & 0 \\ \boldsymbol{d^2} & \boldsymbol{-d} & \boldsymbol{0} & \cdots & \boldsymbol{ad} & \boldsymbol{0} & \cdots & \boldsymbol{0} & \boldsymbol{0} \\ \boldsymbol{d^2} & \boldsymbol{-d^2} & \boldsymbol{0} & \cdots & \boldsymbol{0} & \boldsymbol{-ad} & \cdots & \boldsymbol{0} & \boldsymbol{a} \\ d^2 & -d^2 & 0 & \cdots & 0 & 0 & \cdots & 0 & 0 \\ \vdots & \vdots & \vdots & \vdots & \vdots & \vdots & \vdots & \vdots \\ d^2 & -d^2 & 0 & \cdots & 0 & 0 & \cdots & 0 & 0 \end{array} \right)$$

Consideremos la igualdad (4.1.1)

$$\widetilde{M}Y^{n-1}\widetilde{M} = a\sum_{i=0}^{n} Y^i \widetilde{M} Y^{n-i} - a^2 Y^{n+1}.$$

Evaluando cada término de la ecuación en la entrada $(2n, j)$ obtenemos:

$$(\widetilde{M}Y^{n-1}\widetilde{M})_{2n,j} = \widetilde{M}_{2n,*} \cdot \widetilde{M}_{*+n-1,j} = (d^2, -d, \cdots, ad) \cdot (\widetilde{M}_{n-1,j}, \widetilde{M}_{n,j}, \cdots, \widetilde{M}_{2n,j})$$

$$= d^2 \widetilde{M}_{n-1,j} - d\widetilde{M}_{n,j} + ad\widetilde{M}_{2n,j},$$

$$(Y^i \widetilde{M} Y^{n-i})_{2n,j} = (\widetilde{M} Y^{n-i})_{2n+i,j} = \widetilde{M}_{2n+i, j-n+i}$$

$$(Y^{n+1})_{2n,j} = \delta_{3n+1,j}.$$

Para $j = 0$ en la sumatoria sólo sobrevive el término correspondiente a $i = n$ pues si $i < n$ se tiene $\widetilde{M}_{2n+i,-n+i} = 0$.

De esta manera reemplazando los términos correspondientes a $j = 0$ se obtiene

$$d^2 \widetilde{M}_{n-1,0} - d\widetilde{M}_{n,0} + ad\widetilde{M}_{2n,0} = a\widetilde{M}_{3n,0} - a^2 \delta_{3n+1,0}.$$

Es decir
$$a\widetilde{M}_{3n,0} = d^2 - d^2 + ad^3 = ad^3.$$

Por lo tanto $\boxed{\widetilde{M}_{3n,0} = d^3}$.

Para $j = 1$ en la sumatoria sólo sobreviven los términos correspondientes a $i = n$, $i = n-1$ pues si $i < n-1$ se tiene $\widetilde{M}_{2n+i,1-n+i} = 0$.

De esta manera reemplazando los términos correspondientes a $j = 1$ se obtiene

$$d^2\widetilde{M}_{n-1,1} - d\widetilde{M}_{n,1} + ad\widetilde{M}_{2n,1} = a(\widetilde{M}_{3n-1,0} + \widetilde{M}_{3n,1}) - a^2\delta_{3n+1,1}.$$

Es decir
$$a\widetilde{M}_{3n,1} = \underbrace{-d^2 + d}_{ad^2} - ad^2 - ad^2 = -ad^2$$

Por lo tanto $\boxed{\widetilde{M}_{3n,1} = -d^2}$.

Para $j = n+1$ se tiene

$$(\widetilde{M}Y^{n-1}\widetilde{M})_{2n,n+1} = d^2\widetilde{M}_{n-1,n+1} - d\widetilde{M}_{n,n+1} + ad\widetilde{M}_{2n,n+1} = -ad + a^2d^2.$$

$$(Y^i\widetilde{M}Y^{n-i})_{2n,n+1} = (\widetilde{M}Y^{n-i})_{2n+i,n+1} = \widetilde{M}_{2n+i,1+i}$$

Si $0 \leq i \leq n$ entonces $1 \leq i+1 \leq n+1$ y $2n \leq 2n+i \leq 3n$

Caso 1. Si $1 \leq i \leq n-1$ entonces $2 \leq i+1 \leq n$ y $2n+1 \leq 2n+i \leq 3n-1$.

En este caso $\widetilde{M}_{2n+i,1+i} = 0$

Caso 2. Si $i = 0$ entonces $\widetilde{M}_{2n+i,1+i} = \widetilde{M}_{2n,1} = -d$.

Caso 3. Si $i = n$ entonces $\widetilde{M}_{2n+i,1+i} = \widetilde{M}_{3n,n+1}$.

Por último
$$(Y^{n+1})_{2n,n+1} = \delta_{3n+1,n+1} = 0$$

Reemplazando se obtiene
$$-ad + a^2d^2 = a(-d + \widetilde{M}_{3n,n+1}).$$

Es decir $\boxed{\widetilde{M}_{3n,n+1} = ad^2}$.

Para $2 \leq j \leq n$ en la sumatoria sobreviven los términos correspondientes a $n-j \leq i \leq n$ pues $j - n + i \geq 0$.

Ahora si $n - j \leq i \leq n$ entonces $3n - j \leq 2n + i \leq 3n$.

Caso 1. Si $i = n-j$ entonces $2n \leq 2n+i = 3n-j \leq 3n-2$ y $j-n+i = 0$.
Luego $\widetilde{M}_{2n+i,j-n+i} = \widetilde{M}_{3n-j,0} = d^2$.

Caso 2. Si $i = n-j+1$ entonces $2n+1 \leq 2n+i = 3n-j+1 \leq 3n-1$ y $j-n+i = 1$.
Luego $\widetilde{M}_{2n+i,j-n+i} = \widetilde{M}_{3n-j+1,1} = -d^2$.

Caso 3. Si $n-j+2 \leq i \leq n$ entonces $3n-j+2 \leq 2n+i \leq 3n$.
Pero $2n+2 \leq 3n-j+2 \leq 3n$.

Luego $2n+2 \leq 2n+i \leq 3n$.

Además $2 \leq j-n+i \leq j \leq n$.

Por lo tanto $\widetilde{M}_{2n+i,j-n+i} = \begin{cases} 0, & n-j+2 \leq i \leq n-1, \\ \widetilde{M}_{3n,j}, & i = n. \end{cases}$

Por otro lado
$$(Y^{n+1})_{2n,j} = \delta_{3n+1,j} = 0$$

pues $2 \leq j \leq n$.

Reemplazando en
$$d^2 \widetilde{M}_{n-1,j} - d\widetilde{M}_{n,j} + ad\widetilde{M}_{2n,j} = a\sum_{i=0}^{n}(Y^i \widetilde{M} Y^{n-i})_{2n,j} - a^2(Y^{n+1})_{2n,j}.$$

se obtiene
$$0 = a(d^2 - d^2 + \widetilde{M}_{3n,j})$$

Por lo tanto $\widetilde{M}_{3n,j} = 0$, para $2 \leq j \leq n$.

Si $n+2 \leq j \leq 3n$, entonces
$$(Y^{n+1})_{2n,j} = \delta_{3n+1,j} = 0$$

y
$$(\widetilde{M} Y^{n-1} \widetilde{M})_{2n,j} = d^2 \widetilde{M}_{n-1,j} - d\widetilde{M}_{n,j} + ad\widetilde{M}_{2n,j} = 0,$$

pues $\widetilde{M}_{k,j} = 0$ para $k \leq 2n$ (por ejemplo ver la matriz más arriba). Entonces obtenemos
$$0 = a\left(\sum_{i=0}^{n} Y^i \widetilde{M} Y^{n-i}\right)_{2n,j}.$$

Pero $a \neq 0$ y para $2 \leq i \leq n-1$ se tiene
$$\left(Y^i \widetilde{M} Y^{n-i}\right)_{2n,j} = \widetilde{M}_{2n+i,j+i-n} = 0,$$

pues $2n+2 \leq 2n+i \leq 3n-1$ implica $\widetilde{M}_{2n+i,*} = (d^2, -d^2, 0, \ldots)$ y entonces $j+i-n > 1$ implica $\widetilde{M}_{2n+i,j+i-n} = 0$.

Por lo tanto tenemos
$$0 = \left(\sum_{i=0}^{n} Y^i \widetilde{M} Y^{n-i}\right)_{2n,j} = \left(\widetilde{M} Y^n\right)_{2n,j} + \left(Y\widetilde{M} Y^{n-1}\right)_{2n,j} + \left(Y^n \widetilde{M}\right)_{2n,j}$$
$$= \widetilde{M}_{2n,j-n} + \widetilde{M}_{2n+1,j+1-n} + \widetilde{M}_{3n,j}.$$

Pero $\widetilde{M}_{2n,*} = d^2 E_0 - dE_1 + adE_{n+1}$ y $\widetilde{M}_{2n+1,*} = d^2 E_0 - d^2 E_1 - adE_{n+2} + aE_{2n+2}$.

Entonces, si $n+2 \leq j \leq 3n$, se tiene que $\widetilde{M}_{2n,j-n} = 0$ si $j \neq 2n+1$ y $\widetilde{M}_{2n,2n+1-n} =$

$\widetilde{M}_{2n,n+1} = ad$. También se tiene $\widetilde{M}_{2n+1,j+1-n} = 0$ si $j \neq 2n+1$ y $\widetilde{M}_{2n+1,2n+2-n} = \widetilde{M}_{2n+1,n+2} = -ad$. Por lo tanto, se tiene

$$\widetilde{M}_{3n,j} = -(\widetilde{M}_{2n,j-n} + \widetilde{M}_{2n+1,j+1-n}) = 0, \quad \text{para } n+2 \leq j \leq 3n.$$

Si $j = 3n+1$ entonces la sumatoria de (4.1.1) es

$$\sum_{i=0}^{n}(Y^i \widetilde{M} Y^{n-i})_{2n,3n+1} = \sum_{i=0}^{n} \widetilde{M}_{2n+i,2n+1+i}$$

$$= \underbrace{\widetilde{M}_{2n,2n+1}}_{0} + \underbrace{\widetilde{M}_{2n+1,2n+2}}_{a} + \cdots + \underbrace{\widetilde{M}_{2n+i,2n+1+i}}_{0} + \cdots + \underbrace{\widetilde{M}_{3n-1,3n}}_{0} + \widetilde{M}_{3n,3n+1}.$$

Como $j = 3n+1$ se obtiene $0 = a(a + \widetilde{M}_{3n,3n+1}) - a^2$, es decir $\boxed{\widetilde{M}_{3n,3n+1} = 0}$.

En la matriz del segundo caso de la Proposición 4.2.2, las primeras $3n$ filas de la matriz $\widetilde{M}$ son

$$\begin{array}{r} 0 \to \\ \\ \\ n \to \\ \\ \\ \\ 2n \to \\ 2n+1 \to \\ \\ \\ \\ 3n \to \end{array} \left( \begin{array}{ccccccccccc} 1 & -1 & 0 & \cdots & 0 & 0 & \cdots & 0 & 0 & \cdots & 0 & 0 \\ \vdots & \vdots & \vdots & \vdots & \vdots & \vdots & \vdots & \vdots & \vdots & \cdots & \vdots & \vdots \\ 1 & -1 & 0 & \cdots & 0 & 0 & \cdots & 0 & 0 & \cdots & 0 & 0 \\ \mathbf{d} & \mathbf{-1} & \mathbf{0} & \cdots & \mathbf{a} & \mathbf{0} & \cdots & \mathbf{0} & \mathbf{0} & \cdots & 0 & 0 \\ d & -d & 0 & \cdots & 0 & 0 & \cdots & 0 & 0 & \cdots & \mathbf{0} & \mathbf{0} \\ \vdots & \vdots & \vdots & \cdots & \vdots & \vdots & \vdots & \vdots & \vdots & \cdots & \vdots & \vdots \\ d & -d & 0 & \cdots & 0 & 0 & \cdots & 0 & 0 & \cdots & 0 & 0 \\ d^2 & -d & \mathbf{0} & \cdots & \mathbf{ad} & \mathbf{0} & \cdots & \mathbf{0} & \mathbf{0} & \cdots & \mathbf{0} & \mathbf{0} \\ d^2 & -d^2 & \mathbf{0} & \cdots & \mathbf{0} & -ad & \cdots & \mathbf{0} & a & \cdots & \mathbf{0} & \mathbf{0} \\ d^2 & -d^2 & 0 & \cdots & 0 & 0 & \cdots & 0 & 0 & \cdots & 0 & 0 \\ \vdots & \vdots & \vdots & \vdots & \vdots & \vdots & \vdots & \vdots & \vdots & \cdots & \vdots & \vdots \\ d^2 & -d^2 & 0 & \cdots & 0 & 0 & \cdots & 0 & 0 & \cdots & 0 & 0 \\ \mathbf{d^3} & \mathbf{-d^2} & \mathbf{0} & \cdots & \mathbf{ad^2} & \mathbf{0} & \cdots & \mathbf{0} & \mathbf{0} & \cdots & \mathbf{0} & \mathbf{0} \end{array} \right)$$

Ahora pongamos $t_i = \widetilde{M}_{3n+1,i}$. Es decir la fila $3n+1$ de $\widetilde{M}$ es

$$(t_0, t_1, t_2, \cdots, t_{3n+1}, t_{3n+2}, 0, \cdots).$$

Evaluando cada término de la ecuación (4.1.1) en la entrada $(2n+2, j)$ obtenemos:

$$(\widetilde{M} Y^{n-1} \widetilde{M})_{2n+2,j} = \widetilde{M}_{2n+2,*} \cdot \widetilde{M}_{*+n-1,j} = (d^2, -d^2) \cdot (\widetilde{M}_{n-1,j}, \widetilde{M}_{n,j})$$

$$= d^2 \widetilde{M}_{n-1,j} - d^2 \widetilde{M}_{n,j},$$

$$(Y^i \widetilde{M} Y^{n-i})_{2n+2,j} = (\widetilde{M} Y^{n-i})_{2n+2+i,j} = \widetilde{M}_{2n+2+i,j-n+i},$$

$$(Y^{n+1})_{2n+2,j} = \delta_{3n+3,j}.$$

Para $j = 0$ en la sumatoria sólo sobrevive el término correspondiente a $i = n$ pues si $i < n$ se tiene $\widetilde{M}_{2n+2+i,-n+i} = 0$.

De esta manera reemplazando los términos correspondientes a $j = 0$ se obtiene

$$d^2 \widetilde{M}_{n-1,0} - d^2 \widetilde{M}_{n,0} = a\widetilde{M}_{3n+2,0} - a^2 \delta_{3n+3,0},$$

$$d^2 - d^3 = a\widetilde{M}_{3n+2,0}.$$

Es decir

$$ad^3 = a\widetilde{M}_{3n+2,0}.$$

Por lo tanto $\boxed{\widetilde{M}_{3n+2,0} = d^3}$.

Para $j = 1$ en la sumatoria sólo sobreviven los términos correspondientes a $i = n$, $i = n-1$ pues si $i < n-1$ se tiene $\widetilde{M}_{2n+2+i,1-n+i} = 0$.

De esta manera reemplazando los términos correspondientes a $j = 1$ se obtiene

$$d^2 \widetilde{M}_{n-1,1} - d^2 \widetilde{M}_{n,1} = a(\widetilde{M}_{3n+2,1} + \widetilde{M}_{3n+1,0}) - a^2 \delta_{3n+3,1},$$

$$-d^2 + d^2 = a(\widetilde{M}_{3n+2,1} + t_0).$$

Es decir

$$0 = a(\widetilde{M}_{3n+2,1} + t_0).$$

Por lo tanto $\boxed{\widetilde{M}_{3n+2,1} = -t_0}$.

Para $j = 2$ en la sumatoria sólo sobreviven los términos correspondientes a $i = n$, $i = n-1$, $i = n-2$ pues si $i < n-2$ se tiene $\widetilde{M}_{2n+2+i,2-n+i} = 0$.

De esta manera reemplazando los términos correspondientes a $j = 2$ se obtiene

$$d^2 \widetilde{M}_{n-1,2} - d^2 \widetilde{M}_{n,2} = a(\widetilde{M}_{3n+2,2} + \widetilde{M}_{3n+1,1} + \widetilde{M}_{3n,0}) - a^2 \delta_{3n+3,2},$$

$$0 = a(\widetilde{M}_{3n+2,2} + t_1 + d^3).$$

Por lo tanto $\boxed{\widetilde{M}_{3n+2,2} = -d^3 - t_1}$.

$\boxed{\text{Caso } 3 \leq j \leq n}$

En este caso en la sumatoria sólo sobreviven los términos correspondientes a $n - j \leq i \leq n$, pues si $i < n - j$ se tiene $\widetilde{M}_{2n+2+i,j-n+i} = 0$.

Luego

$$\sum_{i=0}^{n} \widetilde{M}_{2n+2+i,j-n+i} = \sum_{i=n-j}^{n} \widetilde{M}_{2n+2+i,j-n+i}$$

$$= \widetilde{M}_{3n+2-j,0} + \widetilde{M}_{3n+3-j,1} + \cdots + \widetilde{M}_{2n+2+i,k} + \cdots + \widetilde{M}_{3n+1,j-1} + \widetilde{M}_{3n+2,j}.$$

Ahora $0 \leq i \leq n-2 \Rightarrow 2n+2 \leq 2n+2+i \leq 3n$, entonces $\widetilde{M}_{2n+2+i,k} = 0$, para $2 \leq k = j - n + i \leq j - 2$.

Por lo tanto sobreviven los términos $\widetilde{M}_{3n+2-j,0}, \widetilde{M}_{3n+3-j,1}, \widetilde{M}_{3n+1,j-1}, \widetilde{M}_{3n+2,j}$.

De esta manera reemplazando los términos correspondientes en (4.1.1) se obtiene

$$d^2\widetilde{M}_{n-1,j} - d^2\widetilde{M}_{n,j} = a(\widetilde{M}_{3n+2-j,0} + \widetilde{M}_{3n+3-j,1} + \widetilde{M}_{3n+1,j-1} + \widetilde{M}_{3n+2,j}) - a^2\delta_{3n+3,j},$$

$$0 = a(d^2 - d^2 + t_{j-1} + \widetilde{M}_{3n+2,j}).$$

Por lo tanto $\boxed{\widetilde{M}_{3n+2,j} = -t_{j-1}}$.

$\boxed{\text{Caso } j = n+1}$

En este caso

$$\sum_{i=0}^{n} \widetilde{M}_{2n+2+i,i+1} = \widetilde{M}_{2n+2,1} + \widetilde{M}_{2n+3,2} + \cdots + \widetilde{M}_{2n+2+i,i+1} + \cdots + \widetilde{M}_{3n+1,n} + \widetilde{M}_{3n+2,n+1}.$$

Ahora $0 \le i \le n-2 \Rightarrow 2n+2 \le 2n+2+i \le 3n$, entonces $\widetilde{M}_{2n+2+i,i+1} = 0$ para $2 \le i+1 \le n-1$.

Luego

$$\sum_{i=0}^{n} \widetilde{M}_{2n+2+i,i+1} = \widetilde{M}_{2n+2,1} + \widetilde{M}_{3n+1,n} + \widetilde{M}_{3n+2,n+1}$$

y entonces (4.1.1) lleva a

$$-d^2 a = a(\widetilde{M}_{2n+2,1} + \widetilde{M}_{3n+1,n} + \widetilde{M}_{3n+2,n+1}) = a(-d^2 + t_n + \widetilde{M}_{3n+2,n+1}).$$

Por lo tanto $\boxed{\widetilde{M}_{3n+2,n+1} = -t_n}$.

$\boxed{\text{Caso } j = n+2}$

En este caso

$$\sum_{i=0}^{n} \widetilde{M}_{2n+2+i,i+2} = \widetilde{M}_{2n+2,2} + \widetilde{M}_{2n+3,3} + \cdots + \widetilde{M}_{2n+2+i,i+2} + \cdots + \widetilde{M}_{3n+1,n+1} + \widetilde{M}_{3n+2,n+2}.$$

Ahora $0 \le i \le n-2 \Rightarrow 2n+2 \le 2n+2+i \le 3n$, entonces $\widetilde{M}_{2n+2+i,i+2} = 0$ para $2 \le i+2 \le n$.

Luego

$$\sum_{i=0}^{n} \widetilde{M}_{2n+2+i,i+2} = \widetilde{M}_{3n+1,n+1} + \widetilde{M}_{3n+2,n+2}$$

y entonces (4.1.1) lleva a

$$0 = a(\widetilde{M}_{3n+1,n+1} + \widetilde{M}_{3n+2,n+2}) = a(t_{n+1} + \widetilde{M}_{3n+2,n+2}).$$

Por lo tanto $\boxed{\widetilde{M}_{3n+2,n+2} = -t_{n+1}}$.

$\boxed{\text{Caso } j = n+3}$

En este caso

$$\sum_{i=0}^{n} \widetilde{M}_{2n+2+i,i+3} = \widetilde{M}_{2n+2,3} + \widetilde{M}_{2n+3,4} + \cdots + \widetilde{M}_{2n+2+i,i+3} + \cdots + \widetilde{M}_{3n+1,n+2} + \widetilde{M}_{3n+2,n+3}.$$

Ahora $0 \leq i \leq n-2 \Rightarrow 2n+2 \leq 2n+2+i \leq 3n$, entonces $\widetilde{M}_{2n+2+i,i+3} = 0$, para $3 \leq i+3 \leq n$.

Luego
$$\sum_{i=0}^{n} \widetilde{M}_{2n+2+i,i+3} = \widetilde{M}_{3n,n+1} + \widetilde{M}_{3n+1,n+2} + \widetilde{M}_{3n+2,n+3}$$

y entonces (4.1.1) lleva a
$$0 = a(\widetilde{M}_{3n,n+1} + \widetilde{M}_{3n+1,n+2} + \widetilde{M}_{3n+2,n+3}) = a(ad^2 + t_{n+2} + \widetilde{M}_{3n+2,n+3}).$$

Por lo tanto $\boxed{\widetilde{M}_{3n+2,n+3} = -ad^2 - t_{n+2}}$.

$\boxed{\text{Caso } n+4 \leq j \leq 3n+2}$.

En este caso $4 \leq j-n+i \leq 3n+2$, luego

$$\sum_{i=0}^{n} \widetilde{M}_{2n+2+i,j-n+i} = \widetilde{M}_{2n+2,j-n} + \widetilde{M}_{2n+3,j-n+1} + \cdots + \widetilde{M}_{2n+2+i,j-n+i} + \cdots + \widetilde{M}_{3n+1,j-1} + \widetilde{M}_{3n+2,j}$$

En la sumatoria sólo sobreviven los dos últimos términos.

Luego
$$\sum_{i=0}^{n} \widetilde{M}_{2n+2+i,j-n+i} = \widetilde{M}_{3n+1,j-1} + \widetilde{M}_{3n+2,j} = t_{j-1} + \widetilde{M}_{3n+2,j}$$

y entonces (4.1.1) lleva a
$$0 = a(t_{j-1} + \widetilde{M}_{3n+2,j}).$$

Por lo tanto $\boxed{\widetilde{M}_{3n+2,j} = -t_{j-1}}$.

$\boxed{\text{Caso } j = 3n+3}$

En este caso

$$\sum_{i=0}^{n} \widetilde{M}_{2n+2+i,2n+3+i} = \widetilde{M}_{2n+2,2n+3} + \widetilde{M}_{2n+3,2n+4} + \cdots + \widetilde{M}_{2n+2+i,2n+3+i} \cdots + \widetilde{M}_{3n+1,3n+2} + \widetilde{M}_{3n+2,3n+3}.$$

En la sumatoria sólo sobreviven los dos últimos términos.

Luego
$$\sum_{i=0}^{n} \widetilde{M}_{2n+2+i,i+1} = \widetilde{M}_{3n+1,3n+2} + \widetilde{M}_{3n+2,3n+3} = t_{3n+2} + \widetilde{M}_{3n+2,3n+3}.$$

Además $\delta_{3n+3,j} = \delta_{3n+3,3n+3} = 1$. Entonces por (4.1.1) obtenemos

$$0 = a(t_{3n+2} + \widetilde{M}_{3n+2,3n+3}) - a^2$$

Por lo tanto $\boxed{\widetilde{M}_{3n+2,3n+3} = a - t_{3n+2}}$.

En el segundo caso de la Proposición 4.2.2 la matriz $\widetilde{M}$ hasta la fila $3n+2$ se ve como sigue:

$$\begin{array}{r} 0 \to \\ \\ \\ n \to \\ \\ \\ \\ \\ 2n \to \\ 2n+1 \to \\ \\ \\ \\ \\ 3n \to \\ 3n+1 \to \\ 3n+2 \to \end{array} \left( \begin{array}{cccccccccccc} 1 & -1 & 0 & \cdots & 0 & 0 & \cdots & 0 & 0 & \cdots & 0 & 0 \\ \vdots & \vdots & \vdots & \cdots & \vdots & \vdots & \cdots & \vdots & \vdots & \cdots & \vdots & \vdots \\ 1 & -1 & 0 & \cdots & 0 & 0 & \cdots & 0 & 0 & \cdots & 0 & 0 \\ d & -1 & \mathbf{0} & \cdots & a & \mathbf{0} & \cdots & \mathbf{0} & \mathbf{0} & \cdots & 0 & 0 \\ d & -d & 0 & \cdots & 0 & 0 & \cdots & 0 & 0 & \cdots & \mathbf{0} & \mathbf{0} \\ \vdots & \vdots & \vdots & \cdots & \vdots & \vdots & \cdots & \vdots & \vdots & \cdots & \vdots & \vdots \\ d & -d & 0 & \cdots & 0 & 0 & \cdots & 0 & 0 & \cdots & 0 & 0 \\ d^2 & -d & \mathbf{0} & \cdots & ad & \mathbf{0} & \cdots & \mathbf{0} & \mathbf{0} & \cdots & \mathbf{0} & \mathbf{0} \\ d^2 & -d^2 & \mathbf{0} & \cdots & \mathbf{0} & -ad & \cdots & \mathbf{0} & a & \cdots & \mathbf{0} & \mathbf{0} \\ d^2 & -d^2 & 0 & \cdots & 0 & 0 & \cdots & 0 & 0 & \cdots & 0 & 0 \\ \vdots & \vdots & \vdots & \cdots & \vdots & \vdots & \cdots & \vdots & \vdots & \cdots & \vdots & \vdots \\ d^2 & -d^2 & 0 & \cdots & 0 & 0 & \cdots & 0 & 0 & \cdots & 0 & 0 \\ d^3 & -d^2 & \mathbf{0} & \cdots & ad^2 & \mathbf{0} & \cdots & \mathbf{0} & \mathbf{0} & \cdots & \mathbf{0} & \mathbf{0} \\ t_0 & t_1 & t_2 & \cdots & t_{n+1} & t_{n+2} & \cdots & t_{2n} & t_{2n+1} & \cdots & t_{3n+2} & \mathbf{0} \\ d^3 & -t_0 & \mathbf{u} & \cdots & -t_n & -t_{n+1} & \mathbf{v} & \cdots & -t_{2n} & \cdots & -t_{3n+1} & a-t_{3n+2} \end{array} \right)$$

donde $\mathbf{u} = -d^3 - t_1$, $\mathbf{v} = -t_{n+2} - ad^2$.

En la matriz del segundo caso de la Proposición 4.2.4 vemos que $\widetilde{M}_{j,*} = \widetilde{M}_{j-1,*}$ para $n+2 \leq j \leq 2n-1$.

Por el Lema 3.1.2 se sigue que $\widetilde{M} Y^j \widetilde{M} = 0$ para $n+1 \leq j \leq 2n-2$.

Luego

$$0 = (\widetilde{M} Y^j \widetilde{M})_{n,i} = \widetilde{M}_{n,*} \cdot \widetilde{M}_{*+j,i} = (d, -1, 0, \cdots, 0, a) \cdot (\widetilde{M}_{j,i}, \widetilde{M}_{j+1,i}, \cdots, \widetilde{M}_{j+n+1,i})$$

Es decir

$$0 = d\widetilde{M}_{j,i} - \widetilde{M}_{j+1,i} + a\widetilde{M}_{j+n+1,i}.$$

Caso $\boxed{i=0}$

Se tiene $\widetilde{M}_{j,0} = d$, $\widetilde{M}_{j+1,0} = d$ y se obtiene $\widetilde{M}_{j+n+1,0} = \dfrac{d-d^2}{a} = d^2$.

Caso $\boxed{i=1}$

Se tiene $\widetilde{M}_{j,1} = -d$, $\widetilde{M}_{j+1,1} = -d$ y se obtiene $\widetilde{M}_{j+n+1,1} = \dfrac{d^2-d}{a} = -d^2$.

Caso $\boxed{i \geq 2}$

En este caso como $n+1 \leq j \leq 2n-2$ se cumple $\widetilde{M}_{j,i} = \widetilde{M}_{j+1,i} = 0$.

Luego se tiene que $0 = a\widetilde{M}_{j+n+1,i}$ y como $a \neq 0$ entonces $\widetilde{M}_{j+n+1,i} = 0$.

Por lo tanto se cumple que

$$\widetilde{M}_{k,*} = (d^2, -d^2, 0, 0, \cdots)$$

para $2n+2 \leq k = j+n+1 \leq 3n-1$.

De esta manera la matriz $\widetilde{M}$ hasta la fila $3n-1$ se ve como sigue:

$$\begin{array}{c} 0 \to \\ \\ n \to \\ \\ \\ \\ 2n \to \\ 2n+1 \to \\ \\ \\ 3n-1 \to \end{array} \left( \begin{array}{ccccccccc} 1 & -1 & 0 & \cdots & 0 & 0 & \cdots & 0 & 0 \\ \vdots & \vdots & \vdots & \vdots & \vdots & \vdots & \vdots & \vdots & \vdots \\ 1 & -1 & 0 & \cdots & 0 & 0 & \cdots & 0 & 0 \\ d & -1 & 0 & \cdots & 0 & a & \cdots & 0 & 0 \\ d & -d & 0 & \cdots & 0 & 0 & \cdots & 0 & 0 \\ \vdots & \vdots & \vdots & \cdots & \vdots & \vdots & \vdots & \vdots & \vdots \\ d & -d & 0 & \cdots & 0 & 0 & \cdots & 0 & 0 \\ d & -d & 0 & \cdots & -ad & 0 & \cdots & a & 0 \\ d^2 & -d & 0 & \cdots & 0 & ad & \cdots & 0 & 0 \\ d^2 & -d^2 & 0 & \cdots & 0 & 0 & \cdots & 0 & 0 \\ \vdots & \vdots & \vdots & \vdots & \vdots & \vdots & \vdots & \vdots \\ d^2 & -d^2 & 0 & \cdots & 0 & 0 & \cdots & 0 & 0 \end{array} \right)$$

Ahora pongamos $t_i = \widetilde{M}_{3n,i}$. Es decir la fila $3n$ de $\widetilde{M}$ es

$$(t_0, t_1, t_2, \cdots, t_{3n}, t_{3n+1}, 0, \cdots).$$

Evaluando cada término de la ecuación (4.1.1) en la entrada $(2n+1, j)$ obtenemos:

$$(\widetilde{M} Y^{n-1} \widetilde{M})_{2n+1,j} = \widetilde{M}_{2n+1,*} \cdot \widetilde{M}_{*+n-1,j} = (d^2, -d, 0, \cdots, 0, ad) \cdot (\widetilde{M}_{n-1,j}, \widetilde{M}_{n,j}, \cdots, \widetilde{M}_{2n,j})$$

$$= d^2 \widetilde{M}_{n-1,j} - d \widetilde{M}_{n,j} + ad \widetilde{M}_{2n,j},$$

$$(Y^i \widetilde{M} Y^{n-i})_{2n+1,j} = (\widetilde{M} Y^{n-i})_{2n+1+i,j} = \widetilde{M}_{2n+1+i, j-n+i},$$

$$(Y^{n+1})_{2n+1,j} = \delta_{3n+2,j}.$$

Para $j = 0$ en la sumatoria sólo sobrevive el término correspondiente a $i = n$ pues si $i < n$ se tiene $\widetilde{M}_{2n+1+i, -n+i} = 0$.

De esta manera reemplazando los términos correspondientes a $j = 0$ se obtiene

$$d^2 \widetilde{M}_{n-1,0} - d \widetilde{M}_{n,0} + ad \widetilde{M}_{2n,0} = a \widetilde{M}_{3n+1,0} - a^2 \delta_{3n+2,0},$$

$$d^2 - d^2 + ad^2 = a \widetilde{M}_{3n+1,0}.$$

Es decir

$$ad^2 = a \widetilde{M}_{3n+1,0}.$$

Por lo tanto $\boxed{\widetilde{M}_{3n+1,0} = d^2}$.

Para $j = 1$ en la sumatoria sólo sobreviven los términos correspondientes a $i = n$,

$i = n-1$ pues si $i < n-1$ se tiene $\widetilde{M}_{2n+1+i,1-n+i} = 0$.

De esta manera reemplazando los términos correspondientes a $j = 1$ se obtiene

$$d^2\widetilde{M}_{n-1,1} - d\widetilde{M}_{n,1} + ad\widetilde{M}_{2n,1} = a(\widetilde{M}_{3n+1,1} + \widetilde{M}_{3n,0}) - a^2\delta_{3n+2,1},$$

$$-d^2 + d - ad^2 = a(\widetilde{M}_{3n+1,1} + t_0),$$

$$(1-d)d - ad^2 = a(\widetilde{M}_{3n+1,1} + t_0),$$

$$ad^2 - ad^2 = a(\widetilde{M}_{3n+1,1} + t_0.)$$

Es decir

$$0 = a(\widetilde{M}_{3n+1,1} + t_0).$$

Por lo tanto $\boxed{\widetilde{M}_{3n+1,1} = -t_0}$.

Para $j = 2$ en la sumatoria sólo sobreviven los términos correspondientes a $i = n$, $i = n-1$, $i = n-2$ pues si $i < n-2$ se tiene $\widetilde{M}_{2n+1+i,2-n+i} = 0$.

De esta manera reemplazando los términos correspondientes a $j = 2$ se obtiene

$$d^2\widetilde{M}_{n-1,2} - d\widetilde{M}_{n,2} + ad\widetilde{M}_{2n,2} = a(\widetilde{M}_{3n+1,2} + \widetilde{M}_{3n,1} + \widetilde{M}_{3n-1,0}) - a^2\delta_{3n+2,2},$$

$$0 = a(\widetilde{M}_{3n+1,2} + t_1 + d^2).$$

Por lo tanto $\boxed{\widetilde{M}_{3n+1,2} = -t_1 - d^2}$.

$\boxed{\text{Caso } 3 \leq j \leq n}$.

En este caso en la sumatoria sólo sobreviven los términos correspondientes a $n-j \leq i \leq n$, pues si $i < n-j$ se tiene $\widetilde{M}_{2n+1+i,j-n+i} = 0$.

Luego

$$\sum_{i=0}^{n}\widetilde{M}_{2n+1+i,j-n+i} = \sum_{i=n-j}^{n}\widetilde{M}_{2n+1+i,j-n+i}$$

$$= \widetilde{M}_{3n+1-j,0} + \widetilde{M}_{3n+2-j,1} + \cdots + \widetilde{M}_{2n+1+i,k} + \cdots + \widetilde{M}_{3n,j-1} + \widetilde{M}_{3n+1,j}$$

Pero $\widetilde{M}_{2n+1+i,k} = 0$ para $2 \leq k = j-n+i \leq j-2$, pues en este caso $2n+3 \leq 2n+1+i \leq 3n-1$.

Por lo tanto sobreviven los términos $\widetilde{M}_{3n+1-j,0}, \widetilde{M}_{3n+2-j,1}, \widetilde{M}_{3n,j-1}, \widetilde{M}_{3n+1,j}$.

De esta manera reemplazando los términos correspondientes en (4.1.1) se obtiene

$$d^2\widetilde{M}_{n-1,j} - d\widetilde{M}_{n,j} + ad\widetilde{M}_{2n,j} = a(\widetilde{M}_{3n+1-j,0} + \widetilde{M}_{3n+2-j,1} + \widetilde{M}_{3n,j-1} + \widetilde{M}_{3n+1,j}) - a^2\delta_{3n+2,j},$$

Si $3 \leq j \leq n-1$ se obtiene:

$$0 = a(d^2 - d^2 + t_{j-1} + \widetilde{M}_{3n+1,j}).$$

Por lo tanto $\boxed{\widetilde{M}_{3n+1,j} = -t_{j-1}}$.
Si $j = n$ se obtiene:
$$-a^2 d^2 = a(d^2 - d^2 + t_{n-1} + \widetilde{M}_{3n+1,n}).$$
Por lo tanto $\boxed{\widetilde{M}_{3n+1,n} = -ad^2 - t_{n-1}}$.

$\boxed{\text{Caso } j = n+1}$

En este caso
$$\sum_{i=0}^{n} \widetilde{M}_{2n+1+i,i+1} = \widetilde{M}_{2n+1,1} + \widetilde{M}_{2n+2,2} + \cdots + \widetilde{M}_{2n+1+i,i+1} + \cdots + \widetilde{M}_{3n,n} + \widetilde{M}_{3n+1,n+1}.$$

Si $2 \leq i+1 \leq n-1$, entonces $\widetilde{M}_{2n+1+i,i+1} = 0$, pues $2n+2 \leq 2n+1+i \leq 3n-1$.
Luego
$$\sum_{i=0}^{n} \widetilde{M}_{2n+1+i,i+1} = \widetilde{M}_{2n+1,1} + \widetilde{M}_{3n,n} + \widetilde{M}_{3n+1,n+1}.$$

y entonces (4.1.1) lleva a
$$-da = a(-d + t_n + \widetilde{M}_{3n+1,n+1}).$$

Por lo tanto $\boxed{\widetilde{M}_{3n+1,n+1} = -t_n}$.

$\boxed{\text{Caso } n+2 \leq j \leq 2n}$

En este caso $2 \leq j-n+i \leq 2n$. Luego
$$\sum_{i=0}^{n} \widetilde{M}_{2n+1+i,j-n+i} = \widetilde{M}_{2n+1,j-n} + \widetilde{M}_{2n+2,j-n+1} + \cdots + \widetilde{M}_{2n+1+i,j-n+i} + \cdots + \widetilde{M}_{3n,j-1} + \widetilde{M}_{3n+1,j}$$

Consideremos $2 \leq j-n+i \leq j-2$.
Si $1 \leq i \leq n-2$ entonces $2n+2 \leq 2n+1+i \leq 3n-1$.
Si $i = 0$ entonces $2n+1+i = 2n+1$ y $2 \leq j-n+i = j-n \leq n$.
En ambos casos se tiene $\widetilde{M}_{2n+1+i,j-n+i} = 0$
Por lo tanto
$$\sum_{i=0}^{n} \widetilde{M}_{2n+1+i,j-n+i} = \widetilde{M}_{3n,j-1} + \widetilde{M}_{3n+1,j} = t_{j-1} + \widetilde{M}_{3n+1,j}$$

y entonces (4.1.1) lleva a
$$0 = a(t_{j-1} + \widetilde{M}_{3n+1,j}).$$

Por lo tanto $\boxed{\widetilde{M}_{3n+1,j} = -t_{j-1}}$.

$\boxed{\text{Caso } j = 2n+1}$

En este caso

$$\sum_{i=0}^{n}\widetilde{M}_{2n+1+i,i+n+1} = \widetilde{M}_{2n+1,n+1}+\widetilde{M}_{2n+2,n+2}+\cdots+\widetilde{M}_{2n+1+i,i+n+1}+\cdots+\widetilde{M}_{3n,2n}+\widetilde{M}_{3n+1,2n+1}.$$

Si $n+2 \leq i+n+1 \leq 2n-1$, entonces $\widetilde{M}_{2n+1+i,i+n+1} = 0$, pues $2n+2 \leq 2n+1+i \leq 3n-1$.
Luego

$$\sum_{i=0}^{n}\widetilde{M}_{2n+1+i,i+n+1} = \widetilde{M}_{2n+1,n+1} + \widetilde{M}_{3n,2n} + \widetilde{M}_{3n+1,2n+1}.$$

y entonces (4.1.1) lleva a

$$a^2 d = a(ad + t_{2n} + \widetilde{M}_{3n+1,2n+1}) = a^2 d + a(t_{2n} + \widetilde{M}_{3n+1,2n+1}).$$

Por lo tanto $\boxed{\widetilde{M}_{3n+1,2n+1} = -t_{2n}}$.

$\boxed{\text{Caso } 2n+2 \leq j \leq 3n+1}$

En este caso $n+2 \leq j-n+i \leq 3n+1$, luego

$$\sum_{i=0}^{n}\widetilde{M}_{2n+1+i,j-n+i} = \widetilde{M}_{2n+1,j-n}+\widetilde{M}_{2n+2j,j-n+1}+\cdots+\widetilde{M}_{2n+1+i,j-n+i}+\cdots+\widetilde{M}_{3n,j-1}+\widetilde{M}_{3n+1,j}$$

En la sumatoria sólo sobreviven los dos últimos términos, pues $n+2 \leq j-n+i \leq j-2 \Rightarrow 2n+1 \leq 2n+1+i \leq 3n-1$.
Luego

$$\sum_{i=0}^{n}\widetilde{M}_{2n+1+i,j-n+i} = \widetilde{M}_{3n,j-1} + \widetilde{M}_{3n+1,j} = t_{j-1} + \widetilde{M}_{3n+1,j}$$

y entonces (4.1.1) lleva a

$$0 = a(t_{j-1} + \widetilde{M}_{3n+1,j}).$$

Por lo tanto $\boxed{\widetilde{M}_{3n+1,j} = -t_{j-1}}$.

$\boxed{\text{Caso } j = 3n+2}$

En este caso

$$\sum_{i=0}^{n}\widetilde{M}_{2n+1+i,i+2n+2} = \widetilde{M}_{2n+1,2n+2}+\widetilde{M}_{2n+2,2n+3}+\cdots+\widetilde{M}_{2n+1+i,i+2n+2}\cdots+\widetilde{M}_{3n,3n+1}+\widetilde{M}_{3n+1,3n+2}.$$

En la sumatoria sobreviven los dos últimos términos.

Luego
$$\sum_{i=0}^{n} \widetilde{M}_{2n+1+i, i+2n+2} = \widetilde{M}_{3n,3n+1} + \widetilde{M}_{3n+1,3n+2} = t_{3n+1} + \widetilde{M}_{3n+1,3n+2}.$$

Además $\delta_{3n+2,3n+2} = 1$.

Entonces (4.1.1) lleva a
$$0 = a(t_{3n+1} + \widetilde{M}_{3n+1,3n+2}) - a^2.$$

Por lo tanto $\boxed{\widetilde{M}_{3n+1,3n+2} = a - t_{3n+1}}$.

Evaluando cada término de la ecuación (4.1.1) en la entrada $(2n+2, j)$ obtenemos:

$$(\widetilde{M} Y^{n-1} \widetilde{M})_{2n+2,j} = \widetilde{M}_{2n+2,*} \cdot \widetilde{M}_{*+n-1,j} = (d^2, -d^2) \cdot (\widetilde{M}_{n-1,j}, \widetilde{M}_{n,j})$$
$$= d^2 \widetilde{M}_{n-1,j} - d^2 \widetilde{M}_{n,j},$$
$$(Y^i \widetilde{M} Y^{n-i})_{2n+2,j} = (\widetilde{M} Y^{n-i})_{2n+2+i,j} = \widetilde{M}_{2n+2+i, j-n+i},$$
$$(Y^{n+1})_{2n+2,j} = \delta_{3n+3, j}.$$

Para $j = 0$ en la sumatoria sólo sobrevive el término correspondiente a $i = n$ pues si $i < n$ se tiene $\widetilde{M}_{2n+2+i, -n+i} = 0$.

De esta manera reemplazando los términos correspondientes a $j = 0$ se obtiene
$$d^2 \widetilde{M}_{n-1,0} - d^2 \widetilde{M}_{n,0} = a \widetilde{M}_{3n+2,0} - a^2 \delta_{3n+3,0},$$
$$d^2 - d^3 = a \widetilde{M}_{3n+2,0},$$
$$d^2(1 - d) = a \widetilde{M}_{3n+2,0}$$

Es decir
$$a d^3 = a \widetilde{M}_{3n+2,0}.$$

Por lo tanto $\boxed{\widetilde{M}_{3n+2,0} = d^3}$.

Para $j = 1$ en la sumatoria sólo sobreviven los términos correspondientes a $i = n$, $i = n-1$ pues si $i < n-1$ se tiene $\widetilde{M}_{2n+2+i, 1-n+i} = 0$.

De esta manera reemplazando los términos correspondientes a $j = 1$ se obtiene
$$d^2 \widetilde{M}_{n-1,1} - d^2 \widetilde{M}_{n,1} = a(\widetilde{M}_{3n+2,1} + \widetilde{M}_{3n+1,0}) - a^2 \delta_{3n+3,1},$$
$$-d^2 + d^2 = a(d^2 + \widetilde{M}_{3n+2,1}),$$
$$0 = a(\widetilde{M}_{3n+2,1} + d^2),$$

Por lo tanto $\boxed{\widetilde{M}_{3n+2,1} = -d^2}$.

Para $j = 2$ en la sumatoria sólo sobreviven los términos correspondientes a $i = n$,

$i = n-1$, $i = n-2$ pues si $i < n-2$ se tiene $\widetilde{M}_{2n+2+i,2-n+i} = 0$.

De esta manera reemplazando los términos correspondientes a $j = 2$ se obtiene

$$d^2\widetilde{M}_{n-1,2} - d\widetilde{M}_{n,2} = a(\widetilde{M}_{3n+2,2} + \widetilde{M}_{3n+1,1} + \widetilde{M}_{3n,0}) - a^2\delta_{3n+3,2},$$

$$0 = a(\widetilde{M}_{3n+2,2} - t_0 + t_0).$$

Por lo tanto $\boxed{\widetilde{M}_{3n+2,2} = 0}$.

Para $j = 3$ en la sumatoria sólo sobreviven los términos correspondientes a $i = n$, $i = n-1$, $i = n-2$, $i = n-3$ pues si $i < n-3$ se tiene $\widetilde{M}_{2n+2+i,3-n+i} = 0$.

De esta manera reemplazando los términos correspondientes a $j = 3$ se obtiene

$$d^2\widetilde{M}_{n-1,3} - d\widetilde{M}_{n,3} = a(\widetilde{M}_{3n+2,3} + \widetilde{M}_{3n+1,2} + \widetilde{M}_{3n,1} + \widetilde{M}_{3n-1,0}) - a^2\delta_{3n+3,3},$$

$$0 = a(\widetilde{M}_{3n+2,3} - d^2 - t_1 + t_1 + d^2).$$

Por lo tanto $\boxed{\widetilde{M}_{3n+2,3} = 0}$.

$\boxed{\text{Caso } 4 \leq j \leq n}$.

En este caso en la sumatoria sólo sobreviven los términos correspondientes a $n - j \leq i \leq n$, pues si $i < n - j$ se tiene $\widetilde{M}_{2n+2+i,j-n+i} = 0$.

Luego

$$\sum_{i=0}^{n} \widetilde{M}_{2n+2+i,j-n+i} = \sum_{i=n-j}^{n} \widetilde{M}_{2n+2+i,j-n+i}$$

$$= \widetilde{M}_{3n+2-j,0} + \widetilde{M}_{3n+3-j,1} + \cdots + \widetilde{M}_{2n+2+i,k} + \cdots + \widetilde{M}_{3n,j-2} + \widetilde{M}_{3n+1,j-1} + \widetilde{M}_{3n+2,j}$$

Si $2 \leq k = j - n + i \leq j - 3$ entonces $\widetilde{M}_{2n+2+i,k} = 0$, pues en este caso $2n + 4 \leq 2n + 2 + i \leq 3n - 1$

Por lo tanto sobreviven los términos $\widetilde{M}_{3n+2-j,0}, \widetilde{M}_{3n+3-j,1}, \widetilde{M}_{3n,j-2}, \widetilde{M}_{3n+1,j-1}, \widetilde{M}_{3n+2,j}$.

De esta manera reemplazando los términos correspondientes en (4.1.1) se obtiene

$$d^2\widetilde{M}_{n-1,j} - d\widetilde{M}_{n,j} = a(\widetilde{M}_{3n+2-j,0} + \widetilde{M}_{3n+3-j,1} + \widetilde{M}_{3n,j-2} + \widetilde{M}_{3n+1,j-1} + \widetilde{M}_{3n+2,j}) - a^2\delta_{3n+3,j},$$

$$0 = a(d^2 - d^2 + t_{j-2} - t_{j-2} + \widetilde{M}_{3n+2,j}).$$

Por lo tanto $\boxed{\widetilde{M}_{3n+2,j} = 0}$.

$\boxed{\text{Caso } j = n+1}$

En este caso

$$\sum_{i=0}^{n} \widetilde{M}_{2n+2+i,i+1} = \widetilde{M}_{2n+2,1} + \widetilde{M}_{2n+3,2} + \cdots + \widetilde{M}_{2n+2+i,i+1} + \cdots + \widetilde{M}_{3n+1,n} + \widetilde{M}_{3n+2,n+1}.$$

Ahora $0 \leq i \leq n-3 \Rightarrow 2n+2 \leq 2n+2+i \leq 3n-1$, entonces $\widetilde{M}_{2n+2+i,i+1} = 0$, para $2 \leq i+1 \leq n-2$

Luego
$$\sum_{i=0}^{n} \widetilde{M}_{2n+2+i,i+1} = \widetilde{M}_{2n+2,1} + \widetilde{M}_{3n,n-1} + \widetilde{M}_{3n+1,n} + \widetilde{M}_{3n+2,n+1}.$$

y entonces (4.1.1) lleva a
$$-ad^2 = a(-d^2 + t_{n-1} - t_{n-1} - ad^2 + \widetilde{M}_{3n+2,n+1}).$$

Por lo tanto $\boxed{\widetilde{M}_{3n+2,n+1} = ad^2}$.

$\boxed{\text{Caso } n+2 \leq j \leq 3n+2}$

En este caso $2 \leq j-n+i \leq 3n+2$, luego
$$\sum_{i=0}^{n} \widetilde{M}_{2n+2+i,j-n+i} = \widetilde{M}_{2n+2,j-n} + \widetilde{M}_{2n+3,j-n+1} + \cdots + \widetilde{M}_{2n+2+i,j-n+i} + \cdots + \widetilde{M}_{3n,j-2} + \widetilde{M}_{3n+1,j-1} + \widetilde{M}_{3n+2,j}$$

Si $2 \leq j-n+i \leq j-3$ entonces $\widetilde{M}_{2n+2+i,j-n+i} = 0$ pues, en este caso, $2n+2 \leq 2n+2+i \leq 3n-1$

Luego
$$\sum_{i=0}^{n} \widetilde{M}_{2n+2+i,j-n+i} = \widetilde{M}_{3n,j-2} + \widetilde{M}_{3n+1,j-1} + \widetilde{M}_{3n+2,j} = t_{j-2} - t_{j-2} + \widetilde{M}_{3n+2,j}$$

y entonces (4.1.1) lleva a
$$0 = a\widetilde{M}_{3n+2,j}.$$

Por lo tanto $\boxed{\widetilde{M}_{3n+2,j} = 0.}$

$\boxed{\text{Caso } j = 3n+3}$

En este caso
$$\sum_{i=0}^{n} \widetilde{M}_{2n+2+i,i+2n+3} = \widetilde{M}_{2n+2,2n+3} + \widetilde{M}_{2n+3,2n+4} + \cdots + \widetilde{M}_{2n+2+i,i+2n+3} + \cdots + \widetilde{M}_{3n+1,3n+2} + \widetilde{M}_{3n+2,3n+3}.$$

Ahora $0 \leq i \leq n-3 \Rightarrow 2n+2 \leq 2n+2+i \leq 3n-1$, entonces $\widetilde{M}_{2n+2+i,i+2n+3} = 0$, para $2n+3 \leq i+2n+3 \leq 3n$.

Luego
$$\sum_{i=0}^{n} \widetilde{M}_{2n+2+i,i+2n+3} = \widetilde{M}_{3n,3n+1} + \widetilde{M}_{3n+1,3n+2} + \widetilde{M}_{3n+2,3n+3}.$$

y entonces (4.1.1) lleva a

$$0 = a(t_{3n+1} + a - t_{3n+1} + \widetilde{M}_{3n+2,3n+3}) - a^2 \delta_{3n+3,3n+3} = a\widetilde{M}_{3n+2,3n+3}.$$

Por lo tanto $\boxed{\widetilde{M}_{3n+2,3n+3} = 0}$.

En el segundo caso de la Proposición 4.2.4 la matriz $\widetilde{M}$ hasta la fila $3n+2$ se ve como sigue:

$$\begin{array}{c}0\rightarrow\\ \\ n\rightarrow\\ \\ \\ \\ \\ 2n\rightarrow\\ 2n+1\rightarrow\\ \\ \\ \\ 3n\rightarrow\\ 3n+1\rightarrow\\ 3n+2\rightarrow\end{array}\left(\begin{array}{ccccccccccccc} 1 & -1 & 0 & \cdots & 0 & 0 & \cdots & 0 & 0 & \cdots & 0 & 0 \\ \vdots & \vdots & \vdots & \vdots & \vdots & \vdots & \vdots & \vdots & \vdots & \cdots & \vdots & \vdots \\ 1 & -1 & 0 & \cdots & 0 & 0 & \cdots & 0 & 0 & \cdots & 0 & 0 \\ \boldsymbol{d} & \boldsymbol{-1} & \boldsymbol{0} & \cdots & \boldsymbol{0} & \boldsymbol{a} & \cdots & \boldsymbol{0} & \boldsymbol{0} & \cdots & 0 & 0 \\ d & -d & 0 & \cdots & 0 & 0 & \cdots & 0 & 0 & \cdots & 0 & 0 \\ \vdots & \vdots & \vdots & \cdots & \vdots & \vdots & \vdots & \vdots & \vdots & \cdots & 0 & 0 \\ d & -d & 0 & \cdots & 0 & 0 & \cdots & 0 & 0 & \cdots & 0 & 0 \\ \boldsymbol{d} & \boldsymbol{-d} & \boldsymbol{0} & \cdots & \boldsymbol{-ad} & \boldsymbol{0} & \cdots & \boldsymbol{a} & \boldsymbol{0} & \cdots & \boldsymbol{0} & \boldsymbol{0} \\ \boldsymbol{d^2} & \boldsymbol{-d} & \boldsymbol{0} & \cdots & \boldsymbol{0} & \boldsymbol{ad} & \cdots & \boldsymbol{0} & \boldsymbol{0} & \cdots & \boldsymbol{0} & \boldsymbol{0} \\ d^2 & -d^2 & 0 & \cdots & 0 & 0 & \cdots & 0 & 0 & \cdots & 0 & 0 \\ \vdots & \vdots & \vdots & \vdots & \vdots & \vdots & \vdots & \vdots & \vdots & \cdots & \vdots & \vdots \\ d^2 & -d^2 & 0 & \cdots & 0 & 0 & \cdots & 0 & 0 & \cdots & 0 & 0 \\ \boldsymbol{t_0} & \boldsymbol{t_1} & \boldsymbol{t_2} & \cdots & \boldsymbol{t_n} & \boldsymbol{t_{n+1}} & \cdots & \boldsymbol{t_{2n+1}} & \boldsymbol{t_{2n+2}} & \cdots & \boldsymbol{t_{3n+1}} & 0 \\ \boldsymbol{d^2} & \boldsymbol{-t_0} & \boldsymbol{u} & \cdots & \boldsymbol{v} & \boldsymbol{-t_n} & \cdots & \boldsymbol{-t_{2n}} & \boldsymbol{-t_{2n+1}} & \cdots & \boldsymbol{-t_{3n}} & \boldsymbol{a-t_{3n+1}} \\ \boldsymbol{d^3} & \boldsymbol{-d^2} & \boldsymbol{0} & \cdots & \boldsymbol{0} & \boldsymbol{ad^2} & \cdots & \boldsymbol{0} & \boldsymbol{0} & \cdots & \boldsymbol{0} & \boldsymbol{0} \end{array}\right)$$

donde $\mathbf{u} = -\boldsymbol{d^2} - \boldsymbol{t_1}$, $\mathbf{v} = -\boldsymbol{ad^2} - \boldsymbol{t_{n-1}}$.

Cada una de las cuatro matrices a su vez nos da de nuevo cuatro posibilidades y, resolviendo un sistema similar como antes para los $t'_i$s, se determinan completamente las filas $3n$, $3n+1$ and $3n+2$. Así tenemos 16 casos.

Sin embargo, el sistema de ecuaciones es más complicado, especialmente en las dos últimas matrices. La Proposición 4.1.3 puede ser adaptada a estas dos últimas matrices, pero no el Lema 4.1.4. Este lema expresa las potencias de $M$ en términos de $Y^k \widetilde{M} Y^j$, lo cual es la clave para obtener el sistema de ecuaciones. Sin embargo este lema se aplica solamente a los dos primeros casos.

De manera que nuestros métodos no proveen una clasificación completa. En particular, existen restricciones que impiden que todos los casos lleguen a ser una aplicación de torcimiento. Algunas de estas restricciones serán vistas en la siguiente sección donde describiremos una familia de aplicaciones de torcimiento para la cual se aplica una generalización del Lema 4.1.4 .

## 4.3. La familia $B(n, L)$ y sucesiones casi-balanceadas

En esta sección describiremos una familia de productos tensoriales torcidos que surgen en el caso (2) de la Proposición 3.1.6. Así $\sigma$ es una aplicación de torcimiento, $Y$ y $M$ son como en la Proposición 1.4.1 y $\widetilde{M} = M - Y$. Además existe $n \in \mathbb{N}$, con $n \geq 2$, tal que $\widetilde{M}_{k,*} = \widetilde{M}_{0,*}$, para $1 < k < n$ y $\widetilde{M}_{n,*} = dE_0 - E_1 + aE_{n+1}$ para algunos $a, d \in K^\times$ con $d(a+1) = 1$. C Por la Proposición 4.1.2 (2),(3) la aplicación de torcimiento define una sucesión $L = (L_1, L_2, \cdots)$ tal que $\widetilde{M}_{L_i, L_i+1} = a$ y tal que $\widetilde{M}_{k,k+1} = 0$ si $L_i < k < L_{i+1}$ para algún $i$. Además $L_{j+1} - L_j \in \{n, n+1\}$ para todo $j \geq 1$, de manera que $L$ pertenece al conjunto de sucesiones

$$\Delta(n, n+1) = \{L \in \mathbb{N}^\mathbb{N} : L_1 = n \text{ y } L_{j+1} - L_j \in \{n, n+1\} \text{ para todo } j \geq 1\}.$$

Motivados por los casos de las Proposiciones 4.2.2 (1) y 4.2.4 (1) y como una continuación de estos dos casos, para una sucesión dada $L \in \Delta(n, n+1)$, $a \neq 0, -1$ y haciendo $d := \dfrac{1}{1+a}$, definimos la matriz infinita $M(L, a)$ mediante:

$$M(L,a)_{j,*} := \begin{cases} E_0 - E_1, & \text{si } j < L_1, \\ d^k(E_0 - E_1), & \text{si } L_k < j < L_{k+1} \text{ para algún } k, \\ d^k E_0 - d^{k-1} E_1 + aE_{L_k+1}, & \text{si } j = L_k \text{ para algún } k. \end{cases}$$

Probaremos en la Proposición 4.3.5 que si la matriz infinita $\widetilde{M} := M(L, a)$ define una aplicación de torcimiento via la Observación 3.1.1, entonces $L$ es una sucesión casi-balanceada, donde el conjunto de sucesiones casi-balanceadas $\mathscr{L}$ está definido como sigue:

$$\mathscr{L} := \{L \in \Delta(n, n+1) : L_r - 1 \leq L_j + L_{r-j} \leq L_r \text{ para } j, r \in N, \text{tal que } 0 < j < r\}.$$

Primero describiremos algunas propiedades de las sucesiones casi-balanceadas.

**Definición 4.3.1.** *Decimos que una sucesión $L \in \Delta(n, n+1)$ es $r$-balanceada, si $L_r = L_j + L_{r-j}$ para todo $0 < j < r$. Decimos que la sucesión es $r$-casi-balanceada, si $L_r - 1 \leq L_j + L_{r-j} \leq L_r$ para todo $0 < j < r$. También decimos que una sucesión finita $(L_1, \ldots, L_{r_0})$ es casi-balanceada, si $L_r - 1 \leq L_j + L_{r-j} \leq L_r$ para todo $0 < j < r \leq r_0$.*

Notemos que toda sucesión $L \in \Delta(n, n+1)$ es 1-balanceada (en este caso $r = 1$ y no existe $0 < j < r$), así mismo es 2-casi-balanceada.
En efecto en este caso $r = 2$, $j = 1$ y $r - j = 1$.
Luego $L_r - 1 \in \{2n - 1, 2n\}$, $L_j = L_1 = n, L_{r-j} = L_1 = n$ y $L_r = L_2 \in \{2n, 2n+1\}$.
Si $L_r = L_2 = 2n$ se tiene

$$L_r - 1 = 2n - 1 \leq L_1 + L_1 = 2n \leq 2n = L_r.$$

Si $L_r = L_2 = 2n+1$ se tiene

$$L_r - 1 = 2n \leq L_1 + L_1 = 2n \leq 2n+1 = L_r.$$

Cualquier sucesión $L \in \Delta(n, n+1)$ también es 3-casi-balanceada.
En efecto en este caso $r = 3$, $j = 1, j = 2$.
Para $j = 1$ se tiene $r - j = 2$.
Debemos demostrar que $L_3 - 1 \leq L_1 + L_2 \leq L_3$.
En efecto:

$$L_1 = n, L_2 = 2n, L_3 = 3n \Rightarrow 3n - 1 \leq L_1 + L_2 = 3n \leq 3n,$$
$$L_1 = n, L_2 = 2n, L_3 = 3n+1 \Rightarrow 3n \leq L_1 + L_2 = 3n \leq 3n+1,$$
$$L_1 = n, L_2 = 2n+1, L_3 = 3n+1 \Rightarrow 3n \leq L_1 + L_2 = 3n+1 \leq 3n+1,$$
$$L_1 = n, L_2 = 2n+1, L_3 = 3n+2 \Rightarrow 3n+1 \leq L_1 + L_2 = 3n+1 \leq 3n+2.$$

Para $j = 2$ se tiene $r - j = 1$, y es el mismo caso anterior.
Sin embargo no todas la sucesiones son 4-casi-balanceadas. Por ejemplo las sucesiones que empiezan con $(L_1 = n, L_2 = 2n, L_3 = 3n+1, L_4 = 4n+2, \ldots)$ o con $(L_1 = n, L_2 = 2n+1, L_3 = 3n+1, L_4 = 4n+1, \ldots)$ no son 4-casi-balanceadas.
En efecto estas sucesiones deberían cumplir $L_4 - 1 \leq L_2 + L_2 \leq L_4$.
Sin embargo para la primera sucesión se tiene $L_2 + L_2 = 4n < L_4 - 1$ y para la segunda sucesión se tiene $L_2 + L_2 = 4n+2 > L_4$.
Evidentemente una sucesión $L \in \Delta(n, n+1)$ es casi-balanceada (es decir, pertenece a $\mathscr{L}$), si y sólo si es $r$-casi-balanceada para todo $r$, si y sólo si todas las sucesiones $L_{\leq r} := (L_1, \ldots, L_r)$ son casi-balanceadas.
Para una sucesión dada $L$ y $0 < j < r$ definimos $\Delta_{r,j} := L_r - L_j - L_{r-j}$.
El hecho de que $(L_1, \ldots, L_{r_0})$ es casi-balanceada es equivalente al hecho de que $\Delta_{r,j} \in \{0, 1\}$ para todo $0 < j < r \leq r_0$.

**Proposición 4.3.2.** *Supongamos que la sucesión finita $(L_1, \ldots, L_{r_0})$ es casi-balanceada. Entonces $(L_1, \ldots, L_{r_0}, L_{r_0} + n)$ es casi-balanceada o $(L_1, \ldots, L_{r_0}, L_{r_0} + n + 1)$ es casi-balanceada.*

**Demostración.** Sabemos que $L_{r_0+1} \in \{L_{r_0} + n, L_{r_0} + n + 1\}$. Debemos demostrar que $\Delta_{r_0+1, j} = L_{r_0+1} - L_j - L_{r_0+1-j} \in \{0, 1\}$, para todo $0 < j < r_0 + 1$, en alguno de los dos casos: en el caso en que $L_{r_0+1} = L_{r_0} + n$, o en el caso en que $L_{r_0+1} = L_{r_0} + n + 1$.
Es decir debemos demostrar que

$$\Delta^-_{r_0+1, j} := L_{r_0} + n - L_j - L_{r_0+1-j} \in \{0, 1\} \text{ para todo } 0 < j < r_0 + 1,$$

o que

$$\Delta^+_{r_0+1, j} := L_{r_0} + n + 1 - L_j - L_{r_0+1-j} \in \{0, 1\} \text{ para todo } 0 < j < r_0 + 1.$$

Evidentemente $\Delta^+_{r_0+1,j} = \Delta^-_{r_0+1,j} + 1$ y además

$$\Delta^-_{r_0+1,j} = L_{r_0} + n - L_j - L_{r_0+1-j} = L_{r_0} + n - L_j - L_{r_0+1-j} + L_{r_0-j} - L_{r_0-j}$$

$$= \Delta_{r_0,j} + n - (L_{r_0+1-j} - L_{r_0-j}).$$

Pero $(L_{r_0+1-j} - L_{r_0-j}) \in \{n, n+1\}$.

Luego

$$\Delta^-_{r_0+1,j} \in \{\Delta_{r_0,j}, \Delta_{r_0,j} - 1\} \subset \{-1, 0, 1\},$$

por lo tanto $\Delta^+_{r_0+1,j} \in \{0, 1, 2\}$. Ahora, si se cumple

$$\Delta^-_{r_0+1,j} \in \{0, 1\} \text{ para todo } 0 < j < r_0 + 1,$$

ya no hay nada que demostrar. Pero si no se cumple, es decir $\Delta^-_{r_0+1,j_0} = -1$ para algún $0 < j_0 < r_0 + 1$, entonces tenemos que demostrar que $\Delta^+_{r_0+1,j} \in \{0, 1\}$ para todo $0 < j < r_0 + 1$.

Asumamos por contradicción que existe $j_1$ tal que $\Delta^+_{r_0+1,j_1} = 2$, , lo cual significa que $\Delta^-_{r_0+1,j_1} = 1$.

Podemos asumir que $j_1 > j_0$ ya que $\Delta^-_{r_0+1,j} = \Delta^-_{r_0+1,r_0+1-j}$.

Hacemos $s := j_1 - j_0 = (r_0 + 1 - j_0) - (r_0 + 1 - j_1) > 0$ y obtenemos

$$\begin{aligned} 2 &= \Delta^-_{r_0+1,j_1} - \Delta^-_{r_0+1,j_0} \\ &= -L_{r_0+1-j_1} - L_{j_1} + L_{r_0+1-j_0} + L_{j_0} + L_s - L_s \\ &= L_{r_0+1-j_0} - L_s - L_{r_0+1-j_1} - (L_{j_1} - L_{j_0} - L_s) \\ &= \Delta_{r_0+1-j_0,s} - \Delta_{j_1,s}. \end{aligned}$$

Pero $\Delta_{r_0+1-j_0,s}, \Delta_{j_1,s} \in \{0, 1\}$, de manera que obtenemos una contradicción y concluye la demostración. $\square$

**Corolario 4.3.3.** *Si $\tilde{L} = (L_1, \ldots, L_r)$ es casi-balanceada, entonces existe $L \in \mathscr{L}$ tal que $L_{\leq r} = \tilde{L}$.*

**Demostración.** Usamos la Proposición 4.3.2 con $r_0 = r$. Luego haciendo $L_{r+1} = L_r + n$ o $L_{r+1} = L_r + n + 1$ se tiene que la sucesión $(L_1, \ldots, L_r, L_{r+1})$ es casi-balanceada. Volviendo a usar la Proposición 4.3.2 para $r_0 = r + 1$ se construye $L_{r+2}$ tal que la sucesión $(L_1, \ldots, L_r, L_{r+1}, L_{r+2})$ es casi-balanceada. De la misma forma se construyen inductivamente $L_{r+3}, L_{r+4}, \ldots$ $\square$

**Lema 4.3.4.** *Si $L \in \Delta(n, n+1)$, pero $L \notin \mathscr{L}$, existen $k, r \in \mathbb{N}$ tales que se cumple una de las siguientes posibilidades*

(1) $L_{k+r} = L_k + L_r - 1$ y $L_{r-1} + 2 < L_r$,

(2) $L_{k+r} = L_k + L_r + 2$ y $L_r + 2 < L_{r+1}$.

**Demostración.** Sea $m := \min\{i : L_{\leq i} \text{ no es casi-balanceada}\}$. Para $0 < j < m$ se tiene
$$\Delta_{m,j} = L_m - L_j - L_{m-j},$$
$$\Delta_{m-1,j} = L_{m-1} - L_j - L_{m-1-j}.$$

Luego
$$\Delta_{m,j} = \Delta_{m-1,j} + (L_m - L_{m-1}) - (L_{m-j} - L_{m-j-1}).$$

Dado que $\Delta_{m-1,j} \in \{0,1\}$, $(L_m - L_{m-1}) \in \{n, n+1\}$ y $(L_{m-j} - L_{m-j-1}) \in \{n, n+1\}$, necesariamente $\Delta_{m,j} \in \{-1, 0, 1, 2\}$. Pero $L_{\leq m}$ no es casi-balanceada, de manera que existe $j$ tal que $\Delta_{m,j} = -1$ o $\Delta_{m,j} = 2$. En el caso $\Delta_{m,j} = -1$ pongamos $r := \min\{j : \Delta_{m,j} = -1\}$ y $k := m - r$.

Entonces
$$\Delta_{m,r} = L_m - L_r - L_{m-r},$$
$$\Delta_{m,r-1} = L_m - L_{r-1} - L_{m-r+1}.$$
$$\Delta_{m,r-1} - \Delta_{m,r} = (L_r - L_{r-1}) - (L_{m-r+1} - L_{m-r}).$$

Luego $|\Delta_{m,r} - \Delta_{m,r-1}| \leq 1$ y por lo tanto $\Delta_{m,r-1} = 0$.

Pero $\Delta_{m,r-1} = \Delta_{m,r} + (L_r - L_{r-1}) - (L_{m-r+1} - L_{m-r})$, y reemplazando los valores hallados se obtiene $0 = -1 + (L_r - L_{r-1}) - (L_{m-r+1} - L_{m-r})$.

Ahora $(L_r - L_{r-1}), (L_{m-r+1} - L_{m-r}) \in \{n, n+1\}$, entonces $L_r - L_{r-1} = n+1$ y $L_{m-r+1} - L_{m-r} = n$.

De manera que $L_r - L_{r-1} = n + 1 > 2$.

Por otro lado $\Delta_{k+r,r} = \Delta_{m,r} = L_{k+r} - L_k - L_r = -1$.

Esto demuestra que estamos en el caso (1).

En el caso $\Delta_{m,j} = 2$ pongamos $r := \max\{j : \Delta_{m,j} = 2\}$ y $k := m - r$.

Entonces $\Delta_{m,r+1} = 1$, ya que $|\Delta_{m,r} - \Delta_{m,r+1}| \leq 1$.

Pero $\Delta_{m,r+1} = \Delta_{m,r} + (L_{m-r} - L_{m-r-1}) - (L_{r+1} - L_r)$.

Por lo tanto $1 = 2 + (L_{m-r} - L_{m-r-1}) - (L_{r+1} - L_r)$, y como $(L_{m-r} - L_{m-r-1}), (L_{r+1} - L_r) \in \{n, n+1\}$, tenemos $L_{r+1} - L_r = n + 1$ y $L_{m-r} - L_{m-r-1} = n$.

Así $L_{r+1} - L_r = n + 1 > 2$ y $\Delta_{k+r,r} = L_{k+r} - L_k - L_r = 2$.

Esto demuestra que estamos en el caso (2), lo cual concluye la demostración. □

**Proposición 4.3.5.** *Sea $L \in \Delta(n, n+1)$. Si la matriz infinita $\widetilde{M} := M(L, a)$ define una aplicación de torcimiento vía la Observación 3.1.1, entonces $L \in \mathcal{L}$.*

**Demostración.** Supongamos que $\widetilde{M} = M(L, a)$ define una aplicación de torcimiento y supongamos por contradicción que $L \in \Delta(n, n+1) \setminus \mathcal{L}$. En el primer caso del Lema 4.3.4 notemos que $L_{r-1} < L_r - 2 < L_r - 1 < L_r$, por lo tanto por definición $\widetilde{M}_{L_r-2,*} = \widetilde{M}_{L_r-1,*} = d^{r-1}(E_0 - E_1)$, y así, por el Lema 3.1.2(1) tenemos $\widetilde{M} Y^{L_r-2} \widetilde{M} = 0$.
Pero entonces, usando el hecho de que $L_{k+r} = L_k + L_r - 1$ y que $a, d, 1-d \neq 0$ obtenemos

$$\begin{aligned}
0 = (\widetilde{M} Y^{L_r-2} \widetilde{M})_{L_k,0} &= \widetilde{M}_{L_k,*} \cdot \widetilde{M}_{*+L_r-2,0} \\
&= (d^k, -d^{k-1}, 0, \cdots, a) \cdot (\widetilde{M}_{L_r-2,0}, \widetilde{M}_{L_r-1,0}, \cdots, \widetilde{M}_{L_k+L_r-1,0}) \\
&= (d^k, -d^{k-1}, 0, \cdots, a) \cdot (\widetilde{M}_{L_r-2,0}, \widetilde{M}_{L_r-1,0}, \cdots, \widetilde{M}_{L_{k+r},0}) \\
&= d^k \widetilde{M}_{L_r-2,0} - d^{k-1} \widetilde{M}_{L_r-1,0} + a \widetilde{M}_{L_{k+r},0} \\
&= d^{k+r-1} - d^{k+r-2} + a d^{k+r} \\
&= d^{r+k-2}(d - 1 + a d^2) = d^{r+k-2}(-ad + a d^2) \\
&= -a d^{k+r-1}(1-d) \neq 0,
\end{aligned}$$

una contradicción que descarta este caso.

Por otro lado, en el segundo caso del Lema 4.3.4 notamos que $L_r < L_r + 1 < L_r + 2 < L_{r+1}$, por lo tanto por definición $\widetilde{M}_{L_r+1,*} = \widetilde{M}_{L_r+2,*} = d^r(E_0 - E_1)$, y asi por el Lema 3.1.2(1) tenemos $\widetilde{M} Y^{L_r+1} \widetilde{M} = 0$.

Pero entonces

$$0 = (\widetilde{M} Y^{L_r+1} \widetilde{M})_{L_k, L_{r+k}+1} = \widetilde{M}_{L_k,*} \cdot \widetilde{M}_{*+L_r+1, L_{r+k}+1}$$

$$= (d^k, -d^{k-1}, 0, \cdots, a) \cdot (\widetilde{M}_{L_r+1, L_{r+k}+1}, \widetilde{M}_{L_r+2, L_{r+k}+1}, \cdots, \widetilde{M}_{L_k+L_r+2, L_{r+k}+1})$$

$$= (d^k, -d^{k-1}, 0, \cdots, a) \cdot (\widetilde{M}_{L_r+1, L_{r+k}+1}, \widetilde{M}_{L_r+2, L_{r+k}+1}, \cdots, \widetilde{M}_{L_{k+r}, L_{r+k}+1})$$

$$= d^k \widetilde{M}_{L_r+1, L_{r+k}+1} - d^{k-1} \widetilde{M}_{L_r+2, L_{r+k}+1} + a \widetilde{M}_{L_{k+r}, L_{r+k}+1}.$$

Pero $\widetilde{M}_{L_r+1, L_{r+k}+1} = 0 = \widetilde{M}_{L_r+2, L_{r+k}+1}$ y $\widetilde{M}_{L_{k+r}, L_{r+k}+1} = a$, y así llegamos a

$$0 = a^2 \neq 0,$$

una contradicción que descarta el segundo caso del Lema 4.3.4 y concluye la demostración. □

## 4.4. Existencia de aplicaciones de torcimiento para la familia $B(n, L)$

En el Teorema 3.2.3 determinamos la forma de las matrices $\widetilde{M}$ correspondiente a las aplicaciones de torcimiento de la familia $A(n, d, a)$ y en el Corolario 3.2.7 demostramos que cada tal matriz define realmente una aplicación de torcimiento. Similarmente, en la Proposición 4.3.5 demostramos que necesariamente $\widetilde{M} :=$ $M(L,a)$ tiene que satisfacer $L \in \mathscr{L}$, para definir una aplicación de torcimiento via la Observación 3.1.1. Esta última sección está dedicada a la demostración de que la condición $L \in \mathscr{L}$ es suficiente, es decir, que toda $L \in \mathscr{L}$ define una aplicación de torcimiento, poniendo $\widetilde{M} := M(L,a)$ y usando la Observación 3.1.1.

De manera que, en toda esta última sección, asumimos que $L \in \Delta(n, n+1)$ y ponemos $\widetilde{M} := M(L,a)$.

**Lema 4.4.1.** *Supongamos que se cumplen las siguientes condiciones*

(1) $\widetilde{M} Y^k \widetilde{M} = 0$ *para* $k \in \mathbb{N}_0$ *con* $k, k+1 \neq L_t$ *para todo* $t$,

(2) $\widetilde{M} Y^{L_r} \widetilde{M} + Y \widetilde{M} Y^{L_r - 1} \widetilde{M} = a Y^{L_r + 1} \widetilde{M}$ *para todo* $r \in \mathbb{N}$,

(3) $\widetilde{M} Y^{L_r - 1} \widetilde{M} = a \sum_{i=0}^{L_r} Y^i \widetilde{M} Y^{L_r - i} - r a^2 Y^{L_r + 1}$ *para todo* $r \in \mathbb{N}$.

*Pongamos $L_0 := 0$, entonces para todo $r \in \mathbb{N}_0$ y $L_r \leq j < L_{r+1}$, se cumple*

$$d^r M^{j+1} = \sum_{i=0}^{j} Y^i \widetilde{M} Y^{j-i} + (1 - ra) Y^{j+1}, \qquad (4.4.29)$$

*y para todo* $k \geq 0$ *se cumple* (3.1.1)

$$Y^k \widetilde{M} = \sum_{j=0}^{k+1} \widetilde{M}_{k,j} M^{k+1-j} Y^j.$$

**Demostración.** Primero demostraremos (4.4.29) por inducción sobre $j$.

Para $j = L_0 = 0$ se cumple trivialmente.

Ahora supongamos que (4.4.29) es válido para algún $j$, demostraremos que se cumple para $j + 1$. Primero asumamos que $L_r \leq j < L_{r+1} - 1$. Luego multiplicamos (4.4.29) por $M = Y + \widetilde{M}$ y obtenemos:

$$d^r M^{j+2} = \sum_{i=0}^{j} Y^i \widetilde{M} Y^{j-i} (Y + \widetilde{M}) + (1 - ra) Y^{j+1} (Y + \widetilde{M})$$

$$= \sum_{i=0}^{j} Y^i \widetilde{M} Y^{j+1-i} + \sum_{i=0}^{j} Y^i \widetilde{M} Y^{j-i} \widetilde{M} + (1 - ra) Y^{j+2} + (1 - ra) Y^{j+1} \widetilde{M}$$

$$= \sum_{i=0}^{j} Y^i \widetilde{M} Y^{j+1-i} + Y^{j+1} \widetilde{M} - ra Y^{j+1} \widetilde{M} + \sum_{i=0}^{j} Y^i \widetilde{M} Y^{j-i} \widetilde{M} + (1 - ra) Y^{j+2}.$$

Por la condición (1), para $L_t < i < L_{t+1} - 1$ se tiene $\widetilde{M} Y^i \widetilde{M} = 0$ y por la condición (2) se tiene

$$\widetilde{M} Y^{L_t} \widetilde{M} + Y \widetilde{M} Y^{L_t - 1} \widetilde{M} = a Y^{L_t + 1} \widetilde{M}$$

para $1 \leq t \leq r$.

Multiplicando a la izquierda por $Y^{j-L_t}$ se obtiene

$$Y^{j-L_t}\widetilde{M}Y^{L_t}\widetilde{M} + Y^{j-L_t+1}\widetilde{M}Y^{L_t-1}\widetilde{M} = a Y^{j+1}\widetilde{M}. \qquad (4.4.30)$$

Luego

$$\sum_{i=0}^{j} Y^i \widetilde{M} Y^{j-i}\widetilde{M} = \sum_{i=0}^{j} Y^{j-i}\widetilde{M}Y^i\widetilde{M} = \sum_{i=0}^{L_1-1} Y^{j-i}\widetilde{M}Y^i\widetilde{M} + \sum_{i=L_1}^{L_2-1} Y^{j-i}\widetilde{M}Y^i\widetilde{M} + \cdots$$

$$+ \sum_{i=L_{r-1}}^{L_r-1} Y^{j-i}\widetilde{M}Y^i\widetilde{M} + \sum_{i=L_r}^{j} Y^{j-i}\widetilde{M}Y^i\widetilde{M}.$$

Pero

$$\widetilde{M}Y^i\widetilde{M} = 0, \quad 0 \leq i < L_1 - 1, \quad L_t < i < L_{t+1} - 1, \quad t = 1, 2, \cdots, r.$$

Luego los únicos términos que sobreviven en la sumatoria son los correspondientes a

$$i = L_t, i = L_t - 1, \quad 1 \leq t \leq r.$$

Entonces por la igualdad (4.4.30) se tiene

$$\sum_{i=0}^{j} Y^{j-i}\widetilde{M}Y^i\widetilde{M} = \sum_{t=1}^{r} \left( Y^{j-L_t+1}\widetilde{M}Y^{L_t-1}\widetilde{M} + Y^{j-L_t}\widetilde{M}Y^{L_t}\widetilde{M} \right)$$

$$= \sum_{t=1}^{r} a Y^{j+1}\widetilde{M} = r a Y^{j+1}\widetilde{M}.$$

Por lo tanto

$$d^r M^{j+2} = \sum_{i=0}^{j} Y^i \widetilde{M} Y^{j+1-i} + Y^{j+1}\widetilde{M} + (1-ra)Y^{j+2} - ra Y^{j+1}\widetilde{M} + ra Y^{j+1}\widetilde{M}$$

$$= \sum_{i=0}^{j+1} Y^i \widetilde{M} Y^{j+1-i} + (1-ra)Y^{j+2}.$$

Esto demuestra (4.4.29) para $L_r \leq j < L_{r+1} - 1$.

Ahora si $\boxed{j = L_r - 1}$ para algún $r \geq 1$ entonces $L_{r-1} \leq j < L_r$ y usando el mismo argumento anterior para $r-1$ obtenemos

$$d^{r-1}M^{j+2} = \sum_{i=0}^{j} Y^i \widetilde{M} Y^{j+1-i} + Y^{j+1}\widetilde{M}$$

$$+ (1-(r-1)a)Y^{j+2} - (r-1)a Y^{j+1}\widetilde{M} + \sum_{i=0}^{j} Y^{j-i}\widetilde{M}Y^i\widetilde{M}$$

$$= \sum_{i=0}^{j+1} Y^i \widetilde{M} Y^{j+1-i} + (1-(r-1)a)Y^{j+2} - (r-1)a Y^{j+1}\widetilde{M} + \sum_{i=0}^{j} Y^{j-i}\widetilde{M} Y^i \widetilde{M}.$$

Dado que $L_{r-1} \le L_r - 2 = j - 1 < L_r - 1$, obtenemos

$$\sum_{i=0}^{j-1} Y^{j-i}\widetilde{M} Y^i \widetilde{M} = \sum_{i=0}^{L_1-1} Y^{j-i}\widetilde{M} Y^i \widetilde{M} + \sum_{i=L_1}^{L_2-1} Y^{j-i}\widetilde{M} Y^i \widetilde{M} + \cdots$$

$$+ \sum_{i=L_{r-2}}^{L_{r-1}-1} Y^{j-i}\widetilde{M} Y^i \widetilde{M} + \sum_{i=L_{r-1}}^{j-1} Y^{j-i}\widetilde{M} Y^i \widetilde{M}$$

$$= \sum_{t=1}^{r-1} a Y^{j+1}\widetilde{M} = (r-1)a Y^{j+1}\widetilde{M}.$$

Por lo tanto

$$d^{r-1}M^{j+2} = \sum_{i=0}^{j+1} Y^i \widetilde{M} Y^{j+1-i} + \widetilde{M} Y^j \widetilde{M} + [1-(r-1)a]Y^{j+2},$$

lo cual implica

$$d^{r-1}M^{j+2} = \sum_{i=0}^{L_r} Y^i \widetilde{M} Y^{L_r-i} + \widetilde{M} Y^{L_r-1}\widetilde{M} + [1-(r-1))a]Y^{L_r+1}.$$

Reemplazando el valor de $\widetilde{M} Y^{L_r-1}\widetilde{M}$, por la condición (3), se obtiene

$$d^{r-1}M^{j+2} = \sum_{i=0}^{L_r} Y^i \widetilde{M} Y^{L_r-i} + a\sum_{i=0}^{L_r} Y^i \widetilde{M} Y^{L_r-i} - ra^2 Y^{L_r+1} + [1-(r-1)a]Y^{L_r+1},$$

$$d^{r-1}M^{L_r+1} = (1+a)\sum_{i=0}^{L_r} Y^i \widetilde{M} Y^{L_r-i} + (1+a)(1-ra)Y^{L_r+1}.$$

Multiplicando por $d$ y usando $d(a+1) = 1$ obtenemos (4.4.29) para $j+1 = L_r$.
Esto concluye la demostración de (4.4.29).
Ahora demostramos (3.1.1).
De la condición (1) obtenemos

$$M Y^k \widetilde{M} = (\widetilde{M} + Y)Y^k \widetilde{M} = \widetilde{M} Y^k \widetilde{M} + Y^{k+1}\widetilde{M} = Y^{k+1}\widetilde{M},$$

si $L_r < k < L_{r+1} - 1$ para algún $r$ o si $k < n - 1 = L_1 - 1$.
Por lo tanto, si $Y^k \widetilde{M} = d^r M^k \widetilde{M}$, entonces

$$Y^{k+1}\widetilde{M} = M Y^k \widetilde{M} = M d^r M^k \widetilde{M} = d^r M^{k+1}\widetilde{M},$$

para tal $k$.
Es decir

$$Y^k \widetilde{M} = d^r M^k \widetilde{M} \Rightarrow Y^{k+1}\widetilde{M} = d^r M^{k+1}\widetilde{M} \qquad (4.4.31)$$

para tal $k$.

Primero consideremos $0 \leq k < L_1 - 1 = n - 1$.

Para $k = 0$ se cumple trivialmente $Y^0 \widetilde{M} = d^0 M^0 \widetilde{M}$.

Luego por (4.4.31) se obtiene $Y \widetilde{M} = M \widetilde{M}$.

Para $k = 1$ se cumple $Y \widetilde{M} = M \widetilde{M}$.

Luego por (4.4.31) se obtiene $Y^2 \widetilde{M} = M^2 \widetilde{M}$.

Así inductivamente obtenemos

$$Y^k \widetilde{M} = M^k \widetilde{M} \quad \text{para } k \leq n - 1,$$

lo cual demuestra (3.1.1) para $k = 1, \ldots, n - 1$.

Multiplicando $Y^{n-1} \widetilde{M} = M^{n-1} \widetilde{M}$ por $M^2 = \widetilde{M} Y + Y \widetilde{M} + Y^2$ a la izquierda obtenemos

$$(\widetilde{M} Y + Y \widetilde{M} + Y^2) Y^{n-1} \widetilde{M} = M^{n+1} \widetilde{M},$$

$$\widetilde{M} Y^n \widetilde{M} + Y \widetilde{M} Y^{n-1} \widetilde{M} + Y^{n+1} \widetilde{M} = M^{n+1} \widetilde{M},$$

y combinando con $\widetilde{M} Y^n \widetilde{M} + Y \widetilde{M} Y^{n-1} \widetilde{M} = a Y^{n+1} \widetilde{M}$, lo cual es válido por la condición (2) para $L_1 = n$, esto produce

$$M^{n+1} \widetilde{M} = (a + 1) Y^{n+1} \widetilde{M}.$$

Usando $(a + 1)d = 1$ se obtiene

$$Y^{n+1} \widetilde{M} = d M^{n+1} \widetilde{M},$$

lo cual es (3.1.1) para $k = n + 1$.

Por el argumento anterior, obtenemos $Y^k \widetilde{M} = d M^k \widetilde{M}$ para $k = n + 1, \ldots, L_2 - 1$, lo cual es (3.1.1) para tal $k$.

Multiplicando $Y^{L_2 - 1} \widetilde{M} = d M^{L_2 - 1} \widetilde{M}$ por $M^2 = \widetilde{M} Y + Y \widetilde{M} + Y^2$ a la izquierda obtenemos

$$\widetilde{M} Y^{L_2} \widetilde{M} + Y \widetilde{M} Y^{L_2 - 1} \widetilde{M} + Y^{L_2 + 1} \widetilde{M} = d M^{L_2 + 1} \widetilde{M},$$

y combinando con $\widetilde{M} Y^{L_2} \widetilde{M} + Y \widetilde{M} Y^{L_2 - 1} \widetilde{M} = a Y^{L_2 + 1} \widetilde{M}$, lo cual es válido por la condición (2) para $L_2$, esto produce

$$d M^{L_2 + 1} \widetilde{M} = (a + 1) Y^{L_2 + 1} \widetilde{M}.$$

Usando $(a + 1)d = 1$ se obtiene $Y^{L_2 + 1} \widetilde{M} = d^2 M^{L_2 + 1} \widetilde{M}$, lo cual es (3.1.1) para $k = L_2 + 1$.

Repitiendo inductivamente estos argumentos obtenemos (3.1.1) para todo $k$ que no sea igual a algún $L_r$.

Finalmente demostraremos (3.1.1) para $k = L_r$.

Para esto consideremos la igualdad

$$d^r M^{L_r+1} = \sum_{i=0}^{L_r} Y^i \widetilde{M} Y^{L_r-i} + (1-ra)Y^{L_r+1},$$

la cual es la igualdad (4.4.29) para $j = L_r$, y la igualdad

$$d^{r-1} M^{L_r} Y = \sum_{i=0}^{L_r-1} Y^i \widetilde{M} Y^{L_r-i} + (1-(r-1)a)Y^{L_r+1},$$

la cual es la igualdad (4.4.29) para $j = L_r - 1$, multiplicada por $Y$ a la derecha.
Restando la segunda igualdad de la primera, obtenemos

$$\begin{aligned} d^r M^{L_r+1} - d^{r-1} M^{L_r} Y &= \sum_{i=0}^{L_r} Y^i \widetilde{M} Y^{L_r-i} - \sum_{i=0}^{L_r-1} Y^i \widetilde{M} Y^{L_r-i} - a Y^{L_r+1} \\ &= Y^{L_r} \widetilde{M} - a Y^{L_r+1}. \end{aligned}$$

Por lo tanto $Y^{L_r} \widetilde{M} = d^r M^{L_r+1} - d^{r-1} M^{L_r} Y + a Y^{L_r+1}$, lo cual es (3.1.1)
para $k = L_r$ y de esta manera concluye la demostración. □

Para demostrar que para $L \in \mathscr{L}$ la matriz $\widetilde{M} := M(L, a)$ define una aplicación de torcimiento, descomponemos la matriz $\widetilde{M}$ en tres sumandos.
Ponemos

$$m_i := \begin{cases} 1, & \text{si } i \leq n = L_1 \\ d^r, & \text{si } L_r < i \leq L_{r+1} \end{cases}$$

y definimos la matriz $M_1$ mediante $(M_1)_{i,j} := \delta_{0,j} m_i$.
Ahora definimos el conjunto

$$|L| := \{L_t\}_{t \in \mathbb{N}} = \{L_1, L_2, \ldots, L_t, \ldots\},$$

y ponemos

$$n_i = \begin{cases} a, & \text{si } i \in |L| \\ 0, & \text{de lo contrario.} \end{cases}$$

Definimos la matriz $B$ mediante $B_{i,j} := \delta_{i,j} n_i$.
Entonces $\widetilde{M} = BY + YM_1 - M_1 Y$ y

$$BYM_1 + YM_1 - M_1 = 0. \tag{4.4.32}$$

En efecto:
$(BY)_{i,j} = \delta_{i,j-1} n_i, \quad (YM_1)_{i,j} = \delta_{0,j} m_{i+1}, \quad (M_1 Y)_{i,j} = \delta_{0,j-1} m_i$
Si $i < L_1 = n$ entonces $(BY + YM_1 - M_1 Y)_{i,j} = \delta_{0,j} - \delta_{0,j-1}$.

Por lo tanto $(BY+YM_1-M_1Y)_{i,*}=E_0-E_1=\widetilde{M}_{i,*}$.
Si $L_r<i<L_{r+1}$ entonces $(BY+YM_1-M_1Y)_{i,j}=\delta_{0,j}d^r-\delta_{0,j-1}d^r$.
Por lo tanto $(BY+YM_1-M_1Y)_{i,*}=d^r(E_0-E_1)=\widetilde{M}_{i,*}$.
Si $i=L_r$ entonces $(BY+YM_1-M_1Y)_{i,j}=\delta_{L_r,j-1}a+\delta_{0,j}d^r-\delta_{0,j-1}d^{r-1}$.
Por lo tanto $(BY+YM_1-M_1Y)_{i,*}=d^rE_0-d^{r-1}E_1+aE_{L_r+1}=\widetilde{M}_{i,*}$.
Ahora demostraremos (4.4.32).

$$(BYM_1+YM_1-M_1)_{i,j}=(BYM_1)_{i,j}+(YM_1)_{i,j}-(M_1)_{i,j}$$

$$=\delta_{0,j}n_im_{i+1}+\delta_{0,j}m_{i+1}-\delta_{0,j}m_i.$$

Luego $(BYM_1+YM_1-M_1)_{i,j}=0$, si $j\neq 0$.
Para $j=0$ tenemos
$(BYM_1+YM_1-M_1)_{i,0}=n_im_{i+1}+m_{i+1}-m_i$.
Si $i<n$ tenemos $(BYM_1+YM_1-M_1)_{i,0}=0+1-1=0$.
Si $i=n$ tenemos $(BYM_1+YM_1-M_1)_{i,0}=ad+d-1=0$.
Si $L_r<i<L_{r+1}$ tenemos $(BYM_1+YM_1-M_1)_{i,0}=0+d^r-d^r=0$.
Si $i=L_r$ tenemos $(BYM_1+YM_1-M_1)_{i,0}=ad^r+d^r-d^{r-1}=d^{r-1}(ad+d-1)=0$.

**Lema 4.4.2.** *Las siguientes igualdades son válidas para $k\in\mathbb{N}_0$:*

1. $M_1Y^kM_1=m_kM_1$,

2. $M_1Y^kB=n_kM_1Y^k$,

3. $(BY^kM_1)_{i,j}=\delta_{0,j}n_im_{k+i}$,

4. $(BY^kB)_{i,j}=\delta_{i+k,j}n_in_{i+k}$.

**Demostración**

Demostración de 1.
$(M_1Y^kM_1)_{i,j}=\sum_l(M_1)_{i,l}(M_1)_{l+k,j}=\sum_l\delta_{0,l}m_i\delta_{0,j}m_{l+k}=\delta_{0,j}m_im_k$.
Por otro lado $m_k(M_1)_{i,j}=\delta_{0,j}m_im_k$.
Demostración de 2.
$(M_1Y^kB)_{i,j}=\sum_l(M_1)_{i,l}(B)_{l+k,j}=\sum_l\delta_{0,l}m_i\delta_{l+k,j}n_{l+k}=\delta_{k,j}m_in_k$.
Por otro lado $n_k(M_1Y^k)_{i,j}=n_k(M_1)_{i,j-k}=\delta_{0,j-k}m_in_k=\delta_{k,j}m_in_k$.
Demostración de 3.
$(BY^kM_1)_{i,j}=\sum_l(B)_{i,l}(M_1)_{l+k,j}=\sum_l\delta_{i,l}n_i\delta_{0,j}m_{l+k}=\delta_{0,j}n_im_{k+i}$.
Demostración de 4.
$(BY^kB)_{i,j}=\sum_l(B)_{i,l}(B)_{l+k,j}=\sum_l\delta_{i,l}n_i\delta_{l+k,j}n_{l+k}=\delta_{i+k,j}n_in_{i+k}$. $\square$

**Observación 4.4.3.** Notemos que si $L \in \mathscr{L}$ y $k, k+1 \notin |L|$, entonces $BY^{k+1}B = 0$. De hecho, por el Lema 4.4.2(4), es suficiente demostrar que si $n_i n_{i+k+1} \neq 0$ para algún $i$, entonces $k \in |L|$ o $k+1 \in |L|$.

Pero $n_i n_{i+k+1} \neq 0$ implica $i = L_t$ y $i+k+1 = L_{r+t}$ para algún $r \in \mathbb{N}, t \in \mathbb{N}$.
Dado que $L$ es casi-balanceada, entonces $L_{r+t} = L_t + L_r$, lo cual implica $k+1 = L_r$, o $L_{r+t} = L_t + L_r + 1$, lo cual implica $k = L_r$.

**Lema 4.4.4.** *Sea $L \in \mathscr{L}$. Si $i \in |L|$ y $k \notin |L|$, entonces $m_{i+k+1} = d\, m_i m_k$.*

**Demostración.** Sabemos que $i = L_t$ para algún $t \in \mathbb{N}$.
Si $k < n = L_1$, entonces $L_t < i+k+1 \le i+n \le L_t + n \le L_{t+1}$, y así $m_{i+k+1} = d^t$.
Dado que $m_i = d^{t-1}$ y $m_k = 1$, esto demuestra el resultado en este caso.
En cambio si $L_r < k < L_{r+1}$ para algún $r \in \mathbb{N}$, y dado que $L$ es casi-balanceada, tenemos
$$L_r + L_t = L_r + i < i + k < L_{r+1} + i = L_{r+1} + L_t,$$
$$L_{r+t} - 1 < i + k < L_{t+r+1},$$
$$L_{r+t} < i + k + 1 \le L_{t+r+1},$$
y así $m_{i+k+1} = d^{r+t}$.
Dado que $m_i = d^{t-1}$ y $m_k = d^r$, esto concluye la demostración. $\square$

**Lema 4.4.5.** *Supongamos que $L \in \mathscr{L}$, y sean $i, k \in \mathbb{N}_0$ y $r \in \mathbb{N}$. Entonces*

$$n_i m_{i+k+1} + m_k m_{i+1} - m_{k+1} m_i = n_k m_{i+k+1}, \tag{4.4.33}$$

$$n_i n_{L_r+i+1} + n_{i+1} n_{L_r+i+1} = a n_{L_r+i+1}, \tag{4.4.34}$$

$$a \sum_{j=0}^{L_r} n_{i+j} = r a^2 + n_i n_{L_r+i}. \tag{4.4.35}$$

**Demostración.** Primero demostramos (4.4.33) en cada uno de los cuatro casos posibles. Notemos que para $k \in |L|$ tenemos $m_{k+1} = d m_k$, y para $k \notin |L|$ tenemos $m_{k+1} = m_k$.

**Caso $i, k \notin |L|$:**
Aquí $n_i = n_k = 0$, $m_{i+1} = m_i$, $m_{k+1} = m_k$, por lo tanto $m_k m_{i+1} - m_{k+1} m_i = 0$ y ambos lados de (4.4.33) se anulan.

**Caso $i \in |L|, k \notin |L|$:**
Aquí $n_i = a$, $n_k = 0$, $m_{i+1} = d m_i$, $m_{k+1} = m_k$ y por el Lema 4.4.4 tenemos $m_{i+k+1} = d m_i m_k$. Por lo tanto

$$n_i m_{i+k+1} + m_k m_{i+1} - m_{k+1} m_i = m_i m_k (ad + d - 1) = 0,$$

y ambos lados de (4.4.33) se anulan.

**Caso $i \notin |L|, k \in |L|$:**

Aqui $n_i = 0$, $n_k = a$, $m_{i+1} = m_i$, $m_{k+1} = d\, m_k$ y por el Lema 4.4.4 tenemos $d\, m_i m_k = m_{i+k+1}$. Por lo tanto

$$n_i m_{i+k+1} + m_k m_{i+1} - m_{k+1} m_i = (1-d) m_i m_k = a d\, m_i m_k = a m_{i+k+1} = n_k m_{i+k+1},$$

que es lo que se quería.

**Caso $i, k \in |L|$:**

Aquí $n_i = a$, $n_k = a$, $m_{i+1} = d\, m_i$ y $m_{k+1} = d\, m_k$. Por lo tanto

$$n_i m_{i+k+1} + m_k m_{i+1} - m_{k+1} m_i = a m_{i+k+1} = n_k m_{i+k+1},$$

que es lo que se quería, concluyendo la demostración de (4.4.33).

Ahora demostramos (4.4.34). Si $L_r + i + 1 \notin |L|$, entonces ambos lados se anulan. Si $L_r + i + 1 \in |L|$, entonces tenemos que demostrar que $n_i + n_{i+1} = a$.

Existe $t > 0$ tal que $L_r + i + 1 = L_{r+t}$. Dado que $L \in \mathscr{L}$, tenemos

$$L_{r+t} = L_r + L_t \quad \text{o} \quad L_{r+t} = L_r + L_t + 1.$$

En el primer caso $L_t = i+1$, $n_i = 0$ y $n_{i+1} = a$; y en el segundo caso $L_t = i$, $n_i = a$ y $n_{i+1} = 0$. Por lo tanto en ambos casos $n_i + n_{i+1} = a$, lo cual demuestra (4.4.34).

Para demostrar (4.4.35) consideramos tres casos.

Si $\boxed{i < L_1 = n}$, entonces $L_r \leq L_r + i < L_{r+1}$, y así

$$\sum_{j=0}^{L_r} n_{j+i} = \sum_{j=i}^{L_r+i} n_j$$

$$= \sum_{j=i}^{L_1+i} n_j + \sum_{j=L_1+i+1}^{L_2+i} n_j + \cdots + \sum_{j=L_{r-1}+i+1}^{L_r+i} n_j = \sum_{s=1}^{r} n_{L_s} = a r.$$

Dado que $n_i = 0$, obtenemos (4.4.35).

Si $\boxed{L_t < i < L_{t+1}}$, entonces $L_{r+t} \leq L_r + i < L_{r+t+1}$, dado que $L \in \mathscr{L}$ implica $L_{r+t} - 1 \leq L_r + L_t < L_r + i < L_r + L_{t+1} \leq L_{r+t+1}$. Por lo tanto

$$\sum_{j=0}^{L_r} n_{j+i} = \sum_{j=i}^{L_r+i} n_j = \sum_{s=t+1}^{t+r} n_{L_s} = a r,$$

y usando $n_i = 0$, obtenemos (4.4.35) en este caso.

Si $\boxed{i = L_t}$ para algún $t$, entonces $a n_{i+L_r} = n_i n_{L_r+i}$, y así es suficiente demostrar

$$r a = \sum_{j=0}^{L_r-1} n_{j+i} = \sum_{j=L_t}^{L_r+L_t-1} n_j.$$

Pero $L \in \mathscr{L}$ implica

$$L_{r+t-1} \le L_r + L_{t-1} + 1 \le L_r + L_t - 1 < L_{r+t},$$

y así

$$\sum_{j=L_t}^{L_r+L_t-1} n_j = \sum_{s=0}^{r-1} n_{L_{t+s}} = ra,$$

que es lo que se quería.

De esta manera (4.4.35) es válida en todos los casos, concluyendo la demostración. □

Queremos calcular

$$\widetilde{M} Y^k \widetilde{M} = \widetilde{M} Y^k BY + \widetilde{M} Y^{k+1} M_1 - \widetilde{M} Y^k M_1 Y, \qquad (4.4.36)$$

y así tenemos que calcular $\widetilde{M} Y^k M_1$.

**Lema 4.4.6.** *Si $L \in \mathscr{L}$, entonces tenemos*

$$\widetilde{M} Y^k M_1 = n_k Y^{k+1} M_1.$$

**Demostración.** Tenemos

$$\left(n_k Y^{k+1} M_1\right)_{i,j} = n_k (M_1)_{i+k+1,j} = \delta_{0,j} n_k m_{i+k+1}$$

y por el Lema 4.4.2 y la igualdad (4.4.33), también tenemos

$$\left(\widetilde{M} Y^k M_1\right)_{i,j} = \left((BY + YM_1 - M_1 Y) Y^k M_1\right)_{i,j} = \left(BY^{k+1} M_1 + YM_1 Y^k M_1 - M_1 Y^{k+1} M_1\right)_{i,j}$$
$$= \delta_{0,j} (n_i m_{i+k+1} + m_k m_{i+1} - m_{k+1} m_i)$$
$$= \delta_{0,j} n_k m_{i+k+1},$$

como se quería. □

**Proposición 4.4.7.** *Si $L \in \mathscr{L}$, entonces*

$$\widetilde{M} Y^k \widetilde{M} = BY^{k+1} BY + n_k YM_1 Y^{k+1} - n_{k+1} M_1 Y^{k+2} + n_{k+1} Y^{k+2} M_1 - n_k Y^{k+1} M_1 Y.$$

**Demostración.** Por la igualdad (4.4.36) y el Lema 4.4.6 tenemos

$$\widetilde{M} Y^k \widetilde{M} = \widetilde{M} Y^k BY + \widetilde{M} Y^{k+1} M_1 - \widetilde{M} Y^k M_1 Y = \widetilde{M} Y^k BY + n_{k+1} Y^{k+2} M_1 - n_k Y^{k+1} M_1 Y.$$

Del Lema 4.4.2(2) obtenemos

$$\widetilde{M} Y^k BY = (BY + YM_1 - M_1 Y) Y^k BY = BY^{k+1} BY + YM_1 Y^k BY - M_1 Y^{k+1} BY$$
$$= BY^{k+1} BY + n_k YM_1 Y^{k+1} - n_{k+1} M_1 Y^{k+2},$$

lo cual concluye la demostración. □

**Proposición 4.4.8.** *Para cada $a \in K \setminus \{0, -1\}$ y $L \in \mathscr{L}$, la matriz $\widetilde{M} = M(L, a)$ define una aplicación de torcimiento vía la Observación 3.1.1.*

**Demostración.** Dado que $\widetilde{M}_{0,j} = \delta_{0,j} - \delta_{1,j}$ and $\widetilde{M}_{k,j} = 0$ para $j > k+1$, por la Observación 3.1.1 y el Lema 4.4.1 es suficiente demostrar que

1. $\widetilde{M} Y^k \widetilde{M} = 0$ para $k \in \mathbb{N}_0$ con $k, k+1 \neq L_t$ para todo $t$,

2. $\widetilde{M} Y^{L_r} \widetilde{M} + Y \widetilde{M} Y^{L_r - 1} \widetilde{M} = a Y^{L_r + 1} \widetilde{M}$ para todo $r \in \mathbb{N}$,

3. $\widetilde{M} Y^{L_r - 1} \widetilde{M} = a \sum_{i=0}^{L_r} Y^i \widetilde{M} Y^{L_r - i} - r a^2 Y^{L_r + 1}$ para todo $r \in \mathbb{N}$.

Si $k, k+1 \notin |L|$, entonces por la Proposición 4.4.7 sabemos que $\widetilde{M} Y^k \widetilde{M} = B Y^{k+1} B Y$. Por la Observación 4.4.3 también sabemos que $B Y^{k+1} B = 0$ para $k, k+1 \notin |L|$, lo cual demuestra el item (1).

Para demostrar el item (2), usamos la Proposición 4.4.7 y calculamos

$$\widetilde{M} Y^{L_r} \widetilde{M} = B Y^{L_r + 1} B Y + a Y M_1 Y^{L_r + 1} - a Y^{L_r + 1} M_1 Y,$$
$$Y \widetilde{M} Y^{L_r - 1} \widetilde{M} = Y B Y^{L_r} B Y - a Y M_1 Y^{L_r + 1} + a Y^{L_r + 2} M_1,$$
$$a Y^{L_r + 1} \widetilde{M} = a Y^{L_r + 1} B Y + a Y^{L_r + 2} M_1 - a Y^{L_r + 1} M_1 Y.$$

Así tenemos que demostrar

$$B Y^{L_r + 1} B + Y B Y^{L_r} B = a Y^{L_r + 1} B.$$

Tenemos

$$(B Y^{L_r + 1} B)_{i,j} = \sum_k B_{i,k} B_{k + L_r + 1, j} = \sum_k \delta_{i,k} \delta_{k + L_r + 1, j} n_i n_{k + L_r + 1} = \delta_{i + L_r + 1, j} n_i n_{i + L_r + 1},$$

$$(Y B Y^{L_r} B)_{i,j} = \sum_k B_{i+1,k} B_{k + L_r, j} = \sum_k \delta_{i+1,k} \delta_{k + L_r, j} n_{i+1} n_{k + L_r} = \delta_{i + L_r + 1, j} n_{i+1} n_{i + L_r + 1},$$

$$(a Y^{L_r + 1} B)_{i,j} = a B_{i + L_r + 1, j} = a \delta_{i + L_r + 1, j} n_{i + L_r + 1},$$

y así (4.4.34) concluye la demostración del item(2).

Para demostrar el item (3), calculamos

$$a \sum_{i=0}^{L_r} Y^i \widetilde{M} Y^{L_r - i} = a \sum_{i=0}^{L_r} Y^i (B Y + Y M_1 - M_1 Y) Y^{L_r - i}$$

$$= a \sum_{i=0}^{L_r} Y^i B Y^{L_r + 1 - i} + a \sum_{i=0}^{L_r} Y^{i+1} M_1 Y^{L_r - i} - a \sum_{i=0}^{L_r} Y^i M_1 Y^{L_r + 1 - i}. \quad (4.4.37)$$

Dado que

$$\sum_{i=0}^{L_r} Y^{i+1} M_1 Y^{L_r - i} = \sum_{i=1}^{L_r + 1} Y^i M_1 Y^{L_r + 1 - i},$$

se tiene

$$a\sum_{i=0}^{L_r} Y^i \widetilde{M} Y^{L_r-i} = a\sum_{i=0}^{L_r} Y^i B Y^{L_r+1-i} + a Y^{L_r+1} M_1 - a M_1 Y^{L_r+1}.$$

Dado que por la Proposición 4.4.7 sabemos que

$$\widetilde{M} Y^{L_r-1} \widetilde{M} = B Y^{L_r} B Y - a M_1 Y^{L_r+1} + a Y^{L_r+1} M_1,$$

tenemos que demostrar

$$a\sum_{s=0}^{L_r} Y^s B Y^{L_r+1-s} = B Y^{L_r} B Y + r a^2 Y^{L_r+1}.$$

Pero

$$(Y^s B Y^{L_r+1-s})_{i,j} = B_{s+i, j+s-L_r-1} = \delta_{i, j-L_r-1} n_{s+i},$$
$$(B Y^{L_r} B Y)_{i,j} = \sum_k B_{i,k} B_{k+L_r, j-1} = \sum_k \delta_{i,k} \delta_{k+L_r, j-1} n_i n_{k+L_r} = n_i n_{i+L_r} \delta_{i+L_r, j-1},$$
$$r a^2 (Y^{L_r+1})_{i,j} = r a^2 \delta_{i, j-L_r-1},$$

por lo tanto es suficiente demostrar

$$a \sum_{s=0}^{L_r} n_{s+i} = n_i n_{i+L_r} + r a^2,$$

lo cual es válido por (4.4.35). Esto concluye la demostración. □

**Observación 4.4.9.** Nuestra estrategia contiene dos componentes principales. Por un lado, las ecuaciones para matrices infinitas producen condiciones que reducen las posibilidades a muy pocas familias. Incluso en el caso más complicado de la sección 4.2 podemos mostrar cómo se puede lograr la clasificación de todas las aplicaciones de torcimiento posibles hasta un grado arbitrario, con un aumento del trabajo computacional. Por otro lado, probar que una matriz infinita dada produce una aplicación de torcimiento requiere verificar un número infinito de igualdades de matrices infinitas. Logramos culminar esa difícil tarea en el Corolario 3.2.7 y en la Proposición 4.4.8.

En el primer caso sólo tuvimos que demostrar una de las igualdades, dado que aquellas aplicaciones de torcimiento tienen la propiedad de $n$- extensión, es decir, están completamente determinadas por los valores de $\widetilde{M}_{k,*}$ para $k \leq n$. En el caso de la Proposición 4.4.8 logramos descomponer la matriz infinita en tres matrices más simples, y probamos las igualdades de matrices infinitas requeridas, usando propiedades de estas matrices más simples.

Sin embargo, ninguna de las aplicaciones de torcimiento construídas en la Proposición 4.4.8 tiene la propiedad de $m$-extensión para ningún $m$. Esto es una consecuencia directa de la siguiente propiedad de las sucesiones casi-balanceadas:

Sea $L_{\leq r} = (L_1, \ldots, L_r)$ una sucesión casi-balanceada parcial.

Entonces existe una extensión de $L_{\leq r}$ de la forma $(L_1, \ldots, L_r, \ldots, L_{r+k})$ tal que tanto

$$(L_1, \ldots, L_r, \ldots, L_{r+k}, L_{r+k} + n) \quad \text{como} \quad (L_1, \ldots, L_r, \ldots, L_{r+k}, L_{r+k} + n + 1)$$

son sucesiones casi- balancedas parciales.

Existen otras propiedades de estas sucesiones. Por ejemplo, las sucesiones casi-balanceadas muestran una sorprendente conexión con la función Totient de Euler.

# Conclusiones

1. Para todo producto tensorial torcido graduado de $K[x]$ por $K[y]$ (plano torcido graduado) existe una representación en $L(K^{\aleph_0})$ que traduce el problema de clasificar todos los planos torcidos graduados, en un problema de clasificar matrices infinitas con entradas en $K$ que satisfacen ciertas condiciones.

2. Para todo plano torcido graduado existen $a, b, c \in K$ tales que se cumple $yx = ax^2 + bxy + cy^2$. Si $a \neq 0$, entonces se puede asumir que $a = 1$.

3. Si $a = 1$, $c \neq 1$ y $(b, c)$ no es raíz de ningún miembro una cierta familia de polinomios $Q_n(b, c)$, entonces se obtiene un álgebra cuadrática, es decir, la estructura multiplicativa del plano torcido graduado está determinada completamente por $b$ y $c$.

4. Existe un plano torcido graduado con $y^k x = x^{k+1} - x^k y + y^{k+1}$ para todo $k$.

5. Para cada $n > 1$, $a$ y $d$ tales que $(a, d)$ no es raíz de ningún miembro una cierta familia de polinomios $R_n(a, d)$, existe un único plano torcido graduado llamado $A(n, d, a)$, tal que $y^k x = x^{k+1} - x^k y + y^{k+1}$ para todo $k < n$ y $y^n x = dx^{n+1} - dx^n y - axy^n + (a+1)y^{n+1}$.

6. Todos los planos torcidos graduados están clasificados en uno de los casos de los ítems 3., 4., o 5., o de lo contrario existe $n > 1$ tal que $y^k x = x^{k+1} - x^k y + y^{k+1}$ para todo $k < n$ y $y^n x = dx^{n+1} - x^n y + (a+1)y^{n+1}$, donde $d(a+1) = 1$.

7. Los planos torcidos del caso del item 6. no están clasificados completamente, pero si está caracterizada una familia $B(a, L)$ dentro de ese caso, donde $L$ es una sucesión que llamamos "sucesión casi-balanceada".

**Problemas abiertos**

1. ¿Existe algún plano torcido graduado tal que las primeras 2n+1 filas correspondan al caso 2 de la Proposición 4.2.2 o al caso 2 de la Proposición 4.2.4?

2. ¿Algún plano torcido graduado correspondiente al caso 2 de la Proposición 3.1.6 tiene la propiedad de $m$-extensión para algún $m$?

3. Cálculo de la cohomología de Hochschild para los planos torcidos graduados.

# Bibliografía